\documentclass[leqno,openany,draft]{book}

\usepackage{makeidx}
\makeindex

\def\re{\mathop{\rm Re}}
\def\im{\mathop{\rm Im}}
\def\hom{\mathop{\rm Hom}}

\newtheorem{theorem}{Theorem}
\newtheorem{lemma}[theorem]{Lemma}
\newtheorem{proposition}[theorem]{Proposition}
\newtheorem{definition}[theorem]{Definition}
\newtheorem{corollary}[theorem]{Corollary}

\newcommand{\begintheorem}{\addtocounter{equation}{1}\begin{theorem}}
\newcommand{\beginlemma}{\addtocounter{equation}{1}\begin{lemma}}
\newcommand{\beginproposition}{\addtocounter{equation}{1}\begin{proposition}}
\newcommand{\begindefinition}{\addtocounter{equation}{1}\begin{definition}}
\newcommand{\begincorollary}{\addtocounter{equation}{1}\begin{corollary}}

\begin{document}

\frontmatter

\title{Some topics related to real, complex, \\ and Fourier analysis}

\author{Stephen Semmes \\
        Rice University}

\date{}

\maketitle

\chapter*{Preface}

        Here we look at some situations that are like the unit circle
or the real line in some ways, but which can be more complicated or fractal
in other ways.

\tableofcontents

\mainmatter

\chapter{The unit circle}
\label{unit circle}

\section{Integrable functions}
\label{integrable functions}

        Let ${\bf T}$ be the unit circle in the complex plane ${\bf C}$,
consisting of complex numbers $z$ whose modulus $|z|$ is equal to $1$.
Also let $f$ be a complex-valued function on ${\bf T}$ which is
integrable with respect to the usual arc-length measure on ${\bf T}$.
The $j$th \emph{Fourier coefficient}\index{Fourier coefficients}
$\widehat{f}(j)$ of $f$ is defined for each integer $j$ by
\begin{equation}
\label{widehat{f}(j) = frac{1}{2 pi} int_{bf T} f(z) overline{z}^j |dz|}
 \widehat{f}(j) = \frac{1}{2 \, \pi} \, \int_{\bf T} f(z) \, \overline{z}^j
                                                                    \, |dz|,
\end{equation}
where $\overline{z}$ is the complex conjugate of $z \in {\bf C}$, and
$|dz|$ is the element of arc-length measure on ${\bf T}$.  Note that
\begin{equation}
\label{|widehat{f}(j)| le frac{1}{2 pi} int_{bf T} |f(z)| |dz|}
        |\widehat{f}(j)| \le \frac{1}{2 \, \pi} \, \int_{\bf T} |f(z)| \, |dz|
\end{equation}
for every $j \in {\bf Z}$, where ${\bf Z}$ denotes the set of all
integers.

        Let $\exp z = \sum_{n = 0}^\infty z^n/n!$ be the usual complex
exponential function on ${\bf C}$.  It is well known that $\exp(i \, t)$
defines a mapping from the real line ${\bf R}$ onto ${\bf T}$, which
is periodic with period $2 \, \pi$.  The derivative of this mapping
has modulus equal to $1$ at each $t \in {\bf R}$, so that $\exp (i \, t)$
parameterizes ${\bf T}$ by arc length.  This permits the integral of a
function $h(z)$ on ${\bf T}$ with respect to arc-length measure to be
identified with the integral of $h(\exp(i \, t))$ on $[0, 2 \, \pi)$
  with respect to ordinary Lebesgue measure.  In particular, one can
  use this to check that
\begin{equation}
\label{int_{bf T} z^j |dz| = 0}
        \int_{\bf T} z^j \, |dz| = 0
\end{equation}
for each $j \in {\bf Z}$ with $j \ne 0$.

        If $f$ and $g$ are complex-valued square-integrable functions on
${\bf T}$, then put
\begin{equation}
\label{langle f, g rangle = frac{1}{2 pi} int_{bf T} f(z) overline{g(z)} |dz|}
 \langle f, g \rangle = \frac{1}{2 \, \pi} \, \int_{\bf T} \, f(z) \, 
                                                   \overline{g(z)} \, |dz|.
\end{equation}
Of course, the space $L^2({\bf T})$ of complex-valued
square-integrable functions on ${\bf T}$ is a Hilbert space with
respect to this inner product, and $\widehat{f}(j)$ is the same as the
inner product of $f$ with $g(z) = z^j$ for each $j \in {\bf Z}$.
Using (\ref{int_{bf T} z^j |dz| = 0}), we get that
\begin{equation}
\label{langle z^j, z^k rangle = frac{1}{2 pi} int_{bf T} z^{j - k} |dz| = 0}
 \langle z^j, z^k \rangle = \frac{1}{2 \, \pi} \, \int_{\bf T} z^{j - k} \, |dz|
                          = 0
\end{equation}
when $j, k \in {\bf Z}$ and $j \ne k$, because $\overline{z} = 1/z$
when $z \in {\bf C}$ satisfies $|z| = 1$.  Clearly
\begin{equation}
\label{langle z^j, z^j rangle = frac{1}{2 pi} int_{bf T} |dz| = 1}
        \langle z^j, z^j \rangle = \frac{1}{2 \, \pi} \, \int_{\bf T} |dz| = 1
\end{equation}
for every $j \in {\bf Z}$, so that the functions $z^j$ on ${\bf T}$
with $j \in {\bf Z}$ are orthonormal with respect to (\ref{langle f, g
  rangle = frac{1}{2 pi} int_{bf T} f(z) overline{g(z)} |dz|}).  If $f
\in L^2({\bf T})$, then it follows that
\begin{equation}
\label{sum_{j = -infty}^infty |widehat{f}(j)|^2 le langle f, f rangle = ...}
        \sum_{j = -\infty}^\infty |\widehat{f}(j)|^2 \le \langle f, f \rangle
                      = \frac{1}{2 \, \pi} \, \int_{\bf T} |f(z)|^2 \, |dz|,
\end{equation}
by standard arguments about inner product spaces.  In particular,
\begin{equation}
\label{lim_{|j| to infty} |widehat{f}(j)| = 0}
        \lim_{|j| \to \infty} |\widehat{f}(j)| = 0
\end{equation}
for every $f \in L^2({\bf T})$, which also holds for every $f \in
L^1({\bf T})$.  This uses the fact that $L^2({\bf T})$ is dense in
$L^1({\bf T})$, and the simple estimate (\ref{|widehat{f}(j)| le
  frac{1}{2 pi} int_{bf T} |f(z)| |dz|}).

        Let $C({\bf T})$ be the space of continuous complex-valued
functions on ${\bf T}$, and let $\mathcal{E}({\bf T})$ be the linear
subspace $C({\bf T})$ consisting of finite linear combinations of the
functions $z^j$ on ${\bf T}$ with $j \in {\bf Z}$.  The
Stone--Weierstrass theorem implies that $\mathcal{E}({\bf T})$ is
dense in $C({\bf T})$ with respect to the supremum norm on $C({\bf
  T})$, and hence that $\mathcal{E}({\bf T})$ is dense in $L^p({\bf
  T})$ when $p < \infty$, because $C({\bf T})$ is dense in $L^p({\bf
  T})$ when $p < \infty$.  Thus the functions $z^j$ on ${\bf T}$ with
$j \in {\bf Z}$ form an orthonormal basis for $L^2({\bf T})$.  It
follows that the Fourier series\index{Fourier series}
\begin{equation}
\label{sum_{j = -infty}^infty widehat{f}(j) z^j}
        \sum_{j = -\infty}^\infty \widehat{f}(j) \, z^j
\end{equation}
converges to $f$ with respect to the $L^2$ norm for every $f \in
L^2({\bf T})$, in the sense that the partial sums
\begin{equation}
\label{sum_{j = -N}^N widehat{f}(j) z^j}
        \sum_{j = -N}^N \widehat{f}(j) \, z^j
\end{equation}
converge to $f$ with respect to the $L^2$ norm as $N \to \infty$.
In particular,
\begin{equation}
\label{sum_{j = -infty}^infty |widehat{f}(j)|^2 = ...}
        \sum_{j = -\infty}^\infty |\widehat{f}(j)|^2 
                       = \frac{1}{2 \, \pi} \, \int_{\bf T} |f(z)|^2 \, |dz|
\end{equation}
for every $f \in L^2({\bf T})$.

\section{Abel sums and Cauchy products}
\label{abel sums, cauchy products}

        Let $\sum_{j = 0}^\infty a_j$ be an infinite series of real or
complex numbers, and suppose that
\begin{equation}
\label{sum_{j = 0}^infty a_j r^j}
        \sum_{j = 0}^\infty a_j \, r^j
\end{equation}
converges for every nonnegative real number $r$ with $r < 1$.  If the limit
\begin{equation}
\label{lim_{r to 1-} (sum_{j = 0}^infty a_j r^j)}
        \lim_{r \to 1-} \Big(\sum_{j = 0}^\infty a_j \, r^j\Big)
\end{equation}
exists, then $\sum_{j = 0}^\infty a_j$ is said to be \emph{Abel
  summable},\index{Abel summability} with the Abel sum equal to
(\ref{lim_{r to 1-} (sum_{j = 0}^infty a_j r^j)}).  Of course, if
(\ref{sum_{j = 0}^infty a_j r^j}) converges for some $r \ge 0$, then
$a_j \, r^j \to 0$ as $j \to \infty$, which implies that $\{a_j \,
r^j\}_{j = 0}^\infty$ is a bounded sequence.  Conversely, if $\{a_j \,
t^j\}_{j = 0}^\infty$ is a bounded sequence for some $t > 0$, then
(\ref{sum_{j = 0}^infty a_j r^j}) converges absolutely when $0 \le r <
t$, by comparison with the convergent geometric series $\sum_{j =
  0}^\infty (r/t)^j$.

        In particular, if $\sum_{j = 0}^\infty a_j$ converges, then
$\{a_j\}_{j = 0}^\infty$ is a bounded sequence, and hence 
(\ref{sum_{j = 0}^infty a_j r^j}) converges absolutely when $0 \le r < 1$.
In this case, it is well known that $\sum_{j = 0}^\infty a_j$ is Abel
summable, with Abel sum equal to the ordinary sum.  To see this, let
\begin{equation}
\label{s_n = sum_{j = 0}^n a_j}
        s_n = \sum_{j = 0}^n a_j
\end{equation}
be the $n$th partial sum of $\sum_{j = 0}^\infty a_j$ for each $n \ge 0$,
and put $s_{-1} = 0$.  Thus $a_n = s_n - s_{n - 1}$ for each $n \ge 0$,
and
\begin{eqnarray}
 \sum_{j = 0}^\infty a_j \, r^j & = & \sum_{j = 0}^\infty s_j \, r^j 
                                  - \sum_{j = 0}^\infty s_{j - 1} \, r^j \\
 & = & \sum_{j = 0}^\infty s_j \, r^j - \sum_{j = 0}^\infty s_j \, r^{j + 1} 
                                                            \nonumber \\
 & = & (1 - r) \, \sum_{j = 0}^\infty s_j \, r^j.              \nonumber
\end{eqnarray}
when $0 \le r < 1$.  Note that these series all converge when $0 \le r
<1$, because $\{s_j\}_{j = 0}^\infty$ is a bounded sequence too.  It
follows that
\begin{equation}
\label{sum_{j = 0}^infty a_j r^j - sum_{j = 0}^infty a_j = ...}
        \sum_{j = 0}^\infty a_j \, r^j - \sum_{j = 0}^\infty a_j
                = (1 - r) \, \sum_{j = 0}^\infty (s_j - s) \, r^j
\end{equation}
when $0 \le r < 1$, where $s = \sum_{j = 0}^\infty a_j$, and using the
fact that $(1 - r) \, \sum_{j = 0}^\infty r^j = 1$.  The remaining
point is to show that this converges to $0$ as $r \to 1-$, because
$s_j \to s$ as $j \to \infty$.  Note that this argument also works for
convergent series whose terms are elements of a real or complex vector
space equipped with a norm.

        If $\sum_{j = 0}^\infty a_j$ and $\sum_{k = 0}^\infty b_k$ are
any two infinite series of complex numbers, then their \emph{Cauchy
  product}\index{Cauchy products} is the infinite series $\sum_{l =
  0}^\infty c_l$, where
\begin{equation}
\label{c_l = sum_{j = 0}^l a_j b_{l - j}}
        c_l = \sum_{j = 0}^l a_j \, b_{l - j}
\end{equation}
for each $l \ge 0$.  It is easy to see that
\begin{equation}
\label{sum_{l = 0}^infty c_l = (sum_{j = 0}^infty a_j) (sum_{k = 0}^infty b_k)}
        \sum_{l = 0}^\infty c_l = \Big(\sum_{j = 0}^\infty a_j\Big)
                                 \, \Big(\sum_{k = 0}^\infty b_k\Big)
\end{equation}
formally, in the sense that every term $a_j \, b_k$ occurs exactly
once on both sides of the equation.  In particular, if $a_j = 0$ for
all but finitely many $j$ and $b_k = 0$ for all but finitely many $k$,
then $c_l = 0$ for all but finitely many $l$, and (\ref{sum_{l =
    0}^infty c_l = (sum_{j = 0}^infty a_j) (sum_{k = 0}^infty b_k)})
holds.  If $\sum_{j = 0}^\infty a_j$ and $\sum_{k = 0}^\infty b_k$
converge absolutely, then it is well known and not difficult to check
that $\sum_{l = 0}^\infty c_l$ also converges absolutely, and that
(\ref{sum_{l = 0}^infty c_l = (sum_{j = 0}^infty a_j) (sum_{k =
    0}^infty b_k)}) holds.  Note that
\begin{equation}
\label{|c_l| le sum_{j = 0}^l |a_j| |b_{l - j}|}
        |c_l| \le \sum_{j = 0}^l |a_j| \, |b_{l - j}|
\end{equation}
for each $l \ge 0$, where the right side of (\ref{|c_l| le sum_{j =
    0}^l |a_j| |b_{l - j}|}) corresponds exactly to the Cauchy product
of $\sum_{j = 0}^\infty |a_j|$ and $\sum_{k = 0} |b_k|$.

        Similarly,
\begin{equation}
\label{c_l r^l = sum_{j = 0}^l (a_j r^j) (b_{l - j} r^{l - j})}
        c_l \, r^l = \sum_{j = 0}^l (a_j \, r^j) \, (b_{l - j} \, r^{l - j})
\end{equation}
for every $r$, which corresponds to the Cauchy product of
\begin{equation}
\label{sum_{j = 0}^infty a_j r^j and sum_{k = 0}^infty b_k r^k}
        \sum_{j = 0}^\infty a_j \, r^j \quad\hbox{and}\quad
         \sum_{k = 0}^\infty b_k \, r^k.
\end{equation}
If the series in (\ref{sum_{j = 0}^infty a_j r^j and sum_{k = 0}^infty
  b_k r^k}) converge absolutely when $0 \le r <1$, then $\sum_{l =
  0}^\infty c_l \, r^l$ also converges absolutely when $0 \le r < 1$, with
\begin{equation}
\label{sum_{l = 0}^infty c_l r^l = ...}
 \sum_{l = 0}^\infty c_l \, r^l = \Big(\sum_{j = 0}^\infty a_j \, r^j\Big)
                                 \, \Big(\sum_{k = 0}^\infty b_k \, r^k\Big),
\end{equation}
as in the previous paragraph.  If, in addition, $\sum_{j = 0}^\infty
a_j$ and $\sum_{k = 0}^\infty b_k$ are Abel summable, then it follows
that $\sum_{l = 0}^\infty c_l$ is Abel summable as well, and that the
Abel sum of $\sum_{l = 0}^\infty c_l$ is equal to the product of the
Abel sums of $\sum_{j = 0}^\infty a_j$ and $\sum_{k = 0}^\infty b_k$.
There is an analogous statement for the Abel summability of $\sum_{j =
  0}^\infty (a_j + b_j)$, which is much simpler.

\section{The Poisson kernel}
\label{poisson kernel}

        Let $f$ be a complex-valued integrable function on the unit
circle again, and let us consider the Abel sums of the corresponding
Fourier series (\ref{sum_{j = -infty}^infty widehat{f}(j) z^j}).  More
precisely, put
\begin{equation}
\label{A_r(f)(z) = sum_{j = -infty}^infty widehat{f}(j) r^{|j|} z^j}
        A_r(f)(z) = \sum_{j = -\infty}^\infty \widehat{f}(j) \, r^{|j|} \, z^j
\end{equation}
for each $z \in {\bf T}$ and $r \in {\bf R}$ with $0 \le r <1$, where
the absolute convergence of the sum follows from the boundedness of
the Fourier coefficients $\widehat{f}(j)$, as in (\ref{|widehat{f}(j)|
  le frac{1}{2 pi} int_{bf T} |f(z)| |dz|}).  Of course, the
doubly-infinite Fourier series (\ref{sum_{j = -infty}^infty
  widehat{f}(j) z^j}) may be treated as a sum of two ordinary infinite
series, or it can be reduced to an ordinary infinite series by
combining the terms for $j$ and $-j$ when $j \ge 1$.  In both cases,
one gets (\ref{A_r(f)(z) = sum_{j = -infty}^infty widehat{f}(j)
  r^{|j|} z^j}) by taking the Abel sums of the associated ordinary
infinite series.

        Equivalently,
\begin{equation}
\label{A_r(f)(z) = frac{1}{2 pi} int_{bf T} f(w) p_r(w, z) |dw|}
        A_r(f)(z) = \frac{1}{2 \, \pi} \, \int_{\bf T} f(w) \, p_r(w, z) \, |dw|
\end{equation}
for every $z \in {\bf T}$ and $r \in {\bf R}$ with $0 \le r < 1$, where
\begin{equation}
\label{p_r(w, z) = sum_{j = -infty}^infty r^{|j|} z^j overline{w}^j}
        p_r(w, z) = \sum_{j = -\infty}^\infty r^{|j|} \, z^j \, \overline{w}^j.
\end{equation}
As before, this series converges absolutely when $w, z \in {\bf T}$
and $0 \le r < 1$, by comparison with a convergent geometric series.
In fact, the partial sums
\begin{equation}
\label{sum_{j = -N}^N r^{|j|} z^j overline{w}^j}
        \sum_{j = -N}^N r^{|j|} \, z^j \, \overline{w}^j
\end{equation}
converge to (\ref{p_r(w, z) = sum_{j = -infty}^infty r^{|j|} z^j
  overline{w}^j}) as $N \to \infty$ uniformly over $w, z \in {\bf T}$
when $0 \le r < 1$, by the Weierstrass $M$-test.  This permits one to
interchange the order of summation and integration to get
(\ref{A_r(f)(z) = frac{1}{2 pi} int_{bf T} f(w) p_r(w, z) |dw|}).
Similarly, if $0 \le t < 1$, then the partial sums (\ref{sum_{j =
    -N}^N r^{|j|} z^j overline{w}^j}) converge to (\ref{p_r(w, z) =
  sum_{j = -infty}^infty r^{|j|} z^j overline{w}^j}) as $N \to \infty$
uniformly over $w, z \in {\bf T}$ and $0 \le r \le t$.

        Alternatively, put
\begin{equation}
\label{p_r(z) = sum_{j = -infty}^infty r^{|j|} z^j}
        p_r(z) = \sum_{j = -\infty}^\infty r^{|j|} \, z^j
\end{equation}
for $z \in {\bf T}$ and $0 \le r < 1$, where the series converges
absolutely by comparison with a geometric series again.  If $0 \le t
<1$, then the partial sums
\begin{equation}
\label{sum_{j = -N}^N r^{|j|} z^j}
        \sum_{j = -N}^N r^{|j|} \, z^j
\end{equation}
converge to (\ref{p_r(z) = sum_{j = -infty}^infty r^{|j|} z^j}) as $N
\to \infty$ uniformly over $z \in {\bf T}$ and $0 \le r \le t$, for
the same reasons as before.  This implies that $p_r(z)$ is continuous
as a function of $z$ and $r$, with $z \in {\bf T}$ and $0 \le r < 1$.
By construction,
\begin{equation}
\label{p_r(w, z) = p_r(z overline{w})}
        p_r(w, z) = p_r(z \, \overline{w})
\end{equation}
for every $w, z \in {\bf T}$ and $0 \le r < 1$, and so it suffices to
compute $p_r(z)$ to determine $p_r(w, z)$.

        Observe that
\begin{equation}
\label{p_r(z) = ...}
        p_r(z) = \sum_{j = 0}^\infty r^j \, z^j 
                   + \sum_{j = 0}^\infty r^j \, \overline{z}^j - 1
\end{equation}
for every $z \in {\bf T}$ and $0 \le r < 1$, since $1/z = \overline{z}$
when $|z| = 1$.  Hence
\begin{equation}
\label{p_r(z) = 2 re sum_{j = 0}^infty r^j z^j - 1}
        p_r(z) = 2 \, \re \sum_{j = 0}^\infty r^j \, z^j - 1
 = 2 \, \re (1 - r \, z)^{-1} - 1 
\end{equation}
for every $z \in {\bf T}$ and $0 \le r < 1$, where $\re a$ denotes the
real part of a omplex number $a$, and summing the geometric series in
the second step.  It follows that
\begin{eqnarray}
\label{p_r(z) = ... = frac{1 - r^2}{|1 - r z|^2}}
        p_r(z) & = & \re \Big(\frac{2 - (1 - r \, z)}{1 - r \, z}\Big) \\
 & = & \re \Big(\frac{(1 + r \, z) \, (1 - r \, \overline{z})}
            {(1 - r \, z) \, (1 - r \, \overline{z})}\Big) \nonumber \\
 & = & \re \Big(\frac{1 + r \, (z - \overline{z}) - r^2 \, |z|^2}
                                    {|1 - r \, z|^2}\Big)   \nonumber \\
 & = & \frac{1 - r^2}{|1 - r \, z|^2}  \nonumber
\end{eqnarray}
for every $z \in {\bf T}$ and $0 \le r < 1$, because $z -
\overline{z}$ is purely imaginary and $|z| = 1$.  Combining this with
(\ref{p_r(w, z) = p_r(z overline{w})}), we get that
\begin{equation}
\label{p_r(w, z) = ... = frac{1 - r^2}{|w - r z|^2}}
        p_r(w, z) = \frac{1 - r^2}{|1 - r \, z \, \overline{w}|^2}
                  = \frac{1 - r^2}{|w - r \, z|^2}
\end{equation}
for every $w, z \in {\bf T}$ and $0 \le r < 1$.

\section{Continuous functions}
\label{continuous functions}

        It is easy to see that
\begin{equation}
\label{frac{1}{2 pi} int_{bf T} p_r(z) |dz| = 1}
        \frac{1}{2 \, \pi} \, \int_{\bf T} p_r(z) \, |dz| = 1
\end{equation}
for each $0 \le r < 1$, using (\ref{int_{bf T} z^j |dz| = 0}) and the
fact that the partial sums (\ref{sum_{j = -N}^N r^{|j|} z^j}) converge
to $p_r(z)$ as $N \to \infty$ uniformly over $z \in {\bf T}$.
This implies that
\begin{equation}
\label{frac{1}{2 pi} int_{bf T} p_r(w, z) |dw| = 1}
        \frac{1}{2 \, \pi} \, \int_{\bf T} p_r(w, z) \, |dw| = 1
\end{equation}
for every $z \in {\bf T}$ and $0 \le r < 1$, by (\ref{p_r(w, z) =
  p_r(z overline{w})}).  Thus (\ref{A_r(f)(z) = frac{1}{2 pi} int_{bf
    T} f(w) p_r(w, z) |dw|}) expresses $A_r(f)(z)$ as an average of
$f$ on ${\bf T}$, since $p_r(w, z)$ is also real-valued and
nonnegative.  In particular,
\begin{equation}
\label{A_r(f)(z) - f(z) = frac{1}{2 pi} int_{bf T} (f(w) - f(z)) p_r(w, z) |dw|}
        A_r(f)(z) - f(z) = \frac{1}{2 \, \pi} \, \int_{\bf T} (f(w) - f(z))
                                                     \, p_r(w, z) \, |dw|,
\end{equation}
and hence
\begin{equation}
\label{|A_r(f)(z) - f(z)| le ...}
        |A_r(f)(z) - f(z)| \le \frac{1}{2 \, \pi} \, \int_{\bf T} |f(w) - f(z)|
                                                        \, p_r(w, z) \, |dw|.
\end{equation}

        As $r$ approaches $1$, this average is concentrated near
$z$, because $p_r(w, z) \to 0$ as $r \to 1-$ uniformly over $w \in {\bf T}$
that do not get too close to $z$.  More precisely, if $w$ does not get too
close to $z$, then the denominator in (\ref{p_r(w, z) = ... = frac{1 -
    r^2}{|w - r z|^2}}) remains bounded away from $0$, while the numerator
obviously tends to $0$ as $r \to 1$.  If $f$ is continuous at $z$, then
then one can use this, (\ref{frac{1}{2 pi} int_{bf T} p_r(w, z) |dw| = 1}), 
and (\ref{|A_r(f)(z) - f(z)| le ...}) to show that
\begin{equation}
\label{lim_{r to 1-} A_r(f)(z) = f(z)}
        \lim_{r \to 1-} A_r(f)(z) = f(z).
\end{equation}
If $f$ is continuous at every point in ${\bf T}$, then it is well known
that $f$ is uniformly continuous on ${\bf T}$, because ${\bf T}$ is compact.
In this case, an analogous argument shows that $A_r(f)$ converges to $f$
as $r \to 1-$ uniformly on ${\bf T}$.

        Note that the partial sums
\begin{equation}
\label{sum_{j = -N}^N widehat{f}(j) r^{|j|} z^j}
        \sum_{j = -N}^N \widehat{f}(j) \, r^{|j|} \, z^j
\end{equation}
converge to $A_r(f)(z)$ as $N \to \infty$ uniformly over $z \in {\bf
  T}$, for any fixed $r$ with $0 \le r < 1$.  This follows from the
boundedness of the Fourier coefficients $\widehat{f}(j)$ and the
Weierstrass $M$-test, as before.  Thus $A_r(f)(z)$ can be approximated
by finite linear combinations of the $z^j$'s uniformly on ${\bf T}$
for each $r < 1$, and for any integrable function $f$ on ${\bf T}$.
If $f$ is continuous on ${\bf T}$, then $A_r(f)$ converges to $f$
as $r \to 1-$ uniformly on ${\bf T}$, as in the previous paragraph.
This implies that $f$ can be approximated uniformly by finite linear
combinations of the $z^j$'s on ${\bf T}$, by the corresponding statement
for $A_r(f)$ just mentioned.

        Because $p_r(w, z)$ is nonnegative and real-valued,
\begin{equation}
\label{|A_r(f)(z)| le frac{1}{2 pi} int_{bf T} |f(w)| p_r(w, z) |dw|}
 |A_r(f)(z)| \le \frac{1}{2 \, \pi} \int_{\bf T} |f(w)| \, p_r(w, z) \, |dw|
\end{equation}
for any integrable function $f$ on ${\bf T}$, $z \in {\bf T}$, and
$0 \le r < 1$.  In particular, if $f$ is a bounded measurable function
on ${\bf T}$, then
\begin{equation}
\label{|A_r(f)(z)| le sup_{w in {bf T}} |f(w)|}
        |A_r(f)(z)| \le \sup_{w \in {\bf T}} |f(w)|,
\end{equation}
by (\ref{frac{1}{2 pi} int_{bf T} p_r(w, z) |dw| = 1}).  Of course,
the same conclusion holds when $f$ is essentially bounded on ${\bf T}$,
with the supremum of $|f(w)|$ on ${\bf T}$ replaced with the essential
supremum.  This implies that
\begin{equation}
\label{||A_r(f)||_infty le ||f||_infty}
        \|A_r(f)\|_\infty \le \|f\|_\infty
\end{equation}
for every $0 \le r < 1$, where $\|f\|_\infty$ is the $L^\infty$ norm
of $f$ on ${\bf T}$.  The $L^\infty$ norm of $A_r(f)$ is the same as
the supremum norm of $A_r(f)$ for every $r < 1$, because $A_r(f)$ is a
continuous function on ${\bf T}$ for each $r < 1$, by the remarks
about uniform convergence of (\ref{sum_{j = -N}^N widehat{f}(j)
  r^{|j|} z^j}) to $A_r(f)$ in the preceding paragraph.

\section{$L^p$ Functions}
\label{L^p functions}

        Let $p$ be a real number with $p \ge 1$, and suppose that
$f \in L^p({\bf T})$.  Using (\ref{|A_r(f)(z)| le frac{1}{2 pi} 
int_{bf T} |f(w)| p_r(w, z) |dw|}), we get that
\begin{equation}
\label{|A_r(f)(z)|^p le (frac{1}{2 pi} int_{bf T} |f(w)| p_r(w, z) |dw|)^p}
        |A_r(f)(z)|^p \le \Big(\frac{1}{2 \, \pi} \, \int_{\bf T} |f(w)|
                                                \, p_r(w, z) \, |dw|\Big)^p 
\end{equation}
for every $z \in {\bf T}$ and $0 \le r < 1$.  This implies that
\begin{equation}
\label{|A_r(f)(z)|^p le frac{1}{2 pi} int_{bf T} |f(w)|^p p_r(w, z) |dw|}
 |A_r(f)(z)|^p \le \frac{1}{2 \, \pi} \, \int_{\bf T} |f(w)|^p \, p_r(w, z)
                                                               \, |dw|
\end{equation}
for every $z \in {\bf T}$ and $0 \le r < 1$, by Jensen's inequality,
and the fact that $t \mapsto t^p$ is a convex function on the set of
nonnegative real numbers when $p \ge 1$.  This also uses the
nonnegativity of $p_r(w, z)$ and (\ref{frac{1}{2 pi} int_{bf T} p_r(w,
  z) |dw| = 1}), so that the relevant integrals are in fact averages.

        Integrating (\ref{|A_r(f)(z)|^p le frac{1}{2 pi} int_{bf T} 
|f(w)|^p p_r(w, z) |dw|}) over $z$, we get that
\begin{equation}
\label{int_{bf T} |A_r(f)(z)|^p |dz| le ...}
 \int_{\bf T} |A_r(f)(z)|^p \, |dz| \le \frac{1}{2 \, \pi} \, 
               \int_{\bf T} \int_{\bf T} |f(w)|^p \, p_r(w, z) \, |dw| \, |dz|
\end{equation}
for every $0 \le r < 1$.  Equivalently,
\begin{equation}
\label{int_{bf T} |A_r(f)(z)|^p |dz| le ..., 2}
 \int_{\bf T} |A_r(f)(z)|^p \, |dz| \le \frac{1}{2 \, \pi} \,
                 \int_{\bf T} \int_{\bf T} |f(w)|^p \, p_r(w, z) \, |dz| \, |dw|,
\end{equation}
by Fubini's theorem.  Of course,
\begin{equation}
\label{frac{1}{2 pi} int_{bf T} p_r(w, z) |dz| = 1}
        \frac{1}{2 \, \pi} \, \int_{\bf T} p_r(w, z) \, |dz| = 1
\end{equation}
for every $w \in {\bf T}$ and $0 \le r < 1$, by (\ref{p_r(w, z) =
  p_r(z overline{w})}) and (\ref{frac{1}{2 pi} int_{bf T} p_r(z) |dz|
  = 1}).  Thus (\ref{int_{bf T} |A_r(f)(z)|^p |dz| le ..., 2}) reduces
to
\begin{equation}
\label{int_{bf T} |A_r(f)(z)|^p |dz| le int_{bf T} |f(w)|^p |dw|}
        \int_{\bf T} |A_r(f)(z)|^p \, |dz| \le \int_{\bf T} |f(w)|^p \, |dw|,
\end{equation}
which holds for every $0 \le r < 1$.

        Using this, one can check that
\begin{equation}
\label{A_r(f) to f as r to 1-}
        A_r(f) \to f \quad\hbox{as } r \to 1-
\end{equation}
with respect to the $L^p$ norm when $f \in L^p({\bf T})$ and $1 \le p
< \infty$.  More precisely, if $f$ is a continuous function on ${\bf
  T}$, then we have seen that (\ref{A_r(f) to f as r to 1-}) holds
uniformly on ${\bf T}$, and hence with respect to the $L^p$ norm.  If
$f \in L^p({\bf T})$ and $p < \infty$, then it is well known that $f$
can be approximated by continuous functions on ${\bf T}$ with respect
to the $L^p$ norm.  To get (\ref{A_r(f) to f as r to 1-}), one can use
the analogous statement for these approximations to $f$, and
(\ref{int_{bf T} |A_r(f)(z)|^p |dz| le int_{bf T} |f(w)|^p |dw|}) to
estimate the errors.  It is well known that (\ref{lim_{r to 1-}
  A_r(f)(z) = f(z)}) also holds for almost every $z \in {\bf T}$ with
respect to arc-length measure when $f \in L^1({\bf T})$, but we shall
not go into this here.

        In the $p = 2$ case, we have that
\begin{equation}
\label{frac{1}{2 pi} int_{bf T} |A_r(f)(z)|^2 |dz| = ...}
        \frac{1}{2 \, \pi} \, \int_{\bf T} |A_r(f)(z)|^2 \, |dz|
                        = \sum_{j = -\infty}^\infty |\widehat{f}(j)|^2 \, r^{2 j}
\end{equation}
for every $0 \le r < 1$, because of the orthonormality of the
functions $z^j$ on ${\bf T}$ with respect to the usual integral inner
product.  This gives another way to look at (\ref{int_{bf T}
  |A_r(f)(z)|^p |dz| le int_{bf T} |f(w)|^p |dw|}) when $f \in
L^2({\bf T})$, since
\begin{equation}
\label{sum_{j = -infty}^infty |widehat{f}(j)|^2 r^{2 j} le ...}
        \sum_{j = -\infty}^\infty |\widehat{f}(j)|^2 \, r^{2 j}
            \le \sum_{j = -\infty}^\infty |\widehat{f}(j)|^2
\end{equation}
for each $r < 1$.  Using orthonormality of the $z^j$'s again, we get that
\begin{equation}
\label{frac{1}{2 pi} int_{bf T} |f(z) - A_r(f)(z)|^2 |dz| = ...}
        \frac{1}{2 \, \pi} \, \int_{\bf T} |f(z) - A_r(f)(z)|^2 \, |dz|
                  = \sum_{j = -\infty}^\infty |\widehat{f}(j)|^2 \, (1 - r)^2
\end{equation}
for every $f \in L^2({\bf T})$ and $0 \le r < 1$, which gives another
way to look at the convergence (\ref{A_r(f) to f as r to 1-}) with
respect to the $L^2$ norm.

        If $f \in L^\infty({\bf T})$, then (\ref{A_r(f) to f as r to 1-})
may not hold with respect to the $L^\infty$ norm on ${\bf T}$, because
$f$ may not be equal to a continuous function almost everywhere on
${\bf T}$.  Of course, $L^\infty({\bf T}) \subseteq L^p({\bf T})$ for
every $p < \infty$, and so (\ref{A_r(f) to f as r to 1-}) holds with
respect to the $L^p$ norm when $1 \le p < \infty$, as before.
Similarly, (\ref{lim_{r to 1-} A_r(f)(z) = f(z)}) still holds for
almost every $z \in {\bf T}$ with respect to arc-length measure.
Remember that $L^\infty({\bf T})$ can be identified with the dual of
$L^1({\bf T})$ in the usual way.  Let us check that (\ref{A_r(f) to f
  as r to 1-}) also holds with respect to the corresponding weak$^*$
topology on $L^\infty({\bf T})$.

        More precisely, this means that
\begin{equation}
\label{lim_{r to 1-} frac{1}{2 pi} int_{bf T} A_r(f)(z) phi(z) |dz| = ...}
 \lim_{r \to 1-} \frac{1}{2 \, \pi} \, \int_{\bf T} A_r(f)(z) \, \phi(z) \, |dz|
        = \frac{1}{2 \, \pi} \, \int_{\bf T} f(z) \, \phi(z) \, |dz|
\end{equation}
for every $\phi \in L^1({\bf T})$.  If $\phi \in L^\infty({\bf T})$,
then (\ref{lim_{r to 1-} frac{1}{2 pi} int_{bf T} A_r(f)(z) phi(z)
  |dz| = ...})  follows from the fact that (\ref{A_r(f) to f as r to
  1-}) holds with respect to the $L^1$ norm on ${\bf T}$, as before.
One can then use this to show that (\ref{lim_{r to 1-} frac{1}{2 pi}
  int_{bf T} A_r(f)(z) phi(z) |dz| = ...}) holds for every $\phi \in
L^1({\bf T})$, by approximating $\phi$ by bounded measurable functions
on ${\bf T}$ with respect to the $L^1$ norm, and estimating the errors
with (\ref{||A_r(f)||_infty le ||f||_infty}).  One could just as well
start with $\phi \in L^q({\bf T})$ for any $q > 1$, and then use
(\ref{A_r(f) to f as r to 1-}) with respect to the $L^p$ norm on ${\bf
  T}$, where $p < \infty$ is conjugate to $q$.

        Alternatively, one can verify that
\begin{equation}
\label{frac{1}{2 pi} int_{bf T} A_r(f)(z) phi(z) |dz| = ...}
        \frac{1}{2 \, \pi} \, \int_{\bf T} A_r(f)(z) \, \phi(z) \, |dz|
          = \frac{1}{2 \, \pi} \, \int_{\bf T} f(w) \, A_r(\phi)(w) \, |dw|
\end{equation}
for every $\phi \in L^1({\bf T})$ and $0 \le r < 1$, using Fubini's
theorem.  If one starts with the original definition (\ref{A_r(f)(z) =
  sum_{j = -infty}^infty widehat{f}(j) r^{|j|} z^j}) of $A_r(f)(z)$,
then one should also interchange the order of summation and
integration, and use the change of variables $j \mapsto -j$ in the
sum.  Otherwise, one can use the expression (\ref{A_r(f)(z) =
  frac{1}{2 pi} int_{bf T} f(w) p_r(w, z) |dw|}) for $A_r(f)(z)$, and
observe that $p_r(w, z)$ is symmetric in $w$ and $z$.  Because
$A_r(\phi) \to \phi$ as $r \to 1-$ with respect to the $L^1$ norm on
${\bf T}$ and $f \in L^\infty({\bf T})$, we have that
\begin{equation}
\label{lim_{r to 1-} frac{1}{2 pi} int_{bf T} f(w) A_r(phi)(w) |dw| = ...}
 \lim_{r \to 1-} \frac{1}{2 \, \pi} \, \int_{\bf T} f(w) \, A_r(\phi)(w) \, |dw|
   = \frac{1}{2 \, \pi} \, \int_{\bf T} f(w) \, \phi(w) \, |dw|.
\end{equation}
Combining this with (\ref{frac{1}{2 pi} int_{bf T} A_r(f)(z) phi(z)
  |dz| = ...}), we get (\ref{lim_{r to 1-} frac{1}{2 pi} int_{bf T}
  A_r(f)(z) phi(z) |dz| = ...}), as desired.

\section{Harmonic functions}
\label{harmonic functions}

        Let $f$ be a complex-valued integrable function on ${\bf T}$ again, 
and let
\begin{equation}
\label{D = {zeta in {bf C} : |zeta| < 1}}
        D = \{\zeta \in {\bf C} : |\zeta| < 1\}
\end{equation}
be the open unit disk in the complex plane.  Thus the closure
$\overline{D}$ of $D$ in ${\bf C}$ is the closed unit disk in ${\bf C}$, 
\begin{equation}
\label{overline{D} = {zeta in {bf C} : |zeta| le 1}}
        \overline{D} = \{\zeta \in {\bf C} : |\zeta| \le 1\},
\end{equation}
which is the same as the union of $D$ and the unit circle ${\bf T}$.
Put
\begin{equation}
\label{h(zeta) = ...}
        h(\zeta) = \sum_{j = 0}^\infty \widehat{f}(j) \, \zeta^j
                   + \sum_{j = 1}^\infty \widehat{f}(-j) \, \overline{\zeta}^j
\end{equation}
for each $\zeta \in D$, where the series converge absolutely for each
$\zeta \in D$, because the Fourier coefficients of $f$ are bounded, as
in (\ref{|widehat{f}(j)| le frac{1}{2 pi} int_{bf T} |f(z)| |dz|}).
Similarly, if $0 \le t < 1$, then the partial sums of these series
also converge uniformly on the set of $\zeta \in {\bf C}$ such that
$|\zeta| \le t$, by the Weierstrass $M$-test.  In particular, this
implies that $h$ is continuous on $D$.

        More precisely, let us put
\begin{equation}
\label{h_+(zeta) = sum_{j = 0}^infty widehat{f}(j) zeta^j}
        h_+(\zeta) = \sum_{j = 0}^\infty \widehat{f}(j) \, \zeta^j
\end{equation}
and
\begin{equation}
\label{h_-(zeta) = sum_{j = 1}^infty widehat{f}(-j) overline{zeta}^j}
        h_-(\zeta) = \sum_{j = 1}^\infty \widehat{f}(-j) \, \overline{\zeta}^j
\end{equation}
for every $\zeta \in D$.  As before, both of these series converge
absolutely for each $\zeta \in D$, and their partial sums converge
uniformly on compact subsets of $D$, because the Fourier coefficients
of $f$ are bounded, as in (\ref{|widehat{f}(j)| le frac{1}{2 pi}
  int_{bf T} |f(z)| |dz|}).  Of course, $h_+(\zeta)$ is a holomorphic
function of $\zeta$ on $D$, and $h_-(\zeta)$ is the complex-conjugate
of a holomorphic function of $\zeta$ on $D$.  This implies that
$h_+(\zeta)$ and $h_-(\zeta)$ are both harmonic functions of $\zeta$
on $D$, so that
\begin{equation}
\label{h(zeta) = h_+(zeta) + h_-(zeta)}
        h(\zeta) = h_+(\zeta) + h_-(\zeta)
\end{equation}
is a harmonic function of $\zeta$ on $D$ as well.  

        Observe that
\begin{equation}
\label{h(r z) = A_r(f)(z)}
        h(r \, z) = A_r(f)(z)
\end{equation}
for every $z \in {\bf T}$ and $0 \le r < 1$, where $A_r(f)(z)$ is as
in Section \ref{poisson kernel}.  Suppose for the moment that $f$ is a
continuous function on ${\bf T}$, and consider the function on
$\overline{D}$ equal to $h$ on $D$ and to $f$ on ${\bf T}$.  One can
check that this function is continuous on $\overline{D}$, for
essentially the same reasons as in Section \ref{continuous functions}.
We already know that $h$ is continuous on $D$ for any integrable
function $f$ on ${\bf T}$, and so the main point is to show that this
extension of $h$ to $\overline{D}$ is continuous at each point in
${\bf T}$.  Because $f$ is supposed to be continuous on ${\bf T}$,
this means that for each $z \in {\bf T}$,
\begin{equation}
\label{h(zeta) to f(z)}
        h(\zeta) \to f(z)
\end{equation}
as $\zeta \in D$ approaches $z$, which is analogous to (\ref{lim_{r to
    1-} A_r(f)(z) = f(z)}).

        It is well known that a continuous function $g$ on $\overline{D}$
that is harmonic on $D$ is uniquely determined by its restriction to
${\bf T}$, because of the maximum principle.  If $f$ denotes the restriction
of $g$ to ${\bf T}$, and if $h$ is defined on $D$ as in (\ref{h(zeta) = ...}),
then it follows that $g = h$ on $D$.  If $g$ is any harmonic function on
$D$, then similar arguments can be applied to
\begin{equation}
\label{g_t(zeta) = g(t zeta)}
        g_t(\zeta) = g(t \, \zeta)
\end{equation}
for each $0 \le t < 1$.  More precisely, if $f_t(z)$ is the
restriction of $g_t$ to ${\bf T}$, and if $h_t$ is defined on $D$ as
in (\ref{h(zeta) = ...}) with $f$ replaced by $f_t$, then $g_t = h_t$
on $D$.  If $f_t(z)$ converges as $t \to 1-$ to an integrable function
$f(z)$ on ${\bf T}$ with respect to the $L^1$ norm, and if $h$ is
defined on $D$ as in (\ref{h(zeta) = ...}) again, then it is easy to
see that $h_t \to h$ as $t \to 1-$ pointwise on $D$, and hence that $g
= h$ on $D$.

\section{Complex Borel measures}
\label{complex borel measures}

        Let $\mu$ be a complex Borel measure on the unit circle ${\bf T}$,
and let $|\mu|$ be the corresponding total variation measure on ${\bf T}$.
Thus $|\mu|$ is a finite nonnegative Borel measure on ${\bf T}$ such that
\begin{equation}
\label{|mu(E)| le |mu|(E)}
        |\mu(E)| \le |\mu|(E)
\end{equation}
for every Borel set $E \subseteq {\bf T}$, and in fact $|\mu|$ is the
smallest nonnegative Borel measure on ${\bf T}$ that satisfies
(\ref{|mu(E)| le |mu|(E)}).  If $a(z)$ is a bounded complex-valued
Borel measurable function on ${\bf T}$, then the integral of $a(z)$
with respect to $\mu$ on ${\bf T}$ can be defined in a standard way,
and satisfies
\begin{equation}
\label{|int_{bf T} a(z) d mu(z)| le int_{bf T} |a(z)| d |mu|(z)}
 \biggl|\int_{\bf T} a(z) \, d\mu(z)\biggr| \le \int_{\bf T} |a(z)| \, d|\mu|(z).
\end{equation}
Of course, continuous functions on ${\bf T}$ are automatically Borel
measurable, and so the $j$th Fourier coefficient\index{Fourier coefficients}
of $\mu$ may be defined for each integer $j$ by
\begin{equation}
\label{widehat{mu}(j) = int_{bf T} overline{z}^j d mu(z)}
        \widehat{\mu}(j) = \int_{\bf T} \overline{z}^j \, d\mu(z).
\end{equation}
It follows that
\begin{equation}
\label{|widehat{mu}(j)| le |mu|({bf T})}
        |\widehat{\mu}(j)| \le |\mu|({\bf T})
\end{equation}
for each $j \in {\bf Z}$, by taking $a(z) = \overline{z}^j$ in
(\ref{|int_{bf T} a(z) d mu(z)| le int_{bf T} |a(z)| d |mu|(z)}).

        This is compatible with the earlier definition for integrable
functions, in the following sense.  If $f \in L^1({\bf T})$, then it is
well known that
\begin{equation}
\label{mu(E) = frac{1}{2 pi} int_E f(z) |dz|}
        \mu(E) = \frac{1}{2 \, \pi} \, \int_E f(z) \, |dz|
\end{equation}
defines a Borel measure on ${\bf T}$, and that
\begin{equation}
\label{|mu|(E) = frac{1}{2 pi} int_E |f(z)| |dz|}
        |\mu|(E) = \frac{1}{2 \, \pi} \, \int_E |f(z)| \, |dz|
\end{equation}
for every Borel set $E \subseteq {\bf T}$.  Moreover,
\begin{equation}
\label{int_{bf T} a(z) d mu(z) = frac{1}{2 pi} int_{bf T} a(z) f(z) |dz|}
        \int_{\bf T} a(z) \, d\mu(z) 
                  = \frac{1}{2 \, \pi} \, \int_{\bf T} a(z) \, f(z) \, |dz|
\end{equation}
for every bounded complex-valued Borel measurable function $a(z)$ on
${\bf T}$, by approximating $a(z)$ uniformly on ${\bf T}$ by Borel
measurable simple functions on ${\bf T}$.  This implies that
\begin{equation}
\label{widehat{mu}(j) = widehat{f}(j)}
        \widehat{\mu}(j) = \widehat{f}(j)
\end{equation}
for every $j \in {\bf Z}$ in this case, by taking $a(z) = z^j$ in
(\ref{int_{bf T} a(z) d mu(z) = frac{1}{2 pi} int_{bf T} a(z) f(z)
  |dz|}).

        Let $\mu$ be an arbitrary complex Borel measure on ${\bf T}$ again,
and put
\begin{equation}
\label{A_r(mu)(z) = sum_{j = -infty}^infty widehat{mu}(j) r^{|j|} z^j}
        A_r(\mu)(z) = \sum_{j = -\infty}^\infty \widehat{\mu}(j) \, r^{|j|} \, z^j
\end{equation}
for $z \in {\bf T}$ and $0 \le r < 1$, as in Section \ref{poisson kernel}.
As before, this series converges absolutely for every $z \in {\bf T}$
and $0 \le r < 1$, because the Fourier coefficients of $\mu$ are bounded,
as in (\ref{|widehat{mu}(j)| le |mu|({bf T})}).  The partial sums
\begin{equation}
\label{sum_{j = -N}^N widehat{mu}(j) r^{|j|} z^j}
        \sum_{j = -N}^N \widehat{\mu}(j) \, r^{|j|} \, z^j
\end{equation}
also converge to $A_r(\mu)(z)$ as $N \to \infty$ uniformly over $z \in
{\bf T}$ for any fixed $r < 1$, by the Weierstrass $M$-test, which
implies that $A_r(\mu)(z)$ is a continuous function of $z$ on ${\bf
  T}$ for each $r < 1$.  Similarly,
\begin{equation}
\label{h(zeta) = ..., 2}
        h(\zeta) = \sum_{j = 0}^\infty \widehat{\mu}(j) \, \zeta^j
                     + \sum_{j = 1}^\infty \widehat{\mu}(-j) \, \overline{\zeta}^j
\end{equation}
defines a harmonic function on the open unit disk $D$, as in the
previous section.  We also have that
\begin{equation}
\label{h(r z) = A_r(mu)(z)}
        h(r \, z) = A_r(\mu)(z)
\end{equation}
for every $z \in {\bf T}$ and $0 \le r < 1$, as in (\ref{h(r z) =
  A_r(f)(z)}).

        In analogy with Section \ref{poisson kernel},
\begin{equation}
\label{A_r(mu)(z) = int_{bf T} p_r(w, z) d mu(w)}
        A_r(\mu)(z) = \int_{\bf T} p_r(w, z) \, d\mu(w)
\end{equation}
for every $z \in {\bf T}$ and $0 \le r <1$, and hence
\begin{equation}
\label{|A_r(mu)(z)| le int_{bf T} p_r(w, z) d |mu|(w)}
        |A_r(\mu)(z)| \le \int_{\bf T} p_r(w, z) \, d|\mu|(w)
\end{equation}
for every $z \in {\bf T}$ and $0 \le r < 1$.  Integrating over $z$
and interchanging the order of integration, we get that
\begin{eqnarray}
\label{frac{1}{2 pi} int_{bf T} |A_r(mu)(z)| |dz| le ...}
        \frac{1}{2 \, \pi} \, \int_{\bf T} |A_r(\mu)(z)| \, |dz|
           & \le & \frac{1}{2 \, \pi} \, \int_{\bf T} \int_{\bf T} p_r(w, z)
                                                   \, d|\mu|(w) \, |dz| \\
 & = & \frac{1}{2 \, \pi} \, \int_{\bf T} \int_{\bf T} p_r(w, z) \, |dz| \,
                                                         d|\mu|(w) \nonumber
\end{eqnarray}
for every $0 \le r < 1$.  It follows that
\begin{equation}
\label{frac{1}{2 pi} int_{bf T} |A_r(mu)(z)| |dz| le |mu|({bf T})}
 \frac{1}{2 \, \pi} \, \int_{\bf T} |A_r(\mu)(z)| \, |dz| \le |\mu|({\bf T}),
\end{equation}
because of (\ref{frac{1}{2 pi} int_{bf T} p_r(w, z) |dz| = 1}).  This
corresponds exactly to (\ref{int_{bf T} |A_r(f)(z)|^p |dz| le int_{bf
    T} |f(w)|^p |dw|}) with $p = 1$ when $\mu$ is given by an integrable
function $f$ on ${\bf T}$ as in (\ref{mu(E) = frac{1}{2 pi} int_E f(z) |dz|}).

        If $\phi$ is a continuous complex-valued function on the unit circle,
then
\begin{eqnarray}
\label{frac{1}{2 pi} int_{bf T} A_r(mu)(z) phi(z) |dz| = ...}
        \frac{1}{2 \, \pi} \, \int_{\bf T} A_r(\mu)(z) \, \phi(z) \, |dz|
 & = & \sum_{j = -\infty}^\infty \frac{\widehat{\mu}(j) \, r^{|j|}}{2 \, \pi} \, 
                              \int_{\bf T} \phi(z) \, z^j \, |dz| \\
 & = & \sum_{j = -\infty}^\infty \widehat{\mu}(j) \, \widehat{\phi}(-j) \, r^{|j|}
                                                                 \nonumber
\end{eqnarray}
for every $0 \le r < 1$.  This follows by plugging the definition
(\ref{A_r(mu)(z) = sum_{j = -infty}^infty widehat{mu}(j) r^{|j|} z^j})
of $A_r(\mu)(z)$ into the integral on the left side, and then
interchanging the order of summation and integration, which is
possible because of the uniform convergence of the partial sums
(\ref{sum_{j = -N}^N widehat{mu}(j) r^{|j|} z^j}).  Equivalently,
\begin{equation}
\label{frac{1}{2 pi} int_{bf T} A_r(mu)(z) phi(z) |dz| = ..., 2}
 \frac{1}{2 \, \pi} \, \int_{\bf T} A_r(\mu)(z) \, \phi(z) \, |dz|
 = \sum_{j = -\infty}^\infty \widehat{\phi}(j) \, r^{|j|} \, \widehat{\mu}(-j),
\end{equation}
for every $0 \le r < 1$, using the change of variables $j \mapsto -j$.
Plugging the definition (\ref{widehat{mu}(j) = int_{bf T}
  overline{z}^j d mu(z)}) of $\widehat{\mu}(-j)$ into the sum on the
right side, we get that
\begin{eqnarray}
  \sum_{j = -\infty}^\infty \widehat{\phi}(j) \, r^{|j|} \, \widehat{\mu}(-j)
 & = & \sum_{j = -\infty}^\infty \int_{\bf T} \widehat{\phi}(j) \, r^{|j|} \,
                                                           z^j \, d\mu(z) \\
        & = & \int_{\bf T} A_r(\phi)(z) \, d\mu(z) \nonumber
\end{eqnarray}
for every $0 \le r < 1$.  Here $A_r(\phi)(z)$ is defined in the same
way as before, and we can interchange the order of summation and
integration in the second step for the usual reasons of uniform
convergence.

        This shows that
\begin{equation}
\label{frac{1}{2 pi} int_{bf T} A_r(mu)(z) phi(z) |dz| = ..., 3}
        \frac{1}{2 \, \pi} \, \int_{\bf T} A_r(\mu)(z) \, \phi(z) \, |dz|
                                  = \int_{\bf T} A_r(\phi)(z) \, d\mu(z)
\end{equation}
for every continuous function $\phi$ on ${\bf T}$ and $0 \le r < 1$.
One could also get this using (\ref{A_r(mu)(z) = int_{bf T} p_r(w, z)
  d mu(w)}), the analogous expression for $A_r(\phi)$, and the fact
that $p_r(w, z)$ is symmetric in $w$ and $z$.  Note that
\begin{equation}
\label{lim_{r to 1-} int_{bf T} A_r(phi)(z) d mu(z) = int_{bf T} phi(z) d mu(z)}
        \lim_{r \to 1-} \int_{\bf T} A_r(\phi)(z) \, d\mu(z) 
           = \int_{\bf T} \phi(z) \, d\mu(z)
\end{equation}
for every continuous function $\phi$ on ${\bf T}$, because $A_r(\phi)$
converges to $\phi$ as $r \to 1-$ uniformly on ${\bf T}$, as in
Section \ref{continuous functions}.  It follows that
\begin{equation}
\label{lim_{r to 1-} frac{1}{2 pi} int_{bf T} A_r(mu)(z) phi(z) |dz| = ...}
 \lim_{r \to 1-} \frac{1}{2 \, \pi} \, \int_{\bf T} A_r(\mu)(z) \, \phi(z) \, |dz|
                                     = \int_{\bf T} \phi(z) \, d\mu(z)
\end{equation}
for every continuous function $\phi$ on ${\bf T}$.

        Of course,
\begin{equation}
\label{lambda(phi) = int_{bf T} phi(z) d mu(z)}
        \lambda(\phi) = \int_{\bf T} \phi(z) \, d\mu(z)
\end{equation}
defines a linear functional on the vector space $C({\bf T})$ of
continuous complex-valued functions $\phi$ on ${\bf T}$.  More
precisely, this is a bounded linear functional on $C({\bf T})$ with
respect to the supremum norm, because
\begin{equation}
\label{|lambda(phi)| le |mu|({bf T}) sup_{z in {bf T}} |phi(z)|}
        |\lambda(\phi)| \le |\mu|({\bf T}) \, \sup_{z \in {\bf T}} |\phi(z)|
\end{equation}
for every $\phi \in C({\bf T})$, by (\ref{|int_{bf T} a(z) d mu(z)| le
  int_{bf T} |a(z)| d |mu|(z)}).  A version of the Riesz
representation theorem implies that every bounded linear functional on
$C({\bf T})$ corresponds to a unique Borel measure on ${\bf T}$ in
this way.

        Let $\mu_r$ be the Borel measure on ${\bf T}$ defined by
\begin{equation}
\label{mu_r(E) = frac{1}{2 pi} int_E A_r(mu) |dz|}
        \mu_r(E) = \frac{1}{2 \, \pi} \, \int_E A_r(\mu) \, |dz|
\end{equation}
for every Borel set $E \subseteq {\bf T}$, and for each real number
$r$ with $0 \le r < 1$.  Thus
\begin{equation}
\label{lambda_r(phi) = int_{bf T} phi(z) d mu_r(z) = ...}
        \lambda_r(\phi) = \int_{\bf T} \phi(z) \, d\mu_r(z)
          = \frac{1}{2 \, \pi} \, \int_{\bf T} A_r(\mu)(z) \, \phi(z) \, |dz|
\end{equation}
is the bounded linear functional on $C({\bf T})$ that corresponds to
$\mu_r$ as in the preceding paragraph.  In this notation, (\ref{lim_{r
    to 1-} frac{1}{2 pi} int_{bf T} A_r(mu)(z) phi(z) |dz| = ...})
says that
\begin{equation}
\label{lim_{r to 1-} lambda_r(phi) = lambda(phi)}
        \lim_{r \to 1-} \lambda_r(\phi) = \lambda(\phi)
\end{equation}
for every $\phi \in C({\bf T})$.  This is the same as saying that
$\lambda_r \to \lambda$ as $r \to 1-$ with respect to the weak$^*$
topology on the space of bounded linear functionals on $C({\bf T})$.

\section{Holomorphic functions}
\label{holomorphic functions}

        Let $f$ be an integrable function on the unit circle again,
and let $h$, $h_+$, and $h_-$ be defined on the open unit disk as in
Section \ref{harmonic functions}.  If
\begin{equation}
\label{widehat{f}(-j) = 0}
        \widehat{f}(-j) = 0
\end{equation}
for each positive integer $j$, then $h_-(\zeta) = 0$ for every $\zeta
\in D$, and hence
\begin{equation}
\label{h(zeta) = h_+(zeta)}
        h(\zeta) = h_+(\zeta)
\end{equation}
is a holomorphic function of $\zeta$ on $D$.  In this case,
\begin{equation}
\label{A_r(f)(z) = sum_{j = 0}^infty widehat{f}(j) r^j z^j}
        A_r(f)(z) = \sum_{j = 0}^\infty \widehat{f}(j) \, r^j \, z^j
\end{equation}
for each $z \in {\bf T}$ and $0 \le r < 1$, and the partial sums
\begin{equation}
\label{sum_{j = 0}^N widehat{f}(j) r^j z^j}
        \sum_{j = 0}^N \widehat{f}(j) \, r^j \, z^j
\end{equation}
converge to (\ref{A_r(f)(z) = sum_{j = 0}^infty widehat{f}(j) r^j
  z^j}) as $N \to \infty$ uniformly over $z \in {\bf T}$ for any fixed
$r < 1$.  If $f$ is continuous on ${\bf T}$, then we have seen that
$A_r(f)$ converges to $f$ uniformly on ${\bf T}$ as $r \to 1-$, which
implies that $f$ can be approximated uniformly by finite linear
combinations of $z^j$'s with $j \ge 0$.  Similarly, if $f \in L^p({\bf
  T})$ for some $p$, $1 \le p < \infty$, then $A_r(f)$ converges to
$f$ with respect to the $L^p$ norm on ${\bf T}$ as $r \to 1-$, which
implies that $f$ can be approximated by finite linear combinations of
the $z^j$'s with $j \ge 0$ with respect to the $L^p$ norm on ${\bf T}$.

        If $g$ is any holomorphic function on $D$, then
\begin{equation}
\label{int_{bf T} g(t z) z^j dz = 0}
        \int_{\bf T} g(t \, z) \, z^j \, dz = 0
\end{equation}
for every $t \in {\bf R}$ with $0 \le t < 1$ and nonnegative integer
$j$, by Cauchy's theorem.  The complex line integral element $dz$ is
equal to $i \, z \, |dz|$ on the unit circle, so that (\ref{int_{bf T}
  g(t z) z^j dz = 0}) is equivalent to
\begin{equation}
\label{int_{bf T} g(t z) z^{j + 1} |dz| = 0}
        \int_{\bf T} g(t \, z) \, z^{j + 1} \, |dz| = 0
\end{equation}
for each nonnegative integer $j$.  If we put
\begin{equation}
\label{f_t(z) = g(t z)}
        f_t(z) = g(t \, z)
\end{equation}
for each $z \in {\bf T}$ and $0 \le t < 1$, then (\ref{int_{bf T} g(t
  z) z^{j + 1} |dz| = 0}) is the same as saying that
\begin{equation}
\label{widehat{f_t}(-j - 1) = 0}
        \widehat{f_t}(-j - 1) = 0
\end{equation}
for each nonnegative integer $j$.  If $f_t(z)$ converges as $t \to 1-$
to an integrable function $f(z)$ on ${\bf T}$ with respect to the
$L^1$ norm, then it follows that
\begin{equation}
\label{widehat{f}(-j - 1) = 0}
        \widehat{f}(-j - 1) = 0
\end{equation}
for every nonnegative integer $j$.  Under these conditions, the Cauchy
integral formula can be used to express $g(\zeta)$ in terms of an
integral of $f$ on ${\bf T}$ for each $\zeta \in D$, which amounts to
the identification of $g(\zeta)$ with (\ref{h(zeta) = h_+(zeta)}).

        Now let $\mu$ be a complex Borel measure on ${\bf T}$, and consider 
the corresponding harmonic function $h(\zeta)$ on $D$, as in 
(\ref{h(zeta) = ..., 2}).  If $\widehat{\mu}(-j) = 0$ for every
positive integer $j$, then $h(\zeta)$ is holomorphic on $D$, as before.
In this case, a famous theorem of F.\ and M.~Riesz implies that $\mu$
is absolutely continuous with respect to arc-length measure on ${\bf
  T}$, so that $\mu$ can be represented as in (\ref{mu(E) = frac{1}{2
    pi} int_E f(z) |dz|}) for some $f \in L^1({\bf T})$.

\section{Products of functions}
\label{products of functions}

        Remember that $\mathcal{E}({\bf T})$ is the space of continuous
functions on ${\bf T}$ that can be expressed as finite linear combinations
of the functions $z^j$ with $j \in {\bf Z}$.  If $f$, $g$ are elements of
$\mathcal{E}({\bf T})$, then
\begin{equation}
\label{sums for f(z), g(z)}
 f(z) = \sum_{j = -\infty}^\infty \widehat{f}(j) \, z^j, \quad
 g(z) = \sum_{k = -\infty}^\infty \widehat{g}(k) \, z^k
\end{equation}
for every $z \in {\bf T}$, where all but finitely many terms in these sums
are equal to $0$.  This implies that
\begin{eqnarray}
\label{f(z) g(z) = ...}
 f(z) \, g(z) & = & \Big(\sum_{j = -\infty}^\infty \widehat{f}(j) \, z^j\Big) \,
                 \Big(\sum_{k = -\infty}^\infty \widehat{g}(k) \, z^k\Big) \\
 & = & \sum_{k = -\infty}^\infty \, \sum_{j = -\infty}^\infty \widehat{f}(j) \,
                              \widehat{g}(k) \, z^{j + k} \nonumber \\
 & = & \sum_{k = -\infty}^\infty \, \sum_{l = -\infty}^\infty \widehat{f}(l - k) \,
                              \widehat{g}(k) \, z^l \nonumber \\
 & = & \sum_{l = -\infty}^\infty \Big(\sum_{k = -\infty}^\infty \widehat{f}(l - k) \,
                               \widehat{g}(k)\Big) \, z^l \nonumber
\end{eqnarray}
for each $z \in {\bf T}$, using a change of variables in the third
step.  It follows that
\begin{equation}
\label{widehat{(f g)}(l) = sum_k widehat{f}(l - k) widehat{g}(k)}
        \widehat{(f \, g)}(l) = \sum_{k = -\infty}^\infty \widehat{f}(l - k) \,
                                                        \widehat{g}(k)
\end{equation}
for each $l \in {\bf Z}$ in this case.

        Observe that
\begin{eqnarray}
\label{sum_{k = -infty}^infty |widehat{f}(l - k)| |widehat{g}(k)| le ...}
 \quad  \sum_{k = -\infty}^\infty |\widehat{f}(l - k)| \, |\widehat{g}(k)|
 & \le & \Big(\sum_{k = -\infty}^\infty |\widehat{f}(l - k)|^2\Big)^{1/2} \,
          \Big(\sum_{k = -\infty}^\infty |\widehat{g}(k)|^2\Big)^{1/2} \\
 & = & \Big(\sum_{j = -\infty}^\infty |\widehat{f}(j)|^2\Big)^{1/2} \,
        \Big(\sum_{k = -\infty}^\infty |\widehat{g}(k)|^2\Big)^{1/2} \nonumber
\end{eqnarray}
for each $l \in {\bf Z}$, by the Cauchy--Schwarz inequality.  This
implies that
\begin{equation}
\label{sum_{k = -infty}^infty |widehat{f}(l - k)| |widehat{g}(k)| le ..., 2}
 \qquad  \sum_{k = -\infty}^\infty |\widehat{f}(l - k)| \, |\widehat{g}(k)|
 \le \Big(\frac{1}{2 \, \pi} \, \int_{\bf T} |f(z)|^2 \, |dz|\Big)^{1/2} \,
      \Big(\frac{1}{2 \, \pi} \, \int_{\bf T} |g(w)|^2 \, |dw|\Big)^{1/2}
\end{equation}
for every $l \in {\bf Z}$, by (\ref{sum_{j = -infty}^infty
  |widehat{f}(j)|^2 le langle f, f rangle = ...}) in Section
\ref{integrable functions}.  More precisely, this works for all $f, g
\in L^2({\bf T})$, so that the sum on the left side of (\ref{sum_{k =
    -infty}^infty |widehat{f}(l - k)| |widehat{g}(k)| le ..., 2})
converges for every $f, g \in L^2({\bf T})$.

        If $f, g \in L^2({\bf T})$, then their product $f \, g$ is an
integrable function on ${\bf T}$, and hence the corresponding Fourier
coefficient
\begin{equation}
\label{widehat{(f g)}(l) = frac{1}{2 pi} int_{T} f(z) g(z) overline{z}^l |dz|}
        \widehat{(f \, g)}(l) = \frac{1}{2 \, \pi} \, \int_{\bf T} f(z) \, g(z)
                                           \, \overline{z}^l \, |dz|
\end{equation}
is defined for every $l \in {\bf Z}$.  If $g \in \mathcal{E}({\bf
  T})$, then it is easy to get (\ref{widehat{(f g)}(l) = sum_k
  widehat{f}(l - k) widehat{g}(k)}) from (\ref{widehat{(f g)}(l) =
  frac{1}{2 pi} int_{T} f(z) g(z) overline{z}^l |dz|}), where all but
finitely many terms in the sum on the right side of (\ref{widehat{(f
    g)}(l) = sum_k widehat{f}(l - k) widehat{g}(k)}) are equal to $0$.
Using this, one can check that (\ref{widehat{(f g)}(l) = sum_k
  widehat{f}(l - k) widehat{g}(k)}) holds for every $f, g \in L^2({\bf
  T})$ and $l \in {\bf Z}$, by approximating $g$ by elements of
$\mathcal{E}({\bf T})$ with respect to the $L^2$ norm.  This also uses
(\ref{sum_{k = -infty}^infty |widehat{f}(l - k)| |widehat{g}(k)| le
  ..., 2}), which implies in particular that the sum on the right side
of (\ref{widehat{(f g)}(l) = sum_k widehat{f}(l - k) widehat{g}(k)})
converges absolutely for every $f, g \in L^2({\bf T})$ and $l \in {\bf
  Z}$.

        Suppose now that the Fourier coefficients of $g \in L^1({\bf T})$
are absolutely summable, which is to say that
\begin{equation}
\label{sum_{k = -infty}^infty |widehat{g}(k)|}
        \sum_{k = -\infty}^\infty |\widehat{g}(k)|
\end{equation}
converges.  This implies that the corresponding Fourier series
\begin{equation}
\label{sum_{k = -infty}^infty widehat{g}(k) z^k}
        \sum_{k = -\infty}^\infty \widehat{g}(k) \, z^k
\end{equation}
converges absolutely for each $z \in {\bf T}$, and that the partial sums
\begin{equation}
\label{sum_{k = -N}^N widehat{g}(k) z^k}
        \sum_{k = -N}^N \widehat{g}(k) \, z^k
\end{equation}
converge to (\ref{sum_{k = -infty}^infty widehat{g}(k) z^k}) as $N \to
\infty$ uniformly on ${\bf T}$.  In particular, (\ref{sum_{k =
    -infty}^infty widehat{g}(k) z^k}) defines a continuous function on
${\bf T}$.  As in Section \ref{abel sums, cauchy products}, the Abel
sums associated to (\ref{sum_{k = -infty}^infty widehat{g}(k) z^k})
converge to the ordinary sum for each $z \in {\bf T}$, which can also
be verified more directly in the case of absolute convergence.
Similar arguments show that the Abel sums converge uniformly on ${\bf
  T}$ in this situation, which is again simplified by absolute
convergence.  Of course, the Abel sums associated to (\ref{sum_{k =
    -infty}^infty widehat{g}(k) z^k}) should also converge to $g$ with
respect to the $L^1$ norm, as in Section \ref{L^p functions}.  This
implies that $g$ is equal to (\ref{sum_{k = -infty}^infty
  widehat{g}(k) z^k}) almost everywhere on ${\bf T}$, so that we may
as well take $g$ to be equal to (\ref{sum_{k = -infty}^infty
  widehat{g}(k) z^k}) everywhere on ${\bf T}$.

        Equivalently, one might start with a function $g$ on ${\bf T}$
given by
\begin{equation}
\label{g(z) = sum_{k = -infty}^infty b_k z^k}
        g(z) = \sum_{k = -\infty}^\infty b_k \, z^k
\end{equation}
for some complex numbers $b_k$ such that
\begin{equation}
\label{sum_{k = -infty}^infty |b_k|}
        \sum_{k = -\infty}^\infty |b_k|
\end{equation}
converges.  This ensures that the sum in (\ref{g(z) = sum_{k =
    -infty}^infty b_k z^k}) converges absolutely for every $z \in {\bf
  T}$, and that the corresponding partial sums
\begin{equation}
\label{sum_{k = -N}^N b_k z^k}
        \sum_{k = -N}^N b_k \, z^k
\end{equation}
converge to $g$ as $N \to \infty$ uniformly on ${\bf T}$.  It follows
that $g$ is a continuous function on ${\bf T}$, and that
\begin{equation}
        \widehat{g}(k) = b_k
\end{equation}
for every $k \in {\bf Z}$.

        If $f \in L^1({\bf T})$, then the Fourier coefficients of
$f$ are uniformly bounded, as in (\ref{|widehat{f}(j)| le frac{1}{2 pi} 
int_{bf T} |f(z)| |dz|}) in Section \ref{integrable functions}.
This implies that
\begin{eqnarray}
\label{sum_{k = -infty}^infty |widehat{f}(l - k)| |widehat{g}(k)| le ..., 3}
 \sum_{k = -\infty}^\infty |\widehat{f}(l - k)| \, |\widehat{g}(k)|
 & \le & \Big(\sup_{j \in {\bf Z}} |\widehat{f}(j)|\Big) \, 
 \Big(\sum_{k = -\infty}^\infty |\widehat{g}(k)|\Big) \\
 & \le & \Big(\frac{1}{2 \, \pi} \, \int_{\bf T} |f(z)| \, |dz|\Big) \,
 \Big(\sum_{k = -\infty}^\infty |\widehat{g}(k)|\Big) \nonumber
\end{eqnarray}
for each $l \in {\bf Z}$.  If $g \in \mathcal{E}({\bf T})$, then we
can still get (\ref{widehat{(f g)}(l) = sum_k widehat{f}(l - k)
  widehat{g}(k)}) from (\ref{widehat{(f g)}(l) = frac{1}{2 pi} int_{T}
  f(z) g(z) overline{z}^l |dz|}), as before, where all but finitely
many terms on the right side of (\ref{widehat{(f g)}(l) = sum_k
  widehat{f}(l - k) widehat{g}(k)}) are equal to $0$.  The same
conclusion also holds when the Fourier coefficients of $g$ are
absolutely summable, by a simple approximation argument, and using
(\ref{sum_{k = -infty}^infty |widehat{f}(l - k)| |widehat{g}(k)| le
  ..., 3}) to deal with the right side of (\ref{widehat{(f g)}(l) =
  sum_k widehat{f}(l - k) widehat{g}(k)}).  There are analogous
statements when $f$ is replaced by a complex Borel measure on ${\bf
  T}$.

        Suppose that both
\begin{equation}
\label{sum_{j = -infty}^infty |widehat{f}(j)|}
        \sum_{j = -\infty}^\infty |\widehat{f}(j)|
\end{equation}
and (\ref{sum_{k = -infty}^infty |widehat{g}(k)|}) converge, and let
us check that $f \, g$ has the analogous property.  Using
(\ref{widehat{(f g)}(l) = sum_k widehat{f}(l - k) widehat{g}(k)}), we
get that
\begin{equation}
\label{|widehat{(f g)}(l)| le ...}
 |\widehat{(f \, g)}(l)| 
        \le \sum_{k = -\infty}^\infty |\widehat{f}(l - k)| \, |\widehat{g}(k)|
\end{equation}
for each $l \in {\bf Z}$, so that
\begin{equation}
\label{sum_{l = -infty}^infty |widehat{(f g)}(l)| le ...}
 \sum_{l = -\infty}^\infty |\widehat{(f \, g)}(l)|
             \le \sum_{l = -\infty}^\infty \sum_{k = -\infty}^\infty 
                                 |\widehat{f}(l - k)| \, |\widehat{g}(k)|.
\end{equation}
Interchanging the order of summation on the right side, and then
making a change of variables, we get that
\begin{eqnarray}
\label{sum_{l = -infty}^infty |widehat{(f g)}(l)| le ..., 2}
 \sum_{l = -\infty}^\infty |\widehat{(f \, g)}(l)|
 & \le & \sum_{k = - \infty}^\infty \, \sum_{l = -\infty}^\infty |\widehat{f}(l - k)|
                                                         \, |\widehat{g}(k)| \\
 & = & \sum_{k = -\infty}^\infty \, \sum_{j = -\infty}^\infty |\widehat{f}(j)| \,
                                              |\widehat{g}(k)| \nonumber \\
 & = & \Big(\sum_{j = -\infty}^\infty |\widehat{f}(j)|\Big) \,
             \Big(\sum_{k = -\infty}^\infty |\widehat{g}(k)|\Big), \nonumber
\end{eqnarray}
as desired.  Alternatively, one can start with $f$ given on ${\bf T}$ by
\begin{equation}
\label{f(z) = sum_{j = -infty}^infty a_j z^j}
        f(z) = \sum_{j = -\infty}^\infty a_j \, z^j
\end{equation}
for some complex numbers $a_j$ such that
\begin{equation}
\label{sum_{j = -infty}^infty |a_j|}
        \sum_{j = -\infty}^\infty |a_j|
\end{equation}
converges, and similarly for $g$, as in (\ref{g(z) = sum_{k =
    -infty}^infty b_k z^k}) and (\ref{sum_{k = -infty}^infty |b_k|}).
One can then multiply the expansions for $f$ and $g$ directly, as in
(\ref{f(z) g(z) = ...}), and estimate the coefficients in the
corresponding expansion for $f \, g$ as before.

        If $\widehat{f}(-j) = \widehat{g}(-j) = 0$ for each positive
integer $j$, then (\ref{widehat{(f g)}(l) = sum_k widehat{f}(l - k) 
widehat{g}(k)}) reduces to
\begin{equation}
\label{widehat{(f g)}(l) = sum_{k = 0}^l widehat{f}(l - k) widehat{g}(k)}
 \widehat{(f \, g)}(l) = \sum_{k = 0}^l \widehat{f}(l - k) \, \widehat{g}(k)
\end{equation}
when $l \ge 0$, and to $\widehat{(f \, g)}(l) = 0$ when $l < 0$.  This
holds when $f \in L^p({\bf T})$ and $g \in L^q({\bf T})$ for some $1
\le p, q \le \infty$ which are conjugate exponents, in the sense that
$1/p + 1/q = 1$, so that $f \, g$ is an integrable function on ${\bf
  T}$.  More precisely, if $f$ or $g$ is in $\mathcal{E}({\bf T})$,
then (\ref{widehat{(f g)}(l) = sum_k widehat{f}(l - k) widehat{g}(k)})
follows from (\ref{widehat{(f g)}(l) = frac{1}{2 pi} int_{T} f(z) g(z)
  overline{z}^l |dz|}), which implies (\ref{widehat{(f g)}(l) = sum_{k
    = 0}^l widehat{f}(l - k) widehat{g}(k)}) in this case.  If $1 \le
q < \infty$, then we can get the same comclusion when $f \in L^p({\bf
  T})$ and $g \in L^q({\bf T})$ by approximating $g$ by finite linear
combinations of $z^k$'s with $k \ge 0$ with respect to the $L^q$ norm,
as in the previous section.  Otherwise, if $q = \infty$, then $p = 1$,
and one can approximate $f$ by finite linear combinations of $z^j$'s
with $j \ge 0$ with respect to the $L^1$ norm.  In this situation, the
Fourier series for $f \, g$ may be considered as the Cauchy product of
the Fourier series for $f$ and $g$.  This implies that
\begin{equation}
\label{A_r(f g)(z) = A_r(f)(z) A_r(g)(z)}
        A_r(f \, g)(z) = A_r(f)(z) \, A_r(g)(z)
\end{equation}
for every $z \in {\bf T}$ and $0 \le r < 1$, as in Section \ref{abel
  sums, cauchy products}.  Similarly, the power series
\begin{equation}
\label{sum_{l = 0}^infty widehat{(f g)}(l) zeta^l}
        \sum_{l = 0}^\infty \widehat{(f \, g)}(l) \, \zeta^l
\end{equation}
corresponding to $f \, g$ as in the previous section is the same as
the Cauchy product of the power series
\begin{equation}
\label{sum_{j = 0}^infty widehat{f}(j) zeta^j, ...}
        \sum_{j = 0}^\infty \widehat{f}(j) \, \zeta^j, \quad
         \sum_{k = 0}^\infty \widehat{g}(k) \, \zeta^k
\end{equation}
corresponding to $f$ and $g$.  Thus the holomorphic function on $D$
associated to $f \, g$ is the same as the product of the holomorphic
functions on $D$ associated to $f$ and $g$ under these conditions.

\section{Representation theorems}
\label{representation theorems}

        Let $g(\zeta)$ be a complex-valued harmonic function on the
open unit disk $D$, and put $g_t(\zeta) = g(t \, \zeta)$ for $0 \le t
< 1$.  Also put
\begin{equation}
\label{f_t(z) = g_t(z)}
        f_t(z) = g_t(z)
\end{equation}
when $z \in {\bf T}$, and let $h_t$ be the harmonic function on $D$
corresponding to $f_t$, as in Section \ref{harmonic functions}.  Thus
\begin{equation}
\label{g_t(zeta) = h_t(zeta)}
        g_t(\zeta) = h_t(\zeta)
\end{equation}
for every $\zeta \in D$, by the maximum principle, as before.  Suppose
now that the $L^p$ norm of $f_t$ on ${\bf T}$ is uniformly bounded in
$t$ for some $p$, $1 < p \le \infty$.  Let $1 \le q < \infty$ be the
exponent conjugate to $p$, so that $1/p + 1/q = 1$, and remember that
$L^p({\bf T})$ can be identified with the dual of $L^q({\bf T})$ in
the usual way.  The Banach--Alaoglu theorem implies that the closed
unit ball in $L^p({\bf T})$ is compact with respect to the
corresponding weak$^*$ topology on $L^p({\bf T})$.  The separability
of $L^q({\bf T})$ implies that the topology induced on the closed unit
ball in $L^p({\bf T})$ by the weak$^*$ topology is metrizable, and
hence that the closed unit ball in $L^p({\bf T})$ is sequentially
compact.  Similarly, every closed ball in $L^p({\bf T})$ is
sequentially compact with respect to the weak$^*$ topology.  It
follows that there is a sequence $\{t_j\}_{j = 1}^\infty$ of
nonnegative real numbers less than $1$ that converges to $1$ such that
$\{f_{t_j}\}_{j = 1}^\infty$ converges to some $f \in L^p({\bf T})$
with respect to the weak$^*$ topology under these conditions.  

        If $h(\zeta)$ is the harmonic function on $D$ that corresponds 
to $f$ as before, then
\begin{equation}
\label{lim_{j to infty} h_{t_j}(zeta) = h(zeta)}
        \lim_{j \to \infty} h_{t_j}(\zeta) = h(\zeta)
\end{equation}
for every $\zeta \in D$.  This uses the fact that $h_t(\zeta)$ is
equal to the integral of $f_t(z)$ times a continuous function on ${\bf
  T}$ that depends only on $\zeta$, so that (\ref{lim_{j to infty}
  h_{t_j}(zeta) = h(zeta)}) corresponds exactly to the convergence of
$\{f_{t_j}\}_{j = 1}^\infty$ to $f$ with respect to the weak$^*$
topology on $L^p({\bf T})$.  By construction,
\begin{equation}
\label{h_{t_j}(zeta) = g(t_j zeta) to g(zeta) as j to infty}
 h_{t_j}(\zeta) = g(t_j \, \zeta) \to g(\zeta) \quad\hbox{as } j \to \infty
\end{equation}
for each $\zeta \in D$, so that
\begin{equation}
\label{g(zeta) = h(zeta)}
        g(\zeta) = h(\zeta)
\end{equation}
for every $\zeta \in D$.  If $1 < p < \infty$, then we may conclude
that $f_t \to f$ as $t \to 1-$ with respect to the $L^p$ norm on ${\bf
  T}$, as in Section \ref{L^p functions}.  Some variants of this for
$p = \infty$ were also discussed in Section \ref{L^p functions}.

        In the $p = 1$ case, it is better to consider the Borel measure
\begin{equation}
\label{mu_t(E) = frac{1}{2 pi} int_E f_t(z) |dz|}
        \mu_t(E) = \frac{1}{2 \, \pi} \, \int_E f_t(z) \, |dz|
\end{equation}
on ${\bf T}$ corresponding to $f_t$ for each $t$, and the associated
bounded linear functional
\begin{equation}
\label{lambda_t(phi) = int_{bf T} phi(z) d mu_t(z) =  ...}
        \lambda_t(\phi) = \int_{\bf T} \phi(z) \, d\mu_t(z)
                 = \frac{1}{2 \, \pi} \, \int_{\bf T} \phi(z) \, f_t(z) \, |dz|
\end{equation}
on $C({\bf T})$.  The hypothesis that the $L^1$ norm of $f_t$ be
uniformly bounded in $t$ says exactly that the total variation norm of
$\mu_t$ is uniformly bounded, and hence that the dual norm of
$\lambda_t$ with respect to the supremum norm on $C({\bf T})$ is
uniformly bounded.  The Banach--Alaoglu theorem implies that closed
balls in the dual of $C({\bf T})$ are compact with respect to the
weak$^*$ topology again, and in fact they are sequentially compact,
because $C({\bf T})$ is separable.  As before, it follows that there
is a sequence $\{t_j\}_{j = 1}^\infty$ of real numbers less than $1$
that converges to $1$ for which the corresponding sequence
$\{\lambda_{t_j}\}_{j = 1}^\infty$ converges to a bounded linear
functional $\lambda$ on $C({\bf T})$ with respect to the weak$^*$
topology on the dual of $C({\bf T})$.  Remember that there is a
complex Borel measure $\mu$ on ${\bf T}$ corresponding to $\lambda$ as
in (\ref{lambda(phi) = int_{bf T} phi(z) d mu(z)}), by a version of
the Riesz representation theorem.  If $h(\zeta)$ is the harmonic
function on $D$ that corresponds to $\mu$ as in Section \ref{complex
  borel measures}, then (\ref{lim_{j to infty} h_{t_j}(zeta) =
  h(zeta)}) holds for every $\zeta \in D$ again, because
$\lambda_{t_j} \to \lambda$ as $j \to \infty$ with respect to the
weak$^*$ topology on the dual of $C({\bf T})$.  We still have
(\ref{h_{t_j}(zeta) = g(t_j zeta) to g(zeta) as j to infty}) for each
$\zeta \in D$, by construction, and hence that (\ref{g(zeta) =
  h(zeta)}) holds for every $\zeta \in D$.  It follows that $\lambda_t
\to \lambda$ as $t \to 1-$ with respect to the weak$^*$ topology on
the dual of $C({\bf T})$, as in (\ref{lim_{r to 1-} lambda_r(phi) =
  lambda(phi)}).

        If $g(\zeta)$ is a holomorphic function on $D$ such that
the $L^1$ norm of $f_t$ on ${\bf T}$ is uniformly bounded in $t$, then
a famous theorem states that $f_t$ converges to some $f \in L^1({\bf
  T})$ with respect to the $L^1$ norm as $t \to 1-$.

\chapter{The real line}
\label{real line}

\section{Fourier transforms}
\label{fourier transforms}

        The \emph{Fourier transform}\index{Fourier transforms} of an
integrable function $f$ on the real line is defined by
\begin{equation}
\label{widehat{f}(xi) = int_{bf R} f(x) exp(- i x xi) dx}
        \widehat{f}(\xi) = \int_{\bf R} f(x) \, \exp(- i \, x \, \xi) \, dx,
\end{equation}
for every $\xi \in {\bf R}$, where $dx$ refers to ordinary Lebesgue
measure on ${\bf R}$.  Similarly, the Fourier transform of a complex
Borel measure $\mu$ on ${\bf R}$ is defined by
\begin{equation}
\label{widehat{mu}(xi) = int_{bf R} exp(- i x xi) d mu(x)}
        \widehat{\mu}(\xi) = \int_{\bf R} \exp(- i \, x \, \xi) \, d\mu(x)
\end{equation}
for every $\xi \in {\bf R}$.  If
\begin{equation}
\label{mu(E) = int_E f(x) dx}
        \mu(E) = \int_E f(x) \, dx
\end{equation}
for some $f \in L^1({\bf R})$ and every Borel set $E \subseteq {\bf
  R}$, then (\ref{widehat{mu}(xi) = int_{bf R} exp(- i x xi) d mu(x)})
is the same as (\ref{widehat{f}(xi) = int_{bf R} f(x) exp(- i x xi)
  dx}).  In the first case (\ref{widehat{f}(xi) = int_{bf R} f(x)
  exp(- i x xi) dx}), we have that
\begin{equation}
\label{|widehat{f}(xi)| le int_{bf R} |f(x)| dx}
        |\widehat{f}(\xi)| \le \int_{\bf R} |f(x)| \, dx
\end{equation}
for each $\xi \in {\bf R}$, since $|\exp(i \, t)| = 1$ for every $t \in
{\bf R}$.  In the second case (\ref{widehat{mu}(xi) = int_{bf R} exp(-
  i x xi) d mu(x)}), we have that
\begin{equation}
\label{|widehat{mu}(xi)| le |mu|({bf R})}
        |\widehat{\mu}(\xi)| \le |\mu|({\bf R}),
\end{equation}
where $|\mu|$ is the total variation measure on ${\bf R}$ associated
to $\mu$.

        Let $\mu$ be a complex Borel measure on ${\bf R}$, and
consider
\begin{equation}
\label{widehat{mu}_N(xi) = int_{[-N, N]} exp(- i x xi) d mu(x)}
        \widehat{\mu}_N(\xi) = \int_{[-N, N]} \exp(- i \, x \, \xi) \, d\mu(x)
\end{equation}
for each $N \ge 0$ and $\xi \in {\bf R}$.  Observe that
\begin{equation}
\label{|widehat{mu}_N(xi) - widehat{mu}(xi)| le |mu|({x in {bf R} : |x| > N})}
        |\widehat{\mu}_N(\xi) - \widehat{\mu}(\xi)| 
                  \le |\mu|(\{x \in {\bf R} : |x| > N\})
\end{equation}
for every $N \ge 0$ and $\xi \in {\bf R}$, which implies that
$\widehat{\mu}_N(\xi) \to \widehat{\mu}(\xi)$ as $N \to \infty$
uniformly on ${\bf R}$.  One can check that $\widehat{\mu}_N(\xi)$ is
uniformly continuous on ${\bf R}$ for each $N \ge 0$, using the fact
that $\exp(i \, t)$ is uniformly continuous on ${\bf R}$.  This implies
that $\widehat{\mu}(\xi)$ is also uniformly continuous on ${\bf R}$,
since uniform continuity is preserved by uniform convergence.  In
particular, $\widehat{f}(\xi)$ is uniformly continuous on ${\bf R}$
for every $f \in L^1({\bf R})$.

        If $f$ is the characteristic or indicator function associated
to an interval in ${\bf R}$, then it is easy to see that
\begin{equation}
\label{lim_{|xi| to infty} widehat{f}(xi) = 0}
        \lim_{|\xi| \to \infty} \widehat{f}(\xi) = 0,
\end{equation}
by direct computation.  The same conclusion holds when $f$ is a step
function, which is to say a finite linear combination of
characteristic functions of intervals, by linearity.  It is well known
that step functions are dense in $L^1({\bf R})$, which implies that
(\ref{lim_{|xi| to infty} widehat{f}(xi) = 0}) holds for every $f \in
L^1({\bf R})$, using (\ref{|widehat{f}(xi)| le int_{bf R} |f(x)| dx})
to estimate the errors.  Now let $a \in {\bf R}$ be given, and let
$\delta_a$ be the Dirac mass at $a$, which is the Borel measure on
${\bf R}$ defined by $\delta_a(E) = 1$ when $a \in E$, and
$\delta_a(E) = 0$ when $a \in {\bf R} \backslash E$.  In this case, it
is easy to see that
\begin{equation}
\label{widehat{delta_a}(xi) = exp(- i a xi)}
        \widehat{\delta_a}(\xi) = \exp(- i \, a \, \xi)
\end{equation}
for every $\xi \in {\bf R}$, which obviously does not tend to $0$ as
$|\xi| \to \infty$.

        If $f \in L^1({\bf R})$ and $a \in {\bf R}$, then
\begin{equation}
\label{f_a(x) = f(x - a)}
        f_a(x) = f(x - a),
\end{equation}
is also an integrable function on ${\bf R}$, and
\begin{eqnarray}
\label{widehat{f_a}(xi) = ... = widehat{f}(xi) exp(- i a xi)}
 \widehat{f_a}(\xi) & = & \int_{\bf R} f(x - a) \, \exp(- i \, x \, \xi) \, dx \\
 & = & \int_{\bf R} f(x) \, \exp(- i \, (x + a) \, \xi) \, dx \nonumber \\
 & = & \widehat{f}(\xi) \, \exp(- i \, a \, \xi)        \nonumber
\end{eqnarray}
for every $\xi \in {\bf R}$.  Similarly,
\begin{equation}
\label{f^alpha(x) = f(x) exp(- i x alpha)}
        f^\alpha(x) = f(x) \, \exp(- i \, x \, \alpha)
\end{equation}
is an integrable function on ${\bf R}$ for every $\alpha \in {\bf R}$,
and
\begin{equation}
\label{widehat{f^alpha}(xi) = ... = widehat{f}(xi + alpha)}
 \widehat{f^\alpha}(\xi) = \int_{\bf R} f(x) \, \exp(- i \, x \, \alpha) \,
                                               \exp(- i \, x \, \xi) \, dx
                        = \widehat{f}(\xi + \alpha)
\end{equation}
for every $\alpha, \xi \in {\bf R}$.  Of course, there are analogous
statements for complex Borel measures on ${\bf R}$ instead of
integrable functions.  Remember that
\begin{equation}
\label{lim_{a to 0} int_{bf R} |f(x) - f(x - a)| dx = 0}
        \lim_{a \to 0} \int_{\bf R} |f(x) - f(x - a)| \, dx = 0
\end{equation}
for every $f \in L^1({\bf R})$, since this holds when $f$ is a step
function on ${\bf R}$, and step functions are dense in $L^1({\bf R})$.
This implies that $\widehat{f_a}(\xi) \to \widehat{f}(\xi)$ as $a \to
0$ uniformly on ${\bf R}$, which can be used to give another proof of
(\ref{lim_{|xi| to infty} widehat{f}(xi) = 0}).

\section{Holomorphic extensions}
\label{holomorphic extensions}

        If $\mu$ is a complex Borel measure on ${\bf R}$ with compact support,
then
\begin{equation}
\label{widehat{mu}(zeta) = int_{bf R} exp(- i x zeta) d mu(x)}
        \widehat{\mu}(\zeta) = \int_{\bf R} \exp(- i \, x \, \zeta) \, d\mu(x)
\end{equation}
makes sense for every $\zeta \in {\bf C}$, and defines an entire
holomorphic function on the complex plane.  Similarly,
(\ref{widehat{mu}(zeta) = int_{bf R} exp(- i x zeta) d mu(x)}) makes sense
when $\mu$ is a complex Borel measure on ${\bf R}$ such that
\begin{equation}
\label{int_{bf R} |exp(- i x zeta)| d|mu|(x) = ...}
        \int_{\bf R} |\exp(- i \, x \, \zeta)| \, d|\mu|(x)
              = \int_{\bf R} \exp(x \, \im \zeta) \, d|\mu|(x)
\end{equation}
is finite, where $\im \zeta$ denotes the imaginary part of $\zeta$.
In this case, the modulus of (\ref{widehat{mu}(zeta) = int_{bf R}
  exp(- i x zeta) d mu(x)}) is less than or equal to (\ref{int_{bf R}
  |exp(- i x zeta)| d|mu|(x) = ...}).  Of course, (\ref{int_{bf R}
  |exp(- i x zeta)| d|mu|(x) = ...}) is finite when $\zeta$ is real,
so that (\ref{widehat{mu}(zeta) = int_{bf R} exp(- i x zeta) d mu(x)})
reduces to the ordinary Fourier transform of $\mu$.

        Now suppose that
\begin{equation}
\label{|mu|((0, +infty)) = 0}
        |\mu|((0, +\infty)) = 0,
\end{equation}
which is the same as saying that $\mu \equiv 0$ on the set of positive
real numbers.  Let
\begin{equation}
\label{U = {zeta in {bf C} : im zeta > 0}}
        U = \{\zeta \in {\bf C} : \im \zeta > 0\}
\end{equation}
be the open upper half-plane, whose closure
\begin{equation}
\label{overline{U} = {zeta in {bf C} : im zeta ge 0}}
        \overline{U} = \{\zeta \in {\bf C} : \im \zeta \ge 0\}
\end{equation}
in ${\bf C}$ is the closed upper half-plane.  Under these conditions,
(\ref{int_{bf R} |exp(- i x zeta)| d|mu|(x) = ...}) is equal to
\begin{equation}
\label{int_{(-infty, 0]} exp(x im zeta) d|mu|(x)}
        \int_{(-\infty, 0]} \exp(x \, \im \zeta) \, d|\mu|(x)
\end{equation}
for every $\zeta \in {\bf C}$, which is less than or equal to
\begin{equation}
\label{|mu|({bf R}) = |mu|((-infty, 0])}
        |\mu|({\bf R}) = |\mu|((-\infty, 0])
\end{equation}
when $\zeta \in \overline{U}$.  Thus (\ref{widehat{mu}(zeta) = int_{bf
    R} exp(- i x zeta) d mu(x)}) reduces to
\begin{equation}
\label{widehat{mu}(zeta) = int_{(-infty, 0]} exp(- i x zeta) d mu(x)}
 \widehat{\mu}(\zeta) = \int_{(-\infty, 0]} \exp(- i \, x \, \zeta) \, d\mu(x)
\end{equation}
in this situation, which is defined for every $\zeta \in \overline{U}$, 
and satisfies
\begin{equation}
\label{|widehat{mu}(zeta)| le |mu|((-infty, 0])}
        |\widehat{\mu}(\zeta)| \le |\mu|((-\infty, 0]).
\end{equation}
One can also check that $\widehat{\mu}(\zeta)$ is uniformly continuous
on $\overline{U}$, in essentially the same way as in the previous
section, and that $\widehat{\mu}(\zeta)$ is holomorphic on $U$.

        If $f$ is an integrable function on ${\bf R}$ such that
$f(x) = 0$ almost everywhere when $x > 0$, then the remarks in the
previous paragraph can be applied to the measure $\mu$ on ${\bf R}$
that corresponds to $f$ as in (\ref{mu(E) = int_E f(x) dx}).  In particular,
\begin{equation}
\label{widehat{f}(zeta) = int_{-infty}^0 f(x) exp(- i x zeta) dx}
 \widehat{f}(\zeta) = \int_{-\infty}^0 f(x) \, \exp(- i \, x \, \zeta) \, dx
\end{equation}
is defined for every $\zeta \in \overline{U}$, which corresponds to
(\ref{widehat{mu}(zeta) = int_{(-infty, 0]} exp(- i x zeta) d mu(x)}).
  If $f$ is the characteristic function associated to a subinterval of
  $(-\infty, 0]$, then one can check that
\begin{equation}
\label{lim_{|zeta| to infty} f(zeta) = 0}
        \lim_{|\zeta| \to \infty} f(\zeta) = 0,
\end{equation}
by direct computation, and where the limit is taken over $\zeta \in
\overline{U}$.  This is analogous to (\ref{lim_{|xi| to infty}
  widehat{f}(xi) = 0}), and it also works when $f$ is a step function
supported in $(-\infty, 0]$, by linearity.  It follows that
  (\ref{lim_{|zeta| to infty} f(zeta) = 0}) holds for every $f \in
  L^1({\bf R})$ supported in $(-\infty, 0]$, by approximation by step
    functions.

        Similarly, the Fourier transform of an integrable function or
complex Borel measure on ${\bf R}$ supported in $[0, +\infty)$ can be
extended to the closed lower half-plane, with the same type of properties
as before.  The Fourier transform of an integrable function or complex
Borel measure on ${\bf R}$ may also have a natural extension to a
strip in the plane, a half-plane with boundary parallel to the real
line, or to the whole complex plane, depending on the finiteness of
(\ref{int_{bf R} |exp(- i x zeta)| d|mu|(x) = ...}), and without such
conditions on the support.

\section{Some examples}
\label{some examples}

        Let $r$ be a positive real number, and put
\begin{eqnarray}
\label{a_r(x) = exp(r x) when x le 0, = 0  when x > 0}
        a_r(x) & = & \exp(r \, x) \quad\hbox{when } x \le 0 \\
               & = & 0    \qquad\qquad\hbox{when } x > 0. \nonumber
\end{eqnarray}
This is an integrable function on ${\bf R}$, whose Fourier transform
\begin{equation}
\label{widehat{a_r}(zeta) = ... = int_{-infty}^0 exp(r x - i x zeta) dx}
 \widehat{a_r}(\zeta) = \int_{\bf R} a_r(x) \, \exp(- i \, x \, \zeta) \, dx
                      = \int_{-\infty}^0 \exp(r \, x - i \, x \, \zeta) \, dx
\end{equation}
may be defined for all $\zeta \in {\bf C}$ such that $\im \zeta > -
r$, as in the previous section.  In this case, we can integrate the
exponential directly, to get that
\begin{equation}
\label{widehat{a_r}(zeta) = frac{1}{r - i zeta}}
        \widehat{a_r}(\zeta) = \frac{1}{r - i \, \zeta}
\end{equation}
when $\im \zeta > - r$.

         Similarly, put
\begin{eqnarray}
\label{b_r(x) = 0  when x le 0, = exp(- r x) when x > 0}
        b_r(x) & = & 0   \enspace\qquad\qquad\hbox{ when } x \le 0 \\
               & = & \exp(- r \, x) \quad\hbox{when } x > 0. \nonumber
\end{eqnarray}
This is also an integrable function on ${\bf R}$, whose Fourier transform
\begin{equation}
\label{widehat{b_r}(zeta) = ... = int_0^infty exp(- r x - i x zeta) dx}
 \widehat{b_r}(\zeta) = \int_{\bf R} b_r(x) \, \exp(- i \, x \, \zeta) \, dx
                      = \int_0^\infty \exp(- r \, x - i \, x \, \zeta) \, dx
\end{equation}
may be defined for all $\zeta \in {\bf C}$ with $\im \zeta < r$, as
before.  The exponential can be integrated directly again, to get that
\begin{equation}
\label{widehat{b_r}(zeta) = frac{1}{r + i zeta}}
        \widehat{b_r}(\zeta) = \frac{1}{r + i \, \zeta}
\end{equation}
when $\im \zeta < r$.

        Now consider
\begin{equation}
\label{c_r(x) = a_r(x) + b_r(x) = exp(- r |x|)}
        c_r(x) = a_r(x) + b_r(x) = \exp(- r \, |x|).
\end{equation}
This is an integrable function on ${\bf R}$, whose Fourier transform
\begin{equation}
\label{widehat{c_r}(zeta) = ... = int_{-infty}^infty exp(- r |x| - i x zeta) dx}
 \widehat{c_r}(\zeta) = \int_{\bf R} c_r(x) \, \exp(- i \, x \, \zeta) \, dx
              = \int_{-\infty}^\infty \exp(- r \, |x| - i \, x \, \zeta) \, dx
\end{equation}
may be defined for all $\zeta \in {\bf C}$ with $|\im \zeta| < r$.
Combining (\ref{widehat{a_r}(zeta) = frac{1}{r - i zeta}}) and
(\ref{widehat{b_r}(zeta) = frac{1}{r + i zeta}}), we get that
\begin{equation}
\label{widehat{c_r}(zeta) = widehat{a_r}(zeta) + widehat{b_r}(zeta) = ...}
 \widehat{c_r}(\zeta) = \widehat{a_r}(\zeta) + \widehat{b_r}(\zeta)
                      = \frac{1}{r - i \, \zeta} + \frac{1}{r + i \, \zeta}
                      = \frac{2 \, r}{(r - i \, \zeta) \, (r + i \, \zeta)}
\end{equation}
when $|\im \zeta| < r$.

        Let us restrict our attention to $\zeta = \xi \in {\bf R}$,
so that (\ref{widehat{c_r}(zeta) = widehat{a_r}(zeta) +
  widehat{b_r}(zeta) = ...}) reduces to
\begin{equation}
\label{widehat{c_r}(xi) = frac{2 r}{|r - i xi|^2} = frac{2 r}{r^2 + xi^2}}
        \widehat{c_r}(\xi) = \frac{2 \, r}{|r - i \, \xi|^2}
                           = \frac{2 \, r}{r^2 + \xi^2}.
\end{equation}
This is a nonnegative real-valued integrable function on ${\bf R}$ for
each $r > 0$.  Equivalently,
\begin{equation}
\label{widehat{c_r}(xi) = ... = 2 re frac{1}{r + i xi}}
 \widehat{c_r}(\xi) = \frac{1}{r - i \, \xi} + \frac{1}{r + i \, \xi}
                   = 2 \, \re \frac{1}{r + i \, \xi}
\end{equation}
for each $\xi \in {\bf R}$.  Using the principal branch of the logarithm,
we have that
\begin{equation}
\label{frac{d}{d xi} log(r + i xi) = frac{i}{r + i xi}}
        \frac{d}{d\xi} \, \log(r + i \, \xi) = \frac{i}{r + i \, \xi},
\end{equation}
and hence that
\begin{equation}
\label{int_{-N}^N frac{1}{r + i xi} d xi = - i log(r + i N) + i log(r - i N)}
        \int_{-N}^N \frac{1}{r + i \, \xi} \, d\xi 
                       = - i \, \log(r + i \, N) + i \, \log(r - i \, N)
\end{equation}
for each $N \ge 0$.  It follows that
\begin{eqnarray}
\label{int_{-N}^N widehat{c_r}(xi) d xi = ...}
        \int_{-N}^N \widehat{c_r}(\xi) \, d\xi 
                & = & 2 \, \re \int_{-N}^N \frac{1}{r + i \, \xi} \, d\xi \\
          & = & 2 \, \im \log(r + i \, N) - 2 \, \im \log(r - i \, N) \nonumber
\end{eqnarray}
for each $N \ge 0$.  Of course,
\begin{equation}
\label{lim_{N to infty} im log(r + i N) = frac{pi}{2}}
        \lim_{N \to \infty} \im \log(r + i \, N) = \frac{\pi}{2}
\end{equation}
and
\begin{equation}
\label{lim_{N to infty} im log(r - i N) = - frac{pi}{2}}
        \lim_{N \to \infty} \im \log(r - i \, N) = - \frac{\pi}{2}.
\end{equation}
This shows that
\begin{equation}
\label{int_{-infty}^infty widehat{c_r}(xi) d xi = 2 pi}
        \int_{-\infty}^\infty \widehat{c_r}(\xi) \, d\xi = 2 \, \pi,
\end{equation}
by taking the limit as $N \to \infty$ in (\ref{int_{-N}^N
  widehat{c_r}(xi) d xi = ...}).

\section{The multiplication formula}
\label{multiplication formula}

        If $\mu$ and $\nu$ are complex Borel measures on the real line, then
the \emph{multiplication formula}\index{multiplication formula} states that
\begin{equation}
\label{int_{bf R} widehat{nu}(t) d mu(t) = int_{bf R} widehat{mu}(xi) d nu(xi)}
        \int_{\bf R} \widehat{\nu}(t) \, d\mu(t) 
                        = \int_{\bf R} \widehat{\mu}(\xi) \, d\nu(\xi).
\end{equation}
This follows by plugging in the definition of the Fourier transform
and applying Fubini's theorem.  Note that both sides of (\ref{int_{bf
    R} widehat{nu}(t) d mu(t) = int_{bf R} widehat{mu}(xi) d nu(xi)})
make sense, because $\widehat{\mu}$ and $\widehat{\nu}$ are bounded
and continuous on ${\bf R}$.  Similarly,
\begin{equation}
\label{int_{bf R} widehat{nu}(t - w) d mu(t) = ...}
        \int_{\bf R} \widehat{\nu}(t - w) \, d\mu(t)
 = \int_{\bf R} \widehat{\mu}(\xi) \, \exp(i \, w \, \xi) \, d\nu(\xi)
\end{equation}
for every $w \in {\bf R}$, which is the same as (\ref{int_{bf R}
  widehat{nu}(t) d mu(t) = int_{bf R} widehat{mu}(xi) d nu(xi)})
with $\nu$ replaced by the measure $\nu_w$ on ${\bf R}$ defined by
\begin{equation}
\label{nu_w(E) = int_E exp(i w u) d nu(u)}
        \nu_w(E) = \int_E \exp(i \, w \, u) \, d\nu(u).
\end{equation}

        Let $r > 0$ be given, and let us apply the preceding remarks to
\begin{equation}
\label{nu(E) = int_E c_r(u) du = int_E exp(- r |u|) du}
        \nu(E) = \int_E c_r(u) \, du = \int_E \exp(- r \, |u|) \, du.
\end{equation}
Thus (\ref{int_{bf R} widehat{nu}(t - w) d mu(t) = ...}) implies that
\begin{equation}
\label{int_{bf R} frac{2 r}{r^2 + (t - w)^2} d mu(t) = ...}
        \int_{\bf R} \frac{2 \, r}{r^2 + (t - w)^2} \, d\mu(t)
 = \int_{\bf R} \widehat{\mu}(\xi) \, \exp(i \, w \, \xi) \, \exp(- r \, |\xi|)
                                                                  \, d\xi
\end{equation}
for every $w \in {\bf R}$, using the expression for $\widehat{c_r}$ in
(\ref{widehat{c_r}(xi) = frac{2 r}{|r - i xi|^2} = frac{2 r}{r^2 +
    xi^2}}).  If $\widehat{\mu}(\xi)$ is an integrable function on
${\bf R}$, then the right side of (\ref{int_{bf R} frac{2 r}{r^2 + (t
    - w)^2} d mu(t) = ...}) converges to
\begin{equation}
\label{int_{bf R} widehat{mu}(xi) exp(i w xi) d xi}
        \int_{\bf R} \widehat{\mu}(\xi) \, \exp(i \, w \, \xi) \, d\xi
\end{equation}
as $r \to 0$.  This can be derived from the dominated convergence
theorem applied to any sequence of $r$'s converging to $0$.  The same
type of argument can also be used more directly, since
\begin{equation}
\label{exp(- r |xi|) le 1}
        \exp(- r \, |\xi|) \le 1
\end{equation}
for every $r \ge 0$ and $\xi \in {\bf R}$, and $\exp(- r \, |\xi|) \to 1$
as $r \to 0$ uniformly on any bounded set of $\xi$'s.

        Put
\begin{equation}
\label{A_r(mu)(w) = frac{1}{pi} int_{bf R} frac{r}{r^2 + (t - w)^2} d mu(t)}
 A_r(\mu)(w) = \frac{1}{\pi} \, \int_{\bf R} \frac{r}{r^2 + (t - w)^2} \, d\mu(t)
\end{equation}
for every $w \in {\bf R}$ and $r > 0$, which is the same as the common
value of (\ref{int_{bf R} frac{2 r}{r^2 + (t - w)^2} d mu(t) = ...})
divided by $2 \, \pi$.  Thus
\begin{equation}
\label{|A_r(mu)(w)| le frac{1}{pi} int_{bf R} frac{r}{r^2 + (t - w)^2} d|mu|(t)}
 |A_r(\mu)(w)| \le \frac{1}{\pi} \, \int_{\bf R} \frac{r}{r^2 + (t - w)^2}
                                                                \, d|\mu|(t)
\end{equation}
for every $w \in {\bf R}$ and $r > 0$, and one can check that
$A_r(\mu)(w)$ is continuous in $w$.  Integrating in $w$, we get that
\begin{eqnarray}
\label{int_{bf R} |A_r(mu)(w)| dw le ...}
 \int_{\bf R} |A_r(\mu)(w)| \, dw 
 & \le & \frac{1}{\pi} \, \int_{\bf R} \int_{\bf R} \frac{r}{r^2 + (t - w)^2}
                                               \, d|\mu|(t) \, dw \\
 & = & \frac{1}{\pi} \, \int_{\bf R} \int_{\bf R} \frac{r}{r^2 + (t - w)^2}
                                               \, dw \, d|\mu|(t) \nonumber \\
 & = & \frac{1}{\pi} \, \int_{\bf R} \int_{\bf R} \frac{r}{r^2 + w^2}
                                                \, dw \, d|\mu|(t) \nonumber
\end{eqnarray}
for every $r > 0$, using Fubini's theorem in the second step, and then
translation-invariance of Lebesgue measure on ${\bf R}$.  We also know
that
\begin{equation}
\label{frac{1}{pi} int_{bf R} frac{r}{r^2 + w^2} dw = 1}
        \frac{1}{\pi} \, \int_{\bf R} \frac{r}{r^2 + w^2} \, dw = 1,
\end{equation}
by (\ref{int_{-infty}^infty widehat{c_r}(xi) d xi = 2 pi}), so that
(\ref{int_{bf R} |A_r(mu)(w)| dw le ...}) reduces to
\begin{equation}
\label{int_{bf R} |A_r(mu)(w)| dw le |mu|({bf R})}
        \int_{\bf R} |A_r(\mu)(w)| \, dw \le |\mu|({\bf R})
\end{equation}
for every $r > 0$.

        If $f \in L^1({\bf R})$, then
\begin{equation}
\label{A_r(f)(w) = frac{1}{pi} int_{bf R} f(t) frac{r}{r^2 + (t - w)^2} dt}
A_r(f)(w) = \frac{1}{\pi} \, \int_{\bf R} f(t) \, \frac{r}{r^2 + (t - w)^2} \, dt
\end{equation}
is the same as (\ref{A_r(mu)(w) = frac{1}{pi} int_{bf R} frac{r}{r^2 +
    (t - w)^2} d mu(t)}), where $\mu$ corresponds to $f$ as in
(\ref{mu(E) = int_E f(x) dx}).  This integral also makes sense for any
locally-integrable function $f$ on ${\bf R}$ such that
\begin{equation}
\label{|f(t)| frac{1}{1 + t^2}}
        |f(t)| \, \frac{1}{1 + t^2}
\end{equation}
is integrable on ${\bf R}$.  In particular, this holds when $f \in
L^p({\bf R})$ for some $p \ge 1$, by H\"older's inequality.  As before,
\begin{equation}
\label{|A_r(f)(w)| le frac{1}{pi} int_{bf R} |f(t)| frac{r}{r^2 + (t - w)^2} dt}
 |A_r(f)(w)| \le \frac{1}{\pi} \, \int_{\bf R} |f(t)| \,
                                               \frac{r}{r^2 + (t - w)^2} \, dt
\end{equation}
for every $w \in {\bf R}$ and $r > 0$, and one can check that
$A_r(f)(w)$ is continuous in $w$ under these conditions.  Note that
\begin{equation}
\label{frac{1}{pi} int_{bf R} frac{r}{r^2 + (t - w)^2} dt = 1}
        \frac{1}{\pi} \, \int_{\bf R} \frac{r}{r^2 + (t - w)^2} \, dt = 1
\end{equation}
for every $w \in {\bf R}$ and $r > 0$, as in (\ref{int_{-infty}^infty
  widehat{c_r}(xi) d xi = 2 pi}) and (\ref{frac{1}{pi} int_{bf R}
  frac{r}{r^2 + w^2} dw = 1}).  If $f \in L^p({\bf R})$ for some $p$,
$1 \le p < \infty$, then it follows that
\begin{equation}
\label{|A_r(f)(w)|^p le ...}
 |A_r(f)(w)|^p \le \frac{1}{\pi} \, \int_{\bf R} |f(t)|^p \,
                                              \frac{r}{r^2 + (t - w)^2} \, dt
\end{equation}
for every $w \in {\bf R}$ and $r > 0$, by Jensen's inequality.
Integrating this in $w$ and using Fubini's theorem as in (\ref{int_{bf
    R} |A_r(mu)(w)| dw le ...}), we get that
\begin{equation}
\label{int_{bf R} |A_r(f)(w)|^p dw le int_{bf R} |f(t)|^p dt}
        \int_{\bf R} |A_r(f)(w)|^p \, dw \le \int_{\bf R} |f(t)|^p \, dt
\end{equation}
for every $r > 0$.  This shows that $A_r(f) \in L^p({\bf R})$ for
every $r > 0$ when $f \in L^p({\bf R})$ for some $p$, $1 \le p <
\infty$, and that
\begin{equation}
\label{||A_r(f)||_p le ||f||_p}
        \|A_r(f)\|_p \le \|f\|_p
\end{equation}
for every $r > 0$, where $\|f\|_p$ is the usual $L^p$ norm of $f$.  It
is easy to see that this also holds when $p = \infty$, by
(\ref{|A_r(f)(w)| le frac{1}{pi} int_{bf R} |f(t)| frac{r}{r^2 + (t -
    w)^2} dt}) and (\ref{frac{1}{pi} int_{bf R} frac{r}{r^2 + (t -
    w)^2} dt = 1}).

\section{Convergence}
\label{convergence}

        Let $f$ be a locally integrable function on ${\bf R}$ such that
(\ref{|f(t)| frac{1}{1 + t^2}}) is integrable on ${\bf R}$, and let
$A_r(f)(w)$ be as in (\ref{A_r(f)(w) = frac{1}{pi} int_{bf R} f(t) 
frac{r}{r^2 + (t - w)^2} dt}).  This is basically an average of the
values of $f$, because of (\ref{frac{1}{pi} int_{bf R} frac{r}{r^2 + 
(t - w)^2} dt = 1}).  As $r \to 0$, this average becomes more and more
concentrated near $w$.  In particular, if $f$ is continuous at $w$, then
one can check that
\begin{equation}
\label{lim_{r to 0} A_r(f)(w) = f(w)}
        \lim_{r \to 0} A_r(f)(w) = f(w).
\end{equation}
Similarly, if $f$ is bounded and uniformly continuous on ${\bf R}$,
then $A_r(f) \to f$ as $r \to 0$ uniformly on ${\bf R}$.  If $f$ is a
continuous function on ${\bf R}$ with compact support, then $f$ is
bounded and uniformly continuous on ${\bf R}$, and hence $A_r(f) \to
f$ as $r \to 0$ uniformly on ${\bf R}$.  In this case, $|A_r(f)(w)|$
is bounded by a constant times
\begin{equation}
\label{frac{r}{r^2 + w^2}}
        \frac{r}{r^2 + w^2}
\end{equation}
when $|w|$ is large, so that $A_r(f)(w) \to 0$ as $|w| \to \infty$.
Using this estimate when $|w|$ is large and uniform convergence, it is
easy to see that
\begin{equation}
\label{lim_{r to 0} ||A_r(f) - f||_p = 0}
        \lim_{r \to 0} \|A_r(f) - f\|_p = 0
\end{equation}
for every $p \ge 1$.  This also holds for every $f \in L^p({\bf R})$
when $1 \le p < \infty$, by approximating $f$ by continuous functions
on ${\bf R}$ with compact support with respect to the $L^p$ norm.
This uses the uniform bound (\ref{||A_r(f)||_p le ||f||_p}) on the
$L^p$ norms as well, to estimate the errors in the approximation.

        Let $f$ be a locally integrable function on ${\bf R}$ such that 
(\ref{|f(t)| frac{1}{1 + t^2}}) is integrable on ${\bf R}$ again,
let $I$ be a bounded interval in ${\bf R}$, and let $I_1$ be a bounded
interval in ${\bf R}$ that contains $I$ in its integrior.  If $g$ is
the function on ${\bf R}$ equal to $f$ on $I_1$ and to $0$ on ${\bf R}
\backslash I_1$, then $g \in L^1({\bf R})$, and hence $A_r(g) \to g$
as $r \to 0$ with respect to the $L^1$ norm, as in the previous
paragraph.  One can also check that $A_r(f - g) \to 0$ as $r \to 0$
uniformly on $I$, because $f - g \equiv 0$ on $I_1$ and $I$ is
contained in the interior of $I_1$.  Combining these two statements,
we get that $A_r(f) \to f$ as $r \to 0$ with respect to the $L^1$ norm
on $I$.  If $|f(t)|^p$ is locally integrable on ${\bf R}$ for some
$p$, $1 \le p < \infty$, and if (\ref{|f(t)| frac{1}{1 + t^2}}) is
integrable on ${\bf R}$, then an analogous argument shows that
$A_r(f) \to f$ as $r \to 0$ with respect to the $L^p$ norm on bounded
intervals on ${\bf R}$.  Similarly, if $f$ is a continuous function on
${\bf R}$ such that (\ref{|f(t)| frac{1}{1 + t^2}}) is integrable on ${\bf R}$,
then $A_r(f) \to f$ as $r \to 0$ uniformly on compact subsets of ${\bf R}$.
If $f$ is a locally integrable function on ${\bf R}$ such that
(\ref{|f(t)| frac{1}{1 + t^2}}) is integrable on ${\bf R}$, then it is
well known that (\ref{lim_{r to 0} A_r(f)(w) = f(w)}) holds for almost
every $w \in {\bf R}$ with respect to Lebesgue measure, but we shall go
into this here.  Note that it suffices to show this for $f \in L^1({\bf R})$,
by the same type of localization argument.

        Now let $\mu$ be a complex Borel measure on ${\bf R}$, and let
$C_0({\bf R})$ be the space of continuous complex-valued functions on
${\bf R}$ that vanish at infinity, equipped with the supremum norm.
Thus
\begin{equation}
\label{lambda(phi) = int_{bf R} phi(t) d mu(t)}
        \lambda(\phi) = \int_{\bf R} \phi(t) \, d\mu(t)
\end{equation}
defines a bounded linear functional on $C_0({\bf R})$, and it is well
known that every bounded linear functional on $C_0({\bf R})$ is of
this form, by a version of the Riesz representation theorem.  Let
$A_r(\mu)$ be as in (\ref{A_r(mu)(w) = frac{1}{pi} int_{bf R}
  frac{r}{r^2 + (t - w)^2} d mu(t)}) for each $r > 0$, which is an
integrable function on ${\bf R}$ for each $r > 0$, by (\ref{int_{bf R}
  |A_r(mu)(w)| dw le |mu|({bf R})}).  Every integrable function on
${\bf R}$ determines a Borel measure on ${\bf R}$ as in (\ref{mu(E) =
  int_E f(x) dx}), and
\begin{equation}
\label{lambda_r(phi) = int_{bf R} A_r(mu)(w) phi(w) dw}
        \lambda_r(\phi) = \int_{\bf R} A_r(\mu)(w) \, \phi(w) \, dw
\end{equation}
is the bounded linear functional on $C_0({\bf R})$ corresponding to
the Borel measure on ${\bf R}$ associated to $A_r(\mu)$ in this way.
Observe that
\begin{equation}
\label{lambda_r(phi) = int_{bf R} A_r(phi)(t) d mu(t)}
        \lambda_r(\phi) = \int_{\bf R} A_r(\phi)(t) \, d\mu(t)
\end{equation}
for every $\phi \in C_0({\bf R})$ and $r > 0$, where $A_r(\phi)$ is as
in (\ref{A_r(f)(w) = frac{1}{pi} int_{bf R} f(t) frac{r}{r^2 + (t -
    w)^2} dt}), by Fubini's theorem.  As before, $A_r(\phi) \to \phi$
as $r \to 0$ uniformly on ${\bf R}$ for every $\phi \in C_0({\bf R})$,
because elements of $C_0({\bf R})$ are bounded and uniformly continuous
on ${\bf R}$.  It follows that
\begin{equation}
\label{lim_{r to 0} lambda_r(phi) = lambda(phi)}
        \lim_{r \to 0} \lambda_r(\phi) = \lambda(\phi)
\end{equation}
for every $\phi \in C_0({\bf R})$, which means that $\lambda_r \to
\lambda$ as $r \to 0$ with respect to the weak$^*$ topology on the
dual of $C_0({\bf R})$.

        If $f \in L^\infty({\bf R})$, then $A_r(f)$ does not normally
converge to $f$ as $r \to 0$ with respect to the $L^\infty$ norm.
However, we do have that the $L^\infty$ norm of $A_r(f)$ remains
bounded, as in (\ref{||A_r(f)||_p le ||f||_p}).  We also have that
$A_r(f) \to f$ as $r \to 0$ with respect to the $L^p$ norm on bounded
intervals and pointwise almost everywhere on ${\bf R}$, as mentioned
earlier.  Alternatively, we can identify $L^\infty({\bf R})$ with the
dual of $L^1({\bf R})$ in the usual way, and consider convergence with
respect to the corresponding weak$^*$ topology on $L^\infty({\bf R})$.
In order to show that $A_r(f) \to f$ as $r \to 0$ with respect to the
weak$^*$ topology on $L^\infty({\bf R})$ as the dual of $L^1({\bf R})$,
it suffices to show that
\begin{equation}
\label{lim_{r to 0} int_{bf R} A_r(f)(w) phi(w) dw = int_{bf R} f(w) phi(w) dw}
        \lim_{r \to 0} \int_{\bf R} A_r(f)(w) \, \phi(w) \, dw
                     = \int_{\bf R} f(w) \, \phi(w) \, dw
\end{equation}
for each $\phi \in L^1({\bf R})$.  As in the previous paragraph,
one can check that
\begin{equation}
\label{int_{bf R} A_r(f)(w) phi(w) dw = int_{bf R} f(t) A_r(phi)(t) dt}
 \int_{\bf R} A_r(f)(w) \, \phi(w) \, dw = \int_{\bf R} f(t) \, A_r(\phi)(t) \, dt
\end{equation}
for every $\phi \in L^1({\bf R})$ and $r > 0$, using Fubini's theorem.
We have already seen that $A_r(\phi) \to \phi$ as $r \to 0$ with
respect to the $L^1$ norm when $\phi \in L^1({\bf R})$, which implies
that
\begin{equation}
\label{lim_{r to 0} int_{bf R} f(t) A_r(phi)(t) dt = int_{bf R} f(t) phi(t) dt}
        \lim_{r \to 0} \int_{\bf R} f(t) \, A_r(\phi)(t) \, dt
                       = \int_{\bf R} f(t) \, \phi(t) \, dt.
\end{equation}
Combining this with (\ref{int_{bf R} A_r(f)(w) phi(w) dw = int_{bf R}
  f(t) A_r(phi)(t) dt}), we get that (\ref{lim_{r to 0} int_{bf R}
  A_r(f)(w) phi(w) dw = int_{bf R} f(w) phi(w) dw}) holds, as desired.

\section{Harmonic extensions}
\label{harmonic extensions}

        Let $\mu$ be a complex Borel measure, and let $x, y \in {\bf R}$
be given, with $y > 0$, so that $z = x + i \, y$ is an element of the
upper half-plane $U$.  Put
\begin{equation}
\label{h(z) = A_y(mu)(x) = ...}
 h(z) = A_y(\mu)(x) = \frac{1}{\pi} \, \int_{\bf R} \frac{y}{(x - t)^2 + y^2}
                                                                  \, d\mu(t),
\end{equation}
where $A_y(\mu)(x)$ is as in (\ref{A_r(mu)(w) = frac{1}{pi} int_{bf R}
  frac{r}{r^2 + (t - w)^2} d mu(t)}).  Observe that
\begin{equation}
\label{frac{1}{t - z} - frac{1}{t - overline{z}} = ...}
 \frac{1}{t - z} - \frac{1}{t - \overline{z}} = \frac{2 \, i \, y}{|z - t|^2}
                                        = \frac{2 \, i \, y}{(x - t)^2 + y^2}
\end{equation}
for each $t \in {\bf R}$, by partial fractions.  If we also put
\begin{equation}
\label{h_+(z) = frac{1}{2 pi i} int_{bf R} frac{1}{t - z} d mu(t)}
 h_+(z) = \frac{1}{2 \, \pi \, i} \, \int_{\bf R} \frac{1}{t - z} \, d\mu(t)
\end{equation}
and
\begin{equation}
\label{h_-(z) = frac{-1}{2 pi i} int_{bf R} frac{1}{t - overline{z}} d mu(t)}
 h_-(z) = \frac{-1}{2 \, \pi \, i} \, \int_{\bf R} \frac{1}{t - \overline{z}}
                                                                 \, d\mu(t),
\end{equation}
then it follows that
\begin{equation}
\label{h(z) = h_+(z) + h_-(z)}
        h(z) = h_+(z) + h_-(z).
\end{equation}
Of course, $h_+(z)$ is a holomorphic function of $z$ in $U$, and
$h_-(z)$ is a conjugate-holomorphic function of $z$ in $U$, which
implies that $h(z)$ is a harmonic function of $z$ in $U$.

        Equivalently,
\begin{equation}
\label{h(z) = A_y(mu)(x) = ..., 2}
 h(z) = A_y(\mu)(x) = \frac{1}{2 \, \pi} \, \int_{\bf R} \widehat{\mu}(\xi) \, 
                             \exp(i \, x \, \xi) \, \exp(- y \, |\xi|) \, d\xi,
\end{equation}
as in (\ref{int_{bf R} frac{2 r}{r^2 + (t - w)^2} d mu(t) = ...}).
If we put
\begin{eqnarray}
\label{h_+(z) = ..., 2}
 h_+(z) & = & \frac{1}{2 \, \pi} \, \int_0^\infty \widehat{\mu}(\xi) \,
                           \exp(i \, x \, \xi) \, \exp(-y \, \xi) \, d\xi \\
        & = & \frac{1}{2 \, \pi} \, \int_0^\infty \widehat{\mu}(\xi) \,
                                      \exp(i \, z \, \xi) \, d\xi \nonumber
\end{eqnarray}
and
\begin{eqnarray}
\label{h_-(z) = ..., 2}
 h_-(z) & = & \frac{1}{2 \, \pi} \, \int_{-\infty}^0 \widehat{\mu}(\xi) \,
                            \exp(i \, x \, \xi) \, \exp(y \, \xi) \, d\xi \\
        & = & \frac{1}{2 \, \pi} \, \int_{-\infty}^0 \widehat{\mu}(\xi) \,
                             \exp(i \, \overline{z} \, \xi) \, d\xi, \nonumber
\end{eqnarray}
then $h_+(z)$ is a holomorphic function of $z$ in $U$, $h_-(z)$ is a
conjugate-holomorphic function of $z$ in $U$, and $h_+(z) + h_-(z) =
h(z)$, so that $h(z)$ is a harmonic function of $z$ in $U$.  Using the
multiplication formula (\ref{int_{bf R} widehat{nu}(t - w) d mu(t) =
  ...}) with $w = x$, we get that
\begin{equation}
\label{h_+(z) = frac{1}{2 pi} int_{bf R} widehat{b_y}(t - x) d mu(t)}
 h_+(z) = \frac{1}{2 \, \pi} \, \int_{\bf R} \widehat{b_y}(t - x) \, d\mu(t)
\end{equation}
and
\begin{equation}
\label{h_-(z) = frac{1}{2 pi} int_{bf R} widehat{a_y}(t - x) d mu(t)}
 h_-(z) = \frac{1}{2 \, \pi} \, \int_{\bf R} \widehat{a_y}(t - x) \, d\mu(t),
\end{equation}
where $a_y$ and $b_y$ are as in Section \ref{some examples}.  It follows that
\begin{equation}
\label{h_+(z) = frac{1}{2 pi} int_{bf R} frac{1}{y + i (t - x)} d mu(t)}
        h_+(z) = \frac{1}{2 \, \pi} \, \int_{\bf R} \frac{1}{y + i \, (t - x)}
                                                                   \, d\mu(t)
\end{equation}
and
\begin{equation}
\label{h_-(z) = frac{1}{2 pi} int_{bf R} frac{1}{y - i (t - x)} d mu(t)}
        h_-(z) = \frac{1}{2 \, \pi} \, \int_{\bf R} \frac{1}{y - i \, (t - x)}
                                                                   \, d\mu(t),
\end{equation}
by substituting the expressions (\ref{widehat{b_r}(zeta) = frac{1}{r +
    i zeta}}) and (\ref{widehat{a_r}(zeta) = frac{1}{r - i zeta}}) for
the Fourier transforms of $b_y$ and $a_y$ into (\ref{h_+(z) =
  frac{1}{2 pi} int_{bf R} widehat{b_y}(t - x) d mu(t)}) and
(\ref{h_-(z) = frac{1}{2 pi} int_{bf R} widehat{a_y}(t - x) d mu(t)}),
repsectively.  It is easy to see that (\ref{h_+(z) = frac{1}{2 pi}
  int_{bf R} frac{1}{y + i (t - x)} d mu(t)}) and (\ref{h_-(z) =
  frac{1}{2 pi} int_{bf R} frac{1}{y - i (t - x)} d mu(t)}) are
exactly the same as (\ref{h_+(z) = frac{1}{2 pi i} int_{bf R}
  frac{1}{t - z} d mu(t)}) and (\ref{h_-(z) = frac{-1}{2 pi i} int_{bf
    R} frac{1}{t - overline{z}} d mu(t)}), respectively.

        Similarly, if $f$ is a locally-integrable function on ${\bf R}$
such that (\ref{|f(t)| frac{1}{1 + t^2}}) is an integrable function on
${\bf R}$, then we can put
\begin{equation}
\label{h(z) = A_y(f)(x) = ...}
 h(z) = A_y(f)(x) = \frac{1}{\pi} \, \int_{\bf R} f(t) \frac{y}{(x - t)^2 + y^2}
                                                                 \, dt,
\end{equation}
where $A_y(f)(x)$ is as in (\ref{A_r(f)(w) = frac{1}{pi} int_{bf R}
  f(t) frac{r}{r^2 + (t - w)^2} dt}).  If
\begin{equation}
\label{|f(t)| frac{1}{1 + |t|}}
        |f(t)| \, \frac{1}{1 + |t|}
\end{equation}
is also an integrable function on ${\bf R}$, then
\begin{equation}
\label{h_+(z) = frac{1}{2 pi i} int_{bf R} f(t) frac{1}{t - z} dt}
 h_+(z) = \frac{1}{2 \, \pi \, i} \, \int_{\bf R} f(t) \, \frac{1}{t - z} \, dt
\end{equation}
and
\begin{equation}
\label{h_-(z) = frac{-1}{2 pi i} int_{bf R} f(t) frac{1}{t - overline{z}} dt}
 h_-(z) = \frac{-1}{2 \, \pi \, i} \, \int_{\bf R} f(t) \,
                                              \frac{1}{t - \overline{z}} \, dt
\end{equation}
can be defined, in analogy with (\ref{h_+(z) = frac{1}{2 pi i} int_{bf
    R} frac{1}{t - z} d mu(t)}) and (\ref{h_-(z) = frac{-1}{2 pi i}
  int_{bf R} frac{1}{t - overline{z}} d mu(t)}), and satisfy
(\ref{h(z) = h_+(z) + h_-(z)}).  As before, $h_+(z)$ is a holomorphic
function of $z \in U$ under these conditions, $h_-(z)$ is a
conjugate-holomorphic function of $z$ in $U$, and hence $h(z)$ is a
harmonic function of $z$ in $U$.  The integrability of (\ref{|f(t)|
  frac{1}{1 + |t|}}) is more restrictive than the integrability of
(\ref{|f(t)| frac{1}{1 + t^2}}), and it holds when $f \in L^p({\bf
  R})$ for some $p$ with $1 \le p < \infty$, by H\"older's inequality.
One can still check that $h(z)$ is a harmonic function of $z$ in $U$
when (\ref{|f(t)| frac{1}{1 + t^2}}) is integrable, basically because
(\ref{frac{1}{t - z} - frac{1}{t - overline{z}} = ...}) is harmonic in $z$.

        If (\ref{|f(t)| frac{1}{1 + t^2}}) is integrable on ${\bf R}$
and $f$ is continuous at a point $x_0 \in {\bf R}$, then
\begin{equation}
\label{h(z) to f(x_0)}
        h(z) \to f(x_0)
\end{equation}
as $z \in U$ approaches $x_0$.  This is analogous to (\ref{lim_{r to
    0} A_r(f)(w) = f(w)}), since $h(z)$ is given by an average of
values of $f$ which is concentrated near $x_0$ as $z$ approaches
$x_0$.  If $f$ is a continuous function on ${\bf R}$ such that
(\ref{|f(t)| frac{1}{1 + t^2}}) is integrable, then it follows that
the function defined on $\overline{U}$ by taking $h$ on $U$ and $f$ on
${\bf R}$ is continuous.  If $f$ is any bounded measurable function on
${\bf R}$, then (\ref{|f(t)| frac{1}{1 + t^2}}) is integrable on ${\bf
  R}$, and
\begin{equation}
\label{|h(z)| le ||f||_infty}
        |h(z)| \le \|f\|_\infty
\end{equation}
for every $z \in U$.  This is essentially the $p = \infty$ version of
(\ref{||A_r(f)||_p le ||f||_p}), which follows from (\ref{|A_r(f)(w)|
  le frac{1}{pi} int_{bf R} |f(t)| frac{r}{r^2 + (t - w)^2} dt}) and
(\ref{frac{1}{pi} int_{bf R} frac{r}{r^2 + (t - w)^2} dt = 1}), as
mentioned earlier.

\section{Bounded harmonic functions}
\label{bounded harmonic functions}

        Suppose that $g(z)$ is a continuous complex-valued function
on the closed upper half-plane $\overline{U}$ that is harmonic on the
open upper half-plane $U$ and satisfies $g(z) \to 0$ as $|z| \to \infty$.
Under these conditions, the maximum principle implies that the maximum of
$|g(z)|$ on $\overline{U}$ is attained on ${\bf R}$.  In particular,
if $g(x) = 0$ for every $x \in {\bf R}$, then it follows that $g(z) = 0$
for every $z \in \overline{U}$.

        Suppose now that $g(z)$ is a bounded continuous complex-valued 
function on $\overline{U}$ that is harmonic on $U$ and satisfies
\begin{equation}
\label{g(x) = 0}
        g(x) = 0
\end{equation}
for every $x \in {\bf R}$.  If we put
\begin{equation}
\label{g(z) = -g(overline{z})}
        g(z) = -g(\overline{z})
\end{equation}
when $\im z < 0$, then the reflection principle imples that we get a
harmonic function on the complex plane.  This harmonic function is
also bounded and hence constant, by well-known arguments. It follows
that $g(z) = 0$ for every $z$, because of (\ref{g(x) = 0}).  This
shows that a bounded continuous function on $\overline{U}$ that is
harmonic on $U$ is uniquely determined by its restriction to ${\bf R}$.

        As another variant, let $g(z)$ be a bounded continuous function
on $\overline{U}$ that is holomorphic on $U$.  Let $\epsilon > 0$
be given, and put
\begin{equation}
\label{g_epsilon(z) = g(z) (frac{i}{epsilon z + i})}
        g_\epsilon(z) = g(z) \, \Big(\frac{i}{\epsilon \, z + i}\Big)
\end{equation}
for every $z \in \overline{U}$.  This is a bounded continuous function
on $\overline{U}$ that is holomorphic on $U$ and satisfies
$g_\epsilon(z) \to 0$ as $|z| \to \infty$.  The maximum principle
implies that
\begin{equation}
\label{|g_epsilon(z)| le ... le sup_{t in {bf R}} |g(t)|}
        |g_\epsilon(z)| \le \sup_{t \in {\bf R}} |g_\epsilon(t)| 
                        \le \sup_{t \in {\bf R}} |g(t)|
\end{equation}
for every $z \in \overline{U}$, as before.  Taking the limit as
$\epsilon \to 0$, we get that
\begin{equation}
\label{|g(z)| le sup_{t in {bf R}} |g(t)|}
        |g(z)| \le \sup_{t \in {\bf R}} |g(t)|
\end{equation}
for every $z \in \overline{U}$.

        Let $g(z)$ be a bounded continuous function on $\overline{U}$
that is harmonic on $U$ again, and let $f$ be the restriction of $g$
to ${\bf R}$.  Also let $h(z)$ be as in (\ref{h(z) = A_y(f)(x) = ...})
when $z \in U$, and put $h = f$ on ${\bf R}$.  This defines a bounded
continuous function on $\overline{U}$ that is harmonic on $U$, as in
the previous section.  It follows that $g(z) = h(z)$ for every $z \in
\overline{U}$, by the uniqueness principle mentioned earlier.

        Now let $g(z)$ be a harmonic function on $U$, and suppose that
the restriction of $g(z)$ to
\begin{equation}
\label{U_tau = {z in {bf C} : im z > tau}}
        U_\tau = \{z \in {\bf C} : \im z > \tau\}
\end{equation}
is bounded for each $\tau > 0$.  Equivalently, this means that $g(z +
i \, \tau)$ is a bounded harmonic function on $U$ for each $\tau > 0$.
Of course, $g$ is continuous on $U$, and so $g(z + i \, \tau)$ has a
bounded continuous extension to $\overline{U_\tau}$ for each $\tau >
0$.  The discussion in the previous paragraph implies that
\begin{equation}
\label{g(z + i tau) = ...}
        g(z + i \, \tau) = \frac{1}{\pi} \, \int_{\bf R} g(t + i \, \tau) \,
                                             \frac{y}{(x - t)^2 + y^2} \, dt
\end{equation}
for every $z \in U$ and $\tau > 0$.  In particular, this works when
$g(z)$ is a bounded harmonic function on $U$.

        Let a real number $p \ge 1$ be given, and suppose that $g(z)$
is a harmonic function on $U$ such that
\begin{equation}
\label{(int_{bf R} |g(x + i tau)| dx)^{1/p}}
        \Big(\int_{\bf R} |g(x + i \, \tau)| \, dx\Big)^{1/p}
\end{equation}
is finite for each $\tau > 0$, and uniformly bounded over $\tau > 0$.
If $p = \infty$, then one can replace (\ref{(int_{bf R} |g(x + i tau)|
  dx)^{1/p}}) with the supremum of $|g(x + i \, \tau)|$ over $x \in
{\bf R}$, and the resulting condition is the same as saying that
$g(z)$ is bounded on $U$.  If $p < \infty$, then one can still show
that $g(z)$ is bounded on $\overline{U_\tau}$ for each $\tau > 0$,
using the mean value property of harmonic functions.

        Suppose that $1 < p \le \infty$, so that one can identify 
$L^p({\bf R})$ with the dual of $L^q({\bf R})$ in the usual way, where
$1 \le q < \infty$ is the exponent conjugate to $p$.  As in Section
\ref{representation theorems}, closed balls in $L^p({\bf R})$ are 
sequentially compact with respect to the weak$^*$ topology on
$L^p({\bf R})$ as the dual of $L^q({\bf R})$.  This uses the
Banach--Alaoglu theorem, and the separability of $L^q({\bf R})$.  
Thus the boundedness of the $L^p$ norm of $g(x + i \, \tau)$ in $x \in
{\bf R}$ over $\tau > 0$ implies that there is a sequence $\{\tau_j\}_{j
  = 0}^\infty$ of positive real numbers converging to $0$ such that
$g(x + i \, \tau_j)$ converges to some $f(x) \in L^p({\bf R})$ as $j
\to \infty$ with respect to the weak$^*$ topology on $L^p({\bf R})$.
If $h(z)$ is defined as in (\ref{h(z) = A_y(f)(x) = ...}), then it
follows that $g(z) = h(z)$ for every $z \in U$, by taking $\tau =
\tau_j$ in (\ref{g(z + i tau) = ...}), and passing to the limit as $j
\to \infty$.

        If $p = 1$, then the boundedness of the $L^1$ norm of 
$g(x + i \, \tau)$ in $x \in {\bf R}$ over $\tau > 0$ implies that
\begin{equation}
\label{lambda_tau(phi) = int_{bf R} phi(x) g(x + i tau) dx}
        \lambda_\tau(\phi) = \int_{\bf R} \phi(x) \, g(x + i \, \tau) \, dx
\end{equation}
defines a bounded linear functional on $C_0({\bf R})$ for each $\tau >
0$, with bounded dual norm with respect to the supremum norm on
$C_0({\bf R})$.  As before, closed balls in the dual of $C_0({\bf R})$
are sequentially compact with respect to the weak$^*$ topology, by the
Banach--Alaoglu theorem and the separability of $C_0({\bf R})$.  
Hence there is a sequence $\{\tau_j\}_{j = 1}^\infty$ of positive real
numbers converging to $0$ such that $\{\lambda_{\tau_j}\}_{j =
  1}^\infty$ converges with respect to the weak$^*$ topology on the
dual of $C_0({\bf R})$ to a bounded linear functional $\lambda$ on
$C_0({\bf R})$.  A version of the Riesz representation theorem implies
that $\lambda$ corresponds to a complex Borel measure $\mu$ on ${\bf R}$
in the usual way.  If $h(z)$ is defined as in (\ref{h(z) = A_y(mu)(x) = ...}),
then it follows that $g(z) = h(z)$ for every $z \in U$, by taking
$\tau = \tau_j$ in (\ref{g(z + i tau) = ...}), and passing to the limit
as $j \to \infty$.

\section{Cauchy's theorem}
\label{cauchy's theorem}

        Let $g(z)$ be a continuous complex-valued function on the
closed upper half-plane $\overline{U}$ that is holomorphic on the
open upper half-plane $U$.  Also let $R$ be a positive real number,
and let $\Gamma(R)$ be the curve in $\overline{U}$ that goes from $-R$
to $R$ along the real line, and then from $R$ to $-R$ along the top
half of the circle centered at $0$ with radius $R$.  Using Cauchy's
theorem, we get that
\begin{equation}
\label{int_{Gamma(R)} g(z) dz = 0}
        \int_{\Gamma(R)} g(z) \, dz = 0
\end{equation}
for every $R > 0$, where the integral is a complex line integral.  If
$g$ is also integrable on ${\bf R}$ with respect to Lebesgue measure,
and if
\begin{equation}
\label{lim_{|z| to infty} |z| |g(z)| = 0}
        \lim_{|z| \to \infty} |z| \, |g(z)| = 0
\end{equation}
on $\overline{U}$, then one can take the limit as $R \to \infty$ in
(\ref{int_{Gamma(R)} g(z) dz = 0}), to get that
\begin{equation}
\label{int_{-infty}^infty g(x) dx = 0}
        \int_{-\infty}^\infty g(x) \, dx = 0.
\end{equation}
More precisely, (\ref{lim_{|z| to infty} |z| |g(z)| = 0}) ensures that
the integral of $g$ over the semi-circular arc of radius $R$ in
$\Gamma(R)$ converges to $0$ as $R \to \infty$.

        Suppose that $g$ satisfies the same conditions as before,
except that $g$ is bounded on $\overline{U}$, instead of
(\ref{lim_{|z| to infty} |z| |g(z)| = 0}).  Let $\epsilon > 0$ be given,
and observe that
\begin{equation}
\label{g(z) (frac{i}{epsilon z + i})^2}
        g(z) \, \Big(\frac{i}{\epsilon \, z + i}\Big)^2
\end{equation}
defines a continuous function on $\overline{U}$ that is holomorphic on
$U$, integrable on ${\bf R}$, and satisfies the analogue of
(\ref{lim_{|z| to infty} |z| |g(z)| = 0}).  Thus
\begin{equation}
\label{int_{-infty}^infty g(x) (frac{i}{epsilon x + i})^2 dx = 0}
 \int_{-\infty}^\infty g(x) \, \Big(\frac{i}{\epsilon \, x + i}\Big)^2 \, dx = 0
\end{equation}
for every $\epsilon > 0$, for the same reasons as in
(\ref{int_{-infty}^infty g(x) dx = 0}).  It follows that
(\ref{int_{-infty}^infty g(x) dx = 0}) holds in this case as well,
using the dominated convergence theorem to pass to the limit as
$\epsilon \to 0$.

        Now let $g(z)$ be a holomorphic function on $U$ such that
$g(x + i \, \tau)$ is an integrable function of $x$ on ${\bf R}$
for each $\tau > 0$, and
\begin{equation}
\label{int_{-infty}^infty |g(x + i tau)| dx}
        \int_{-\infty}^\infty |g(x + i \, \tau)| \, dx
\end{equation}
is uniformly bounded over $\tau > 0$.  As mentioned in the previous
section, one can show that the restriction of $g(z)$ to
$\overline{U_\tau}$ is bounded for each $\tau > 0$, where $U_\tau$ is
as in (\ref{U_tau = {z in {bf C} : im z > tau}}).  This permits us to
apply the earlier discussion to $g(z + i \, \tau)$ as a function of $z
\in \overline{U}$ for each $\tau > 0$, to get that
\begin{equation}
\label{int_{-infty}^infty g(x + i tau) dx = 0}
        \int_{-\infty}^\infty g(x + i \, \tau) \, dx = 0
\end{equation}
for every $\tau > 0$.  If there is a function $g(x) \in L^1({\bf R})$
such that
\begin{equation}
\label{lim_{tau to 0} int_{-infty}^infty |g(x + i tau) - g(x)| dx = 0}
 \lim_{\tau \to 0} \int_{-\infty}^\infty |g(x + i \, \tau) - g(x)| \, dx = 0,
\end{equation}
then it follows that (\ref{int_{-infty}^infty g(x) dx = 0}) holds.  A
famous theorem states that there is always a function $g(x)$ in
$L^1({\bf R})$ that satisfies (\ref{lim_{tau to 0} int_{-infty}^infty
  |g(x + i tau) - g(x)| dx = 0}) under these conditions.

        If $g(z)$ is as in the previous paragraph, then
\begin{equation}
\label{g(z) exp(-i z xi)}
        g(z) \, \exp(-i \, z \, \xi)
\end{equation}
has the same properties for every $\xi \in {\bf R}$ with $\xi \le 0$, because
\begin{equation}
\label{|exp(-i z xi)| = exp(y xi) le 1}
        |\exp(-i \, z \, \xi)| = \exp(y \, \xi) \le 1
\end{equation}
when $y \ge 0$.  This implies that
\begin{equation}
\label{int_{-infty}^infty g(x) exp(-i x xi) dx = 0}
         \int_{-\infty}^\infty g(x) \, \exp(-i \, x \, \xi) \, dx = 0
\end{equation}
for every $\xi \in {\bf R}$ with $\xi \le 0$, by applying
(\ref{int_{-infty}^infty g(x) dx = 0}) to (\ref{g(z) exp(-i z xi)}).
Conversely, let $\mu$ be a complex Borel measure on ${\bf R}$ such
that
\begin{equation}
\label{widehat{mu}(xi) = 0}
        \widehat{\mu}(\xi) = 0
\end{equation}
for every $\xi \in {\bf R}$ with $\xi \le 0$.  This implies that
$h_-(z) = 0$ for every $z \in U$, where $h_-(z)$ is as in (\ref{h_-(z)
  = ..., 2}), which is the same as (\ref{h_-(z) = frac{-1}{2 pi i}
  int_{bf R} frac{1}{t - overline{z}} d mu(t)}).  Thus the harmonic
function $h(z)$ on $U$ in (\ref{h(z) = A_y(mu)(x) = ...}) is equal to
the holomorphic function $h_+(z)$ on $U$ in (\ref{h_+(z) = frac{1}{2
    pi i} int_{bf R} frac{1}{t - z} d mu(t)}).  This harmonic function
also satisfies the analogue of the boundedness of the integrals
(\ref{int_{-infty}^infty |g(x + i tau)| dx}), because of (\ref{int_{bf
    R} |A_r(mu)(w)| dw le |mu|({bf R})}).  It follows that there is an
integrable function $f(x)$ on ${\bf R}$ such that $h(x + i \, \tau)$
converges to $f(x)$ as $\tau \to 0$ with respect to the $L^1$ norm on
${\bf R}$, as in the preceding paragraph.  However, $h(x + i \, \tau)$
also converges in a weak sense to $\mu$ as $\tau \to 0$, by the
remarks at the end of Section \ref{convergence}.  This means that
$\mu$ is absolutely continuous on ${\bf R}$ with respect to Lebesgue
measure and corresponds to $f$ as in (\ref{mu(E) = int_E f(x) dx}),
which is another version of the F.\ and M.~Riesz theorem.

\section{Cauchy's integral formula}
\label{cauchy's integral formula}

        Let $g(z)$ be a continuous complex-valued function on the closed
upper half-plane $\overline{U}$ that is holomorphic on the open upper
half-plane $U$ again, and let $\Gamma(R)$ be the curve in
$\overline{U}$ corresponding to a positive real number $R$ a described
at the beginning of the previous section.  If $w \in U$ and $|w| < R$,
then the Cauchy integral formula implies that
\begin{equation}
\label{frac{1}{2 pi i} int_{Gamma(R)} g(z) frac{1}{z - w} dz = g(w)}
 \frac{1}{2 \, \pi \, i} \, \int_{\Gamma(R)} g(z) \, \frac{1}{z - w} \, dz
                                                    = g(w).
\end{equation}
If we also ask that
\begin{equation}
\label{|g(x)| frac{1}{1 + |x|}}
        |g(x)| \, \frac{1}{1 + |x|}
\end{equation}
be an integrable function on ${\bf R}$, and that
\begin{equation}
\label{lim_{|z| to infty} |g(z)| = 0}
        \lim_{|z| \to \infty} |g(z)| = 0,
\end{equation}
then one can take the limit as $R \to \infty$ in (\ref{frac{1}{2 pi i}
  int_{Gamma(R)} g(z) frac{1}{z - w} dz = g(w)}), to get that
\begin{equation}
\label{frac{1}{2 pi i} int_{-infty}^infty g(x) frac{1}{x - w} dx = g(w)}
 \frac{1}{2 \, \pi \, i} \, \int_{-\infty}^\infty g(x) \, \frac{1}{x - w} \, dx 
                                                                       = g(w)
\end{equation}
for every $w \in U$.  As before, one can get the same conclusion when
$g(z)$ is bounded on $\overline{U}$ instead of (\ref{lim_{|z| to
    infty} |g(z)| = 0}), by approximating $g(z)$ by $g(z) \,
(i/(\epsilon \, z + i))$ with $\epsilon > 0$.

        Under the same conditions on $g(z)$, we have that for each
$w \in U$,
\begin{equation}
\label{g(z) frac{1}{z - overline{w}}}
        g(z) \, \frac{1}{z - \overline{w}}
\end{equation}
is a bounded continuous function on $\overline{U}$ that is holomorphic
on $U$ and integrable with respect to Lebesgue measure on ${\bf R}$.
Thus
\begin{equation}
\label{int_{-infty}^infty g(x) frac{1}{x - overline{w}} dx = 0}
        \int_{-\infty}^\infty g(x) \, \frac{1}{x - \overline{w}} \, dx = 0
\end{equation}
for every $w \in U$, by applying (\ref{int_{-infty}^infty g(x) dx =
  0}) to (\ref{g(z) frac{1}{z - overline{w}}}).  Combining this with
(\ref{frac{1}{2 pi i} int_{-infty}^infty g(x) frac{1}{x - w} dx =
  g(w)}), we get that
\begin{equation}
\label{frac{1}{2 pi i} int_{-infty}^infty g(x) (...) dx = g(w)}
        \frac{1}{2 \, \pi \, i} \, \int_{-\infty}^\infty g(x) \,
         \Big(\frac{1}{x - w} - \frac{1}{x - \overline{w}}\Big) \, dx = g(w)
\end{equation}
for every $w \in U$.  If $g$ is a bounded continuous function on
$\overline{U}$ that is harmonic on $U$, $f$ is the restriction of $g$
to ${\bf R}$, and $h$ is the function on $U$ corresponding to $f$ as
in (\ref{h(z) = A_y(f)(x) = ...}), then we have seen in Section
\ref{bounded harmonic functions} that $g = h$ on $U$.  This
representation for $g$ on $U$ is the same as (\ref{frac{1}{2 pi i}
  int_{-infty}^infty g(x) (...) dx = g(w)}), but with slightly
different notation.

        Suppose now that $g(z)$ is a holomorphic function on $U$
such that $|g(x + i \, \tau)|^p$ is an integrable function of $x$
on ${\bf R}$ for some $p \in [1, \infty)$ and every $\tau > 0$,
and that
\begin{equation}
\label{int_{-infty}^infty |g(x + i tau)|^p dx}
        \int_{-\infty}^\infty |g(x + i \, \tau)|^p \, dx
\end{equation}
is uniformly bounded over $\tau > 0$.  This implies that the
restriction of $g(z)$ to $\overline{U_\tau}$ is bounded for every
$\tau > 0$, as in Section \ref{bounded harmonic functions}, and where
$U_\tau$ is as in (\ref{U_tau = {z in {bf C} : im z > tau}}).  We also
have that
\begin{equation}
\label{|g(x + i tau)| frac{1}{1 + |x|}}
        |g(x + i \, \tau)| \, \frac{1}{1 + |x|}
\end{equation}
is integrable as a function of $x$ on ${\bf R}$ for each $\tau > 0$,
by H\"older's inequality.  This permits us to apply the previous
discussion to $g(z + i \, \tau)$ as a function of $z \in \overline{U}$
for each $\tau > 0$, to get that
\begin{equation}
\label{frac{1}{2 pi i} int_{-infty}^infty g(x + i tau) ... dx = g(w + i tau)}
 \frac{1}{2 \, \pi \, i} \, \int_{-\infty}^\infty g(x + i \, \tau) \,
                             \frac{1}{x - w} \, dx = g(w + i \, \tau)
\end{equation}
and
\begin{equation}
\label{int_{-infty}^infty g(x + i tau) frac{1}{x - overline{w}} dx = 0}
 \int_{-\infty}^\infty g(x + i \, \tau) \, \frac{1}{x - \overline{w}} \, dx = 0
\end{equation}
for every $w \in U$, and hence that
\begin{equation}
\label{frac{1}{2 pi i} int_{-infty}^infty g(x + i tau) (...) dx = g(w + i tau)}
 \frac{1}{2 \, \pi \, i} \, \int_{-\infty}^\infty g(x + i \, \tau) \,
\Big(\frac{1}{x - w} - \frac{1}{x - \overline{w}}\Big) \, dx = g(w + i \, \tau).
\end{equation}
As before, (\ref{frac{1}{2 pi i} int_{-infty}^infty g(x + i tau) (...)
  dx = g(w + i tau)}) is the same as (\ref{g(z + i tau) = ...}), but
with slightly different notation.

        If $1 < p < \infty$, then we can identify $L^p({\bf R})$
with the dual of $L^q({\bf R})$, where $1 < q < \infty$ and $1/p + 1/q
= 1$.  As in Section \ref{bounded harmonic functions}, one can use
the Banach--Alaoglu theorem and the separability of $L^q({\bf R})$
to get a sequence $\{\tau_j\}_{j = 1}^\infty$ of positive real numbers
converging to $0$ such that $g(x + i \, \tau_j)$ converges to some
$g(x) \in L^p({\bf R})$ as $j \to \infty$ with respect to the weak$^*$
topology on $L^p({\bf R})$ as the dual of $L^q({\bf R})$.  Applying
(\ref{frac{1}{2 pi i} int_{-infty}^infty g(x + i tau) ... dx = g(w + i
  tau)}), (\ref{int_{-infty}^infty g(x + i tau) frac{1}{x -
    overline{w}} dx = 0}), and (\ref{frac{1}{2 pi i}
  int_{-infty}^infty g(x + i tau) (...) dx = g(w + i tau)}) to
$\tau_j$ and taking the limit as $j \to \infty$, we get that
(\ref{frac{1}{2 pi i} int_{-infty}^infty g(x) frac{1}{x - w} dx =
  g(w)}), (\ref{int_{-infty}^infty g(x) frac{1}{x - overline{w}} dx =
  0}), and (\ref{frac{1}{2 pi i} int_{-infty}^infty g(x) (...) dx =
  g(w)}) hold for every $w \in U$.  This is analogous to the
discussion of harmonic functions in Section \ref{bounded harmonic
  functions}, and indeed (\ref{frac{1}{2 pi i} int_{-infty}^infty g(x)
  (...) dx = g(w)}) corresponds exactly to representing $g$ as the
harmonic function $h$ on $U$ associated to some $f \in L^p({\bf R})$
as in (\ref{h(z) = A_y(f)(x) = ...}).  Under these conditions, one can
also conclude that $g(x + i \, \tau)$ converges to $g(x)$ as $\tau \to
0$ with respect to the $L^p$ norm on ${\bf R}$, as in Section
\ref{convergence}.  Similarly, if $p = 1$, then there is a sequence
$\{\tau_j\}_{j = 1}^\infty$ of positive real numbers converging to $0$
such that $g(x + i \, \tau_j)$ converges to a complex Borel measure
$\mu$ on ${\bf R}$ in a suitable weak sense as $j \to \infty$, as in
Section \ref{bounded harmonic functions}.  This leads to equations like
(\ref{frac{1}{2 pi i} int_{-infty}^infty g(x) frac{1}{x - w} dx = g(w)}),
(\ref{int_{-infty}^infty g(x) frac{1}{x - overline{w}} dx = 0}), and
(\ref{frac{1}{2 pi i} int_{-infty}^infty g(x) (...) dx = g(w)}) for each
$w \in U$, but where $g(x)$ and $dx$ are replaced with $d\mu(x)$.
However, as mentioned in the previous section, a famous theorem
implies that $g(x + i \, \tau)$ actually converges to an integrable
function $g(x)$ on ${\bf R}$ with respect to the $L^1$ norm as $\tau
\to 0$ in this situation.  In particular, (\ref{frac{1}{2 pi i}
  int_{-infty}^infty g(x) frac{1}{x - w} dx = g(w)}),
(\ref{int_{-infty}^infty g(x) frac{1}{x - overline{w}} dx = 0}), and
(\ref{frac{1}{2 pi i} int_{-infty}^infty g(x) (...) dx = g(w)}) still
hold, using this function $g(x)$.

\section{Some variants}
\label{some variants}

        Let $g(z)$ be a continuous complex-valued function on the
closed upper half-plane $\overline{U}$ that is holomorphic on the open
upper half-plane $U$, and let $w, w_1 \in U$ be given.  Thus
\begin{equation}
\label{g(z) frac{1}{z - overline{w_1}}}
        g(z) \, \frac{1}{z - \overline{w_1}}
\end{equation}
is also a continuous function of $z$ on $\overline{U}$ that is
holomorphic on $U$.  If $\Gamma(R)$ is the curve in $\overline{U}$
corresponding to a positive real number $R$ as in Section
\ref{cauchy's theorem}, then
\begin{equation}
\label{int_{Gamma(R)} g(z) frac{1}{z - overline{w_1}} dz = 0}
        \int_{\Gamma(R)} g(z) \, \frac{1}{z - \overline{w_1}} \, dz = 0,
\end{equation}
by Cauchy's theorem.  Combining this with (\ref{frac{1}{2 pi i}
  int_{Gamma(R)} g(z) frac{1}{z - w} dz = g(w)}), we get that
\begin{equation}
\label{frac{1}{2 pi i} int_{Gamma(R)} g(z) (...) dz = g(w)}
 \frac{1}{2 \, \pi \, i} \, \int_{\Gamma(R)} g(z) \,
    \Big(\frac{1}{z - w} - \frac{1}{z - \overline{w_1}}\Big) \, dz = g(w)
\end{equation}
when $|w_1| < R$.  Of course,
\begin{equation}
\label{frac{1}{z - w} - frac{1}{z - overline{w_1}} = ...}
        \frac{1}{z - w} - \frac{1}{z - \overline{w_1}} 
 = \frac{(z - \overline{w_1}) - (z - w)}{(z - w) \, (z - \overline{w_1})}
           = \frac{w - \overline{w_1}}{(z - w) \, (z - \overline{w_1})}
\end{equation}
for each $z \in \overline{U}$, the modulus of which is $O(1/|z|^2)$
for $|z|$ large.  If $g(z)$ is bounded on $\overline{U}$, in addition
to the other conditions already mentioned, then we can take the limit
as $R \to \infty$ in (\ref{frac{1}{2 pi i} int_{Gamma(R)} g(z) (...)
  dz = g(w)}), to get that
\begin{equation}
\label{frac{1}{2 pi i} int_{-infty}^infty g(x) (...) dx = g(w), 2}
 \frac{1}{2 \, \pi \, i} \, \int_{-\infty}^\infty g(x) \,
  \Big(\frac{1}{x - w} - \frac{1}{x - \overline{w_1}}\Big) \, dx = g(w).
\end{equation}
In particular, the integrand is an integrable function on ${\bf R}$
when $g(x)$ is a bounded function on ${\bf R}$, because of
(\ref{frac{1}{z - w} - frac{1}{z - overline{w_1}} = ...}).  Note that
(\ref{frac{1}{2 pi i} int_{-infty}^infty g(x) (...) dx = g(w), 2}) is
the same as (\ref{frac{1}{2 pi i} int_{-infty}^infty g(x) (...) dx =
  g(w)}) when $w_1 = w$, which is the same as the representation of
$g(z)$ as $h(z)$ in (\ref{h(z) = A_y(f)(x) = ...}), with $f$ equal to
the restriction of $g$ to ${\bf R}$.  Remember that this
representation works for bounded continuous functions on
$\overline{U}$ that are harmonic on $U$, as in Section \ref{bounded
  harmonic functions}.

        If $w_1, w_2 \in U$, then it follows from (\ref{frac{1}{2 pi i} 
int_{-infty}^infty g(x) (...) dx = g(w), 2}) that
\begin{equation}
\label{int_{-infty}^infty g(x) (...) dx = 0}
        \int_{-\infty}^\infty g(x) \, \Big(\frac{1}{x - \overline{w_1}}
                              - \frac{1}{x - \overline{w_2}}\Big) \, dx = 0,
\end{equation}
because the right side of (\ref{frac{1}{2 pi i} int_{-infty}^infty
  g(x) (...) dx = g(w), 2}) does not depend on $w_1$.  Alternatively,
one can get this by applying (\ref{int_{-infty}^infty g(x) dx = 0}) to
\begin{equation}
\label{g(z) (frac{1}{z - overline{w_1}} - frac{1}{z - overline{w_2}})}
 g(z) \, \Big(\frac{1}{z - \overline{w_1}} - \frac{1}{z - \overline{w_2}}\Big)
\end{equation}
as a function of $z \in \overline{U}$ when $g(z)$ is bounded on
$\overline{U}$.  This uses the fact that
\begin{equation}
\label{frac{1}{z - overline{w_1}} - frac{1}{z - overline{w_2}} = ...}
        \frac{1}{z - \overline{w_1}} - \frac{1}{z - \overline{w_2}}
           = \frac{(z - \overline{w_2}) - (z - \overline{w_1})}
                   {(z - \overline{w_1}) \, (z - \overline{w_2})}
           = \frac{\overline{w_1} - \overline{w_2}}
                   {(z - \overline{w_1}) \, (z - \overline{w_2})}
\end{equation}
for every $z \in \overline{U}$, the modulus of which is $O(1/|z|^2)$
for $|z|$ large, as before.  In the other direction, if one has
(\ref{frac{1}{2 pi i} int_{-infty}^infty g(x) (...) dx = g(w), 2}) in
the special case where $w_1 = w$, and if one knows
(\ref{int_{-infty}^infty g(x) (...) dx = 0}) for every $w_1, w_2 \in
U$, then one can get (\ref{frac{1}{2 pi i} int_{-infty}^infty g(x)
  (...) dx = g(w), 2}) for every $w, w_1 \in U$.

        We also have that
\begin{equation}
\label{frac{1}{2 pi i} int_{-infty}^infty g(x) (...) dx = g(w_1) - g(w_2)}
 \frac{1}{2 \, \pi \, i} \, \int_{-\infty}^\infty g(x) \,
        \Big(\frac{1}{x - w_1} - \frac{1}{x - w_2}\Big) \, dx = g(w_1) - g(w_2)
\end{equation}
for every $w_1, w_2 \in U$ in this situation.  This can be derived
easily from (\ref{frac{1}{2 pi i} int_{-infty}^infty g(x) (...) dx =
  g(w), 2}), with suitable relabelling of the variables.  As a more
direct approch, one can use the Cauchy integral formula as in
(\ref{frac{1}{2 pi i} int_{Gamma(R)} g(z) frac{1}{z - w} dz = g(w)})
to get that
\begin{equation}
\label{frac{1}{2 pi i} int_{Gamma(R)} g(z) (...) dz = g(w_1) - g(w_2)}
 \frac{1}{2 \, \pi \, i} \, \int_{\Gamma(R)} g(z) \,
        \Big(\frac{1}{z - w_1} - \frac{1}{z - w_2}\Big) \, dz = g(w_1) - g(w_2)
\end{equation}
when $|w_1|, |w_2| < R$.  As usual,
\begin{equation}
\label{frac{1}{z - w_1} - frac{1}{z - w_2} = ...}
        \frac{1}{z - w_1} - \frac{1}{z - w_2}
           = \frac{(z - w_2) - (z - w_1)}{(z - w_1) \, (z - w_2)}
           = \frac{w_1 - w_2}{(z - w_1) \, (z - w_2)}
\end{equation}
for every $z \in \overline{U}$, the modulus of which is $O(1/|z|^2)$
for $|z|$ large.  This permits one to take the limit as $R \to \infty$
in (\ref{frac{1}{2 pi i} int_{Gamma(R)} g(z) (...) dz = g(w_1) -
  g(w_2)}), to get (\ref{frac{1}{2 pi i} int_{-infty}^infty g(x) (...)
  dx = g(w_1) - g(w_2)}).

        If $g(z)$ is any bounded holomorphic function on $U$, then
$g(z + i \, \tau)$ is a bounded continuous function of $z$ on $\overline{U}$
for each $\tau > 0$.  Thus (\ref{frac{1}{2 pi i} int_{-infty}^infty g(x) 
(...) dx = g(w), 2}) implies that
\begin{equation}
\label{frac{1}{2 pi i} int_{-infty}^infty ... dx = g(w + i tau)}
 \frac{1}{2 \, \pi \, i} \, \int_{-\infty}^\infty g(x + i \, \tau) \,
            \Big(\frac{1}{x - w} - \frac{1}{x - \overline{w_1}}\Big) \, dx 
         = g(w + i \, \tau)
\end{equation}
for every $w, w_1 \in U$.  Similarly, (\ref{int_{-infty}^infty g(x)
  (...) dx = 0}) and (\ref{frac{1}{2 pi i} int_{-infty}^infty g(x)
  (...) dx = g(w_1) - g(w_2)}) imply that
\begin{equation}
\label{int_{-infty}^infty g(x + i tau) (...) dx = 0}
 \int_{-\infty}^\infty g(x + i \, \tau) \, \Big(\frac{1}{x - \overline{w_1}}
                           - \frac{1}{x - \overline{w_2}}\Big) \, dx = 0
\end{equation}
and
\begin{eqnarray}
\label{frac{1}{2 pi i} int_{-infty}^infty ... dx = ....}
\lefteqn{\frac{1}{2 \, \pi \, i} \, \int_{-\infty}^\infty g(x + i \, \tau) \,
                       \Big(\frac{1}{x - w_1} - \frac{1}{x - w_2}\Big) \, dx} \\
         & = & g(w_1 + i \, \tau) - g(w_2 + i \, \tau) \nonumber
\end{eqnarray}
for every $w_1, w_2 \in U$.  If $w_1 = w$, then (\ref{frac{1}{2 pi i}
  int_{-infty}^infty ... dx = g(w + i tau)}) reduces to
(\ref{frac{1}{2 pi i} int_{-infty}^infty g(x + i tau) (...) dx = g(w +
  i tau)}), which is the same as (\ref{g(z + i tau) = ...}), and which
works for any bounded harmonic function on $U$.  Remember that
$L^\infty({\bf R})$ can be identified with the dual of $L^1({\bf R})$
in the usual way.  As in Section \ref{bounded harmonic functions}, one
can use the Banach--Alaoglu theorem and the separability of $L^1({\bf
  R})$ to get a sequence $\{\tau_j\}_{j = 1}^\infty$ of positive real
numbers that converges to $0$ such that $g(x + i \, \tau_j)$ converges
to some $g(x) \in L^\infty({\bf R})$ as $j \to \infty$ with respect to
the corresponding weak$^*$ topology on $L^\infty({\bf R})$.  Applying
(\ref{frac{1}{2 pi i} int_{-infty}^infty ... dx = g(w + i tau)}) to
$\tau_j$ and taking the limit as $j \to \infty$, we get that
(\ref{frac{1}{2 pi i} int_{-infty}^infty g(x) (...) dx = g(w), 2})
holds for every $w, w_1 \in U$, using this function $g(x)$ on ${\bf
  R}$ in the integral.  Similarly, we can apply
(\ref{int_{-infty}^infty g(x + i tau) (...) dx = 0}) and
(\ref{frac{1}{2 pi i} int_{-infty}^infty ... dx = ....}) to $\tau_j$
and take the limit as $j \to \infty$, to get that
(\ref{int_{-infty}^infty g(x) (...) dx = 0}) and (\ref{frac{1}{2 pi i}
  int_{-infty}^infty g(x) (...) dx = g(w_1) - g(w_2)}) hold for every
$w_1, w_2 \in U$, using this function $g(x)$ on ${\bf R}$ in the
integrals again.  If $g$ is a bounded harmonic function on $U$, then
(\ref{frac{1}{2 pi i} int_{-infty}^infty ... dx = g(w + i tau)}) holds
for $w_1 = w$, which implies that (\ref{frac{1}{2 pi i}
  int_{-infty}^infty g(x) (...)  dx = g(w), 2}) holds for $w_1 = w$,
as in Section \ref{bounded harmonic functions}.

\chapter{Commutative topological groups}
\label{commutative topological groups}

\section{Basic notions}
\label{basic notions}

        Let $A$ be a commutative group, in which the group operations
are expressed additively.  Suppose that $A$ is also equipped with a
topology, such that the group operations are continuous.  More precisely,
this means that
\begin{equation}
\label{(x, y) mapsto x + y}
        (x, y) \mapsto x + y
\end{equation}
is continuous as a mapping from $A \times A$ into $A$, using the
product topology on $A \times A$ associated to the given topology on
$A$, and that
\begin{equation}
\label{x mapsto -x}
        x \mapsto -x
\end{equation}
is continuous as a mapping from $A$ onto itself.  If $\{0\}$ is a
closed set in $A$, then we say that $A$ is a \emph{topological
  group}.\index{topological groups}\index{commutative topological
  groups} Of course, the definition of a topological group that is not
necessarily commutative is analogous, but we shall be primarily
concerned with the commutative case here.

        Continuity of addition on $A$ implies in particular that
\begin{equation}
\label{x mapsto x + a}
        x \mapsto x + a
\end{equation}
is continuous as a mapping from $A$ onto itself for every $a \in A$.
More precisely, these translation mappings on $A$ are homeomorphisms
from $A$ onto itself, since the inverse of (\ref{x mapsto x + a}) is
given by translation by $-a$.  Thus the hypothesis that $\{0\}$ be a
closed set in $A$ implies that $\{a\}$ is a closed set in $A$ for
every $a \in A$.  One can show that $A$ is also Hausdorff, using
continuity of addition at $0$.  A refinement of this argument shows
that $A$ is regular as a topological space.

        Any group is a topological group with respect to the discrete
topology.  The real line ${\bf R}$ is a commutative topological group
with respect to addition and the standard topology, and the unit
circle ${\bf T}$ is a commutative topological group with respect to
multiplication of complex numbers and the topology induced on ${\bf
  T}$ by the standard topology on the complex plane.  The set ${\bf
  R}_+$ of positive real numbers forms a commutative group with
respect to multiplication, which is a topological group with respect
to the standard topology.  If $A$ and $B$ are topological groups,
$\phi$ is a group isomorphism from $A$ onto $B$, and $\phi$ is also a
homeomorphism from $A$ onto $B$, then $\phi$ is considered to be an
isomorphism from $A$ onto $B$ as topological groups.  The exponential
function defines an isomorphism from ${\bf R}$ as a topological group
with respect to addition onto ${\bf R}_+$ as a topoloigcal group with
respect to multiplication, for instance.

        Let $V$ be a vector space over the real numbers, which is
a commutative group with respect to addition in particular.  If $V$ is
equipped with a topology such that the vector space operations are
continuous, and if $\{0\}$ is a closed set in $V$, then we say that
$V$ is a \emph{topological vector space}.\index{topological vector
  spaces} More precisely, the continuity of the vector space
operations means that addition on $V$ is continuous as a mapping from
$V \times V$ into $V$, and that scalar multiplication on $V$ is
continuous as a mapping from ${\bf R} \times V$ into $V$.  Here $V
\times V$ is equipped with the product topology corresponding to the
given topology on $V$, and ${\bf R} \times V$ is equipped with the
product topology corresponding to the standard topology on ${\bf R}$
and the given topology on $V$.  Continuity of scalar multiplication
implies in particular that multiplication by $-1$ is continuous on
$V$, so that $V$ is a topological group with respect to addition.

        Let $(M, d(x, y))$ be a metric space, and let $p$ be an element
of $M$.  The collection of open balls $B(p, 1/j)$ in $M$ with respect
to $d(x, y)$ centered at $p$ and with radius $1/j$ for each positive
integer $j$ forms a local base for the topology of $M$ at $p$.  This
local base for the topology of $M$ at $p$ has only finitely or
countably many elements, and one could just as well have used any
sequence of positive radii that converges to $0$.  If $A$ is a
commutative topological group, then a local base for the topology of
$A$ at any point leads to local bases for the topology of $A$ at every
other point in $A$, using translations.  If there is a local base for
the topology of $A$ at $0$ with only finitely or countably many
elements, then a famous theorem states that there is a metric on $A$
that is invariant under translations and determines the same topology
on $A$.

        Remember that a topological space $X$ is said to be
\emph{separable}\index{separable topological spaces} if there is a
dense set in $X$ with only finitely or countably many elements.
If there is a base for the topology of $X$ with only finitely or
countably many elements, then $X$ is separable.  As a partial converse,
if the topology on $X$ is determined by a metric, and if $X$ is separable,
then there is a base for the topology of $X$ with only finitely or countably
many elements.  In particular, if $A$ is a commutative topological group
with a local base for its topology at $0$ with only finitely or countably
many elements, and if $A$ is separable, then there is a base for the topology
of $A$ with only finitely or countably many elements.  This can be shown
directly, in much the same way as for metric spaces, or by using the
metrization theorem mentioned in the previous paragraph.

        A subset $E$ of a metric space $M$ is said to be \emph{totally 
bounded}\index{totally bounded sets} if for each $r > 0$, $E$
is contained in the union of finitely many balls of radius $r$.
Similarly, a subset $E$ of a commutative topological group $A$ is said
to be totally bounded if for each open set $U \subseteq A$ with $0 \in
U$, $E$ is contained in the union of finitely many translates of $U$.
Of course, compact subsets of metric spaces are totally bounded, and
compact subsets of commutative topological groups are totally bounded
as well.  If there is a translation-invariant metric $d(x, y)$ on a
commutative topological group $A$ that determines the same topology on
$A$, then it is easy to see that $E \subseteq A$ is totally bounded as
a subset of $A$ as a commutative topological group if and only if $E$
is totally bounded with respect to $d(x, y)$.

        A subset of a topological space is said to be
\emph{$\sigma$-compact}\index{sigma-compact@$\sigma$-compact sets}
if it can be expressed as the union of finitely or countably many
compact sets.  It is well known that compact metric spaces are totally
bounded, and hence separable.  Similarly, $\sigma$-compact metric spaces
are separable.  If $A$ is a commutative topological group with a local
base for its topology at $0$ with only finitely or countably many elements,
and if $A$ is $\sigma$-compact, then $A$ is separable.  As before, this
can be shown directly, in much the same way as for metric spaces, or
using the metrization theorem mentioned earlier.

\section{Haar measure}
\label{haar measure}

        Remember that a topological space $X$ is said to be
\emph{locally compact}\index{locally compact topological spaces}
if for each point $p \in X$ there is an open set $U \subseteq X$ and a
compact set $K \subseteq X$ such that $p \in U$ and $U \subseteq K$.
If $X$ is Hausdorff, then compact subsets of $X$ are closed, and $X$
is locally compact if and only if for each $p \in X$ there is an open
set $U \subseteq X$ such that $p \in U$ and the closure $\overline{U}$
of $U$ in $X$ is compact.  If $A$ is a commutative topological group,
then $A$ is locally compact if and only if there is an open set $U
\subseteq A$ such that $0 \in U$ and $\overline{U}$ is compact,
because of continuity of translations.  If $V$ is a topological vector
space over ${\bf R}$ that is locally compact, then it is well known
that $V$ has finite dimension.  In this case, it is also well known
that $V$ is either isomorphic to ${\bf R}^n$ as a topological vector
space for some positive integer $n$, or $V = \{0\}$.

        If $A$ is a locally compact commutative topological group,
then it is well known that there is a nonnegative Borel measure $H$ on
$A$, known as \emph{Haar measure},\index{Haar measure} that is
invariant under translations and has some other nice properties.
To say that $H$ is invariant under translations means that
\begin{equation}
\label{H(E + a) = H(E)}
        H(E + a) = H(E)
\end{equation}
for every Borel set $E \subseteq A$ and $a \in A$, where
\begin{equation}
\label{E + a = {x + a : x in E}}
        E + a = \{x + a : x \in E\}.
\end{equation}
Note that $E + a$ is a Borel set in $A$ for every Borel set $E
\subseteq A$ and $a \in A$, because the translation mappings (\ref{x
  mapsto x + a}) are homeomorphisms on $A$.  This measure $H$ also
satisfies $H(U) > 0$ when $U$ is a nonempty open set in $A$, $H(K) <
\infty$ when $K \subseteq A$ is compact, and some standard regularity
conditions.  Haar measure is unique in the sense that if $H'$ is
another measure on $A$ with the same properties, then $H'$ is equal to
a positive real number times $H$.  Lebesgue measure on ${\bf R}$
satisfies the requirements of Haar measure, as does arc-length measure
on ${\bf T}$.  If $A$ is equipped with the discrete topology, then
counting measure on $A$ satisfies the requirements of Haar measure.

        Let $A$ be a locally compact commutative topological group, and
let $H$ be a Haar measure on $A$.  We would like to check that
\begin{equation}
\label{H(-E) = H(E)}
        H(-E) = H(E)
\end{equation}
for every Borel set $E \subseteq A$, where
\begin{equation}
\label{-E = {-x : x in E}}
        -E = \{-x : x \in E\}.
\end{equation}
Note that $-E$ is a Borel set in $A$ when $E$ is a Borel set, because
(\ref{x mapsto -x}) is a homeomorphism on $A$.  It easy easy to see
that $H(-E)$ also satisfies the requirements of Haar measure on $A$.
The uniqueness of Haar measure implies that there is a positive real
number $c$ such that
\begin{equation}
\label{H(-E) = c H(E)}
        H(-E) = c \, H(E)
\end{equation}
for every Borel set $E \subseteq A$.  Because $A$ is locally compact,
there is an open set $U \subseteq A$ such that $0 \in U$ and
$\overline{U}$ is compact, as before.  It follows that
\begin{equation}
        W = U \cap (-U)
\end{equation}
is also an open set in $A$ that contains $0$ and has compact closure,
so that $H(W)$ is positive and finite.  By construction,
\begin{equation}
\label{-W = W}
        -W = W,
\end{equation}
 which implies that $c = 1$, as desired, by applying (\ref{H(-E) = c
   H(E)}) to $W$.

        Let $A$ be a locally compact commutative topological group again,
and let $H$ be a Haar measure on $A$.  If $f(x)$ is a nonnegative real-valued
Borel measurable function on $A$, then
\begin{equation}
\label{f_a(x) = f(x - a), 2}
        f_a(x) = f(x - a)
\end{equation}
has the same properties for each $a \in A$.  Using
translation-invariance of Haar measure, one can check that
\begin{equation}
\label{int_A f_a(x) dH(x) = int_A f(x) dH(x)}
        \int_A f_a(x) \, dH(x) = \int_A f(x) \, dH(x)
\end{equation}
for every $a \in A$.  This is the same as (\ref{H(E + a) = H(E)}) when
$f$ is the characteristic or indicator function on $A$ associated to a
Borel set $E \subseteq A$, which is equal to $1$ on $E$ and to $0$ on
$A \backslash E$.  If $f$ is a measurable simple function on $A$,
which is a linear combination of indicator functions associated to
measurable subsets of $A$, then (\ref{int_A f_a(x) dH(x) = int_A f(x)
  dH(x)}) follows from the previous case by linearity.  One can then
get (\ref{int_A f_a(x) dH(x) = int_A f(x) dH(x)}) for arbitrary
nonnegative measurable functions on $A$, by approximating such
functions by simple functions.  Similarly, $f(-x)$ is a nonnegative
real-valued Borel measurable function on $A$ when $f$ is, because
(\ref{x mapsto -x}) is a homeomorphism on $A$, and
\begin{equation}
\label{int_A f(-x) dH(x) = int_A f(x) dH(x)}
        \int_A f(-x) \, dH(x) = \int_A f(x) \, dH(x),
\end{equation}
by (\ref{H(-E) = H(E)}).  If $f$ is a real or complex-valued
integrable function on $A$ with respect to $H$, then it follows that
(\ref{f_a(x) = f(x - a), 2}) is also integrable on $H$ for every $a
\in A$, and that (\ref{int_A f_a(x) dH(x) = int_A f(x) dH(x)}) still
holds.  In this case, $f(-x)$ is an integrable function on $A$ with
respect to $H$ as well, and satisfies (\ref{int_A f(-x) dH(x) = int_A
  f(x) dH(x)}).

        Let $C_{com}(A)$ be the space of continuous complex-valued functions
on $A$ with compact support.  If $f \in C_{com}(A)$, then $f$ is Borel
measurable and integrable with respect to $H$, because the Haar measure
of a compact subset of $A$ is finite.  Thus
\begin{equation}
\label{L(f) = int_A f(x) dH(x)}
        L(f) = \int_A f(x) \, dH(x)
\end{equation}
defines a linear functional on $C_{com}(A)$, as a vector space with
respect to pointwise addition and scalar multiplication.  Note that
$f_a \in C_{com}(A)$ for every $a \in A$ when $f \in C_{com}(A)$,
and that
\begin{equation}
\label{L(f_a) = L(f)}
        L(f_a) = L(f)
\end{equation}
for every $a \in A$, by (\ref{int_A f_a(x) dH(x) = int_A f(x) dH(x)}).
If $f \in C_{com}(A)$ is real-valued and nonnegative on $A$, then
$L(f)$ is a nonnegative real number, because $H$ is a nonnegative
measure on $A$.  If $f \in C_{com}(A)$ is real-valued and nonnegative
on $A$, and if $f(x) > 0$ for some $x \in A$, then $L(f) > 0$, because
$H(U) > 0$ when $U$ is a nonempty open set in $A$.  Similarly,
$C_{com}(A)$ and $L$ are invariant under (\ref{x mapsto -x}).

        Conversely, if $L$ is any nonnegative linear functional on
$C_{com}(A)$, then the Riesz representation theorem implies that $L$
can be expressed as in (\ref{L(f) = int_A f(x) dH(x)}) for a unique
nonnegative Borel measure $H$ on $A$ that is finite on compact sets
and satisfies certain other regularity conditions.  If $L$ is
invariant under translations, in the sense that (\ref{L(f_a) = L(f)})
holds for every $f \in C_{com}(A)$ and $a \in A$, then the uniqueness
of $H$ implies that $H$ is invariant under translations too.  If $L(f)
> 0$ when $f \in C_{com}(A)$ is real-valued, nonnegative, and $f(x) >
0$ for some $x \in A$, then $H(U) > 0$ for every nonempty open set $U
\subseteq A$.  This uses Urysohn's lemma, to get such a function $f$
with compact support contained in $U$.

        A linear functional on $C_{com}(A)$ with the properties just
described is known as a \emph{Haar integral}\index{Haar integral} on
$A$.  Thus the existence and uniqueness of Haar measure can be
reformulated in terms of Haar integrals, and are often established in
that way.  Note that the invariance of a Haar integral under (\ref{x
  mapsto -x}) can be derived from its uniqueness, in essentially the
same way as for Haar measure.

        If $A$ is a compact commutative topological group and $H$
is a Haar measure on $A$, then $H(A) < \infty$, and $H$ is often
normalized so that $H(A) = 1$.  If $A$ is a locally compact
commutative topological group which is also $\sigma$-compact, then
Haar measure on $A$ is $\sigma$-finite.  If there is a base for the
topology of $A$ with only finitely or countably many elements, then
any open covering of a subset of $A$ can be reduced to a subcovering
with only finitely or countably many elements, by Lindel\"of's
theorem.  In this case, it follows that $A$ is $\sigma$-compact when
$A$ is locally compact.  An analogous argument implies that every open
set in $A$ is $\sigma$-compact under these conditions, using the fact
that $A$ is regular as a topological space.  Note that closed subsets
of $A$ are $\sigma$-compact when $A$ is $\sigma$-compact.  If $A$ is
metrizable, then it is well known that every open set in $A$ can be
expressed as the union of finitely or countably many closed sets.
This implies that every open set in $A$ is $\sigma$-compact when $A$
is $\sigma$-compact and metrizable.

\section{Direct sums and products}
\label{direct sums, products}

        Let $I$ be a nonempty set, and suppose that $A_j$ is a
commutative group for each $j \in I$.  The  \emph{direct product}\index{direct 
products} of the $A_j$'s is the Cartesian product
\begin{equation}
\label{prod_{j in I} A_j}
        \prod_{j \in I} A_j,
\end{equation}
where the group operations are defined coordinatewise.  The
\emph{direct sum}\index{direct sums} of the $A_j$'s is denoted
\begin{equation}
\label{sum_{j in I} A_j}
        \sum_{j \in I} A_j,
\end{equation}
and is the subgroup of (\ref{prod_{j in I} A_j}) consisting of
elements of the product for which all but finitely many coordinates
are equal to $0$.  The direct product (\ref{prod_{j in I} A_j}) is
also known as the \emph{complete direct sum}\index{complete direct
  sums} of the $A_j$'s, which is the same as the direct sum when $I$
has only finitely many elements.

        If $A_j$ is a topological group for each $j \in I$, then
the direct product (\ref{prod_{j in I} A_j}) is a topological group
too, with respect to the corresponding product topology.
Alternatively, one could use the \emph{strong product
  topology}\index{strong product topology} on (\ref{prod_{j in I}
  A_j}), which is the topology generated by arbitrary products of open
subsets of the $A_j$'s.  Of course, this is the same as the product
topology on (\ref{prod_{j in I} A_j}) when $I$ has only finitely many
elements.  It is easy to see that the direct sum (\ref{sum_{j in I}
  A_j}) is dense in the direct product (\ref{prod_{j in I} A_j}) with
respect to the corresponding product topology.  However, one can check
that the direct sum (\ref{sum_{j in I} A_j}) is a closed subgroup of
the direct product (\ref{prod_{j in I} A_j}) with respect to the
strong product topology.

        Let $I$ be a nonempty set again, and suppose that $V_j$ is a
vector space over the real numbers for each $j \in I$.  As before, the
direct product\index{direct products} of the $V_j$'s is the Cartesian
product
\begin{equation}
\label{prod_{j in I} V_j}
        \prod_{j \in I} V_j,
\end{equation}
where addition and scalar multiplication are defined coordinatewise,
and the direct sum\index{direct sums}
\begin{equation}
\label{sum_{j in I} V_j}
        \sum_{j \in I} V_j
\end{equation}
of the $V_j$'s is the linear subspace of (\ref{prod_{j in I} V_j})
consisting of elements of the product for which all but finitely many
coordinates are equal to $0$.  If $V_j$ is a topological vector space
for each $j \in I$, then the direct product (\ref{prod_{j in I} V_j})
is a topological vector space with respect to the product topology.
This does not work when $I$ is infinite and the direct product is
equipped with the strong product topology, because scalar
multiplication is not continuous as a function of the scalar at
vectors with infinitely many nonzero coordinates.  However, one can
check that the direct sum (\ref{sum_{j in I} V_j}) is a topological
vector space with respect to the topology induced by the strong
product topology on (\ref{prod_{j in I} V_j}).

        Let $I$ be a nonempty set with only finitely or countably many
elements, and suppose that $A_j$ is a commutative topological group
for each $j \in I$.  If there is a local base for the topology of
$A_j$ at $0$ with only finitely or countably many elements for each $j
\in I$, then there is a local base for the product topology on the
direct product (\ref{prod_{j in I} A_j}) at $0$ with only finitely or
countably many elements.  This is a special case of an analogous
statement for arbitrary topological spaces.  In particular, if there
is a translation-invariant metric on $A_j$ that determines the same
topology on $A_j$ for each $j \in I$, then there is a
translation-invariant metric on the direct product (\ref{prod_{j in I}
  A_j}) for which the corresponding topology is the product topology.
This can be shown directly, using the same type of constructions
as for products of finitely or countably many metric spaces.

        Suppose that $\mathcal{B}_j$ is a base for the topology of
$A_j$ for each $j \in I$.  Let $\mathcal{B}$ be the collection of
subsets of the direct product (\ref{prod_{j in I} A_j}) of the form
$\prod_{j \in I} U_j$, where $U_j = A_j$ for all but finitely many $j
\in I$, and $U_j \in \mathcal{B}_j$ when $U_j \ne A_j$.  Thus $\mathcal{B}$
is a base for the product topology on (\ref{prod_{j in I} A_j}).
If $\mathcal{B}_j$ has only finitely or countably many elements
for each $j \in I$, and $I$ has only finitely or countably many elements,
then it is well known that $\mathcal{B}$ has only finitely or countably
many elements as well.  In this case, it follows that every open set
in (\ref{prod_{j in I} A_j}) can be expressed as the union of finitely
or countably many sets of the form $\prod_{j \in I} U_j$, where $U_j$
is an open set in $A_j$ for each $j \in I$, and $U_j = A_j$ for all
but finitely many $j$.

        If $A_1, \ldots, A_n$ are finitely many locally compact
commutative topological groups, then their direct product $\prod_{j =
  1}^n A_j$ is a locally compact commutative topological group with
respect to the product topology.  In this case, Haar measure on
$\prod_{j = 1}^n A_j$ basically corresponds to a product of Haar
measures on the $A_j$'s.  One way to deal with this is to define a
Haar integral as a linear functional on $C_{com}(\prod_{j = 1}^n
A_j)$, using Haar integrals on the $A_j$'s to integrate each variable
separately.  Alternatively, suppose that there is a base for the
topology of $A_j$ with only finitely or countably many elements for
each $j$.  This implies that $A_j$ is $\sigma$-compact for each $j$,
and hence that Haar measure on $A_j$ is $\sigma$-finite, as in the
previous section.  This also implies that every open set in $\prod_{j
  = 1}^n A_j$ can be expressed as the union of finitely or countably
many sets of the form $\prod_{j = 1}^n U_j$, where $U_j$ is an open
set in $A_j$ for each $j = 1, \ldots, n$, as in the preceding
paragraph.  Under these conditions, one can define a product measure
on $\prod_{j = 1}^n A_j$ corresponding to Haar measures on the
$A_j$'s, where open subsets of $\prod_{j = 1}^n A_j$ are measurable
with respect to the usual product measure construction, starting with
Borel sets in $A_j$ for each $j$.

        Now let $I$ be a nonempty set, and suppose that $A_j$ is a
compact commutative topological group for each $j \in I$, so that the
direct product (\ref{prod_{j in I} A_j}) is also a compact commutative
topological group with respect to the product topology.  In this case,
we can choose Haar measure $H_j$ on $A_j$ such that $H_j(A_j) = 1$ for
each $j \in I$, and Haar measure on (\ref{prod_{j in I} A_j})
basically corresponds to a product of the $H_j$'s.  As before, one way
to deal with this is to define a Haar integral on (\ref{prod_{j in I}
  A_j}), as a linear functional on the space of continuous
complex-valued functions on (\ref{prod_{j in I} A_j}), using Haar
integrals on the $A_j$'s to integrate each variable separately.  If
$I$ has infinitely many elements, then it is helpful to observe that
continuous functions on (\ref{prod_{j in I} A_j}) can be approximated
uniformly by continuous functions that depend only on finitely many
variables.  If $I$ has only finitely or countably many elements, and
if there is a base for the topology of $A_j$ with only finitely or
countably many elements for each $j \in I$, then one can use standard
product measure constructions instead.

\section{Continuous homomorphisms}
\label{continuous homomorphisms}

        If $A$ is a nonempty set and $B$ is a commutative group, then
the collection of all mappings from $A$ into $B$ is a commutative
group with respect to pointwise addition.  If $A$ and $B$ are both
commutative groups, then the collection of all homomorphisms from $A$
into $B$ is a subgroup of the group of all mappings from $A$ into $B$.
Similarly, if $A$ is a nonempty topological space and $B$ is a
commutative topological group, then the collection of continuous
mappings from $A$ into $B$ is a subgroup of the group of all mappings
from $A$ into $B$.  If $A$ and $B$ are both commutative topological
groups, then the collection
\begin{equation}
\label{hom(A, B)}\index{Hom(A, B)@$\hom(A, B)$}
        \hom(A, B)
\end{equation}
of all continuous homomorphisms from $A$ into $B$ is a subgroup of the
group of all mappings from $A$ into $B$, which is the same as the
intersection of the two previous subgroups.  Of course, if $A$ is
equipped with the discrete topology, then $\hom(A, B)$ is the same as
the group of all homomorphisms from $A$ into $B$.

        Suppose that $A$ and $B$ are commutative topological groups,
and that $\phi$ is a homomorphism from $A$ into $B$.  If $\phi$ is
continuous at $0$, then it is easy to see that $\phi$ is continuous
at every point in $A$, by continuity of translations.  In this case,
the kernel of $\phi$ is a closed subgroup of $A$.

        If $B$ is any commutative group, then $j \cdot b \in B$ can be
defined for each $j \in {\bf Z}$ and $b \in B$ in the usual way, and
\begin{equation}
\label{j mapsto j cdot b}
        j \mapsto j \cdot b
\end{equation}
defines a homomorphism from ${\bf Z}$ into $B$.  Conversely, if $\phi$ is
any homomorphism from ${\bf Z}$ into $B$, then
\begin{equation}
\label{phi(j) = j cdot phi(1)}
        \phi(j) = j \cdot \phi(1)
\end{equation}
for every $j \in {\bf Z}$.

        Similarly,
\begin{equation}
\label{x mapsto a x}
        x \mapsto a \, x
\end{equation}
defines a continuous homomorphism from ${\bf R}$ into itself for each
$a \in {\bf R}$.  If $\phi$ is any continuous homomorphism from ${\bf
  R}$ into itself, then one can check that
\begin{equation}
\label{phi(x) = x phi(1)}
        \phi(x) = x \, \phi(1)
\end{equation}
for every $x \in {\bf R}$.  More precisely, one can first show that
this holds when $x$ is rational, and then for all $x \in {\bf R}$, by
continuity.

        Using the complex exponential function, we get a continuous
homomorphism
\begin{equation}
\label{x mapsto exp(i a x)}
        x \mapsto \exp(i \, a \, x)
\end{equation}
from ${\bf R}$ into ${\bf T}$ for every $a \in {\bf R}$.  If $\phi$
is any continuous mapping from ${\bf R}$ into ${\bf T}$ such that
$\phi(0) = 1$, then there is a unique continuous mapping $\psi$
from ${\bf R}$ into itself such that $\psi(0) = 0$ and
\begin{equation}
        \phi(x) = \exp(i \, \psi(x))
\end{equation}
for every $x \in {\bf R}$.  If $\phi$ is a continuous homomorphism
from ${\bf R}$ into ${\bf T}$, then one can show that $\psi$ is a
continuous homomorphism from ${\bf R}$ into itself.  This implies that
$\psi$ can be expressed as (\ref{x mapsto a x}) for some $a \in {\bf
  R}$, and hence that $\phi$ can be expressed as (\ref{x mapsto exp(i
  a x)}).

        Of course,
\begin{equation}
\label{z mapsto z^j}
        z \mapsto z^j
\end{equation}
defines a continuous homomorphism from ${\bf T}$ into itself for every
integer $j$.  If $\phi$ is any continuous homomorphism from ${\bf T}$
into itself, then $\phi(\exp(i \, t))$ defines a continuous
homomorphism from ${\bf R}$ into ${\bf T}$.  It follows that there is
an $a \in {\bf R}$ such that
\begin{equation}
\label{phi(exp(i t)) = exp(i a t)}
        \phi(\exp(i \, t)) = \exp(i \, a \, t)
\end{equation}
for every $t \in {\bf R}$, as in the previous paragraph.  This also has to
be equal to $1$ when $t = 2 \, \pi$, which implies that $a \in {\bf Z}$.
Thus every continuous homomorphism from ${\bf T}$ into itself is of the
form (\ref{z mapsto z^j}) for some $j \in {\bf Z}$.

        Now let $n$ be a positive integer, and consider the quotient group
${\bf Z} / n \, {\bf Z}$, which is a cyclic group of order $n$.  A
homomorphism from ${\bf Z} / n \, {\bf Z}$ into a commutative group $B$
is basically the same as a homomorphism from ${\bf Z}$ into $B$ whose
kernel contains $n \, {\bf Z}$.  This is the same as a homomorphism from
${\bf Z}$ into $B$ that sends $n \in {\bf Z}$ to the identity element of $B$.

        Let $w \in {\bf T}$ be an $n$th root of unity, so that $w^n = 1$.
Thus
\begin{equation}
\label{j mapsto w^j}
        j \mapsto w^j
\end{equation}
defines a homomorphism from ${\bf Z}$ into ${\bf T}$ whose kernel
contains $n \, {\bf Z}$, and which therefore corresponds to a
homomorphism from ${\bf Z} / n \, {\bf Z}$ into ${\bf T}$.
Conversely, every homomorphism from ${\bf Z} / n \, {\bf Z}$ into
${\bf T}$ is of this form.

        Suppose that $V$ is a topological vector space over the real
numbers, which is a commutative topological group with respect to
addition in particular.  Every continuous linear mapping from $V$ into
${\bf R}$ is a continuous group homomorphism with respect to addition.
Conversely, every continuous group homomorphism from $V$ into ${\bf
  R}$ with respect to addition is linear as a mapping between $V$ and
${\bf R}$ as vector spaces over ${\bf R}$, because of the earlier
characterization of continuous homomorphisms from ${\bf R}$ into
itself.  This shows that $\hom(V, {\bf R})$ can be identified as a
commutative group with respect to addition with the dual space $V'$ of
continuous linear functionals on $V$, which is also a vector space
with respect to pointwise scalar multiplication.  However, it is well
known that there are nontrivial topological vector spaces $V$ over
${\bf R}$ such that $V' = \{0\}$.  If $\lambda$ is a continuous linear
functional on $V$, then
\begin{equation}
\label{v mapsto exp(i lambda(v))}
        v \mapsto \exp(i \, \lambda(v))
\end{equation}
defines a continuous homomorphism from $V$ as a topological group with
respect to addition into ${\bf T}$.  One can check that every
continuous homomorphism from $V$ into ${\bf T}$ is of this form, using
the analogous statement for the real line.

        Note that the only bounded subgroup of ${\bf R}$ is the trivial
subgroup $\{0\}$.  This implies that every continuous homomorphism
from a compact commutative topological group into ${\bf R}$ is
trivial.  Similarly,
\begin{equation}
\label{{z in {bf T} : re z > 0}}
        \{z \in {\bf T} : \re z > 0\}
\end{equation}
is a relatively open set in ${\bf T}$ that contains $1$, and any
subgroup of ${\bf T}$ contained in (\ref{{z in {bf T} : re z > 0}}) is
trivial.  Thus every homomorphism from a commutative group into ${\bf
  T}$ with values in (\ref{{z in {bf T} : re z > 0}}) is trivial.

        Let $A$ be a commutative topological group, and let $\phi$ be
a homomorphism from $A$ into the real or complex numbers, as a
commutative topological group with respect to addition and the
standard topology.  Suppose that there is an open set $U_0 \subseteq A$
such that $0 \in U_0$ and
\begin{equation}
\label{|phi(x)| < 1}
        |\phi(x)| < 1
\end{equation}
for every $x \in U_0$.  Using continuity of addition on $A$ at $0$, one
can find a sequence $U_1, U_2, U_3, \ldots$ of open subsets of $A$
such that $0 \in U_j$ and
\begin{equation}
\label{U_j + U_j subseteq U_{j - 1}}
        U_j + U_j \subseteq U_{j - 1}
\end{equation}
for each $j \in {\bf Z}_+$.  If $x \in U_j$ for some $j \ge 0$, then
$2^j \cdot x \in U_0$, and hence
\begin{equation}
\label{2^j |phi(x)| = |phi(2^j cdot x)| < 1}
        2^j \, |\phi(x)| = |\phi(2^j \cdot x)| < 1,
\end{equation}
which is to say that $|\phi(x)| < 2^{-j}$.  This implies that $\phi$
is continuous at $0$, so that $\phi$ is continuous on $A$.

        Now let $\phi$ be a homomorphism from $A$ into ${\bf T}$,
and suppose that there is an open set $U_0$ in $A$ such that $0 \in
U_0$ and
\begin{equation}
\label{re phi(x) > 0}
        \re \phi(x) > 0
\end{equation}
for every $x \in U_0$.  As before, there is a sequece $U_1, U_2, U_3,
\ldots$ of open subsets of $A$ that satisfy $0 \in U_j$ and (\ref{U_j
  + U_j subseteq U_{j - 1}}) for each $j \ge 1$.  Note that $U_j
\subseteq U_{j - 1}$ for each $j \ge 1$, because $0 \in U_j$, so that
$U_j \subseteq U_0$ for each $j$.  If $x \in U_k$ for some $k \ge 0$,
then $2^j \, x \in U_{k - j} \subseteq U_0$ when $0 \le j \le k$, and
hence
\begin{equation}
\label{re phi(x)^{2^j} = re phi(2^j cdot x) > 0}
        \re \phi(x)^{2^j} = \re \phi(2^j \cdot x) > 0
\end{equation}
for $j = 0, \ldots, k$. by (\ref{re phi(x) > 0}).  This implies that
\begin{equation}
        \sup_{x \in U_k} |\phi(x) - 1| \to 0
\end{equation}
as $k \to \infty$, so that $\phi$ is continuous at $0$, and thus
everywhere on $A$.

\section{Sums and products, continued}
\label{sums, products, continued}

        Let $A_1, \ldots, A_n$ be finitely many commutative groups,
and let $B$ be another commutative group.  If $\phi_j$ is a
homomorphism from $A_j$ into $B$ for each $j = 1, \ldots, n$, then
\begin{equation}
\label{phi(x) = sum_{j = 1}^n phi_j(x_j)}
        \phi(x) = \sum_{j = 1}^n \phi_j(x_j)
\end{equation}
defines a homomorphism from $\prod_{j = 1}^n A_j$ into $B$.
Conversely, it is easy to see that every homomorphism from $\prod_{j =
  1}^n A_j$ into $B$ is of this form.  Suppose now that $A_1, \ldots,
A_n$ and $B$ are topological groups, and let $\prod_{j = 1}^n A_j$ be
equipped with the product topology.  If $\phi_j$ is a continuous
homomorphism from $A_j$ into $B$ for $j = 1, \ldots, n$, then
(\ref{phi(x) = sum_{j = 1}^n phi_j(x_j)}) is a continuous homomorphism
from $\prod_{j = 1}^n A_j$ into $B$, and every continuous homomorphism
from $\prod_{j = 1}^n A_j$ into $B$ is of this form.

        Let $I$ be a nonempty set, let $A_j$ be a commutative group
for each $j \in I$, and let $B$ be another commutative group.  Also
let $j_1, \ldots, j_l$ be finitely many elements of $I$, and let
$\phi_{j_k}$ be a homomorphism from $A_{j_k}$ into $B$ for $k = 1,
\ldots, l$.  Thus
\begin{equation}
\label{phi(x) = sum_{k = 1}^l phi_{j_k}(x_{j_k})}
        \phi(x) = \sum_{k = 1}^l \phi_{j_k}(x_{j_k})
\end{equation}
defines a homomorphism from the direct product $\prod_{j \in I} A_j$
into $B$, where $x_j \in A_j$ denotes the $j$th component of $x \in
\prod_{j \in I} A_j$ for each $j \in I$.  If $A_j$ is a commutative
topological group for each $j \in I$, then we have seen that $\prod_{j
  \in I} A_j$ is a commutative topological group with respect to the
corresponding product topology.  If $B$ is also a commutative
topological group, and $\phi_{j_k}$ is continuous as a mapping from
$A_{j_k}$ into $B$ for each $k$, then it is easy to see that $\phi$ is
continuous with respect to the product topology on $\prod_{j \in I}
A_j$.

        Suppose in addition that there is an open set $V$ in $B$
that contains the identity element, but does not contain any
nontrivial subgroup of $B$.  This property is satisfied by both ${\bf
  R}$ and ${\bf T}$, as mentioned in the previous section.  Let $\phi$
be a continuous homomorphism from $\prod_{j \in I} A_j$ into $B$,
where $\prod_{j \in I} A_j$ is equipped with the product topology.
Thus $\phi^{-1}(V)$ is an open set in $\prod_{j \in I} A_j$ that
contains $0$.  It follows that there are open subsets $U_j$ of $A_j$
for each $j \in I$, such that $0 \in U_j$ for each $j \in I$,
$U_j = A_j$ for all but finitely many $j \in I$, and
\begin{equation}
\label{prod_{j in I} U_j subseteq phi^{-1}(V)}
        \prod_{j \in I} U_j \subseteq \phi^{-1}(V).
\end{equation}
If we put $C_j = A_j$ when $U_j = A_j$, and $C_j = \{0\}$ when $U_j
\ne A_j$, then
\begin{equation}
\label{prod_{j in I} C_j}
        \prod_{j \in I} C_j
\end{equation}
is a subgroup of $\prod_{j \in I} A_j$, and
\begin{equation}
\label{phi(prod_{j in I} C_j) subseteq V}
        \phi\Big(\prod_{j \in I} C_j\Big) \subseteq V,
\end{equation}
by (\ref{prod_{j in I} U_j subseteq phi^{-1}(V)}).  Because
$\phi\Big(\prod_{j \in I} C_j\Big)$ is a subgroup of $B$, our
hypothesis on $V$ implies that $\phi\Big(\prod_{j \in I} C_j\Big)$ is
the trivial subgroup of $B$.  This means that $\phi(x)$ depends only
on the $j$th coordinate $x_j$ of $x$ for finitely many $j \in I$,
since $C_j = A_j$ for all but finitely many $j \in I$.  Using this, it
is easy to see that $\phi$ can be expressed as in the preceding
paragraph under these conditions.

        Let $I$ be a nonempty set again, let $A_j$ be a commutative
group for each $j \in I$, and let $B$ be another commutative group.
If $\phi_j$ is a homomorphism from $A_j$ into $B$ for each $j \in I$,
then
\begin{equation}
\label{phi(x) = sum_{j in I} phi_j(x_j)}
        \phi(x) = \sum_{j \in I} \phi_j(x_j)
\end{equation}
defines a homomorphism from the direct sum $\sum_{j \in I} A_j$ into
$B$.  Of course, $x_j = 0$ for all but finitely many $j \in I$ when $x
\in \sum_{j \in I} A_j$, which implies that $\phi_j(x_j) = 0$ for all
but finitely many $j$, so that the sum in (\ref{phi(x) = sum_{j in I}
  phi_j(x_j)}) makes sense.  Conversely, it is easy to see that every
homomorphism from $\sum_{j \in I} A_j$ into $B$ is of this form.

        Now suppose that $A_j$ is a commutative topological group
for each $j \in I$, and that $B$ is a commutative topological group as
well.  If $\phi$ is a homomorphism from $\sum_{j \in I} A_j$ into $B$
such that $\phi(x)$ is continuous as a function of $x_j \in A_j$ for
each $j \in I$, then $\phi$ can be expressed as in (\ref{phi(x) =
  sum_{j in I} phi_j(x_j)}), where $\phi_j$ is a continuous
homomorphism from $A_j$ into $B$ for each $j$.

        Conversely, let $\phi_j$ be a continuous homomorphism from $A_j$
into $B$ for each $j \in I$, and let $\phi$ be the corresponding
homomorphism from $\sum_{j \in I} A_j$ into $B$, as in (\ref{phi(x) =
  sum_{j in I} phi_j(x_j)}).  Of course, if $I$ has only finitely many
elements, then $\sum_{j \in I} A_j$ is the same as $\prod_{j \in I}
A_j$, and $\phi$ is continuous with respect to the product topology on
$\prod_{j \in I} A_j$.  Otherwise, one still has that $\phi(x)$ is
continuous as a function of $x_j \in A_j$ for each $j \in I$, and
that $\phi(x)$ is continuous as a function of any finite collection
of these variables, with respect to the product topology on the
corresponding product of finitely many $A_j$'s.

        If $I$ is countably infinite, then one can show that $\phi$ is
continuous with respect to the topology induced on $\sum_{j \in I}
A_j$ by the strong product topology on $\prod_{j \in I} A_j$.  To do
this, it is convenient to reduce to the case where $I = {\bf Z}_+$.
Let $V$ be an open set in $B$ that contains $0$, and let $V_1, V_2,
V_3, \ldots$ be a sequence of open subsets of $B$ such that $0 \in
V_j$ and
\begin{equation}
\label{V_j + V_j subseteq V_{j - 1}}
        V_j + V_j \subseteq V_{j - 1}
\end{equation}
for each $j \ge 1$, which exists by continuity of addition on $B$ at
$0$.  If we put
\begin{equation}
\label{W_n = V_1 + V_2 + cdots + V_n}
        W_n = V_1 + V_2 + \cdots + V_n
\end{equation}
for each $n \in {\bf Z}_+$, then we get that
\begin{equation}
\label{W_n + V_n subseteq W_{n - 1} + V_{n - 1}}
        W_n + V_n \subseteq W_{n - 1} + V_{n - 1}
\end{equation}
when $n \ge 2$, by (\ref{V_j + V_j subseteq V_{j - 1}}).  This implies
that
\begin{equation}
\label{W_n subseteq W_n + V_n subseteq V_1 + V_1 subseteq V}
        W_n \subseteq W_n + V_n \subseteq V_1 + V_1 \subseteq V
\end{equation}
for every $n \ge 1$.  Suppose that $\phi_j$ is a continuous
homomorphism from $A_j$ into $B$ for each $j \in {\bf Z}_+$, and let
$U_j$ be an open set in $A_j$ such that $0 \in U_j$ and
\begin{equation}
\label{phi_j(U_j) subseteq V_j}
        \phi_j(U_j) \subseteq V_j
\end{equation}
for each $j \in {\bf Z}_+$.  Put
\begin{equation}
\label{U = (prod_{j = 1}^infty U_j) cap (sum_{j = 1}^infty A_j)}
 U = \Big(\prod_{j = 1}^\infty U_j\Big) \cap \Big(\sum_{j = 1}^\infty A_j\Big),
\end{equation}
which is the set of points in the direct sum $\sum_{j = 1}^\infty A_j$
whose $j$th component is in $U_j$ for each $j \in {\bf Z}_+$.  This is
a relatively open set in $\sum_{j = 1}^\infty A_j$ with respect to the
topology induced by the strong product topology on $\prod_{j =
  1}^\infty A_j$, by construction.  If $\phi$ is as in (\ref{phi(x) =
  sum_{j in I} phi_j(x_j)}), then it is easy to see that
\begin{equation}
\label{phi(U) subseteq V}
        \phi(U) \subseteq V,
\end{equation}
by (\ref{W_n subseteq W_n + V_n subseteq V_1 + V_1 subseteq V}) and
(\ref{phi_j(U_j) subseteq V_j}).  This implies that $\phi$ is
continuous on $\sum_{j = 1}^\infty A_j$ with respect to the topology
induced by the strong product topology on $\prod_{j = 1}^\infty A_j$,
as desired.

        If $I$ is uncountable, then this does not work, even when
$A_j = {\bf R}$ for each $j \in I$, and $B = {\bf R}$.  However, if
$A_j$ is equipped with the discrete topology for each $j \in I$,
then the corresponding strong product topology on $\prod_{j \in I} A_j$
is the same as the discrete topology, and there is no problem.

\section{Spaces of continuous functions}
\label{spaces of continuous functions}

        Let $X$ be a (nonempty) topological space, and let $C(X) = 
C(X, {\bf C})$ be the space of continuous complex-valued functions on $X$.
This is a vector space over the complex numbers with respect to
pointwise addition and scalar multiplication, and a commutative
algebra with respect to pointwise multiplication.  If $f$ is an
element of $C(X)$ and $K \subseteq X$ is compact, then $f(K)$ is a
compact subset of ${\bf C}$, and in particular the restriction of $f$
to $K$ is bounded.  Thus we can put
\begin{equation}
\label{||f||_K = sup_{x in K} |f(x)|}
        \|f\|_K = \sup_{x \in K} |f(x)|
\end{equation}
when $f \in C(X)$ and $K$ is a nonempty compact subset of $X$, and in fact
the supremum is attained.  It is easy to see that
\begin{equation}
\label{||f + g||_K le ||f||_K + ||g||_K}
        \|f + g\|_K \le \|f\|_K + \|g\|_K
\end{equation}
and
\begin{equation}
\label{||f g||_K le ||f||_K ||g||_K}
        \|f \, g\|_K \le \|f\|_K \, \|g\|_K
\end{equation}
for every $f, g \in C(X)$, and that
\begin{equation}
\label{||t f||_K = |t| ||f||_K}
        \|t \, f\|_K = |t| \, \|f\|_K
\end{equation}
for every $f \in C(X)$ and $t \in {\bf C}$.  This shows that $\|f\|_K$
defines a seminorm on $C(X)$, which is known as the \emph{supremum
  seminorm}\index{supremum seminorms} associated to $K$.  If $X$ is
compact, then we can take $K = X$, and the supremum seminorm $\|f\|_X$
is the same as the supremum norm on $C(X)$.

        If $K \subseteq X$ is nonempty and compact, then the open ball
in $C(X)$ centered at $f \in C(X)$ and with radius $r > 0$ with respect
to the supremum seminorm associated to $K$ is defined by
\begin{equation}
\label{B_K(f, r) = {g in C(X) : ||f - g||_K < r}}
        B_K(f, r) = \{g \in C(X) : \|f - g\|_K < r\}.
\end{equation}
Let us say that a set $U \subseteq C(X)$ is an open set if for each $f
\in U$ there is an $r > 0$ and a nonempty compact set $K \subseteq X$
such that
\begin{equation}
\label{B_K(f, r) subseteq U}
        B_K(f, r) \subseteq U.
\end{equation}
It is easy to see that this defines a topology on $C(X)$, which makes
$C(X)$ into a locally convex topological vector space over the complex
numbers.  Of course, finite subsets of $X$ are compact, which implies
that this topology is Hausdorff.  Note that the open balls
(\ref{B_K(f, r) = {g in C(X) : ||f - g||_K < r}}) are open sets with
respect to this topology, by the usual argument with the triangle
inequality.  One can also check that multiplication on $C(X)$ defines
a continuous mapping from $C(X) \times C(X)$ into $C(X)$, using the
corresponding product topology on $C(X) \times C(X)$.  Thus $C(X)$ is
actually a topological algebra with respect to this topology.  If $X$
is compact, then this is the same as the topology on $C(X)$ determined
by the supremum norm.  If $X$ is equipped with the discrete topology,
then every function on $X$ is continuous, and $C(X)$ is the same as
the Cartesian product of copies of ${\bf C}$ indexed by $X$.  Every
compact subset of $X$ has only finitely many elements in this case,
and the corresponding topology on $C(X)$ is the same as the product
topology.

        If $X$ is $\sigma$-compact, then there is a sequence $K_1, K_2, K_3, 
\ldots$ of nonempty compact subsets of $X$ such that $K_j \subseteq K_{j + 1}$
for each $j \ge 1$ and
\begin{equation}
\label{bigcup_{j = 1}^infty K_j = X}
        \bigcup_{j = 1}^\infty K_j = X.
\end{equation}
If $X$ is locally compact as well, then one can enlarge the $K_j$'s,
if necessary, to arrange that $K_j$ is contained in the interior of
$K_{j + 1}$ for each $j$.  This implies that each compact set $K
\subseteq X$ is contained in $K_j$ for some $j$, because the interiors
of the $K_j$'s form an open covering of $K$.  It follows that the
topology on $C(X)$ determined by the supremum seminorms associated to
compact subsets of $X$ can be described by the supremum seminorms
associated to the sequence of compact sets $K_j$.  If we put
\begin{equation}
\label{d_j(f, g) = min(||f - g||_{K_j}, 1/j)}
        d_j(f, g) = \min(\|f - g\|_{K_j}, 1/j)
\end{equation}
and
\begin{equation}
\label{d(f, g) = max_{j ge 1} d_j(f, g)}
        d(f, g) = \max_{j \ge 1} d_j(f, g)
\end{equation}
for each $f, g \in C(X)$, then one can check that $d(f, g)$ is a
translation-invariant metric on $C(X)$ that determines the same
topology on $C(X)$.

        Let $X$ be any nonempty topological space again, and let
$C(X, {\bf T})$ be the space of continuous mapping from $X$ into ${\bf T}$.
This is a commutative group with respect to pointwise multiplication
of functions, and a topological group with respect to the topology
induced by the one on $C(X)$ defined earlier.  Note that $C(X, {\bf
  T})$ is a closed set in $C(X)$ with respect to this topology,
because finite subsets of $X$ are compact.  Of course,
\begin{equation}
\label{||a f - a g||_K = ||f - g||_K}
        \|a \, f - a \, g\|_K = \|f - g\|_K
\end{equation}
for every $a \in C(X, {\bf T})$, $f, g \in C(X)$, and nonempty compact
set $K \subseteq X$.  If $X$ is compact, then we can take $K = X$, to
get that the supremum metric on $C(X)$ is invariant under
multiplication by elements of $C(X, {\bf T})$.  Similarly, if $X$ is
locally compact and $\sigma$-compact, then (\ref{d(f, g) = max_{j ge
    1} d_j(f, g)}) is invariant under multiplication by elements of
$C(X, {\bf T})$.  If $X$ is equipped with the discrete topology, then
$C(X, {\bf T})$ is the same as the direct product of copies of ${\bf
  T}$ indexed by $X$, which is compact with respect to the
corresponding product topology.

        Of course, $C(X)$ can also be considered as a vector space
over the real numbers, by forgetting about multiplication by $i$.  The
space $C(X, {\bf R})$ of real-valued continuous functions on $X$ may
be considered as a real-linear subspace of $C(X)$, and a topological
vector space over the real numbers with respect to the topology
induced by the one on $C(X)$ discussed earlier.  If $f(x)$ is a
real-valued continuous function on $X$, then $\exp(i \, f(x))$ is a
continuous function on $X$ with values in ${\bf T}$.  It is easy to
see that
\begin{equation}
\label{f(x) mapsto exp(i f(x))}
        f(x) \mapsto \exp(i \, f(x))
\end{equation}
defines a homomorphism from $C(X, {\bf R})$, as a commutative group
with respect to addition, into $C(X, {\bf T})$, as a commutative group
with respect to multiplication.  One can also check that (\ref{f(x)
  mapsto exp(i f(x))}) is continuous with respect to the corresponding
topologies induced by the one on $C(X)$.

        Remember that a subset $E$ of a commutative topological group
$A$ is totally bounded if for each open set $U$ in $A$ that contains
the identity element, $E$ can be covered by finitely many translates
of $U$.  If $A = C(X)$, then it suffices to consider open sets $U$ of
the form $B_K(0, r)$, where $K$ is a nonempty compact subset of $X$
and $r > 0$.  If $E \subseteq C(X, {\bf T})$, then one can check that
$E$ is totally bounded as a subset of $C(X)$, as a commutative
topological group with respect to addition, if and only if $E$ is
totally bounded as a subset of $C(X, {\bf T})$, as a commutative
topological group with respect to multiplication.

        Note that a sequence or net $\{f_j\}_j$ of elements of $C(X)$
converges to $f \in C(X)$ with respect to the topology on $C(X)$
considered in this section if and only if $\{f_j\}_j$ converges to $f$
uniformly on compact subsets of $X$.  If $\{f_j\}_j$ is a sequence or
net of continuous complex-valued functions on $X$ that converges to a
complex-valued function $f$ on $X$ uniformly on compact subsets of
$X$, then the restriction of $f$ to every compact subset of $X$ is
continuous, by standard arguments.  If $X$ is locally compact, then it
follows that $f$ is continuous on $X$.

        If $\{x_l\}_{l = 1}^\infty$ is a sequence of elements of any
topological space $X$ that converges to an element $x$ of $X$, then
the set $K \subseteq X$ consisting of the $x_l$'s, $l \in {\bf Z}_+$,
and $x$ is compact.  If $f$ is a complex-valued function on $X$ whose
restriction to each compact subset of $X$ is continuous, then it
follows that $f$ is sequentially continuous on $X$.  If $f$ is
sequentially continuous at a point $p \in X$, and if there is a local
base for the topology of $X$ at $p$ with only finitely or countably
many elements, then $f$ is continuous at $p$.  This gives another
criterion for the continuity of limits of sequences and nets of
functions on $X$ as in the preceding paragraph.

\section{The dual group}
\label{dual group}

        If $A$ is a commutative topological group, then the corresponding
\emph{dual group}\index{dual groups} $\widehat{A}$ is defined to be the
group $\hom(A, {\bf T})$ of continuous homomorphisms from $A$ into
${\bf T}$.  This is a subgroup of the group $C(A, {\bf T})$ of
continuous mappings from $A$ into ${\bf T}$, and hence a topological
group with respect to the topology induced on $\widehat{A}$ by the one
described in the previous section.  It is easy to see that
$\widehat{A}$ is also a closed set in $C(A)$ with respect to this
topology, because finite subsets of $A$ are compact.  In particular,
if $A$ is equipped with the discrete topology, then we have seen that
$C(A, {\bf T})$ is compact, which implies that $\widehat{A}$ is compact
as well.
        
        Suppose that $A = {\bf R}$, as a commutative group with respect
to addition, and equipped with the standard topology.  If $a \in {\bf R}$,
then (\ref{x mapsto exp(i a x)}) defines a continuous homomorphism from
${\bf R}$ into ${\bf T}$, and every continuous homomorphism from ${\bf R}$
into ${\bf T}$ is of this form, as in Section \ref{continuous homomorphisms}.
This defines a group isomorphism between ${\bf R}$ and its dual group.
One can check that this isomorphism is also a homeomorphism with respect
to the topology on the dual group induced by the one on $C({\bf R})$
described in the previous section.

        If $A$ is a compact commutative topological group, then the 
topology induced on $\widehat{A}$ by the one on $C(A)$ described in
the previous section is the same as the discrete topology.  To see
this, let $\phi_1, \phi_2 \in \widehat{A}$ be given, and put
\begin{equation}
\label{psi(x) = phi_1(x) phi_2(x)^{-1}}
        \psi(x) = \phi_1(x) \, \phi_2(x)^{-1},
\end{equation}
which defines another element of $\widehat{A}$.  If
\begin{equation}
\label{|phi_1(x) - phi_2(x)| < sqrt{2}}
        |\phi_1(x) - \phi_2(x)| < \sqrt{2}
\end{equation}
for every $x \in A$, then
\begin{equation}
\label{|psi(x) - 1| < sqrt{2}}
        |\psi(x) - 1| < \sqrt{2}
\end{equation}
for every $x \in A$, and hence $\re \psi(x) > 0$ for every $x \in A$.
Thus $\psi(A)$ is a subgroup of ${\bf T}$ contained in (\ref{{z in {bf
      T} : re z > 0}}), which implies that $\psi(A)$ is the trivial
subgroup of ${\bf T}$, as before.  It follows that $\psi(x) = 1$ for
every $x \in A$, so that $\phi_1(x) = \phi_2(x)$ for every $x \in A$
under these conditions.

        Let $A$ be any commutative topological group again, and let $U$
be an open subset of $A$ that contains $0$.  It is easy to see that
\begin{equation}
\label{{phi in widehat{A} : |phi(x) - 1| le 1 for every x in overline{U}}}
        \{\phi \in \widehat{A} : |\phi(x) - 1| \le 1
                                   \hbox{ for every } x \in \overline{U}\}
\end{equation}
is a closed set in $\widehat{A}$ with respect to the topology induced
by the usual one on $C(A)$.  One can also show that (\ref{{phi in
    widehat{A} : |phi(x) - 1| le 1 for every x in overline{U}}}) is a
compact set in $\widehat{A}$.  To do this, one of the main points is
to verify that elements of (\ref{{phi in widehat{A} : |phi(x) - 1| le
    1 for every x in overline{U}}}) are equicontinuous at $0$.  If $A$
is locally compact, then one can choose $U$ so that $\overline{U}$ is
a compact set in $A$.  In this case,
\begin{equation}
\label{{phi in widehat{A} : sup_{x in overline{U}} |phi(x) - 1| < 1}}
 \Big\{\phi \in \widehat{A} : \sup_{x \in \overline{U}} |\phi(x) - 1| < 1\Big\}
\end{equation}
is an open set in $\widehat{A}$, which is contained in (\ref{{phi in
    widehat{A} : |phi(x) - 1| le 1 for every x in overline{U}}}).
This implies that $\widehat{A}$ is locally compact when $A$ is locally
compact.

        If $A$ is any commutative topological group, then the group
$\hom(A, {\bf R})$ may be considered as a subgroup of $C(A)$ with
respect to addition.  In particular, $\hom(A, {\bf R})$ is a
topological group with respect to the topology induced by the usual one
on $C(A)$.  If $\psi \in \hom(A, {\bf R})$, then $\exp(i \, \psi(a))$
defines a continuous homomorphism from $A$ into ${\bf T}$.  As in the
previous section,
\begin{equation}
\label{psi(a) mapsto exp(i psi(a))}
        \psi(a) \mapsto \exp(i \, \psi(a))
\end{equation}
defines a homomorphism from $\hom(A, {\bf R})$ into $\hom(A, {\bf
  T})$.  This mapping is also continuous with respect to the
corresponding topologies induced by the one on $C(A)$, as before.

        Now let $V$ be a topological vector space over the real numbers.
The dual space $V'$ of continuous linear functionals on $V$ is a
subset of $C(V)$, although it is a bit nicer to think of $V'$ as a
linear subspace of the space $C(A, {\bf R})$ of real-valued continuous
functions on $V$.  At any rate, $V'$ is a topological vector space
with respect to the topology induced by the usual topology on $C(V)$.
In particular, $V'$ is a commutative topological group with respect to
addition and this topology.  Remember that every continuous homomorphism
from $V$ as a commutative topological group with respect to addition
into ${\bf R}$ is linear, as in Section \ref{continuous homomorphisms}.

        If $\lambda \in V'$, then $\exp(i \, \lambda(v))$ defines a
continuous homomorphism from $V$ as a commutative topological group
with respect to addition into ${\bf T}$.  In this situation, every
continuous homomorphism from $V$ into ${\bf T}$ is of this form, as in
Section \ref{continuous homomorphisms}.  It is easy to see that this
representation of continuous homomorphisms from $V$ into ${\bf T}$ is
also unique.  Thus
\begin{equation}
\label{lambda(v) mapsto exp(i lambda(v))}
        \lambda(v) \mapsto \exp(i \, \lambda(v))
\end{equation}
defines a group isomorphism from $V' = \hom(V, {\bf R})$, as a
commutative group with respect to addition, onto the dual group
$\widehat{V} = \hom(V, {\bf T})$ associated to $V$ as a commutative
topological group with respect to addition.  Of course, this mapping
is continuous with respect to the topologies induced by the usual one
on $C(V)$, as before.

        Let $K$ be a nonempty compact subset of $V$, and put
\begin{equation}
\label{K_1 = {t v : t in {bf R}, 0 le t le 1, v in K}}
        K_1 = \{t \, v : t \in {\bf R}, \ 0 \le t \le 1, \ v \in K\}.
\end{equation}
It is easy to see that this is a compact subset of $V$ as well,
because of continuity of scalar multiplication.  If $\lambda \in V'$
and $\exp(i \, \lambda(v))$ is uniformly close to $1$ on $K_1$, then
$\lambda(v)$ is uniformly close to $0$ modulo $2 \, \pi \, {\bf Z}$ on
$K_1$.  This implies that $\lambda(v)$ is uniformly close to $0$ on
$K_1$, because $\lambda$ is linear and $K_1$ is star-like about $0$,
by construction.  Using this, one can check that (\ref{lambda(v)
  mapsto exp(i lambda(v))}) is a homeomorphism from $V'$ onto
$\widehat{V}$, with respect to the corresponding topologies induced by
the usual one on $C(V)$.

\section{Compact and discrete groups}
\label{compact, discrete groups}

        Let $A$ be a compact commutative topological group, and let
$H$ be Haar measure on $A$, normalized so that $H(A) = 1$.  If
$\phi \in \widehat{A}$, then
\begin{equation}
\label{int_A phi(x) dH(x) = ... = phi(a) int_A phi(x) dH(x)}
        \int_A \phi(x) \, dH(x) = \int_A \phi(x + a) \, dH(x)
                                = \phi(a) \, \int_A \phi(x) \, dH(x)
\end{equation}
for every $a \in A$, using translation-invariance of Haar measure
in the first step.  This implies that
\begin{equation}
\label{int_A phi(x) dH(x) = 0}
        \int_A \phi(x) \, dH(x) = 0
\end{equation}
when $\phi(a) \ne 1$ for some $a \in A$, which is to say that $\phi$
is nontrivial on $A$.  If $\phi$, $\psi$ are distinct elements of
$\widehat{A}$, then
\begin{equation}
\label{phi(a) overline{psi(a)} = phi(a) psi(a)^{-1}}
        \phi(a) \, \overline{\psi(a)} = \phi(a) \, \psi(a)^{-1}
\end{equation}
is a nontrivial element of $\widehat{A}$.  In this case, it follows that
\begin{equation}
\label{int_A phi(x) overline{psi(x)} dH(x) = 0}
        \int_A \phi(x) \, \overline{\psi(x)} \, dH(x) = 0,
\end{equation}
by applying (\ref{int_A phi(x) dH(x) = 0}) to (\ref{phi(a)
  overline{psi(a)} = phi(a) psi(a)^{-1}}).

        Let $L^2(A)$ be the usual Hilbert space of complex-valued
square-integrable functions on $A$ with respect to $H$, with the
inner product
\begin{equation}
\label{langle f, g rangle = int_A f(x) overline{g(x)} dH(x)}
        \langle f, g \rangle = \int_A f(x) \, \overline{g(x)} \, dH(x).
\end{equation}
Thus distinct elements of $\widehat{A}$ are orthogonal with respect to
this inner product, by the remarks in the previous paragraph.  Note
that
\begin{equation}
\label{langle phi, phi rangle = int_A |phi(x)|^2 dH(x) = 1}
        \langle \phi, \phi \rangle = \int_A |\phi(x)|^2 \, dH(x) = 1
\end{equation}
for every $\phi \in \widehat{A}$, because of the normalization $H(A) =
1$.  If the linear span of the elements of $\widehat{A}$ is dense in
$L^2(A)$, then it follows that the elements of $\widehat{A}$ form an
orthonormal basis for $L^2(A)$.

        Let $\mathcal{E}(A)$ be the linear span of the elements of
$\widehat{A}$, as a linear subspace of $C(A)$.  More precisely,
$\mathcal{E}(A)$ is a subalgebra of $C(A)$ that contains the constant
functions and is invariant under complex conjugation.  If the elements
of $\widehat{A}$ separate points in $A$, then $\mathcal{E}(A)$ separates
points in $A$ as well, and hence $\mathcal{E}(A)$ is dense in $C(A)$
with respect to the supremum norm, by the Stone--Weierstrass theorem.
This implies that $\mathcal{E}(A)$ is dense in $L^2(A)$, because $C(A)$
is dense in $L^2(A)$, since Haar measure is a regular Borel measure on $A$.
It is well known that $\widehat{A}$ separates points on $A$ when $A$ is
compact, or even locally compact.

        Let $B$ be a subset of $\widehat{A}$, and let $\mathcal{E}_B(A)$
be the linear span of the elements of $B$, as a linear subspace of $C(A)$.
If $B$ is a subgroup of $\widehat{A}$, then it is easy to see that
$\mathcal{E}_B(A)$ is a subalgebra of $C(A)$ that contains the constant
functions and is invariant under complex conjugation.  If $B$ is a
subgroup of $\widehat{A}$ that separates points on $A$, then it follows
that $\mathcal{E}_B(A)$ is dense in $C(A)$, by the Stone--Weierstrass theorem.
As before, this implies that $\mathcal{E}_B(A)$ is dense in $L^2(A)$, so that
the elements of $B$ form an orthonormal basis for $L^2(A)$.  Under these
conditions, we get that
\begin{equation}
\label{B = widehat{A}}
        B = \widehat{A},
\end{equation}
since any element of $\widehat{A}$ not in $B$ would have to be orthogonal
to every element of $B$, and hence to every element of $\mathcal{E}_B(A)$.

        Now let $A$ be any commutative group, equipped with the discrete
topology.  Also let $A_0$ be a subgroup of $A$, let $a_1$ be an
element of $A$ not in $A_0$, and let $A_1$ be the subgroup of $A$
generated by $A_0$ and $a_1$.  Thus every element of $A_1$ can be
expressed as
\begin{equation}
\label{j a_1 + x}
        j \, a_1 + x
\end{equation}
for some $j \in {\bf Z}$ and $x \in A_0$.  Observe that
\begin{equation}
\label{{j in {bf Z} : j a_1 in A_0}}
        \{j \in {\bf Z} : j \, a_1 \in A_0\}
\end{equation}
is a subgroup of ${\bf Z}$, which is either trivial or generated by a
positive integer $j_1$.  If (\ref{{j in {bf Z} : j a_1 in A_0}}) is
trivial, then every element of $A_1$ can be expressed as (\ref{j a_1 +
  x}) for a unique $j \in {\bf Z}$ and $x \in A_0$.

        Let $\phi$ be a homomorphism from $A_0$ into ${\bf T}$, which
we would like to extend to a homomorphism from $A_1$ into ${\bf T}$.
The basic idea is to first choose $\phi(a_1) \in {\bf T}$, and then to
try to define $\phi$ on $A_1$ by putting
\begin{equation}
\label{phi(j a_1 + x) = phi(a_1)^j phi(x)}
        \phi(j \, a_1 + x) = \phi(a_1)^j \, \phi(x)
\end{equation}
for every $j \in {\bf Z}$ and $x \in A_0$.  This is easy to do when
(\ref{{j in {bf Z} : j a_1 in A_0}}) is trivial, in which case one can
take $\phi(a_1)$ to be any element of ${\bf T}$.  Otherwise, if
$j_1 \in {\bf Z}_+$ generates (\ref{{j in {bf Z} : j a_1 in A_0}}),
then one should choose $\phi(a_1)$ so that
\begin{equation}
\label{phi(a_1)^{j_1} = phi(j_1 a_1)}
        \phi(a_1)^{j_1} = \phi(j_1 \, a_1),
\end{equation}
where $\phi(j_1 \, a_1)$ is defined because $j_1 \, a_1 \in A_0$.
Under these conditions, one can check that (\ref{phi(j a_1 + x) =
  phi(a_1)^j phi(x)}) determines a well-defined extension of $\phi$
to a homomorphism from $A_1$ into ${\bf T}$, as desired.

        If $A$ is generated by $A_0$ and finitely or countably many
additional elements of $A$, then one can repeat the process to extend
$\phi$ to a homomorphism from $A$ into ${\bf T}$.  Otherwise, the
existence of such an extension can be obtained using Zorn's lemma or
the Hausdorff maximality principle.  In particular, if $a_0 \in A$ and
$a_0 \ne 0$, and if $A_0$ is the subgroup of $A$ generated by $a_0$,
then it is easy to see that there is a homomorphism $\phi$ from $A_0$
into ${\bf T}$ such that $\phi(a_0) \ne 1$, which has an extension to
a homomorphism from $A$ into ${\bf T}$.  If $a, b \in A$ and $a \ne b$,
then one can apply this to $a_0 = a - b$, to get a homomorphism $\phi$
from $A$ into ${\bf T}$ such that $\phi(a) \ne \phi(b)$.

\section{The second dual}
\label{second dual}

        Let $X$ be a (nonempty) topological space, and let $C(X)$ be the
space of complex-valued continuous functions on $X$, equipped with the
topology discussed in Section \ref{spaces of continuous functions}.
If $x \in X$, then
\begin{equation}
\label{Psi_x(f) = f(x)}
        \Psi_x(f) = f(x)
\end{equation}
defines a linear functional on $C(X)$, which is an algebra homomorphism
with respect to multiplication.  It is easy to see that this is also a
continuous function on $C(X)$ for each $x \in X$, because $K = \{x\}$
is a compact subset of $X$.  Thus
\begin{equation}
\label{x mapsto Psi_x}
        x \mapsto \Psi_x
\end{equation}
may be considered as a mapping from $X$ into the space $C(C(X))$ of
continuous complex-valued functions on $C(X)$.  More precisely, this
is a mapping from $X$ into the dual space $C(X)'$ of continuous linear
functionals on $C(X)$, as a topological vector space over the complex
numbers.

        Using the topology already defined on $C(X)$, we can get a
topology on $C(C(X))$ in the same way.  We would like to look at the
continuity properties of (\ref{x mapsto Psi_x}), as a mapping into
$C(C(X))$ equipped with this topology.  To do this, let $E$ be a
nonempty compact subset of $C(X)$, and let $C(E)$ be the space of
continuous complex-valued functions on $E$, with the topology
determined by the supremum norm.  We can also think of (\ref{x mapsto
  Psi_x}) as a mapping into $C(E)$, by restricting $\Psi_x(f)$ to $f
\in E$.  Continuity of (\ref{x mapsto Psi_x}) as a mapping into
$C(C(X))$ is equivalent to continuity of (\ref{x mapsto Psi_x}) as a
mapping into $C(E)$ for each nonempty compact set $E \subseteq C(X)$,
by the definition of this topology on $C(C(X))$.

        Let $E$ be a nonempty compact subset of $C(X)$ again, and let
$K$ be a nonempty compact subset of $X$.  Also let $C(K)$ be the space
of complex-valued continuous functions on $K$, with the topology
determined by the supremum norm.  There is a natural mapping $R_K$
from $C(X)$ into $C(K)$, which sends a continuous function $f$ on $X$
to its restriction to $K$.  It is easy to see that $R_K$ is continuous
with respect to the corresponding topologies on $C(X)$ and $C(K)$, by
construction.  This implies that $R_K(E)$ is a compact subset of
$C(K)$.  In particular, $R_K(E)$ is totally bounded in $C(K)$, with
respect to the supremum norm.  This holds when $E$ is totally bounded
in $C(X)$ as well.

        Under these conditions, (\ref{x mapsto Psi_x}) defines a
continuous mapping from $x \in K$ into $C(E)$, using the topology
induced on $K$ by the given topology on $X$.  This is easy to see when
$E$ has only finitely many elements, because each of those elements is
a continuous function on $X$.  Otherwise, one can basically reduce to
this case when $R_K(E)$ is totally bounded in $C(K)$.  It follows that
(\ref{x mapsto Psi_x}) is continuous as a mapping from $x \in K$ into
$C(C(X))$, since this works for every nonempty compact set $E
\subseteq C(X)$.  If $X$ is locally compact, then we get that (\ref{x
  mapsto Psi_x}) is continuous as a mapping from $X$ into $C(C(X))$.
If $X$ is any topological space, then the continuity of (\ref{x mapsto
  Psi_x}) on compact subsets of $X$ implies that (\ref{x mapsto
  Psi_x}) is sequentially continuous, because any convergent sequence
in $X$ and limit of that sequence is contained in a compact set in
$X$.  This implies that (\ref{x mapsto Psi_x}) is continuous as a
mapping from $X$ into $C(C(X))$ at a point $p \in X$ when there is a
local base for the topology of $X$ at $p$ with only finitely or
countably many elements.

        Now let $A$ be a commutative topological group, and remember
that the dual group $\widehat{A}$ is also a commutative topological
group with respect to the topology on $\widehat{A}$ induced by the
usual topology on $C(A)$.  This permits us to define the second dual
$\widehat{\widehat{A}}$ to be the dual of $\widehat{A}$, consisting
of the continuous homomorphisms from $\widehat{A}$ into ${\bf T}$.
Put
\begin{equation}
\label{Psi_a(phi) = phi(a)}
        \Psi_a(\phi) = \phi(a)
\end{equation}
for each $a \in A$ and $\phi \in \widehat{A}$, which is the same as
the function $\Psi_a$ on $C(A)$ discussed earlier, restricted to
$\widehat{A} \subseteq C(A)$.  This defines a homomorphism from
$\widehat{A}$ into ${\bf T}$ for each $a \in A$, where the group
operations on $\widehat{A}$ and ${\bf T}$ are defined by
multiplication.  As before, $\Psi_a$ is continuous as a mapping from
$\widehat{A}$ into ${\bf T}$ for each $a \in A$, because $K = \{a\}$
is a compact subset of $A$.

        Thus $\Psi_a$ may be considered as an element of the second
dual group $\widehat{\widehat{A}}$ for each $a \in A$.  Observe that
\begin{equation}
\label{Psi_{a + b}(phi) = phi(a + b) = phi(a) phi(b) = Psi_a(phi) Psi_b(phi)}
        \Psi_{a + b}(\phi) = \phi(a + b) = \phi(a) \, \phi(b)
                                        = \Psi_a(\phi) \, \Psi_b(\phi)
\end{equation}
for every $a, b \in A$ and $\phi \in \widehat{A}$.  This implies that
\begin{equation}
\label{a mapsto Psi_a}
        a \mapsto \Psi_a
\end{equation}
defines a homomorphism from $A$ into $\widehat{\widehat{A}}$, where
the group operation on $\widehat{\widehat{A}}$ is defined in terms of
multiplication, as usual.  Remember that $\widehat{\widehat{A}}$ is a
commutative topological group with respect to the topology induced on
$\widehat{\widehat{A}}$ by the usual one on $C(\widehat{A})$, which
involves the nonempty compact subsets of $\widehat{A}$.  We would like
to consider the continuity properties of (\ref{a mapsto Psi_a}) with
respect to this topology on $\widehat{\widehat{A}}$.

        Of course, $\widehat{A} \subseteq C(A)$, which leads to a natural
restriction mapping from $C(C(A))$ into $C(\widehat{A})$.  It is easy
to see that this restriction mapping is continuous with respect to the
usual topologies, since compact subsets of $\widehat{A}$ are compact
in $C(A)$ too.  If $K$ is a compact subset of $A$, then it follows
from the earlier discussion of arbitrary topological spaces that
(\ref{a mapsto Psi_a}) is continuous as a mapping from $a \in K$ into
$\widehat{\widehat{A}}$, with respect to the topology on $K$ induced
by the given topology on $A$.  In particular, (\ref{a mapsto Psi_a})
is continuous as a mapping from $A$ into $\widehat{\widehat{A}}$ when
$A$ is locally compact, and it is sequentially continuous for any $A$.
If $A$ is locally compact, then it is well known that (\ref{a mapsto
  Psi_a}) is actually a group isomorphism from $A$ onto
$\widehat{\widehat{A}}$, and a homeomorphism with respect to the
corresponding topologies.

        Note that (\ref{a mapsto Psi_a}) is injective as a mapping
from $A$ into $\widehat{\widehat{A}}$ if and only if $\widehat{A}$
separates points in $A$.  If $\phi$ and $\psi$ are distinct elements
of $\widehat{A}$, then $\phi(a) \ne \psi(a)$ for some $a \in A$,
and hence
\begin{equation}
\label{Psi_a(phi) ne Psi_a(psi)}
        \Psi_a(\phi) \ne \Psi_a(\psi).
\end{equation}
This implies that
\begin{equation}
\label{{Psi_a : a in A}}
        \{\Psi_a : a \in A\}
\end{equation}
separates points in $\widehat{A}$, and (\ref{{Psi_a : a in A}}) is
also a subgroup of $\widehat{\widehat{A}}$.  Suppose that $A$ is
equipped with the discrete topology, so that $\widehat{A}$ is compact.
In this case, the fact that (\ref{{Psi_a : a in A}}) is equal to
$\widehat{\widehat{A}}$ may be obtained as in (\ref{B = widehat{A}}).
We also saw that $\widehat{A}$ separates points in $A$ in the previous
section, so that (\ref{a mapsto Psi_a}) is injective.  Because
$\widehat{A}$ is compact, the usual topology on
$\widehat{\widehat{A}}$ is discrete as well, which implies that
(\ref{a mapsto Psi_a}) and its inverse are automatically continuous in
this case.

\section{Quotient groups}
\label{quotient groups}

        Let $A$ be a commutative topological group, and let $B$ be a
subgroup of $A$.  Thus the quotient group $A / B$ can be defined in the
usual way, with the canonical quotient homomorphism $q$ from $A$ onto
$A / B$, whose kernel is equal to $B$.  The corresponding
\emph{quotient topology}\index{quotient topology} on $A / B$ is
defined by saying that $V \subseteq A / B$ is an open set if and only if
$q^{-1}(V)$ is an open set in $A$.  Equivalently, $E \subseteq A / B$
is a closed set if and only if $q^{-1}(E)$ is a closed set in $A$.
It is easy to see that this does define a topology on $A / B$,
and $q$ is automatically continuous as a mapping from $A$ onto $A / B$
with respect to this topology on $A / B$.  If $U$ is any subset of $A$,
then
\begin{equation}
\label{q^{-1}(q(U)) = bigcup_{b in B} (U + b)}
        q^{-1}(q(U)) = \bigcup_{b \in B} (U + b),
\end{equation}
which is also denoted $U + B$.  If $U$ is an open set in $A$, then $U
+ b$ is an open set in $A$ for every $b \in B$, and hence
(\ref{q^{-1}(q(U)) = bigcup_{b in B} (U + b)}) is an open set in $A$
too.  This implies that $q(U)$ is an open set in $A / B$ with respect
to the quotient topology for every open set $U \subseteq A$, so that
$q$ is an open mapping.

        It is easy to see that translations on $A / B$ are continuous
with respect to the quotient topology, because of the corresponding
property of $A$.  One can also check that addition on $A / B$ is
continuous as a mapping from $(A / B) \times (A / B)$ into $A / B$
with respect to the quotient topology, and using the associated
product topology on $(A / B) \times (A / B)$.  Because of continuity
of translations, it suffices to verify continuity of addition on $A /
B$ at $0$, again using the analogous property of $A$.  Similarly,
$x \mapsto -x$ is continuous on $A$, which implies that $A / B$
has the same property.

        In order for $A / B$ to be considered a topological group,
$\{0\}$ should be a closed set in $A / B$ with respect to the quotient
topology.  In this situation, this is equivalent to $B$ being a closed
subgroup of $A$, and we suppose now that this is the case.  As in
Section \ref{basic notions}, this is implies that $A / B$ is
Hausdorff, and also regular as a topological space, although one could
look these conditions more directly in terms of $A$ as well.  If $A$
is locally compact, then there is an open set $U \subseteq A$ and a
compact set $K \subseteq U$ such that $0 \in U$ and $U \subseteq K$.
This implies that $q(U)$ is an open set in $A / B$ that contains $0$,
and that $q(K)$ is a compact set in $A / B$ that contains $q(U)$,
so that $A / B$ is locally compact too.

        Suppose that $C$ is another commutative topological group,
and that $\phi$ is a continuous homomorphism from $A / B$ into $C$,
with respect to the quotient topology on $A / B$.  This implies that
\begin{equation}
\label{phi' = phi circ q}
        \phi' = \phi \circ q,
\end{equation}
is a continuous homomorphism from $A$ into $C$ whose kernel contains
$B$.  If $\phi'$ is any homomorphism from $A$ into $C$ whose kernel
contains $B$, then $\phi'$ can be expressed as (\ref{phi' = phi circ
  q}) for some homomorphism $\phi$ from $A / B$ into $C$.  If $\phi'$
is continuous, then it is easy to see that $\phi$ is continuous with
respect to the quotient topology on $A / B$.

        Similarly, let $A_1$ be another commutative topological group,
and let $h$ be a continuous homomorphism from $A$ into $A_1$.  If $\phi$
is a continuous homomorphism from $A_1$ into $C$, then
\begin{equation}
\label{phi' = phi circ h}
        \phi' = \phi \circ h
\end{equation}
is a continuous homomorphism from $A$ into $C$, as before.  Thus
\begin{equation}
\label{phi mapsto phi'}
        \phi \mapsto \phi'
\end{equation}
defines a mapping from $\hom(A_1, C)$ into $\hom(A, C)$, which is a
homomorphism.  Note that the image of $\hom(A_1, C)$ under (\ref{phi
  mapsto phi'}) is contained in the subgroup of $\hom(A, C)$
consisting of continuous homomorphisms from $A$ into $C$ whose kernels
contain the kernel of $h$ in $A$.

        It is easy to see that (\ref{phi mapsto phi'}) is injective
when $h(A)$ is dense in $A_1$, because a continuous function is
uniquely determined by its restriction to a dense set.  Otherwise,
the kernel of (\ref{phi mapsto phi'}) consists of $\phi \in \hom(A_1, C)$
whose restriction to $h(A)$ is trivial, which is equivalent to asking
that the restriction of $\phi$ to the closure $\overline{h(A)}$ in $A_1$
be trivial, by continuity.  As before, this is the same as saying that
$\phi$ can be expressed as the composition of the quotient mapping
from $A_1$ onto $A_1 / \overline{h(A)}$ with a continuous homomorphism
from $A_1 / \overline{h(A)}$ into $C$.

        In particular, one can take $C = {\bf T}$, so that
(\ref{phi mapsto phi'}) defines a homomorphism $\widehat{h}$
from $\widehat{A_1}$ into $\widehat{A}$.  This is known as the \emph{dual
  homomorphism}\index{dual homomorphisms} associated to $h$.  One can
check that $\widehat{h}$ is continuous with respect to the usual
topologies on the dual groups, basically because $h$ sends compact
subsets of $A$ to compact subsets of $A_1$.

        Suppose that $A$ is locally compact, so that there is an open
set $U \subseteq A$ and a compact set $K \subseteq A$ such that $0 \in
U$ and $U \subseteq K$.  If $E \subseteq A / B$ is compact, then $E$
is contained in the union of finitely many translates of $q(U)$, and
hence $E$ is contained in the union of finitely many translates of
$q(K)$.  This implies that $E$ is contained in the image of the union
of finitely many translates of $K$ in $A$ under $q$.  In particular,
$E$ is contained in the image of a compact subset of $A$ under $q$.
Using this, one can check that $\widehat{q}$ is a homeomorphism from
$\widehat{(A / B)}$ onto its image in $\widehat{A}$, with respect to
the topology induced on the image of $\widehat{(A / B)}$ by the usual
topology on $\widehat{A}$.

        We can also think of $B$ as a topological group, with respect
to the topology induced by the given topology on $A$.  The natural
inclusion mapping from $B$ into $A$ is a continuous homomorphism,
which leads to a dual homomorphism from $\widehat{A}$ into
$\widehat{B}$.  This dual homomorphism is the same as the restriction
mapping, which sends a continuous homomorphism from $A$ into ${\bf T}$
to its restriction to $B$.

        If $A$ is equipped with the discrete topology, then every
homomorphism from $B$ into ${\bf T}$ can be extended to a homomorphism
from $A$ into ${\bf T}$, as in Section \ref{compact, discrete groups}.
Similarly, if $B$ is compact, and if the elements of $\widehat{A}$ separate
points in $B$, then it follows that every element of $\widehat{B}$ can be
expressed as the restriction to $B$ of an element of $\widehat{A}$,
by another argument from Section \ref{compact, discrete groups}.
It is well known that if $A$ is locally compact and $B$ is a closed subgroup
of $A$, then every element of $\widehat{B}$ can be expressed as the
restriction to $B$ of an element of $\widehat{A}$.  In this case, the
restriction mapping from $\widehat{A}$ onto $\widehat{B}$ is topologically
equivalent to a quotient mapping.

        Suppose that $A$ is compact, and that $\phi'$ is a continuous
homomorphism from $A$ onto a commutative topological group $C$ with
kernel equal to $B$.  As before, $\phi'$ can be expressed as
(\ref{phi' = phi circ q}), where $\phi$ is a continuous homomorphism
from $A / B$ onto $C$.  In this case, the kernel of $\phi$ is trivial,
which means that $\phi$ is injective.  It is well known that a
one-to-one continuous mapping from a compact topological space onto a
Hausdorff topological space is always a homeomorphism.  It follows
that $\phi$ is a homeomorphism under these conditions, because $A / B$
is compact when $A$ is compact.

\section{Uniform continuity and equicontinuity}
\label{uniform continuity, equicontinuity}

        Let $A$ be a commutative topological group, and let $f$ be a
complex-valued function on $A$.  Let us say that $f$ is
\emph{uniformly continuous}\index{uniform continuity} along a set $E
\subseteq A$ if for each $\epsilon > 0$ there is an open set $U
\subseteq A$ such that $0 \in U$ and
\begin{equation}
\label{|f(x) - f(y)| < epsilon}
        |f(x) - f(y)| < \epsilon
\end{equation}
for every $x \in E$ and $y \in A$ with $y - x \in U$.  If $E = A$,
then we simply say that $f$ is uniformly continuous on $A$.  If
$d(\cdot, \cdot)$ is a translation-invariant metric on $A$ that
determines the same topology on $A$, then this type of uniform
continuity condition is equivalent to the usual one for metric spaces.

        It is well known that continuous functions on compact metric spaces
are uniformly continuous.  Similarly, if $E \subseteq A$ is compact,
and if $f$ is continuous as a function on $A$ at every point in $E$,
then $f$ is uniformly continuous along $E$.  To see this, let
$\epsilon > 0$ be given, and for each $p \in E$, let $V(p)$ be an open
subset of $A$ such that $0 \in V(p)$ and
\begin{equation}
\label{|f(z) - f(p)| < epsilon/2}
        |f(z) - f(p)| < \epsilon/2
\end{equation}
for every $z \in p + V(p)$.  The continuity of addition on $A$ at $0$
implies that for each $p \in E$ there is an open set $U(p) \subseteq A$
that contains $0$ and satisfies
\begin{equation}
\label{U(p) + U(p) subseteq V(p)}
        U(p) + U(p) \subseteq V(p).
\end{equation}
The collection of open sets $p + U(p)$ with $p \in E$ covers $E$, and
hence there are finitely many elements $p_1, \ldots, p_n$ of $E$ such that
\begin{equation}
\label{E subseteq bigcup_{j = 1}^n (p_j + U(p_j))}
        E \subseteq \bigcup_{j = 1}^n (p_j + U(p_j)),
\end{equation}
because $E$ is compact.  Put
\begin{equation}
\label{U = bigcap_{j = 1}^n U(p_j)}
        U = \bigcap_{j = 1}^n U(p_j),
\end{equation}
so that $U$ is also an open set in $A$ that contains $0$, and let us
check that (\ref{|f(x) - f(y)| < epsilon}) holds for every $x \in E$
and $y \in x + U$.  If $x \in E$, then $x \in p_j + U(p_j)$ for some
$j = 1, \ldots, n$, by (\ref{E subseteq bigcup_{j = 1}^n (p_j +
  U(p_j))}).  If $y \in x + U$, then
\begin{equation}
\label{y in (p_j + U(p_j)) + U subseteq ... subseteq p_j + V(p_j)}
        y \in (p_j + U(p_j)) + U \subseteq p_j + U(p_j) + U(p_j)
                                  \subseteq p_j + V(p_j),
\end{equation}
using the fact that $U \subseteq U(p_j)$ in the second step, and
(\ref{U(p) + U(p) subseteq V(p)}) in the third step.  Of course,
$U(p_j) \subseteq V(p_j)$, because of (\ref{U(p) + U(p) subseteq
  V(p)}), so that $x \in p_j + V(p_j)$ as well.  Thus we can apply
(\ref{|f(z) - f(p)| < epsilon/2}) with $p = p_j$ to $z = x$ and $z =
y$, to get that
\begin{equation}
\label{|f(x) - f(y)| le |f(x) - f(p_j)| + |f(p_j) - f(y)| < ... = epsilon}
        |f(x) - f(y)| \le |f(x) - f(p_j)| + |f(p_j) - f(y)|
                         < \epsilon/2 + \epsilon/2 = \epsilon,
\end{equation}
as desired.

        In particular, if $A$ is compact, then it follows that every
continuous function on $A$ is uniformly continuous.  Similarly, if $A$
is locally compact, then continuous functions on $A$ with compact
support are uniformly continuous on $A$.  Remember that a continuous
function $f$ on $A$ is said to \emph{vanish at
  infinity}\index{vanishing at infinity} if for each $\epsilon > 0$
there is a compact set $K(\epsilon) \subseteq A$ such that
\begin{equation}
\label{|f(x)| < epsilon}
        |f(x)| < \epsilon
\end{equation}
for every $x \in A \setminus K(\epsilon)$.  It is easy to see that $f$
is uniformly continuous on $A$ in this case too, because $f$ is
uniformly continuous along $K(\epsilon)$ for every $\epsilon > 0$.

        A collection $\mathcal{E}$ of complex-valued functions on $A$
is said to be \emph{equicontinuous}\index{equicontinuity} at a point
$x \in A$ if for each $\epsilon > 0$ there is an open set $U \subseteq
A$ such that $0 \in U$ and (\ref{|f(x) - f(y)| < epsilon}) holds for
every $f \in \mathcal{E}$ and $y \in x + U$.  Similarly, $\mathcal{E}$
is \emph{uniformly equicontinuous}\index{uniform equicontinuity} along
a set $E \subseteq A$ if for every $\epsilon > 0$ there is an open set
$U \subseteq A$ such that $0 \in U$ and (\ref{|f(x) - f(y)| <
  epsilon}) holds for every $f \in \mathcal{E}$, $x \in E$, and $y \in
A$ with $y - x \in U$.  If $\mathcal{E}$ is equicontinuous at every
point $x \in E$ and $E \subseteq A$ is compact, then $\mathcal{E}$ is
uniformly equicontinuous along $E$, by essentially the same argument
as before.

        Suppose that $f$ is a homomorphism from $A$ into ${\bf C}$,
where ${\bf C}$ is considered as a commutative group with respect to
addition.  If $f$ is continuous at $0$, then it is easy to see that
$f$ is uniformly continuous on $A$.  If $\mathcal{E}$ is a collection
of homomorphisms from $A$ into ${\bf C}$ that is equicontinuous at
$0$, then $\mathcal{E}$ is uniformly equicontinuous on $A$.  In fact,
it suffices to ask that there be an open set $U \subseteq A$ such that
$0 \in U$ and
\begin{equation}
\label{|f(x)| < 1}
        |f(x)| < 1
\end{equation}
for every $x \in U$ and $f \in \mathcal{E}$.  This implies that each
$f \in \mathcal{E}$ is continuous at $0$, as in Section
\ref{continuous homomorphisms}, and the same argument shows that
$\mathcal{E}$ is equicontinuous at $0$ under these conditions.

        Suppose now that $f$ is a homomorphism from $A$ into ${\bf T}$,
where ${\bf T}$ is considered as a commutative group with respect to
multiplication.  If $f$ is continuous at $0$, then one can check that
$f$ is uniformly continuous as a complex-valeud function on $A$.
As before, if $\mathcal{E}$ is a collection of homomorphisms from $A$
into ${\bf T}$ that is equicontinuous at $0$, then $\mathcal{E}$ is
uniformly equicontinuous on $A$.  In this situation, it suffices to
ask that there be an open set $U \subseteq A$ such that $0 \in U$ and
\begin{equation}
\label{re f(x) > 0}
        \re f(x) > 0
\end{equation}
for every $x \in U$ and $f \in \mathcal{E}$.  This implies that each
$f \in \mathcal{E}$ is continuous at $0$, as in Section
\ref{continuous homomorphisms}, and essentially the same argument
shows that $\mathcal{E}$ is equicontinuous at $0$.

        Let $X$ be a (nonempty) topological space, and consider $A = C(X)$
as a commutative topological group with respect to pointwise addition
of functions and the usual topology.  Also let $K$ be a nonempty
compact subset of $X$, and put
\begin{equation}
\label{mathcal{E}_K = {Psi_x : x in K}}
        \mathcal{E}_K = \{\Psi_x : x \in K\},
\end{equation}
where $\Psi_x(f) = f(x)$ for each $f \in C(X)$, as in Section
\ref{second dual}.  Thus $\mathcal{E}_K$ is a collection of continuous
linear functionals on $C(X)$, and it is easy to see that
$\mathcal{E}_K$ is equicontinuous at $0$ on $C(X)$, by definition of
the usual topology on $C(X)$.  Alternatively, one can consider $A =
C(X, {\bf T})$ as a commutative topological group with respect to
pointwise multiplication of functions and the topology induced by the
usual one on $C(X)$.  In this case, $\mathcal{E}_K$ corresponds to a
collection of continuous homomorphisms from $C(X, {\bf T})$ into ${\bf
  T}$, which is automatically equicontinuous at the identity element.

\section{Direct products, revisited}
\label{direct products, revisited}

        Let $A_1, \ldots, A_n$ be finitely many commutative topological
groups, and consider their Cartesian product
\begin{equation}
\label{A = prod_{j = 1}^n A_j}
        A = \prod_{j = 1}^n A_j,
\end{equation}
as a commutative topological group with respect to coordinatewise
addition and the product topology.  As in Section \ref{sums, products,
  continued}, the dual group $\widehat{A}$ can be identified with the
product of the dual groups $\widehat{A_j}$, $j = 1, \ldots, n$.  It is
easy to see that the usual topology on $\widehat{A}$, induced by the
topology defined earlier on $C(A)$, corresponds exactly to the product
topology on $\prod_{j = 1}^n \widehat{A_j}$, where $\widehat{A_j}$ has
the analogous topology for each $j$.  More precisely, if $K_j$ is a
compact subset of $A_j$ for $j = 1, \ldots, n$, then
\begin{equation}
\label{K = prod_{j = 1}^n K_j}
        K = \prod_{j = 1}^n K_j
\end{equation}
is a compact subset of $A$.  Conversely, if $E$ is any compact subset
of $A$, then the projection of $E$ in $A_j$ is compact for each $j$,
which implies that $E$ is contained in a product of compact sets.
Thus the usual topology on $C(A)$ can be defined in terms of nonempty
compact subsets of $A$ of the form (\ref{K = prod_{j = 1}^n K_j}).
Using this, one can check that the induced topology on $\widehat{A}$
corresponds to the product topology on $\prod_{j = 1}^n
\widehat{A_j}$.  There are analogous statements for $\hom(A, {\bf R})$
instead of $\widehat{A}$.

        Now let $I$ be a nonempty set, let $A_j$ be a commutative
topological group for each $j \in I$, and let
\begin{equation}
\label{A = prod_{j in I} A_j}
        A = \prod_{j \in I} A_j
\end{equation}
be the direct product of the $A_j$'s, equipped with the product
topology.  If $K_j$ is a compact subset of $A_j$ for each $j \in I$,
then
\begin{equation}
\label{K = prod_{j in I} K_j}
        K = \prod_{j \in I} K_j
\end{equation}
is a compact subset of $A$, by Tychonoff's theorem.  As before, every
compact subset of $A$ is contained in a product of compact subsets of
the $A_j$'s, so that the usual topology on $C(A)$ can be defined in
terms of nonempty compact subsets of $A$ of the form (\ref{K = prod_{j
    in I} K_j}).  We may as well suppose that $0 \in K_j$ for each $j
\in I$, since otherwise we can replace $K_j$ with $K_j \cup \{0\}$ for
each $j$.  Similarly, we may as well ask that
\begin{equation}
\label{K_j = -K_j}
        K_j = -K_j
\end{equation}
for each $j \in I$, since otherwise we can replace $K_j$ with $K_j
\cup (-K_j)$ for each $j$.

        As in Section \ref{sums, products, continued}, $\hom(A, {\bf R})$
can be identified with the direct sum
\begin{equation}
\label{sum_{j in I} hom(A_j, {bf R}}
        \sum_{j \in I} \hom(A_j, {\bf R}),
\end{equation}
because any subgroup of ${\bf R}$ contained in a bounded neighborhood
of $0$ is trivial.  More precisely, if $\phi$ is an element of
$\hom(A, {\bf R})$, then $\phi$ can be expressed as
\begin{equation}
\label{phi(x) = sum_{l = 1}^n phi_{j_l}(x_{j_l})}
        \phi(x) = \sum_{l = 1}^n \phi_{j_l}(x_{j_l}),
\end{equation}
where $j_1, \ldots, j_n$ are finitely many elements of $I$, and
$\phi_{j_l}$ is an element of $\hom(A_{j_l}, {\bf R})$ for $l = 1,
\ldots, n$.  Note that $\hom(A_j, {\bf R})$ may be considered as a
vector space over the real numbers with respect to pointwise addition
and scalar multiplication for each $j \in I$, which is a real-linear
subspace of $C(A_j)$.  Similarly, $\hom(A, {\bf R})$ is a real-linear
subspace of $C(A)$, and $\hom(A, {\bf R})$ can be identified with the
direct sum (\ref{sum_{j in I} hom(A_j, {bf R}}) as a sum of vector
spaces.  We can also identify $C(A_j)$ with the subspace of $C(A)$
consisting of functions of $x \in A$ that depend only on $x_j \in
A_j$ for each $j \in I$, so that the direct sum
\begin{equation}
\label{sum_{j in I} C(A_j)}
        \sum_{j \in I} C(A_j)
\end{equation}
corresponds to a linear subspace of $C(A)$ as well.

        Let $K_j$ be a compact subset of $A_j$ with $0 \in K_j$ for each
$j \in I$, and let $K$ be the product of the $K_j$'s, as in
(\ref{K = prod_{j in I} K_j}).  If $\|f\|_K$ is the supremum seminorm
on $C(A)$ associated to $K$ and $\phi$ is as in (\ref{phi(x) = sum_{l
    = 1}^n phi_{j_l}(x_{j_l})}), then
\begin{equation}
\label{||phi||_K = ... = sup_{x in K} |sum_{l = 1}^n phi_{j_l}(x_{j_l})|}
        \|\phi\|_K = \sup_{x \in K} |\phi(x)|
             = \sup_{x \in K} \, \biggl|\sum_{l = 1}^n \phi_{j_l}(x_{j_l})\biggr|.
\end{equation}
Observe that
\begin{equation}
\label{widetilde{K} = {x in K : x_j ne 0 for at most one j in I}}
 \widetilde{K} = \{x \in K : x_j \ne 0 \hbox{ for at most one } j \in I\},
\end{equation}
is also a compact subset of $A$, because it can be expressed as the
intersection of $K$ with a closed set in $A$.  Of course, $\|\phi\|_K$
is greater than or equal to
\begin{equation}
\label{||phi||_{widetilde{K}} = ...}
        \|\phi\|_{\widetilde{K}} 
            = \max_{1 \le l \le n} \, \sup_{x_{j_l} \in K_{j_l}} |\phi_{j_l}(x_{j_l})|,
\end{equation}
because $\widetilde{K} \subseteq K$.  If $K_j = -K_j$ for each $j \in
I$, then (\ref{||phi||_K = ... = sup_{x in K} |sum_{l = 1}^n
  phi_{j_l}(x_{j_l})|}) reduces to
\begin{equation}
\label{||phi||_K = sup_{x in K} sum_{l = 1}^n |phi(x_{j_l})|}
        \|\phi\|_K = \sup_{x \in K} \, \sum_{l = 1}^n |\phi(x_{j_l})|,
\end{equation}
since we can replace $x_{j_l}$ with $-x_{j_l}$ for any $l = 1, \ldots,
n$, to get $\phi(x_{j_l}) \ge 0$.

        If $B$ is a commutative topological group and $r \in {\bf Z}_+$,
then it is easy to see that $b \mapsto r \cdot b$ defines a continuous
mapping on $B$, because of contnuity of addition on $B$.  Let $r_j$ be
a positive integer for each $j \in I$, and let $r_j \cdot K_j$ be the
subset of $A_j$ consisting of multiples of elements of $K_j$ by $r_j$.
Thus $r_j \cdot K_j$ is a compact subset of $A_j$ for each $j \in I$,
so that
\begin{equation}
\label{K' = prod_{j in I} (r_j cdot K_j)}
        K' = \prod_{j \in I} (r_j \cdot K_j)
\end{equation}
is a compact set in $A$.  If $\widetilde{K}'$ corresponds to $K'$ as
in (\ref{widetilde{K} = {x in K : x_j ne 0 for at most one j in I}}),
then
\begin{equation}
\label{||phi||_{widetilde{K}'} = ...}
        \|\phi\|_{\widetilde{K}'}
 = \max_{1 \le l \le n} \, \sup_{x_{j_l} \in K_{j_l}} r_{j_l} \, |\phi_{j_l}(x_{j_l})|,
\end{equation}
as in (\ref{||phi||_{widetilde{K}} = ...}).  If $K_j = -K_j$ for each
$j \in I$, then we get that
\begin{equation}
\label{||phi||_{K'} = sup_{y in K'} sum_{l = 1}^n |phi_{j_l}(y_{j_l})| = ...}
        \|\phi\|_{K'} = \sup_{y \in K'} \, \sum_{l = 1}^n |\phi_{j_l}(y_{j_l})|
          = \sup_{x \in K} \, \sum_{l = 1}^n r_{j_l} \, |\phi_{j_l}(x_{j_l})|,
\end{equation}
as in (\ref{||phi||_K = sup_{x in K} sum_{l = 1}^n |phi(x_{j_l})|}).

        The topology on $\hom(A, {\bf R})$ induced by the usual topology
on $C(A)$ may be described by the seminorms (\ref{||phi||_K = ... =
  sup_{x in K} |sum_{l = 1}^n phi_{j_l}(x_{j_l})|}), or equivalently
by the seminorms (\ref{||phi||_K = sup_{x in K} sum_{l = 1}^n
  |phi(x_{j_l})|}), where $K \subseteq A$ is as before.  The seminorms
(\ref{||phi||_{widetilde{K}} = ...}) also determine a topology on
$\hom(A, {\bf R})$, and every open set in $\hom(A, {\bf R})$ with
respect to the seminorms (\ref{||phi||_{widetilde{K}} = ...}) is an
open set with respect to the topology on $\hom(A, {\bf R})$ determined
by the seminorms (\ref{||phi||_K = ... = sup_{x in K} |sum_{l = 1}^n
  phi_{j_l}(x_{j_l})|}), because (\ref{||phi||_{widetilde{K}} = ...})
is less than or equal to (\ref{||phi||_K = ... = sup_{x in K} |sum_{l
    = 1}^n phi_{j_l}(x_{j_l})|}).  Of course, these topologies are the
same when $I$ has only finitely many elements, and one can show that
they are the same when $I$ is countably infinite as well.  More
precisely, if $I$ is countable, then one can choose positive integers
$r_j$ for $j \in I$ such that
\begin{equation}
\label{sum_{j in I} 1/r_j < infty}
        \sum_{j \in I} 1/r_j < \infty.
\end{equation}
Using this, one can estimate (\ref{||phi||_K = sup_{x in K} sum_{l =
    1}^n |phi(x_{j_l})|}) in terms of (\ref{||phi||_{widetilde{K}'} =
  ...}), to show that the two topologies on $\hom(A, {\bf R})$ are the
same in this case.

        Remember that we can identify $\hom(A, {\bf R})$ with the direct sum
(\ref{sum_{j in I} hom(A_j, {bf R}}) for any $I$, which may be considered as
a linear subspace of the Cartesian product
\begin{equation}
\label{prod_{j in I} hom(A_j, {bf R})}
        \prod_{j \in I} \hom(A_j, {\bf R}).
\end{equation}
Using (\ref{||phi||_{widetilde{K}'} = ...}), one can check that the
topology determined on $\hom(A, {\bf R})$ by the seminorms
(\ref{||phi||_{widetilde{K}} = ...}) is the same as the topology
induced on $\hom(A, {\bf R})$ by the strong product topology on
(\ref{prod_{j in I} hom(A_j, {bf R})}) associated to the topology on
$\hom(A_j, {\bf R})$ induced by the usual topology on $C(A_j)$ for
each $j \in I$.

        There is an analogous discussion for $\widehat{A} = \hom(A, {\bf T})$,
with some additional complications.  One of the main points is that
instead of using subsets of $A_j$ of the form $r_j \cdot K_j$ as
before, it is better to take finite unions of multiples of $K_j$.
This is to ensure that if a homomorphism $\phi_j$ from $A_j$ into
${\bf T}$ is reasonably close to $1$ on the larger set, then it will
be as close to $1$ as we want on $K_j$.

        Of course, if $A_j$ is compact, then $\hom(A_j, {\bf R})$ is
trivial, and the topology induced on $\widehat{A_j}$ by the usual
topology on $C(A_j)$ is the discrete topology.  If $A_j$ is compact
for each $j \in I$, then $A$ is compact with respect to the product
topology, $\hom(A, {\bf R})$ is trivial, and the topology induced on
$\widehat{A}$ by the usual topology on $C(A)$ is the discrete
topology.  The strong product topology on
\begin{equation}
\label{prod_{j in I} widehat{A_j}}
        \prod_{j \in I} \widehat{A_j}
\end{equation}
corresponding to the discrete topology on $\widehat{A_j}$ for each $j
\in I$ is the same as the discrete topology on (\ref{prod_{j in I}
  widehat{A_j}}).  Thus the topology induced on
\begin{equation}
\label{sum_{j in I} widehat{A_j}}
        \sum_{j \in I} \widehat{A_j}
\end{equation}
as a subgroup of (\ref{prod_{j in I} widehat{A_j}}) by the strong
product topology on (\ref{prod_{j in I} widehat{A_j}}) is the discrete
topology as well.  In this case, we do not need to be concerned with
other topologies on (\ref{sum_{j in I} widehat{A_j}}).

\section{Direct sums, revisited}
\label{direct sums, revisited}

        Let $I$ be a nonempty set again, and let $A_j$ be a commutative
topological group for each $j \in I$.  Put
\begin{equation}
\label{A = sum_{j in I} A_j}
        A = \sum_{j \in I} A_j,
\end{equation}
which may be considered as a subgroup of the Cartesian product
\begin{equation}
\label{prod_{j in I} A_j, 2}
        \prod_{j \in I} A_j.
\end{equation}
If $I$ has only finitely many elements, then $A$ is the same as
(\ref{prod_{j in I} A_j}), and we are back to the situation discussed
in the previous section.  In this section, we shall be especially
interested in the topology induced on $A$ by the strong product
topology on (\ref{prod_{j in I} A_j, 2}) corresponding to the given
topology on $A_j$ for each $j \in I$, and perhaps some stronger
topologies on $A$ as well.  Remember that the strong product topology
on (\ref{prod_{j in I} A_j, 2}) is the same as the discrete topology
when $A_j$ is equipped with the discrete topology for each $j$.

        Let $E$ be a subset of $A$, and put
\begin{equation}
\label{I(E) = {j in I : there is an x in E such that x_j ne 0}}
 I(E) = \{j \in I : \hbox{there is an } x \in E \hbox{ such that } x_j \ne 0\}.
\end{equation}
If $j \in I(E)$, then let $U_j$ be an open set in $A_j$ such that $0
\in U_j$ and there is an $x \in E$ with $x_j \not\in U_j$.  Otherwise,
put $U_j = A_j$ when $j \in I \setminus I(E)$.  Thus
\begin{equation}
\label{U = prod_{j in I} U_j}
        U = \prod_{j \in I} U_j
\end{equation}
is an open set in (\ref{prod_{j in I} A_j, 2}) with respect to the
strong product topology, which implies that
\begin{equation}
\label{V = A cap U}
        V = A \cap U
\end{equation}
is a relatively open set in $A$ with respect to the topology induced
on $A$ by the strong product topology on (\ref{prod_{j in I} A_j, 2}).
Note that $0 \in U_j$ for every $j \in I$, by construction, so that $0
\in V \subseteq U$.

        We also have that
\begin{equation}
\label{a + U = prod_{j in I} (a_j + U_j)}
        a + U = \prod_{j \in I} (a_j + U_j)
\end{equation}
for every element $a$ of (\ref{prod_{j in I} A_j, 2}), and that
\begin{equation}
\label{a + V = A cap (a + U)}
        a + V = A \cap (a + U)
\end{equation}
when $a \in A$.  If $a_j = 0$ and $j \in I(E)$, then there is an $x
\in E$ such that $x_j \not\in U_j = a_j + U_j$, so that $x \not\in a +
U$, and hence $x \not\in a + V$.  Remember that $a_j = 0$ for all but
finitely many $j \in I$ when $a \in A$, by the definition of $A$.  If
$I(E)$ has infinitely many elements, then it follows that $E$ cannot
be covered by finitely many translates of $V$ in $A$.

        If $E$ is totally bounded in $A$ with respect to the topology
induced by the strong product topology on (\ref{prod_{j in I} A_j,
  2}), then $E$ can be covered by finitely many translates of $V$ in
$A$.  Thus $I(E)$ can have only finitely many elements when $E$ is
totally bounded in $A$ with respect to this topology, by the argument
in the preceding paragraph.  In particular, this holds when $E$ is
compact with respect to this topology on $A$.  Of course, this also
holds when $E$ is totally bounded or compact with respect to any
stronger topology on $A$.

        If $I_1$ is a nonempty subset of $I$, then we can identify
\begin{equation}
\label{prod_{j in I_1} A_j}
        \prod_{j \in I_1} A_j
\end{equation}
with the subgroup of (\ref{prod_{j in I} A_j, 2}) consisting of points
whose $j$th coordinate is equal to $0$ when $j \in I \setminus I_1$.
The strong product topology on (\ref{prod_{j in I_1} A_j}) associated
to the given topologies on the $A_j$'s corresponds exactly to the
topology induced on this subgroup of (\ref{prod_{j in I} A_j, 2}) by
the strong product topology there.  Similarly,
\begin{equation}
\label{sum_{j in I_1} A_j}
        \sum_{j \in I_1} A_j
\end{equation}
can be identified with a subgroup of (\ref{A = sum_{j in I} A_j}), and
the topology induced on (\ref{sum_{j in I_1} A_j}) by the strong
product topology on (\ref{prod_{j in I_1} A_j}) corresponds to the one
induced by the topology on (\ref{A = sum_{j in I} A_j}) that is
induced by the strong product topology on (\ref{prod_{j in I} A_j,
  2}).  If $I_1$ has only finitely many elements, then (\ref{sum_{j in
    I_1} A_j}) is the same as (\ref{prod_{j in I_1} A_j}), and the
strong product topology on (\ref{prod_{j in I_1} A_j}) reduces to the
product topology.

        In particular, we can identify $E \subseteq A$ with a subset of
(\ref{sum_{j in I_1} A_j}) when $I(E) \subseteq I_1$.  If $E$ is totally
bounded or compact with respect to the topology induced on $A$ by the
strong product topology on (\ref{A = sum_{j in I} A_j}), then the
corresponding subset of (\ref{sum_{j in I_1} A_j}) will have the analogous
property, since we have already seen that the topologies match up.
Conversely, if a subset of (\ref{sum_{j in I_1} A_j}) is totally
bounded or compact with respect to the topology induced by the strong
product topology on (\ref{prod_{j in I_1} A_j}), then it will correspond
to a subset of $A$ with the analogous property.

        The Cartesian product of totally bounded or compact subsets of the
$A_j$'s with $j \in I_1$ is totally bounded or compact with respect to the
product topology on (\ref{prod_{j in I_1} A_j}), respectively.
Conversely, every totally bounded or compact subset of (\ref{prod_{j
    in I_1} A_j}) with respect to the product topology is contained in
a Cartesian product of totally bounded or compact subsets of the
$A_j$'s, as in the previous section.  If $I_1$ has only finitely many
elements, then the product topology on (\ref{prod_{j in I_1} A_j}) is
the same as the strong product topology, and (\ref{prod_{j in I_1}
  A_j}) is the same as (\ref{sum_{j in I_1} A_j}).

        The usual topology on $C(A)$ is defined in terms of the supremum
seminorms associated to nonempty compact subsets of $A$, where we
continue to take $A$ to be equipped with the topology induced by the
strong product topology on (\ref{prod_{j in I} A_j, 2}).  The
preceding discussion implies that compact subsets of $A$ with respect
to this topology correspond to compact subsets of sub-sums of the form
(\ref{sum_{j in I_1} A_j}), where $I_1 \subseteq I$ has only finitely
many elements.  Thus one can get the same topology on $C(A)$ by
considering nonempty compact subsets of $A$ that correspond to
Cartesian products of compact subsets of $A_j$ for finitely many $j
\in I$.

\chapter{Some harmonic analysis}
\label{some harmonic analysis}

\section{Fourier transforms}
\label{fourier transforms, 2}

        Let $A$ be a commutative topological group, and let $\mu$
be a complex Borel measure on $A$.  The \emph{Fourier
  transform}\index{Fourier transforms} of $\mu$ is the function
$\widehat{\mu}$ defined on the dual group $\widehat{A} = \hom(A, {\bf
  T})$ by
\begin{equation}
\label{widehat{mu}(phi) = int_A overline{phi(x)} d mu(x)}
        \widehat{\mu}(\phi) = \int_A \overline{\phi(x)} \, d\mu(x).
\end{equation}
Observe that
\begin{equation}
\label{|widehat{mu}(phi)| le int_A |phi(x)| d|mu|(x) = |mu|(A)}
        |\widehat{\mu}(\phi)| \le \int_A |\phi(x)| \, d|\mu|(x) = |\mu|(A)
\end{equation}
for every $\phi \in \widehat{A}$, where $|\mu|$ is the total variation
measure on $A$ associated to $\mu$.

        Suppose for the moment that $\mu$ is supported on a compact
set $K \subseteq A$, so that $|\mu|(A \setminus K) = 0$, and
\begin{equation}
\label{widehat{mu}(phi) = int_K overline{phi(x)} d mu(x)}
        \widehat{\mu}(\phi) = \int_K \overline{\phi(x)} \, d\mu(x)
\end{equation}
for every $\phi \in \widehat{A}$.  This implies that
\begin{eqnarray}
\label{|widehat{mu}(phi) - widehat{mu}(psi)| le ...}
        |\widehat{\mu}(\phi) - \widehat{\mu}(\psi)|
           & \le & \int_K |\phi(x) - \psi(x)| \, d|\mu|(x)    \\
  & \le & \Big(\sup_{x \in K} |\phi(x) - \psi(x)|\Big) \, |\mu|(A) \nonumber
\end{eqnarray}
for every $\phi, \psi \in \widehat{A}$.  It follows that
$\widehat{\mu}$ is continuous on $\widehat{A}$ with respect to the
topology induced by the usual one on $C(A)$.  Equivalently,
\begin{equation}
\label{|widehat{mu}(phi) - widehat{mu}(psi)| le ..., 2}
        |\widehat{\mu}(\phi) - \widehat{\mu}(\psi)|
 \le \Big(\sup_{x \in K} |\phi(x) \, \psi(x)^{-1} - 1|\Big) \, |\mu|(A)
\end{equation}
for every $\phi, \psi \in \widehat{A}$, which shows more explicitly
that $\widehat{\mu}$ is uniformly continuous on $\widehat{A}$ as a
commutative topological group with respect to pointwise multiplication
of functions and the topology induced by the usual topology on $C(A)$.

        Now suppose that for each $\epsilon > 0$ there is a compact set
$K(\epsilon) \subseteq A$ such that
\begin{equation}
\label{|mu|(A setminus K(epsilon)) < epsilon}
        |\mu|(A \setminus K(\epsilon)) < \epsilon,
\end{equation}
so that $\mu$ can be approximated with respect to the total variation
norm by complex Borel measures supported on compact subsets of $A$.
This implies that $\widehat{\mu}$ can be approximated uniformly on
$\widehat{A}$ by the Fourier transforms of complex Borel measures with
compact support on $A$, because of (\ref{|widehat{mu}(phi)| le int_A
  |phi(x)| d|mu|(x) = |mu|(A)}).  It follows that $\widehat{\mu}$ is
uniformly continuous on $\widehat{A}$ in this case too.  In
particular, this condition holds automatically when $A$ is
$\sigma$-compact.  This condition also holds when $A$ is locally
compact and $\mu$ is a regular Borel measure on $A$, which means that
$|\mu|$ is a regular Borel measure on $A$.

        Let us restrict our attention from now on in this section to
the case where $A$ is locally compact, and let $H$ be a choice of Haar
measure on $A$.  If $f(x)$ is a complex-valued integrable function on
$A$ with respect to $H$, then the Fourier transform\index{Fourier transforms}
of $f$ is the function $\widehat{f}$ defined on $\widehat{A}$ by
\begin{equation}
\label{widehat{f}(phi) = int_A f(x) overline{phi(x)} dH(x)}
        \widehat{f}(\phi) = \int_A f(x) \, \overline{\phi(x)} \, dH(x).
\end{equation}
As before,
\begin{equation}
\label{|widehat{f}(phi)| le int_A |f(x)| |phi(x)| dH(x) = int_A |f(x)| dH(x)}
        |\widehat{f}(\phi)| \le \int_A |f(x)| \, |\phi(x)| \, dH(x)
                                = \int_A |f(x)| \, dH(x)
\end{equation}
for every $\phi \in \widehat{A}$.  If we put
\begin{equation}
\label{mu(E) = int_E f(x) dH(x)}
        \mu(E) = \int_E f(x) \, dH(x)
\end{equation}
for every Borel set $E \subseteq A$, then $\mu$ is a complex Borel
measure on $A$, and the Fourier transform of $f$ on $A$ is the same as
the Fourier transform of $\mu$.  More precisely, $\mu$ is a regular
Borel measure on $A$ under these conditions, and hence $\widehat{f}$
is uniformly continuous on $\widehat{A}$.

        Let $L^1(A)$ be the usual vector space of complex-valued
integrable functions on $A$ with respect to $H$.  If $f \in L^1(A)$
and $a \in A$, then
\begin{equation}
\label{f_a(x) = f(x - a), 3}
        f_a(x) = f(x - a)
\end{equation}
also defines an element of $L^1(A)$.  Observe that
\begin{eqnarray}
\label{widehat{f_a}(phi) = ... = widehat{f}(phi) overline{phi(a)}}
 \widehat{f_a}(\phi) & = & \int_A f(x - a) \, \overline{\phi(x)} \, dH(x) \\
            & = & \int_A f(x) \, \overline{\phi(x + a)} \, dH(x) \nonumber \\
            & = & \widehat{f}(\phi) \, \overline{\phi(a)} \nonumber
\end{eqnarray}
for every $\phi \in \widehat{A}$.  This uses the translation-invariance
of Haar measure in the second step, and the fact that $\phi : A \to {\bf T}$
is a homomorphism in the third step.

        Suppose for the moment that $f$ is a continuous complex-valued
function on $A$ supported on a compact set $K \subseteq A$.  Also let
$U$ be an open set in $A$ such that $0 \in U$ and $\overline{U}$ is compact,
which exists because $A$ is supposed to be locally compact.  Thus
$K \times \overline{U}$ is a compact subset of $A \times A$ with respect
to the product topology associated to the given topology on $A$.
This implies that $K + \overline{U}$ is a compact set in $A$, because
addition on $A$ defines a continuous mapping from $A \times A$ into $A$.
Note that the support of $f_a$ is contained in $K + a$ for every $a \in A$,
which is contained in $K + \overline{U}$ when $a \in \overline{U}$.

        It follows that
\begin{eqnarray}
\lefteqn{\int_A |f(x) - f(x - a)| \, dH(x)} \\
         & \le & \Big(\sup_{x \in K + \overline{U}} |f(x) - f(x - a)|\Big)
                                      \, H(K + \overline{U}) \nonumber
\end{eqnarray}
for every $a \in \overline{U}$.  Of course, $H(K + \overline{U}) <
\infty$, because $K + \overline{U}$ is compact.  Using uniform
continuity as in Section \ref{uniform continuity, equicontinuity}, we
get that
\begin{equation}
\label{lim_{a to 0} int_A |f(x) - f(x - a)| dH(x) = 0}
        \lim_{a \to 0} \int_A |f(x) - f(x - a)| \, dH(x) = 0
\end{equation}
where the limit is taken over $a \in A$ converging to $0$ with respect
to the given topology on $A$.  It is well known that continuous
functions with compact support on $A$ are dense in $L^1(A)$, because
of the regularity properties of Haar measure.  This permits one to
extend (\ref{lim_{a to 0} int_A |f(x) - f(x - a)| dH(x) = 0}) to all
$f \in L^1(A)$.

        Let $f \in L^1(A)$ and $\epsilon > 0$ be given, and let $V$
be an open set in $A$ such that $0 \in V$ and
\begin{equation}
\label{int_A |f(x) - f(x - a)| dH(x) le epsilon}
        \int_A |f(x) - f(x - a)| \, dH(x) \le \epsilon
\end{equation}
for every $a \in \overline{V}$.  This implies that
\begin{equation}
\label{|widehat{f}(phi) - widehat{f_a}(phi)| le ... le epsilon}
 |\widehat{f}(\phi) - \widehat{f_a}(\phi)| \le \int_A |f(x) - f_a(x)| \, dH(x) 
                                            \le \epsilon
\end{equation}
for every $a \in \overline{V}$ and $\phi \in \widehat{A}$, as in
(\ref{|widehat{f}(phi)| le int_A |f(x)| |phi(x)| dH(x) = int_A |f(x)|
  dH(x)}).  Combining this with (\ref{widehat{f_a}(phi) = ... =
  widehat{f}(phi) overline{phi(a)}}), we get that
\begin{equation}
\label{|widehat{f}(phi)| |phi(a) - 1| le epsilon}
        |\widehat{f}(\phi)| \, |\phi(a) - 1| \le \epsilon
\end{equation}
for every $\phi \in \widehat{A}$ and $a \in \overline{V}$.  As in
Section \ref{dual group},
\begin{equation}
\label{{phi in widehat{A} : |phi(a) - 1| le 1 for every a in overline{V}}}
        \{\phi \in \widehat{A} : |\phi(a) - 1| \le 1 \hbox{ for every }
                                                         a \in \overline{V}\}
\end{equation}
is a compact set in $\widehat{A}$.  If $\phi \in \widehat{A}$ is not
an element of (\ref{{phi in widehat{A} : |phi(a) - 1| le 1 for every a
    in overline{V}}}), then
\begin{equation}
\label{|phi(a) - 1| > 1}
        |\phi(a) - 1| > 1
\end{equation}
for some $a \in \overline{V}$, and hence
\begin{equation}
\label{|widehat{f}(phi)| < epsilon}
        |\widehat{f}(\phi)| < \epsilon,
\end{equation}
by (\ref{|widehat{f}(phi)| |phi(a) - 1| le epsilon}).  This shows that
$\widehat{f}$ vanishes at infinity on $\widehat{A}$ for every $f \in
L^1(A)$.  Of course, this is vacuous when $A$ is equipped with the
discrete topology, so that $\widehat{A}$ is compact.

        As in Section \ref{fourier transforms}, the Dirac mass on $A$
associated to a point $a \in A$ is the Borel measure $\delta_a$
defined by $\delta_a(E) = 1$ when $a \in E$ and $\delta_a(E) = 0$
when $a \in A \setminus E$.  In this case,
\begin{equation}
\label{widehat{delta_a}(phi) = overline{phi(a)}}
        \widehat{\delta_a}(\phi) = \overline{\phi(a)}
\end{equation}
for every $\phi \in \widehat{A}$, and in particular
$|\widehat{\delta_a}(\phi)| = 1$.

\section{Convolution of measures}
\label{convolution of measures}

        Let $A$ be a commutative topological group again, and put
\begin{equation}
\label{E_2 = {(x, y) in A times A : x + y in E}}
        E_2 = \{(x, y) \in A \times A : x + y \in E\}
\end{equation}
for each subset $E$ of $A$, which is the same as the inverse image of
$E$ under the mapping from $A \times A$ into $A$ corresponding to
addition in $A$.  Thus $E_2$ is an open set in $A \times A$ with
respect to the product topology when $E$ is an open set in $A$, by
continuity of addition on $A$.  Similarly, $E_2$ is a closed set in $A
\times A$ when $A$ is a closed set in $A$.  It follows that $E_2$ is a
Borel set in $A \times A$ when $E$ is a Borel set in $A$.

        If $\mu$ and $\nu$ are complex Borel measures on $A$, then
we would like to define their \emph{convolution}\index{convolutions}
$\mu * \nu$ to be the complex Borel measure on $A$ given by
\begin{equation}
\label{(mu * nu)(E) = (mu times nu)(E_2)}
        (\mu * \nu)(E) = (\mu \times \nu)(E_2)
\end{equation}
for every Borel set $E \subseteq A$, where $\mu \times \nu$ is the
corresponding product measure on $A \times A$.  More precisely, in
order to use the standard product measure construction here, $E_2$
should be in the $\sigma$-algebra of subsets of $A \times A$ generated
by products of Borel subsets of $A$.  If there is a countable base for
the topology of $A$, then there is a countable base for the product
topology on $A \times A$, consisting of product of basic open subsets
of $A$.  In particular, this implies that every open subset of $A
\times A$ can be expressed as a countable union of products of open
subsets of $A$, so that open subsets of $A \times A$ are measurable
with respect to the standard product measure construction.  In this
case, Borel subsets of $A \times A$ are measurable with respect to the
standard product measure construction, as desired.

        Alternatively, suppose that $A$ is locally compact, and that
$\mu$ and $\nu$ are regular Borel measures on $A$.  Under these conditions,
one can use another version of the product measure construction, in which
$\mu \times \nu$ is a regular Borel measure on $A \times A$.  Thus the
Borel measurability of $E_2$ is sufficient in this situation, but one would
also like $\mu * \nu$ to be regular on $A$.  If $\mu$ and $\nu$ are
real-valued and nonnegative, then one can first check that $\mu * \nu$ is
inner regular, and then use that to get that $\mu * \nu$ is outer regular.
Otherwise, one can reduce to that case, by expressing $\mu$ and $\nu$
as linear combinations of nonnegative regular measures.

        At any rate, once one has $\mu * \nu$ as in
(\ref{(mu * nu)(E) = (mu times nu)(E_2)}), one gets that
\begin{equation}
\label{int_A f(z) d(mu * nu)(z) = int_{A times A} f(x + y) d(mu times nu)(x, y)}
        \int_A f(z) \, d(\mu * \nu)(z)
            = \int_{A \times A} f(x + y) \, d(\mu \times \nu)(x, y)
\end{equation}
for every bounded complex-valued Borel measurable function $f$ on $A$.
More precisely, this is the same as (\ref{(mu * nu)(E) = (mu times
  nu)(E_2)}) when $f$ is the characteristic or indicator function
associated to $E$ on $A$.  Thus (\ref{int_A f(z) d(mu * nu)(z) =
  int_{A times A} f(x + y) d(mu times nu)(x, y)}) holds when $f$ is a
Borel measurable simple function on $A \times A$, by linearity.  This
implies that (\ref{int_A f(z) d(mu * nu)(z) = int_{A times A} f(x + y)
  d(mu times nu)(x, y)}) holds for all bounded Borel measurable
functions $f$ on $A \times A$, by standard approximation arguments.
Equivalently,
\begin{eqnarray}
\label{int_A f(z) d(mu * nu)(z) = ...}
        \int_A f(z) \, d(\mu * \nu)(z) 
 & = & \int_A \Big(\int_A f(x + y) \, d\mu(x)\Big) \, d\nu(y) \\
 & = & \int_A \Big(\int_A f(x + y) \, d\nu(y)\Big) \, d\mu(x), \nonumber
\end{eqnarray}
by the appropriate version of Fubini's theorem.

        If $A$ is locally compact, then regular complex Borel measures on
$A$ correspond to bounded linear functionals on $C_0(A)$ with respect to
the supremum norm, by a version of the Riesz representation theorem.
In order to get $\mu \times \nu$, it suffices to define the
appropriate bounded linear functional on $C_0(A \times A)$.  One can
also look at $\mu * \nu$ in terms of the corresponding bounded linear
functional on $C_0(A)$.  This is much simpler when $A$ is compact, so
that $C_0(A) = C(A)$.  Otherwise, $f(x + y)$ need not vanish at
infinity on $A \times A$ even when $f$ has compact support in $A$,
and so it is better to at least be able to integrate bounded continuous
functions on $A \times A$.  Even if $A$ is not compact, there are some
simplifications when $\mu$ or $\nu$ has compact support in $A$.
One can then reduce to this case using suitable approximation arguments.

        Note that $A_2 = A \times A$, and
\begin{equation}
\label{(mu * nu)(A) = mu(A) nu(A)}
        (\mu * \nu)(A) = \mu(A) \, \nu(A).
\end{equation}
One can check that
\begin{equation}
\label{|mu * nu|(E) le (|mu| * |nu|)(E)}
        |\mu * \nu|(E) \le (|\mu| * |\nu|)(E)
\end{equation}
for every Borel set $E \subseteq A$, where $|\mu|$, $|\nu|$, and $|\mu
* \nu|$ are the total variation measures on $A$ associated to $\mu$,
$\nu$, and $\mu * \nu$, respectively, and $|\mu| * |\nu|$ is the
convoution of $|\mu|$ and $|\nu|$ on $A$.  In particular,
\begin{equation}
\label{|mu * nu|(A) le (|mu| * |nu|)(A) = |mu|(A) |nu|(A)}
        |\mu * \nu|(A) \le (|\mu| * |\nu|)(A) = |\mu|(A) \, |\nu|(A),
\end{equation}
so that the total variation norm of $\mu * \nu$ on $A$ is less than or
equal to the product of the total variation norms of $\mu$ and $\nu$
on $A$.  Remember that the total variation norm of a regular complex
Borel measure on $A$ is equal to the dual norm of the corresponding
bounded linear functional on $C_0(A)$, with respect to the supremum
norm on $C_0(A)$.

        It is easy to see that convolution of measures is commutative,
in the sense that
\begin{equation}
\label{mu * nu = nu * mu}
        \mu * \nu = \nu * \mu,
\end{equation}
basically because of commutativity of addition on $A$.  Similarly, one
can check that convolution of measures is associative.  If $\phi$ is a
continuous homomorphism from $A$ into ${\bf T}$, then we can apply
(\ref{int_A f(z) d(mu * nu)(z) = int_{A times A} f(x + y) d(mu times
  nu)(x, y)}) to $f = \overline{\phi}$ to get that
\begin{equation}
\label{widehat{(mu * nu)}(phi) = widehat{mu}(phi) widehat{nu}(phi)}
 \widehat{(\mu * \nu)}(\phi) = \widehat{\mu}(\phi) \, \widehat{\nu}(\phi).
\end{equation}

\section{Convolution of functions}
\label{convolution of functions}

        Let $A$ be a locally compact commutative topological group,
and let $H$ be a Haar measure on $A$.  If $f$, $g$ are real or
complex-valued Borel measurable functions on $A$, then we would like
to define their \emph{convolution}\index{convolutions} on $A$ by
\begin{equation}
\label{(f * g)(x) = int_A f(x - y) g(y) dH(y)}
        (f * g)(x) = \int_A f(x - y) \, g(y) \, dH(y).
\end{equation}
If $f$, $g$ are real-valued and nonnegative on $A$, then this integral
makes sense as a nonnegative extended real number for each $x \in A$.
Otherwise, (\ref{(f * g)(x) = int_A f(x - y) g(y) dH(y)}) is defined
when
\begin{equation}
\label{int_A |f(x - y)| |g(y)| dH(y)}
        \int_A |f(x - y)| \, |g(y)| \, dH(y)
\end{equation}
is finite, in which case the absolute value or modulus of (\ref{(f *
  g)(x) = int_A f(x - y) g(y) dH(y)}) is less than or (\ref{int_A |f(x
  - y)| |g(y)| dH(y)}).  It is easy to see that (\ref{(f * g)(x) =
  int_A f(x - y) g(y) dH(y)}) and (\ref{int_A |f(x - y)| |g(y)|
  dH(y)}) are symmetric in $f$ and $g$, using the change of variables
$y \mapsto x - y$.

        If $f$ and $g$ are real-valued and nonnegative again, then
one can check that
\begin{equation}
\label{int_A (f * g)(x) dH(x) = (int_A f(x) dH(x)) (int_A g(y) dH(y))}
        \int_A (f * g)(x) \, dH(x) = \Big(\int_A f(x) \, dH(x)\Big) \,
                                      \Big(\int_A g(y) \, dH(y)\Big),
\end{equation}
using Fubini's theorem.  More precisely, the Borel measurability of
$f(x - y) \, g(y)$ on $A \times A$ can be derived from the Borel
measurability of $f$, $g$ on $A$ and the continuity of the group
operations on $A$.  If there is a countable base for the topology of
$A$, then $A$ is $\sigma$-compact, which implies that $A$ is
$\sigma$-finite with respect to $H$.  Borel subsets of $A \times A$
are also measurable with respect to the usual product measure
construction in this situation, which permits one to use the stanard
version of Fubini's theorem.  Otherwise, one should use a version for
Borel measures with suitable regularity properties.

        If $f$ and $g$ are real or complex-valued functions on $A$
that are integrable with respect to $H$, then we can apply (\ref{int_A
  (f * g)(x) dH(x) = (int_A f(x) dH(x)) (int_A g(y) dH(y))}) to $|f|$
and $|g|$.  This implies in particular that (\ref{int_A |f(x - y)|
  |g(y)| dH(y)}) is finite for almost every $x \in A$ with respect to
$H$, so that (\ref{(f * g)(x) = int_A f(x - y) g(y) dH(y)}) is defined
for almost every $x \in A$ with respect to $H$.  We also get that
\begin{equation}
\label{int_A |(f * g)(x)| dH(x) le (int_A |f(x)| dH(x)) (int_A |g(y)| dH(y))}
 \int_A |(f * g)(x)| \, dH(x) \le \Big(\int_A |f(x)| \, dH(x)\Big) \,
                                   \Big(\int_A |g(y)| \, dH(y)\Big),
\end{equation}
because the absolute value or modulus of (\ref{(f * g)(x) = int_A f(x
  - y) g(y) dH(y)}) is less than or equal to (\ref{int_A |f(x - y)|
  |g(y)| dH(y)}), and using (\ref{int_A (f * g)(x) dH(x) = (int_A f(x)
  dH(x)) (int_A g(y) dH(y))}) applied to $|f|$ and $|g|$.  Thus $f *
g$ is also integrable on $A$ under these conditions, where
measurability comes from Fubini's theorem.

        If $\mu$ is the Borel measure on $A$ associated to $f$ as in
(\ref{mu(E) = int_E f(x) dH(x)}), and if $\nu$ is associated to $g$ in
the same way, then it is easy to see that $\mu * \nu$ is the Borel
measure corresponding to $f * g$.  It follows that the Fourier
transform of $f * g$ is equal to the product of the Fourier transforms
of $f$ and $g$, as in (\ref{widehat{(mu * nu)}(phi) = widehat{mu}(phi)
  widehat{nu}(phi)}).  Similarly, convolution is associative on
$L^1(A)$, which can also be verified more directly.

        Let $f$ and $g$ be real-valued and nonnegative again, and let
$1 \le p < \infty$ be given.  One can check that
\begin{eqnarray}
\label{(int_A (f * g)(x)^p dH(x))^{1/p} le ...}
\lefteqn{\Big(\int_A (f * g)(x)^p \, dH(x)\Big)^{1/p}} \\
        & \le & \Big(\int_A f(x)^p \, dH(x)\Big)^{1/p} \, 
                 \Big(\int_A g(y) \, dH(y)\Big), \nonumber
\end{eqnarray}
using the integral version of Minkowski's inequality, or Jensen's
inequality and Fubini's theorem.  If $f$ and $g$ are real or
complex-valued functions on $A$ such that $|f|^p$ and $|g|$ are
integrable, then one can apply (\ref{(int_A (f * g)(x)^p dH(x))^{1/p}
  le ...}) to $|f|$ and $|g|$, to get that (\ref{int_A |f(x - y)|
  |g(y)| dH(y)}) is finite for almost every $x \in A$ with respect to
$H$.  Thus (\ref{(f * g)(x) = int_A f(x - y) g(y) dH(y)}) is defined
for almost every $x \in A$ with respect to $H$, and
\begin{eqnarray}
\label{(int_A |(f * g)(x)|^p dH(x))^{1/p} le ...}
\lefteqn{\Big(\int_A |(f * g)(x)|^p \, dH(x)\Big)^{1/p}} \\
        & \le & \Big(\int_A |f(x)|^p \, dH(x)\Big)^{1/p} \,
                 \Big(\int_A |g(y)| \,dH(y)\Big). \nonumber
\end{eqnarray}
This shows that $f * g \in L^p(A)$ when $f \in L^p(A)$ and $g \in L^1(A)$.

        If $f \in L^\infty(A)$ and $g \in L^1(A)$, then
(\ref{int_A |f(x - y)| |g(y)| dH(y)}) is obviously finite for every
$x \in A$, so that (\ref{(f * g)(x) = int_A f(x - y) g(y) dH(y)})
is defined for every $x \in A$.  In this case, we get that
\begin{equation}
\label{|(f * g)(x)| le ||f||_infty int_A |g(y)| dH(y)}
        |(f * g)(x)| \le \|f\|_\infty \, \int_A |g(y)| \, dH(y)
\end{equation}
for every $x \in A$, where $\|f\|_\infty$ denotes the $L^\infty$ norm
of $f$ with respect to $H$.  This implies that $f * g$ is bounded on
$A$, and in fact $f * g$ is uniformly continuous on $A$ under these
conditions.  To see this, it is better to use the analogous expression
for $g * f$, which is the same as $f * g$, as mentioned earlier.  This
permits one to use the analogue of (\ref{lim_{a to 0} int_A |f(x) -
  f(x - a)| dH(x) = 0}) for $g$, to show that the convolution is
uniformly continuous.

        Suppose now that $f \in L^p(A)$ for some $p$, $1 < p < \infty$,
and that $g \in L^q(A)$, where $1 < q < \infty$ is the exponent
conjugate to $p$, so that $1/p + 1/q = 1$.  Thus (\ref{int_A |f(x -
  y)| |g(y)| dH(y)}) is finite for every $x \in A$, by H\"older's
inequality, which implies that (\ref{(f * g)(x) = int_A f(x - y) g(y)
  dH(y)}) is defined for every $x \in A$.  As before,
\begin{equation}
\label{|(f * g)(x)| le ||f|||_p ||g||_q}
        |(f * g)(x)| \le \|f\|_p \, \|g\|_q
\end{equation}
for every $x \in A$, where $\|f\|_p$ denotes the $L^p$ norm of $f$
with respect to $H$, and similarly for $\|g\|_q$.  It follows that $f
* g$ is bounded on $A$, and one can also show that $f * g$ is
uniformly continuous on $A$.  This uses the analogue of (\ref{lim_{a
    to 0} int_A |f(x) - f(x - a)| dH(x) = 0}) for the $L^p$ norm,
which can be obtained from the density of continuous functions on $A$
with compact support in $L^p(A)$ when $p < \infty$.  One could just as
well use the corresponding property for $g$, because $q < \infty$, and
the commutativity of convolution.  By approximating both $f$ and $g$
by functions with compact support in $A$, one can check that $f * g$
vanishes at infinity in this situation.

\section{Functions and measures}
\label{functions, measures}

        Let $A$ be a commutative topological group, and let $\nu$ be
a complex Borel measure on $A$.  If $f$ is a complex-valued Borel
measurable function on $A$, then we would like to define the
\emph{convolution}\index{convolutions} of $f$ and $\nu$ by
\begin{equation}
\label{(f * nu)(x) = int_A f(x - y) d nu(y)}
        (f * \nu)(x) = \int_A f(x - y) \, d\nu(y).
\end{equation}
If $f$ and $\nu$ are real-valued and nonnegative, then this integral
makes sense as a nonnegative extended real number for every $x \in A$.
Otherwise, (\ref{(f * nu)(x) = int_A f(x - y) d nu(y)}) is defined when
\begin{equation}
\label{int_A |f(x - y)| d|nu|(y)}
        \int_A |f(x - y)| \, d|\nu|(y)
\end{equation}
is finite, where $|\nu|$ is the total variation measure associated to
$\nu$ on $A$.  In this case, the modulus of (\ref{(f * nu)(x) = int_A
  f(x - y) d nu(y)}) is less than or equal to (\ref{int_A |f(x - y)|
  d|nu|(y)}).

        As a basic class of examples, let $\delta_a$ be the Borel measure
on $A$ which is the Dirac mass associated to an element $a$ of $A$, as
mentioned in Section \ref{fourier transforms, 2}.  If we take $\nu =
\delta_a$ in (\ref{(f * nu)(x) = int_A f(x - y) d nu(y)}), then we get
that
\begin{equation}
\label{(f * delta_a)(x) = f(x - a)}
        (f * \delta_a)(x) = f(x - a)
\end{equation}
for every $x \in A$.

        If $\nu$ is any complex Borel measure on $A$, and $f$ is a bounded
Borel measurable function on $A$, then (\ref{int_A |f(x - y)|
  d|nu|(y)}) is less than or equal to the supremum norm of $f$ on $A$
times $|\nu|(A)$ for every $x \in A$.  Thus (\ref{(f * nu)(x) =
  int_A f(x - y) d nu(y)}) is defined and satisfies
\begin{equation}
\label{|(f * nu)(x)| le (sup_{w in A} |f(w)|) |nu|(A)}
        |(f * \nu)(x)| \le \Big(\sup_{w \in A} |f(w)|\Big) \, |\nu|(A)
\end{equation}
for every $x \in A$.  If $f$ is also uniformly continuous on $A$, then
it is easy to see that (\ref{(f * nu)(x) = int_A f(x - y) d nu(y)}) is
uniformly continuous on $A$ as well.  If $\phi$ is a continuous
homomorphism from $A$ into ${\bf T}$, then
\begin{equation}
\label{(phi * nu)(x) = widehat{nu}(phi) phi(x)}
        (\phi * \nu)(x) = \widehat{\nu}(\phi) \, \phi(x)
\end{equation}
for every $x \in A$.

        If $f$ is continuous on $A$ and $\nu$ has compact support on $A$,
then (\ref{int_A |f(x - y)| d|nu|(y)}) is finite for every $x \in A$,
because $f$ is bounded on compact subsets of $A$.  Remember that $f$
is uniformly continuous along compact subsets of $A$, which implies
that (\ref{(f * nu)(x) = int_A f(x - y) d nu(y)}) is continuous on $A$
in this situation.  If $f$ is uniformly continuous on $A$, then
(\ref{(f * nu)(x) = int_A f(x - y) d nu(y)}) is uniformly continuous
on $A$ too, as before.

        Suppose now that $A$ is locally compact, and let $H$ be a
Haar measure on $A$.  If $f$ and $\nu$ are real-valued and nonnegative,
then Fubini's theorem implies that
\begin{equation}
\label{int_A (f * nu)(x) dH(x) = (int_A f(x) dH(x)) nu(A)}
 \int_A (f * \nu)(x) \, dH(x) = \Big(\int_A f(x) \, dH(x)\Big) \, \nu(A),
\end{equation}
at least under suitable conditions.  As before, this is a bit simpler
when there is a countable base for the topology of $A$, and otherwise
one can ask that $\nu$ satisfy appropriate regularity conditions.  At
any rate, (\ref{int_A (f * nu)(x) dH(x) = (int_A f(x) dH(x)) nu(A)})
implies that (\ref{(f * nu)(x) = int_A f(x - y) d nu(y)}) is finite
for almost every $x \in A$ with respect to $H$.

        If $f$ and $\nu$ are complex-valued, and if $f$ is integrable
with respect to $H$, then one can apply the remarks in the previous
paragraph to $|f|$ and $|\nu|$.  It follows that (\ref{int_A |f(x -
  y)| d|nu|(y)}) is finite for almost every $x \in A$ with respect to
$H$, so that (\ref{(f * nu)(x) = int_A f(x - y) d nu(y)}) is defined
for almost every $x \in A$ with respect to $H$.  We also get that
\begin{equation}
\label{int_A |(f * nu)(x)| dH(x) le (int_A |f(x)| dH(x)) |nu|(A)}
        \int_A |(f * \nu)(x)| \, dH(x) \le \Big(\int_A |f(x)| \, dH(x)\Big)
                                                               \, |\nu|(A),
\end{equation}
because the modulus of (\ref{(f * nu)(x) = int_A f(x - y) d nu(y)}) is
less than or equal to (\ref{(f * nu)(x) = int_A f(x - y) d nu(y)}).
If $\mu$ is the Borel measure on $A$ corresponding to $f$ as in
(\ref{mu(E) = int_E f(x) dH(x)}), then $\mu * \nu$ corresponds to
(\ref{(f * nu)(x) = int_A f(x - y) d nu(y)}) as an integrable function
on $A$ with respect to $H$ in the same way.

        If $f$ and $\nu$ are real-valued and nonnegative, and 
$1 \le p < \infty$, then
\begin{equation}
\label{(int_A (f * nu)(x)^p dH(x))^{1/p} le ...}
        \Big(\int_A (f * \nu)(x)^p \, dH(x)\Big)^{1/p}
              \le \Big(\int_A f(x)^p \, dH(x)\Big)^{1/p} \, \nu(A).
\end{equation}
As in the preceding section, this can be derived from the integral
version of Minkowski's inequality, or Jensen's inequality and Fubini's
theorem, at least under suitable conditions.  If $f$ and $\nu$ are
complex-valued, and if $|f|^p$ is integrable with respect to $H$,
then we can apply (\ref{(int_A (f * nu)(x)^p dH(x))^{1/p} le ...})
to $|f|$ and $|\nu|$, to get that (\ref{int_A |f(x - y)| d|nu|(y)})
is finite for almost every $x \in A$ with respect to $H$.  This implies
that (\ref{(f * nu)(x) = int_A f(x - y) d nu(y)}) is defined for
almost every $x \in A$ with respect to $H$, and that
\begin{equation}
\label{(int_A |(f * nu)(x)|^p dH(x))^{1/p} le ...}
        \Big(\int_A |(f * \nu)(x)|^p \, dH(x)\Big)^{1/p}
              \le \Big(\int_A |f(x)|^p \, dH(x)\Big)^{1/p} \, |\nu|(A).
\end{equation}

        If $g$ is a complex-valued integrable function on $A$, then
\begin{equation}
\label{nu(E) = int_E g(y) dH(y)}
        \nu(E) = \int_E g(y) \, dH(y)
\end{equation}
defines a complex-valued regular Borel measure on $A$, with
\begin{equation}
\label{|nu|(E) = int_E |g(y)| dH(y)}
        |\nu|(E) = \int_E |g(y)| \, dH(y)
\end{equation}
for every Borel set $E \subseteq A$.  In this case, convolution of a
function $f$ on $A$ with $\nu$ is the same as convolution of $f$ with
$g$, as discussed in the previous section, and with similar estimates
in terms of $|\nu|$.

        Let $f$ be a continuous complex-valued function on $A$.
If $\nu$ is a complex Borel measure on $A$ with compact support, then
we have seen that $f * \nu$ is a continuous function on $A$ as well.
If $f$ also has compact support on $A$, then it is easy to see that $f
* \nu$ has compact support on $A$ too.  Similarly, if $f$ vanishes at
infinity on $A$, then $f * \nu$ vanishes at infinity on $A$.  This
also works when $\nu$ can be approximated by complex Borel measures
with compact support on $A$ with respect to the total variation norm,
and in particular when $\nu$ is a regular Borel measure on $A$.

\section{Some simple estimates}
\label{some simple estimates}

        Let $A$ be a commutative topological group, and fix a real
number $C \ge 1$.  Also let $\nu$ be a complex Borel measure on $A$
such that
\begin{equation}
\label{nu(A) = 1}
        \nu(A) = 1
\end{equation}
and
\begin{equation}
\label{|nu|(A) le C}
        |\nu|(A) \le C,
\end{equation}
where $|\nu|$ is the total variation measure associated to $\nu$ on
$A$.  Of course, if $\nu$ is real-valued and nonnegative, then $|\nu|
= \nu$, and (\ref{nu(A) = 1}) implies that (\ref{|nu|(A) le C}) holds
with $C = 1$.  If $f$ is a complex-valued Borel measurable function on
$A$ such that (\ref{int_A |f(x - y)| d|nu|(y)}) is finite, then $(f *
\nu)(x)$ can be defined as in (\ref{(f * nu)(x) = int_A f(x - y) d
  nu(y)}), and
\begin{equation}
\label{(f * nu)(x) - f(x) = int_A (f(x - y) - f(x)) d nu(y)}
        (f * \nu)(x) - f(x) = \int_A (f(x - y) - f(x)) \, d\nu(y),
\end{equation}
by (\ref{nu(A) = 1}).  This implies that
\begin{equation}
\label{|(f * nu)(x) - f(x)| le int_A |f(x - y) - f(x)| d|nu|(y)}
        |(f * \nu)(x) - f(x)| \le \int_A |f(x - y) - f(x)| \, d|\nu|(y).
\end{equation}

        Let $U$ be an open set in $A$ with $0 \in U$, and suppose that
\begin{equation}
\label{|nu|(A setminus U) = 0}
        |\nu|(A \setminus U) = 0.
\end{equation}
In this case, (\ref{|(f * nu)(x) - f(x)| le int_A |f(x - y) - f(x)|
  d|nu|(y)}) reduces to
\begin{eqnarray}
\label{|(f * nu)(x) - f(x)| le ... le C sup_{y in U} |f(x - y) - f(x)|}
        |(f * \nu)(x) - f(x)| & \le & \int_U |f(x - y) - f(x)| \, d|\nu|(y) \\
         & \le & C \, \sup_{y \in U} |f(x - y) - f(x)|, \nonumber
\end{eqnarray}
using (\ref{|nu|(A) le C}) in the second step.  If $f$ is continuous
at $x$, then the right side of (\ref{|(f * nu)(x) - f(x)| le ... le C
  sup_{y in U} |f(x - y) - f(x)|}) is small when $U$ is sufficiently small.

        Now let $\eta$ be a nonnegative real number, and suppose that
\begin{equation}
\label{|nu|(A setminus U) le eta}
        |\nu|(A \setminus U) \le \eta
\end{equation}
instead of (\ref{|nu|(A setminus U) = 0}), and which reduces to
(\ref{|nu|(A setminus U) = 0}) when $\eta = 0$.  If $f$ is bounded on
$A$, then (\ref{|(f * nu)(x) - f(x)| le int_A |f(x - y) - f(x)|
  d|nu|(y)}) implies that
\begin{eqnarray}
\label{|(f * nu)(x) - f(x)| le ...}
\lefteqn{|(f * \nu)(x) - f(x)|} \\
 & \le & \int_U |f(x - y) - f(x)| \, d|\nu|(y)
           + \int_{A \setminus U} |f(x - y) - f(x)| \, d|\nu|(y) \nonumber \\
 & \le & C \, \sup_{y \in U} |f(x - y) - f(x)|
                           + 2 \, \eta \, \sup_{z \in A} |f(z)|, \nonumber
\end{eqnarray}
using (\ref{|nu|(A) le C}) and (\ref{|nu|(A setminus U) le eta}) in
the second step.  This is small when $f$ is continuous at $x$, and $U$
and $\eta$ are sufficiently small.

        Let us apply this to $x = 0$ and a continuous homomorphism $\phi$
from $A$ into ${\bf T}$, so that $\phi(0) = 1$ and
\begin{equation}
\label{(phi * nu)(0) = widehat{nu}(phi)}
        (\phi * \nu)(0) = \widehat{\nu}(\phi),
\end{equation}
as in (\ref{(phi * nu)(x) = widehat{nu}(phi) phi(x)}).  Plugging this
into (\ref{|(f * nu)(x) - f(x)| le ...}), we get that
\begin{equation}
\label{|widehat{nu}(phi) - 1| le C sup_{y in U} |phi(x - y) - phi(x)| + 2 eta}
 |\widehat{\nu}(\phi) - 1| \le C \, \sup_{y \in U} |\phi(x - y) - \phi(x)|
                                                              + 2 \, \eta,
\end{equation}
because $|\phi(z)| = 1$ for every $z \in A$.  Note that
\begin{equation}
\label{widehat{nu}({bf 1}_A) = 1}
        \widehat{\nu}({\bf 1}_A) = 1,
\end{equation}
by (\ref{nu(A) = 1}), where ${\bf 1}_A$ is the constant function on
$A$ equal to $1$ at every point, which is the identity element of
$\widehat{A}$.  If $\mathcal{E} \subseteq \widehat{A}$ is
equicontinuous at $0$, then (\ref{|widehat{nu}(phi) - 1| le C sup_{y
    in U} |phi(x - y) - phi(x)| + 2 eta}) is uniformly small over
$\phi \in \mathcal{E}$ when $U$ and $\eta$ are sufficiently small.  In
particular, if $A$ is locally compact, and $\mathcal{E} \subseteq
\widehat{A}$ is compact or totally bounded with respect to the usual
topology on $\widehat{A}$, then $\mathcal{E}$ is equicontinuous at $0$
on $A$, by standard arguments.

\section{Some more estimates}
\label{some more estimates}

        Let $A$ be a locally compact commutative topological group,
and let $U$ be an open set in $A$ that contains $0$.  We may as well
ask that the closure $\overline{U}$ of $U$ in $A$ be compact, since
this can always be arranged by taking the intersection of $U$ with a
fixed neighborhood of $0$ in $A$ with this property.  Also let $C \ge
1$ be fixed, and let $\nu$ be a regular complex Borel measure on $A$
that satisfies (\ref{nu(A) = 1}), (\ref{|nu|(A) le C}), and
(\ref{|nu|(A setminus U) = 0}).  Of course, (\ref{|nu|(A setminus U) =
  0}) implies that the support of $\nu$ is contained in
$\overline{U}$, and hence that the support of $\nu$ is compact.

        Thus $(f * \nu)(x)$ is defined for every $x \in A$ when $f$
is a continuous complex-valued function on $A$.  If $f$ is uniformly
continuous along a set $E \subseteq A$, then $f * \nu$ will be uniformly
close to $f * \nu$ on $E$ when $U$ is sufficiently small, as in
(\ref{|(f * nu)(x) - f(x)| le ... le C sup_{y in U} |f(x - y) - f(x)|}).
In particular, $f * \nu$ is uniformly close to $f$ on $A$ when $f$ is
uniformly continuous on $A$ and $U$ is sufficiently small.  If $f$ is
any continuous function on $A$, then $f$ is uniformly continuous on
compact subsets of $A$, and $f * \nu$ is uniformly close to $f$
on compact subsets of $A$ when $U$ is sufficiently small, depending
on the compact set.  If $f$ is a continuous function on $A$ that
vanishes at infinity, then $f$ is uniformly continuous on $A$, so
that $f * \nu$ is uniformly close to $f$ when $U$ is sufficiently small,
and $f * \nu$ also vanishes at infinity on $A$, because $\nu$ has compact
support on $A$.

        Suppose that $f$ has support contained in a compact set
$K_1 \subseteq A$, and that the support of $\nu$ is contained in
a fixed compact set $K_2 \subseteq A$.  Remember that the support of
$\nu$ is automatically contained in $\overline{U}$, and so it suffices 
to ask that
\begin{equation}
\label{overline{U} subseteq K_2}
        \overline{U} \subseteq K_2
\end{equation}
to get the support of $\nu$ contained in $K_2$.  In particular, if $0$
is in the interior of $K_2$, then this holds for all sufficiently
small $U$.  At any rate, it is convenient to ask that $0 \in K_2$, so
that $K_1 \subseteq K_1 + K_2$.  This can always be arranged by
replacing $K_2$ with $K_2 \cup \{0\}$, if necessary.  If $(f * \nu)(x)
\ne 0$ for some $x \in A$, then $f(x - y) \ne 0$ for some $y$ in the
support of $\nu$, which means that $x - y \in K_1$ and $y \in K_2$.
This implies that $x$ is an element of $K_1 + K_2$, which is the image
of $K_1 \times K_2$ under addition as a mapping from $A \times A$ into
$A$.  Note that $K_1 + K_2$ is compact in $A$, because $K_1 \times
K_2$ is a compact subset of $A \times A$ with respect to the product
topology, and addition is continuous as a mapping from $A \times A$
into $A$.

          Let $H$ be a Haar measure on $A$, and let $1 \le p < \infty$
be given.  Because $f$ and $f * \nu$ have supports contained in $K_1 +
K_2$, we have that
\begin{eqnarray}
\label{(int_A |(f * nu)(x) - f(x)|^p dH(x))^{1/p} le ...}
\lefteqn{\Big(\int_A |(f * \nu)(x) - f(x)|^p \, dH(x)\Big)^{1/p}} \\
         & \le & \Big(\sup_{x \in K_1 + K_2} |(f * \nu)(x) - f(x)|\Big)
                                              \, H(K_1 + K_2)^{1/p}. \nonumber
\end{eqnarray}
Of course, $H(K_1 + K_2) < \infty$, because $K_1 + K_2$ is compact.
We already know that $f * \nu$ is uniformly close to $f$ on $K_1 +
K_2$ when $U$ is sufficiently small, by (\ref{|(f * nu)(x) - f(x)| le
  ... le C sup_{y in U} |f(x - y) - f(x)|}) and uniform continuity, so
that (\ref{(int_A |(f * nu)(x) - f(x)|^p dH(x))^{1/p} le ...})  is
small when $U$ is sufficiently small.  If $f$ is any function on $A$
in $L^p(A)$, then $f * \nu$ is also close to $f$ on $A$ with respect
to the $L^p$ norm associated to $H$ when $U$ is sufficiently small,
depending on $f$.  This follows from the preceding argument when $f$
is a continuous function on $A$ with compact support, and otherwise
one can get the same conclusion by approximating $f$ by continuous
functions on $A$ with compact support with respect to the $L^p$ norm.
This also uses (\ref{(int_A |(f * nu)(x)|^p dH(x))^{1/p} le ...}) to
estimate the errors in the approximation.

        Now let $\eta \ge 0$ be given, and suppose that $\nu$ satisfies
(\ref{|nu|(A setminus U) le eta}) instead of (\ref{|nu|(A setminus U) = 0}).
If $f$ is a bounded function on $A$ that is also uniformly continuous on $A$,
then $f * \nu$ is uniformly close to $f$ on $A$ when $U$ and $\eta$ are
sufficiently small, by (\ref{|(f * nu)(x) - f(x)| le ...}).  In order
to deal with $L^p$ norms as before, it will be helpful to approximate
$\nu$ by a measure that satisfies (\ref{|nu|(A setminus U) = 0}).

        Put
\begin{equation}
\label{nu'(E) = nu(E cap U) + nu(A setminus U) delta_0(E)}
        \nu'(E) = \nu(E \cap U) + \nu(A \setminus U) \, \delta_0(E)
\end{equation}
for each Borel set $E \subseteq A$, where $\delta_0$ is the usual
Dirac mass at $0$.  This is a regular complex Borel measure on $A$,
because of the corresponding properties of $\nu$, and
\begin{equation}
\label{nu'(A) = nu(U) + nu(A setminus U) = nu(A) = 1}
        \nu'(A) = \nu(U) + \nu(A \setminus U) = \nu(A) = 1,
\end{equation}
by (\ref{nu(A) = 1}).  It is easy to see that
\begin{equation}
\label{|nu'|(E) le |nu|(E cap U) + |nu(A setminus U)| delta_0(E)}
        |\nu'|(E) \le |\nu|(E \cap U) + |\nu(A \setminus U)| \, \delta_0(E)
\end{equation}
for every Borel set $E \subseteq A$, and hence that
\begin{equation}
\label{|nu'|(A) le |nu|(U) + |nu(A setminus U)| le ... = |nu|(A) le C}
        |\nu'|(A) \le |\nu|(U) + |\nu(A \setminus U)|
                   \le |\nu|(U) + |\nu|(A \setminus U) = |\nu|(A) \le C,
\end{equation}
by (\ref{|nu|(A) le C}).  Note that
\begin{equation}
\label{|nu'|(A setminus U) = 0}
        |\nu'|(A \setminus U) = 0,
\end{equation}
because $0 \in U$, so that $\nu'$ satisfies the analogues of
(\ref{nu(A) = 1}), (\ref{|nu|(A) le C}), and (\ref{|nu|(A setminus U)
  = 0}).

        Similarly, put
\begin{equation}
\label{nu''(E) = nu(E) - nu'(E) = ...}
        \nu''(E) = \nu(E) - \nu'(E)
                 = \nu(E \setminus U) - \nu(A \setminus U) \, \delta_0(E)
\end{equation}
for each Borel set $E \subseteq A$, which is another regular complex
Borel measure on $A$.  As before,
\begin{equation}
\label{|nu''|(E) le |nu|(E setminus U) + |nu(A setminus U)| delta_0(E)}
 |\nu''|(E) \le |\nu|(E \setminus U) + |\nu(A \setminus U)| \, \delta_0(E)
\end{equation}
for every Borel set $E \subseteq A$, so that
\begin{equation}
\label{|nu''|(A) le |nu|(A setminus U) + |nu(A setminus U)| le ... le 2 eta}
        |\nu''|(A) \le |\nu|(A \setminus U) + |\nu(A \setminus U)|
                    \le 2 \, |\nu|(A \setminus U) \le 2 \, \eta,
\end{equation}
by (\ref{|nu|(A setminus U) le eta}).  If $1 \le p < \infty$ and $f
\in L^p(A)$, then $f * \nu$ is close to $f$ with respect to the $L^p$
norm on $A$ when $U$ is sufficiently small, as discussed earlier.  The
$L^p$ norm of $f * \nu''$ is bounded by (\ref{|nu''|(A) le |nu|(A
  setminus U) + |nu(A setminus U)| le ... le 2 eta}) times the $L^p$
norm of $f$, as in (\ref{(int_A |(f * nu)(x)|^p dH(x))^{1/p} le ...}).
This implies that $f * \nu = f * \nu' + f * \nu''$ is close to $f$
with respect to the $L^p$ norm when $U$ and $\eta$ are sufficiently
small, depending on $f$.

\section{Another variant}
\label{another variant}

        Let $A$ be a locally compact commutative topological group, and
let $\mu$, $\nu$ be regular complex Borel measures on $A$.  Also let
$U$ be an open set in $A$ with $0 \in U$, let $C \ge 1$ be fixed, and
let $\eta$ be a nonnegative real number.  As before, we suppose that
$\nu$ satisfies (\ref{nu(A) = 1}), (\ref{|nu|(A) le C}), and
(\ref{|nu|(A setminus U) le eta}).  We would like to look at the
convolution $\mu * \nu$ of $\mu$ and $\nu$ as an approximation to
$\mu$ when $U$ and $\eta$ are sufficiently small.  Note that
\begin{equation}
\label{|mu * nu|(A) le |mu|(A) |nu|(A) le C |mu|(A)}
        |\mu * \nu|(A) \le |\mu|(A) \, |\nu|(A) \le C \, |\mu|(A),
\end{equation}
by (\ref{|mu * nu|(A) le (|mu| * |nu|)(A) = |mu|(A) |nu|(A)}) and
(\ref{|nu|(A) le C}).

        If $f$ is a bounded complex-valued Borel measurable function
on $A$, then
\begin{equation}
\label{int_A f(z) d(mu * nu)(z) = int_A (int_A f(x + y) d nu(y)) d mu(x)}
        \int_A f(z) \, d(\mu * \nu)(z)
                   = \int_A \Big(\int_A f(x + y) \, d\nu(y)\Big) \, d\mu(x),
\end{equation}
as in (\ref{int_A f(z) d(mu * nu)(z) = ...}).  Put
\begin{equation}
\label{widetilde{nu}(E) = nu(-E)}
        \widetilde{\nu}(E) = \nu(-E)
\end{equation}
for every Borel set $E \subseteq A$, which defines another regular
complex Borel measure on $A$.  By construction,
\begin{equation}
\label{(f * widetilde{nu})(x) = ... = int_A f(x + y) d nu(y)}
        (f * \widetilde{\nu})(x) = \int_A f(x - y) \, d\widetilde{\nu}(y)
                                 = \int_A f(x + y) \, d\nu(y),
\end{equation}
so that (\ref{int_A f(z) d(mu * nu)(z) = int_A (int_A f(x + y) d
  nu(y)) d mu(x)}) becomes
\begin{equation}
\label{int_A f(z) d(mu * nu)(z) = int_A (f * widetilde{nu})(x) d mu(x)}
 \int_A f(z) \, d(\mu * \nu)(z) = \int_A (f * \widetilde{\nu})(x) \, d\mu(x).
\end{equation}
It follows that
\begin{equation}
\label{int_A f(z) d(mu * nu)(z) - int_A f(x) d mu(x) = ...}
        \int_A f(z) \, d(\mu * \nu)(z) - \int_A f(x) \, d\mu(x)
                   = \int_A ((f * \widetilde{\nu})(x) - f(x)) \, d\mu(x),
\end{equation}
and hence
\begin{eqnarray}
\label{|int_A f(z) d(mu * nu)(z) - int_A f(x) d mu(x)| le ...}
\lefteqn{\biggl|\int_A f(z) \, d(\mu * \nu)(z) - \int_A f(x) \, d\mu(x)\biggr|}
                                                                           \\
 & \le & \int_A |(f * \widetilde{\nu})(x) - f(x)| \, d|\mu|(x). \nonumber
\end{eqnarray}
In particular,
\begin{eqnarray}
\label{|int_A f(z) d(mu * nu)(z) - int_A f(x) d mu(x)| le ..., 2}
\lefteqn{\biggl|\int_A f(z) \, d(\mu * \nu)(z) - \int_A f(x) \, d\mu(x)\biggr|}
                                                                            \\
 & \le & \Big(\sup_{x \in A} |(f * \widetilde{\nu})(x) - f(x)|\Big) \, |\mu|(A).
                                                                    \nonumber
\end{eqnarray}

        Of course, $\widetilde{\nu}$ satisfies the same conditions
as $\nu$, but with $U$ replaced with $-U$.  If $f$ is bounded and
uniformly continuous on $A$, then $f * \widetilde{\nu}$ is uniformly
close to $f$ on $A$ when $U$ and $\eta$ are sufficiently small, for
the same reasons as before.  This implies that (\ref{int_A f(z) d(mu *
  nu)(z) - int_A f(x) d mu(x) = ...}) is small under these conditions,
by (\ref{|int_A f(z) d(mu * nu)(z) - int_A f(x) d mu(x)| le ..., 2}).

        Remember that
\begin{equation}
\label{lambda_mu(f) = int_A f(x) d mu(x)}
        \lambda_\mu(f) = \int_A f(x) \, d\mu(x)
\end{equation}
and
\begin{equation}
\label{lambda_{mu * nu}(f) = int_A f(z) d(mu * nu)(z)}
        \lambda_{\mu * \nu}(f) = \int_A f(z) \, d(\mu * \nu)(z)
\end{equation}
define bounded linear functionals on $C_0(A)$, with respect to the
supremum norm on $C_0(A)$.  If $f \in C_0(A)$, then $f$ is bounded and
uniformly continuous on $A$, and the preceding discussion shows that
\begin{equation}
\label{lambda_{mu * nu}(f) - lambda_mu(f)}
        \lambda_{\mu * \nu}(f) - \lambda_\mu(f)
\end{equation}
is small when $U$ and $\eta$ are sufficiently small.  This implies
that $\lambda_{\mu * \nu}$ is close to $\lambda_\mu$ with respect to
the weak$^*$ topology on the dual of $C_0(A)$ when $U$ and $\eta$
are sufficiently small.

\section{Discrete groups}
\label{discrete groups}

        Let $A$ be a commutative group equipped with the discrete topology.
Remember that a complex-valued function $f(x)$ on $A$ is said to be
\emph{summable}\index{summable functions} if the sums
\begin{equation}
\label{sum_{x in E} |f(x)|}
        \sum_{x \in E} |f(x)|
\end{equation}
over all nonempty finite subsets $E$ of $A$ are uniformly bounded.
In this case, the sum
\begin{equation}
\label{sum_{x in A} |f(x)|}
        \sum_{x \in A} |f(x)|
\end{equation}
may be defined as the supremum of (\ref{sum_{x in E} |f(x)|}) over all
nonempty finite subsets $E$ of $A$.  The space of summable functions
on $A$ is denoted $\ell^1(A)$, and is a vector space with respect to
pointwise addition and multiplication.  It is easy to see that
(\ref{sum_{x in A} |f(x)|}) defines a norm on $\ell^1(A)$, which may
be denoted $\|f\|_1$ or $\|f\|_{\ell^1(A)}$.  Of course, $f(x)$ is
summable on $A$ if and only if $f(x)$ is integrable with respect to
counting measure on $A$, so that $\ell^1(A)$ is the same as the usual
space $L^1(A)$ associated to counting measure on $A$.  If $f(x)$ is
summable on $A$, then the sum
\begin{equation}
\label{sum_{x in A} f(x)}
        \sum_{x \in A} f(x)
\end{equation}
can be defined in various equivalent ways, and is the same as the
integral of $f(x)$ with respect to counting measure on $A$.

        If $1 \le p < \infty$, then $\ell^p(A)$ is defined to be the
space of complex-valued functions $f(x)$ on $A$ such that $|f(x)|^p$
is summable on $A$.  It is well known that this is also a vector space
with respect to pointwise addition and scalar multiplication, and that
\begin{equation}
\label{||f||_p = ||f||_{ell^p(A)} = (sum_{x in A} |f(x)|^p)^{1/p}}
        \|f\|_p = \|f\|_{\ell^p(A)} = \Big(\sum_{x \in A} |f(x)|^p\Big)^{1/p}
\end{equation}
defines a norm on $\ell^p(A)$.  Similarly, $\ell^\infty(A)$ is defined
to be the space of all bounded complex-valued functions on $A$,
equipped with the supremum norm
\begin{equation}
\label{||f||_infty = ||f||_{ell^infty(A)} = sup_{x in A} |f(x)|}
        \|f\|_\infty = \|f\|_{\ell^\infty(A)} = \sup_{x \in A} |f(x)|.
\end{equation}
These are the same as the usual $L^p$ spaces associated to counting
measure on $A$.  In particular, $\ell^p(A)$ is complete with respect
to the metric corresponding to the $\ell^p$ norm for each $p \ge 1$,
so that $\ell^p(A)$ is a Banach space.

        If $f, g \in \ell^2(A)$, then $|f(x)| \, |g(x)|$ is summable on
$A$, and
\begin{equation}
\label{langle f, g rangle = sum_{x in A} f(x) overline{g(x)}}
        \langle f, g \rangle = \sum_{x \in A} f(x) \, \overline{g(x)}
\end{equation}
defines an inner product on $\ell^2(A)$.  By construction,
\begin{equation}
\label{langle f, f rangle = sum_{x in A} |f(x)|^2 = ||f||_2^2}
        \langle f, f \rangle = \sum_{x \in A} |f(x)|^2 = \|f\|_2^2,
\end{equation}
so that $\ell^2(A)$ is a Hilbert space with respect to this inner
product.

        A complex-valued function $f(x)$ on $A$ is said to
\emph{vanish at infinity}\index{vanishing at infinity} if
\begin{equation}
\label{{x in A : |f(x)| ge epsilon}}
        \{x \in A : |f(x)| \ge \epsilon\}
\end{equation}
has only finitely many elements for each $\epsilon > 0$.  This is
equivalent to saying that $f(x)$ vanishes at infinity on $A$ as a
locally compact Hausdorff topological space, equipped with the
discrete topology.  The space of complex-valued functions $f(x)$
vanishing at infinity on $A$ may be denoted $c_0(A)$, which is a
closed linear subspace of $\ell^\infty(A)$.  Of course, every function
on $A$ is continuous with respect to the discrete topology, so that
$c_0(A)$ is the same as the space $C_0(A)$ defined previously for any
locally compact Hausdorff topological space.  If $f \in \ell^p(A)$ for
some $p < \infty$, then it is easy to see that $f \in c_0(A)$.

        In this situation, counting measure can be used as Haar measure
on $A$, and $\widehat{A}$ consists of all homomorphisms from $A$ into
${\bf T}$.  The Fourier transform of $f \in \ell^1(A)$ can be
expressed as
\begin{equation}
\label{widehat{f}(phi) = sum_{x in A} f(x) overline{phi(x)}}
        \widehat{f}(\phi) = \sum_{x \in A} f(x) \, \overline{\phi(x)}
\end{equation}
for every $\phi \in \widehat{A}$.  Note that
\begin{equation}
\label{|widehat{f}(phi) - sum_{x in E} f(x) overline{phi(x)}| = ...}
 \biggl|\widehat{f}(\phi) - \sum_{x \in E} f(x) \, \overline{\phi(x)}\biggr|
   = \biggl|\sum_{x \in A \setminus E} f(x) \, \overline{\phi(x)}\biggr|
    \le \sum_{x \in A \setminus E} |f(x)|
\end{equation}
for every subset $E$ of $A$ and $\phi \in \widehat{A}$.  The
summability of $f(x)$ on $A$ implies that for each $\epsilon > 0$
there is a finite subset $E_1(\epsilon)$ of $A$ such that
\begin{equation}
\label{sum_{x in A setminus E_1(epsilon)} |f(x)| < epsilon}
        \sum_{x \in A \setminus E_1(\epsilon)} |f(x)| < \epsilon.
\end{equation}
If $E$ is a subset of $A$ that contains $E_1(\epsilon)$, then it
follows that
\begin{equation}
\label{|widehat{f}(phi) - sum_{x in E} f(x) overline{phi(x)}| < epsilon}
 \biggl|\widehat{f}(\phi) - \sum_{x \in E} f(x) \, \overline{\phi(x)}\biggr|
                                                      < \epsilon
\end{equation}
for every $\phi \in \widehat{A}$, by (\ref{|widehat{f}(phi) - sum_{x
    in E} f(x) overline{phi(x)}| = ...}).

        The space $C(A)$ of complex-valued continuous functions on $A$
is the same as the space of all complex-valued functions on $A$,
because $A$ is equipped with the discrete topology.  The usual topology
on $C(A)$ is defined by the collection of supremum seminorms associated
to nonempty compact subsets of $A$, and here the compact subsets of $A$
have only finitely many elements.  Put
\begin{equation}
\label{Psi_x(phi) = phi(x)}
        \Psi_x(\phi) = \phi(x)
\end{equation}
for each $x \in A$ and $\phi \in \widehat{A}$, as in Section
\ref{second dual}.  This defines a homomorphism from $\widehat{A}$
into ${\bf T}$ for each $x \in A$, and this homomorphism is continuous
with respect to the topology induced on $\widehat{A}$ by the usual
topology on $C(A)$.  Thus $\Psi_x$ is an element of the dual
$\widehat{\widehat{A}}$ of $\widehat{A}$ for each $x \in A$, and we
have seen that
\begin{equation}
\label{x mapsto Psi_x, 2}
        x \mapsto \Psi_x
\end{equation}
defines a homomorphism from $A$ into $\widehat{\widehat{A}}$.  We have
also seen that the elements of $\widehat{A}$ separate points in $A$,
which means that (\ref{x mapsto Psi_x, 2}) is injective.  In fact,
(\ref{x mapsto Psi_x, 2}) maps $A$ onto $\widehat{\widehat{A}}$,
as in Section \ref{second dual}.  Note that
\begin{equation}
\label{Psi_{-x}(phi) = phi(-x) = 1/phi(x) = overline{phi(x)}}
        \Psi_{-x}(\phi) = \phi(-x) = 1/\phi(x) = \overline{\phi(x)}
\end{equation}
for every $x \in A$ and $\phi \in \widehat{A}$, using the fact that
$\phi(x) \in {\bf T}$ in the last step.

        Remember that $\widehat{A}$ is compact with respect to the
topology induced by the usual topology on $C(A)$ when $A$ is discrete.
Let $H_{\widehat{A}}$ be Haar measure on $\widehat{A}$, normalized so
that $H_{\widehat{A}}(\widehat{A}) = 1$.  The elements of
$\widehat{\widehat{A}}$ are orthonormal in $L^2(\widehat{A})$ with
respect to the standard integral inner product associated to
$H_{\widehat{A}}$, as in Section \ref{compact, discrete groups}.  Thus
the functions $\Psi_x$ on $\widehat{A}$ with $x \in A$ are orthonormal
in $L^2(\widehat{A})$, which is the same as saying that the functions
$\Psi_{-x}$ with $x \in A$ are orthonormal in $L^2(\widehat{A})$.  If
$f \in \ell^2(A)$, then it follows that
\begin{equation}
\label{sum_{x in A} f(x) Psi_{-x}}
        \sum_{x \in A} f(x) \, \Psi_{-x}
\end{equation}
defines an element of $L^2(\widehat{A})$, by standard Hilbert space
arguments.  In particular, the norm of (\ref{sum_{x in A} f(x)
  Psi_{-x}}) in $L^2(\widehat{A})$ is equal to the norm of $f$ in
$\ell^2(A)$.  Similarly,
\begin{eqnarray}
\label{||sum_{x in A} f(x) Psi_{-x} - sum_{x in E} f(x) Psi_{-x}||_{L^2} = ...}
        \biggl\|\sum_{x \in A} f(x) \, \Psi_{-x} - 
                     \sum_{x \in E} f(x) \, \Psi_{-x}\biggr\|_{L^2(\widehat{A})}
& = & \biggl\|\sum_{x \in A \setminus E} f(x) \, \Psi_{-x}\biggr\|_{L^2(\widehat{A})}
                                                                            \\
& = & \Big(\sum_{x \in A \setminus E} |f(x)|^2\Big)^{1/2} \nonumber
\end{eqnarray}
for each $E \subseteq A$, where $\|\cdot\|_{L^2(\widehat{A})}$ denotes
the usual norm on $L^2(\widehat{A})$.  As before, for each $\epsilon >
0$ there is a finite subset $E_2(\epsilon)$ of $A$ such that
\begin{equation}
\label{(sum_{x in A setminus E_2(epsilon)} |f(x)|^2)^{1/2} < epsilon}
        \Big(\sum_{x \in A \setminus E_2(\epsilon)} |f(x)|^2\Big)^{1/2} < \epsilon,
\end{equation}
because $f \in \ell^2(A)$.  This implies that
\begin{equation}
\label{|sum_{x in A} f(x) Psi_{-x} - sum_{x in E} f(x) Psi_{-x}||_{L^2} < eps}
 \biggl\|\sum_{x \in A} f(x) \, \Psi_{-x} - 
             \sum_{x \in E} f(x) \, \Psi_{-x}\biggr\|_{L^2(\widehat{A})} < \epsilon
\end{equation}
when $E_2(\epsilon) \subseteq E \subseteq A$, by (\ref{||sum_{x in A}
  f(x) Psi_{-x} - sum_{x in E} f(x) Psi_{-x}||_{L^2} = ...}).

        If $f \in \ell^1(A)$, then it is well known that $f \in \ell^2(A)$,
and that
\begin{equation}
\label{||f||_{ell^2(A)} le ||f||_{ell^1(A)}}
        \|f\|_{\ell^2(A)} \le \|f\|_{\ell^1(A)}.
\end{equation}
In this case, $\widehat{f}$ is a continuous function on $\widehat{A}$,
which determines an element of $L^2(\widehat{A})$.  One can check that
$\widehat{f}$ is the same as (\ref{sum_{x in A} f(x) Psi_{-x}}) as an
element of $L^2(\widehat{A})$, because the uniform approximation of
(\ref{widehat{f}(phi) = sum_{x in A} f(x) overline{phi(x)}}) by finite
subsums discussed earlier implies approximation with respect to the
$L^2$ norm on $\widehat{A}$.  If $f \in \ell^2(A)$, then (\ref{sum_{x
    in A} f(x) Psi_{-x}}) may be considered as the definition of
$\widehat{f}$, as an element of $L^2(\widehat{A})$.  Note that the
collection of functions $\Psi_{-x}$ with $x \in A$ is actually an
orthonormal basis for $L^2(\widehat{A})$, by an argument in Section
\ref{compact, discrete groups}.  This uses the fact that the functions
$\Psi_x$ with $x \in A$ automatically separate points in
$\widehat{A}$, by construction.  It follows that every element of
$L^2(\widehat{A})$ can be expressed as $\widehat{f}$ for some $f \in
\ell^2(A)$.

\section{Compact groups}
\label{compact groups}

        Let $A$ be a compact commutative topological group, and let $H$
be Haar measure on $A$, normalized so that $H(A) = 1$.  Thus
continuous homomorphisms from $A$ into ${\bf T}$ are orthonormal with
respect to the standard inner product on $L^2(A)$, as in Section
\ref{compact, discrete groups}.  If $f \in L^2(A)$ and $\phi \in
\widehat{A}$, then $\widehat{f}(\phi)$ is the same as the inner
product of $f$ and $\phi$.  Using standard properties of orthonormal
sets in inner product spaces, we get that
\begin{equation}
\label{sum_{phi in E} |widehat{f}(phi)|^2 le ||f||_2^2}
        \sum_{\phi \in E} |\widehat{f}(\phi)|^2 \le \|f\|_2^2
\end{equation}
for every nonempty finite set $E \subseteq \widehat{A}$,
where $\|f\|_2$ denotes the $L^2$ norm of $f$ on $A$ with respect to
$H$.  This implies that
\begin{equation}
\label{sum_{phi in widehat{A}} |widehat{f}(phi)|^2 le ||f||_2^2}
        \sum_{\phi \in \widehat{A}} |\widehat{f}(\phi)|^2 \le \|f\|_2^2,
\end{equation}
where the sum on the left is defined to be the supremum of the
corresponding subsums on the left side of (\ref{sum_{phi in E}
  |widehat{f}(phi)|^2 le ||f||_2^2}).

        As mentioned previously, it is well known that the elements
of $\widehat{A}$ separate points in $A$ when $A$ is compact, which
implies that the linear span $\mathcal{E}(A)$ of $\widehat{A}$ is
dense in $C(A)$ with respect to the supremum norm.  In particular,
$\mathcal{E}(A)$ is dense in $L^2(A)$, so that $\widehat{A}$ forms an
orthonormal basis for $L^2(A)$.  It follows that
\begin{equation}
\label{f = sum_{phi in widehat{A}} widehat{f}(phi) phi}
        f = \sum_{\phi \in \widehat{A}} \widehat{f}(\phi) \, \phi
\end{equation}
for every $f \in L^2(A)$, where the sum converges with respect to the
$L^2$ norm $\|\cdot \|_2$.  More precisely, this means that for each
$\epsilon > 0$ there is a finite subset $E(\epsilon)$ of $\widehat{A}$
such that
\begin{equation}
\label{||f(x) - sum_{phi in E} widehat{f}(phi) phi||_2 < epsilon}
 \biggl\|f(x) - \sum_{\phi \in E} \widehat{f}(\phi) \, \phi\biggr\|_2 < \epsilon
\end{equation}
for every nonempty finite set $E \subseteq \widehat{A}$ that contains
$E(\epsilon)$.  This also implies that
\begin{equation}
\label{sum_{phi in widehat{A}} |widehat{f}(phi)|^2 = ||f||_2^2}
        \sum_{\phi \in \widehat{A}} |\widehat{f}(\phi)|^2 = \|f\|_2^2.
\end{equation}

        Let $\mu$ be a regular complex Borel measure on $A$, and suppose
that
\begin{equation}
\label{widehat{mu}(phi) = 0}
        \widehat{\mu}(\phi) = 0
\end{equation}
for every $\phi \in \widehat{A}$.  This implies that
\begin{equation}
\label{int_A g(x) d mu(x) = 0}
        \int_A g(x) \, d\mu(x) = 0
\end{equation}
for every $g \in \mathcal{E}(A)$, by linearity of the integral, and
the fact that $\mathcal{E}(A)$ is invariant under complex conjugation.
It follows that (\ref{int_A g(x) d mu(x) = 0}) holds for every
continuous complex-valued function $g$ on $A$, because
$\mathcal{E}(A)$ is dense in $C(A)$ with respect to the supremum norm.
Using the hypothesis that $\mu$ be a regular Borel measure, we get
that $\mu = 0$ on $A$ under these conditions.  In particular, if $f$
is an integrable function on $A$ with respect to $H$ such that
$\widehat{f}(\phi) = 0$ for every $\phi \in \widehat{A}$, then $f =
0$ almost everywhere on $A$ with respect to $H$.

        Let $\ell^2(\widehat{A})$ be the usual space of complex-valued
functions $F(\phi)$ on $\widehat{A}$ that are square-summable, as in the
previous section.  If $F \in \ell^2(\widehat{A})$, then
\begin{equation}
\label{f = sum_{phi in widehat{A}} F(phi) phi}
        f = \sum_{\phi \in \widehat{A}} F(\phi) \, \phi
\end{equation}
defines an element of $L^2(A)$, because of the orthonormality of the
elements of $\widehat{A}$ in $L^2(A)$.  The norm of $f$ in $L^2(A)$ is
the same as the norm of $F$ in $\ell^2(\widehat{A})$, and similarly
\begin{eqnarray}
\label{||f - sum_{phi in E} F(phi) phi||_{L^2(A)} = ...}
        \biggl\|f - \sum_{\phi \in E} F(\phi) \, \phi\biggr\|_{L^2(A)} & = &
 \biggl\|\sum_{\phi \in \widehat{A} \setminus E} F(\phi) \, \phi\biggr\|_{L^2(A)} \\
 & = & \Big(\sum_{\phi \in \widehat{A} \setminus E} |F(\phi)|^2\Big)^{1/2} \nonumber
\end{eqnarray}
for every $E \subseteq \widehat{A}$, where $\|\cdot \|_{L^2(A)}$
denotes the usual norm on $L^2(A)$.  Let $\epsilon > 0$ be given, and
let $E_2(\epsilon)$ be a finite subset of $\widehat{A}$ such that
\begin{equation}
\label{(sum_{phi in widehat{A} setminus E_2(epsilon)} |F(phi)|^2)^{1/2} < eps}
 \Big(\sum_{\phi \in \widehat{A} \setminus E_2(\epsilon)} |F(\phi)|^2\Big)^{1/2}
                                                                < \epsilon,
\end{equation}
which exists because $F \in \ell^2(\widehat{A})$.  Thus
\begin{equation}
\label{||f - sum_{phi in E} F(phi) phi||_{L^2(A)} < epsilon}
 \biggl\|f - \sum_{\phi \in E} F(\phi) \, \phi\biggr\|_{L^2(A)} < \epsilon
\end{equation}
when $E_2(\epsilon) \subseteq E \subseteq \widehat{A}$, by (\ref{||f -
  sum_{phi in E} F(phi) phi||_{L^2(A)} = ...}).

        Similarly, let $\ell^1(\widehat{A})$ be the space of complex-valued
functions $F(\phi)$ on $\widehat{A}$ that are summable, as in the previous
section.  If $F \in \ell^1(\widehat{A})$, then
\begin{equation}
\label{f(x) = sum_{phi in widehat{A}} F(phi) phi(x)}
        f(x) = \sum_{\phi \in \widehat{A}} F(\phi) \, \phi(x)
\end{equation}
can be defined for every $x \in A$, and has modulus less than or equal
to the $\ell^1$ norm of $F$.  If $E$ is any subset of $\widehat{A}$,
then we get that
\begin{equation}
\label{|f(x) - sum_{phi in E} F(phi) phi(x)| le ...}
        \biggl|f(x) - \sum_{\phi \in E} F(\phi) \, \phi(x)\biggr|
  = \biggl|\sum_{\phi \in \widehat{A} \setminus E} F(\phi) \, \phi(x)\biggr|
         \le \sum_{\phi \in \widehat{A} \setminus E} |F(\phi)|
\end{equation}
for every $x \in A$.  The summability of $F$ on $\widehat{A}$ also
implies that for each $\epsilon > 0$ there is a finite subset
$E_1(\epsilon)$ of $\widehat{A}$ such that
\begin{equation}
\label{sum_{phi in widehat{A} setminus E_1(epsilon)} |F(phi)| < epsilon}
        \sum_{\phi \in \widehat{A} \setminus E_1(\epsilon)} |F(\phi)| < \epsilon.
\end{equation}
If $E_1(\epsilon) \subseteq E \subseteq \widehat{A}$, then it follows
that
\begin{equation}
\label{|f(x) - sum_{phi in E} F(phi) phi(x)| < epsilon}
 \biggl|f(x) - \sum_{\phi \in E} F(\phi) \, \phi(x)\biggr| < \epsilon
\end{equation}
for every $x \in A$, by (\ref{|f(x) - sum_{phi in E} F(phi) phi(x)| le
  ...}).  In particular, this implies that $f$ is continuous on $A$,
because $f$ can be approximated by continuous functions uniformly on
$A$.  Remember that $\ell^1(\widehat{A}) \subseteq
\ell^2(\widehat{A})$, as in the previous section.  It is easy to see
that (\ref{f(x) = sum_{phi in widehat{A}} F(phi) phi(x)}) as the same
as (\ref{f = sum_{phi in widehat{A}} F(phi) phi}) as an element of
$L^2(A)$ when $F \in \ell^1(\widehat{A})$, because this type of
uniform convergence on $A$ implies convergence with respect to the
$L^2$ norm on $A$.

        If $F \in \ell^2(\widehat{A})$ and $f$ is as in
(\ref{f = sum_{phi in widehat{A}} F(phi) phi}), then
\begin{equation}
\label{widehat{f}(phi) = langle f, phi rangle_{L^2(A)} = F(phi)}
        \widehat{f}(\phi) = \langle f, \phi \rangle_{L^2(A)} = F(\phi)
\end{equation}
for every $\phi \in \widehat{A}$, where $\langle \cdot, \cdot
\rangle_{L^2(A)}$ is the standard integral inner product on $L^2(A)$.
If $f$ is any element of $L^2(A)$, then $\widehat{f}(\phi)$ is in
$\ell^2(A)$, and (\ref{f = sum_{phi in widehat{A}} widehat{f}(phi)
  phi}) is the same as (\ref{f = sum_{phi in widehat{A}} F(phi) phi}),
with $F(\phi) = \widehat{f}(\phi)$.  If $F = \widehat{f} \in
\ell^1(\widehat{A})$, then $f$ is the same as (\ref{f(x) = sum_{phi in
    widehat{A}} F(phi) phi(x)}) as an element of $L^2(A)$, so that
$f(x)$ is equal to (\ref{f(x) = sum_{phi in widehat{A}} F(phi)
  phi(x)}) for almost every $x \in A$ with respect to $H$.  In
particular, this implies that $f$ is equal to a continuous function on
$A$ almost everywhere on $A$ with respect to $H$.

        Suppose that $\mu$ is a regular complex Borel measure on $A$
such that
\begin{equation}
\label{F(phi) = widehat{mu}(phi)}
        F(\phi) = \widehat{\mu}(\phi)
\end{equation}
is in $\ell^2(\widehat{A})$.  Thus we can define $f \in L^2(A)$ as in
(\ref{f = sum_{phi in widehat{A}} F(phi) phi}), which satisfies
(\ref{widehat{f}(phi) = langle f, phi rangle_{L^2(A)} = F(phi)}).
If we put
\begin{equation}
\label{mu_f(E) = int_E f(x) dH(x)}
        \mu_f(E) = \int_E f(x) \, dH(x)
\end{equation}
for every Borel set $E \subseteq A$, then $\mu_f$ is also a regular
complex-valued Borel measure on $A$, and
\begin{equation}
\label{widehat{mu_f}(phi) = widehat{f}(phi) = F(phi) = widehat{mu}(phi)}
 \widehat{\mu_f}(\phi) = \widehat{f}(\phi) = F(\phi) = \widehat{\mu}(\phi)
\end{equation}
for every $\phi \in \widehat{A}$, by construction.  This implies that
the Fourier transform of $\mu - \mu_f$ is equal to $0$ on
$\widehat{A}$, and hence that $\mu - \mu_f = 0$, as discussed earlier.
Similarly, if (\ref{F(phi) = widehat{mu}(phi)}) is in $\ell^1(A)$,
then $f$ can be defined as in (\ref{f(x) = sum_{phi in widehat{A}}
  F(phi) phi(x)}), so that $f$ is a continuous function on $A$ in this
case.  If $g$ is an integrable function on $A$ with respect to $H$,
then we can apply these remarks to the Borel measure defined by
\begin{equation}
\label{mu(E) = mu_g(E) = int_E g(x) dH(x)}
        \mu(E) = \mu_g(E) = \int_E g(x) \, dH(x)
\end{equation}
for each Borel set $E \subseteq A$.  It follows that $g \in L^2(A)$
when $\widehat{g} \in \ell^2(\widehat{A})$, and that $g$ is equal to a
continuous function on $A$ almost everywhere with respect to $H$ when
$\widehat{g} \in \ell^1(A)$.

\section{Compact groups, continued}
\label{compact groups, continued}

        Let $A$ be a compact commutative topological group again,
with normalized Haar measure $H$.  If $f, g \in L^2(A)$, then
we have seen that their convolution $f * g$ is a continuous function
on $A$, and that
\begin{equation}
\label{widehat{(f * g)}(phi) = widehat{f}(phi) widehat{g}(phi)}
        \widehat{(f * g)}(\phi) = \widehat{f}(\phi) \, \widehat{g}(\phi)
\end{equation}
for every $\phi \in \widehat{A}$.  We also have that $\widehat{f},
\widehat{g} \in \ell^2(\widehat{A})$, as in the previous section,
which implies that
\begin{equation}
\label{sum_{phi in widehat{A}} |widehat{f}(phi)| |widehat{g}(phi)| le ...}
 \sum_{\phi \in \widehat{A}} |\widehat{f}(\phi)| \, |\widehat{g}(\phi)|
 \le \Big(\sum_{\phi \in \widehat{A}} |f(\phi)|^2\Big)^{1/2}
      \, \Big(\sum_{\phi \in \widehat{A}} |g(\phi)|^2\Big)^{1/2}
\end{equation}
is finite, by the Cauchy--Schwarz inequality.  In particular,
$(f * g)(x)$ can be approximated uniformly on $A$ by sums of the form
\begin{equation}
\label{sum_{phi in E} widehat{f}(phi) widehat{g}(phi) phi(x)}
        \sum_{\phi \in E} \widehat{f}(\phi) \, \widehat{g}(\phi) \, \phi(x),
\end{equation}
where $E$ is a finite subset of $\widehat{A}$.

        Alternatively, if $f \in L^1(A)$ and $g \in L^2(A)$, then
$f * g \in L^2(A)$, and hence $f * g$ can be approximated by finite
sums of the form (\ref{sum_{phi in E} widehat{f}(phi) widehat{g}(phi)
  phi(x)}) with respect to the $L^2$ norm on $A$.  If $\widehat{g} \in
\ell^1(\widehat{A})$, then $\widehat{(f * g)} \in \ell^1(\widehat{A})$
too, because of (\ref{widehat{(f * g)}(phi) = widehat{f}(phi) widehat{g}(phi)})
and the fact that $\widehat{f}$ is bounded on $\widehat{A}$.  In this
case, $f * g$ can be approximated by finite sums of the form 
(\ref{sum_{phi in E} widehat{f}(phi) widehat{g}(phi) phi(x)}) uniformly
on $A$, as before.  There are analogous statements in which $f$
is replace by a regular complex Borel measure on $A$, using essentially
the same arguments.

        Let $U$ be an open set in $A$ that contains $0$, and let $g$
be a nonnegative real-valued function on $A$ such that
\begin{equation}
\label{int_A g(x) dH(x) = 1}
        \int_A g(x) \, dH(x) = 1
\end{equation}
and $g(x) = 0$ for every $x \in A \setminus U$.  We may as well take $g$
to be a continuous function on $A$, although it would suffice to have
$g \in L^2(A)$ for some of the arguments that follow.  We can also
choose $g$ so that $\widehat{g} \in \ell^1(\widehat{A})$, by taking
$g$ to be a convolution of functions in $L^2(A)$.  More precisely,
the continuity of addition on $A$ at $0$ implies that there are
open sets $U_1, U_2 \subseteq A$ that contain $0$ and satisfy
\begin{equation}
\label{U_1 + U_2 subseteq U}
        U_1 + U_2 \subseteq U.
\end{equation}
Let $g_1$, $g_2$ be nonnegative real-valued functions on $A$ in
$L^2(A)$ such that
\begin{equation}
\label{int_A g_1(x) dH(x) = int_A g_2(x) dH(x) = 1}
        \int_A g_1(x) \, dH(x) = \int_A g_2(x) \, dH(x) = 1
\end{equation}
and $g_j(x) = 0$ when $x \in A \setminus U_j$ for $j = 1, 2$.  If we
take $g = g_1 * g_2$, then $g$ is a nonnegative real-valued continuous
function on $A$, as before.  The integral of $g$ over $A$ with respect
to $H$ is equal to the product of the integrals of $g_1$ and $g_2$, as
in Section \ref{convolution of functions}, and hence is equal to $1$,
by (\ref{int_A g_1(x) dH(x) = int_A g_2(x) dH(x) = 1}).  It is easy to
see that $g = 0$ on $A \setminus U$ under these conditions, because of
(\ref{U_1 + U_2 subseteq U}) and the analogous properties of $g_1$ and
$g_2$.

        Let $f \in L^p(A)$ be given, for some $p \in [1, \infty)$.
As in Section \ref{some more estimates}, $f * g$ is close to $f$
with respect to the $L^p$ norm on $A$ when $U$ is sufficiently small.
More precisely, this corresponds to taking $\nu$ to be the Borel measure
on $A$ defined by
\begin{equation}
\label{nu(E) = int_E g(x) dH(x)}
        \nu(E) = \int_E g(x) \, dH(x)
\end{equation}
for every Borel set $E \subseteq A$ in the earlier discussion.  If
$\widehat{g} \in \ell^1(\widehat{A})$, then we have seen that $f * g$
can be approximated by finite sums of the form (\ref{sum_{phi in E}
  widehat{f}(phi) widehat{g}(phi) phi(x)}) uniformly on $A$.  This
implies that $f$ is approximated by finite sums of the form
(\ref{sum_{phi in E} widehat{f}(phi) widehat{g}(phi) phi(x)}) with
respect to the $L^p$ norm when $U$ is sufficiently small.

        If $f \in L^2(A)$, then $f$ can already be approximated by finite
subsums of (\ref{f = sum_{phi in widehat{A}} widehat{f}(phi) phi})
with respect to the $L^2$ norm, as in (\ref{||f(x) - sum_{phi in E}
  widehat{f}(phi) phi||_2 < epsilon}).  If $f \in L^p(A)$ and $1 \le p
< 2$, then it suffices to have $g \in L^2(A)$ in the argument in the
preceding paragraph, instead of $\widehat{g} \in \ell^1(\widehat{A})$.
This is because $f * g$ can be approximated by finite sums of the form
(\ref{sum_{phi in E} widehat{f}(phi) widehat{g}(phi) phi(x)}) with
respect to the $L^2$ norm, as mentioned earlier, which implies
approximation with respect to the $L^p$ norm when $p \le 2$.  If $f
\in L^p(A)$ and $2 < p < \infty$, then it also suffices to have $g \in
L^2(A)$ in the argument in the preceding paragraph, instead of
$\widehat{g} \in \ell^1(\widehat{A})$.  In this case, we have that $f
\in L^2(A)$ in particular, so that (\ref{widehat{(f * g)}(phi) =
  widehat{f}(phi) widehat{g}(phi)}) is in $\ell^1(\widehat{A})$,
and the rest of the argument is basically the same as before.

        If $f$ is a continuous function on $A$, then $f * g$ is close
to $f$ uniformly on $A$ when $U$ is sufficiently small.  This implies
that $f$ can be approximated by finite sums of the form (\ref{sum_{phi
    in E} widehat{f}(phi) widehat{g}(phi) phi(x)}) uniformly on $A$,
because $f * g$ can be approximated by finite sums of the form
(\ref{sum_{phi in E} widehat{f}(phi) widehat{g}(phi) phi(x)})
uniformly on $A$ when $g \in L^2(A)$.

        If $\widehat{f} \in \ell^1(\widehat{A})$, then we can already 
express $f$ as in (\ref{f(x) = sum_{phi in widehat{A}} F(phi)
  phi(x)}), with $F = \widehat{f}$.  In this case, $f$ can be
approximated by finite subsums of (\ref{f(x) = sum_{phi in widehat{A}}
  F(phi) phi(x)}) uniformly on $A$, as in (\ref{|f(x) - sum_{phi in E}
  F(phi) phi(x)| < epsilon}).  One can also check that
(\ref{widehat{(f * g)}(phi) = widehat{f}(phi) widehat{g}(phi)}) is
close to $\widehat{f}(\phi)$ with respect to the $\ell^1$ norm on
$\widehat{A}$ when $U$ is sufficiently small.  Indeed, for each $\phi
\in \widehat{A}$, $\widehat{g}(\phi)$ is close to $1$ when $U$ is
sufficiently small, because $\phi(0) = 1$ and $\phi$ is continuous at
$0$.  We also know that $|\widehat{g}(\phi)| \le 1$ for every $\phi
\in \widehat{A}$, since the $L^1$ norm of $g$ is equal to $1$, which
permits one to reduce to approximating $\widehat{f}(\phi)$ by
(\ref{widehat{(f * g)}(phi) = widehat{f}(phi) widehat{g}(phi)})
on finite subsets of $\widehat{A}$.

\section{Products of functions}
\label{products of functions, 2}

        Let us continue to take $A$ to be a compact commutative topological
group, with normalized Haar measure $H$.  Also let $\mathcal{E}(A)$ be
the linear span of $\widehat{A}$ in $C(A)$, as usual.  If $f, g \in
\mathcal{E}(A)$, then
\begin{equation}
\label{sums for f(x), g(x)}
        f(x) = \sum_{\phi \in \widehat{A}} \widehat{f}(\phi) \, \phi(x), \quad
        g(x) = \sum_{\psi \in \widehat{A}} \widehat{g}(\psi) \, \psi(x)
\end{equation}
for every $x \in A$, where all but finitely many terms in the sums
are equal to $0$.  This implies that
\begin{eqnarray}
\label{f(x) g(x) = ...}
 f(x) \, g(x) & = & \sum_{\psi \in \widehat{A}} \, \sum_{\phi \in \widehat{A}}
               \widehat{f}(\phi) \, \widehat{g}(\psi) \, \phi(x) \, \psi(x) \\
 & = & \sum_{\psi \in \widehat{A}} \, \sum_{\phi \in \widehat{A}}
 \widehat{f}(\phi \, \overline{\psi}) \, \widehat{g}(\psi) \, \phi(x) \nonumber
\end{eqnarray}
for every $x \in A$.  More precisely, this uses the change of
variables $\phi \mapsto \phi \, \overline{\psi}$ in the second step,
which is a translation in the group $\widehat{A}$.  Interchanging the
order of summation, we get that
\begin{equation}
\label{f(x) g(x) = ..., 2}
 f(x) \, g(x) = \sum_{\phi \in \widehat{A}} \Big(\sum_{\psi \in \widehat{A}}
     \widehat{f}(\phi \, \overline{\psi}) \, \widehat{g}(\psi)\Big) \, \phi(x)
\end{equation}
for every $x \in A$.  It follows that
\begin{equation}
\label{widehat{(f g)}(phi) = ...}
 \widehat{(f \, g)}(\phi) = \sum_{\psi \in \widehat{A}}
                \widehat{f}(\phi \, \overline{\psi}) \, \widehat{g}(\psi)
\end{equation}
for every $\phi \in \widehat{A}$.

        Suppose now that $f, g \in L^2(A)$, so that $\widehat{f}, \widehat{g}
\in \ell^2(\widehat{A})$, as in Section \ref{compact groups}.
Using the Cauchy--Schwarz inequality, we get that
\begin{equation}
\label{sum_{psi} |widehat{f}(phi overline{psi})| |widehat{g}(psi)| le ...}
 \sum_{\psi \in \widehat{A}} |\widehat{f}(\phi \, \overline{\psi})| \,
                                                   |\widehat{g}(\psi)|
 \le \Big(\sum_{\psi \in \widehat{A}}
          |\widehat{f}(\phi \, \overline{\psi})|^2\Big)^{1/2} \,
 \Big(\sum_{\psi \in \widehat{A}} |\widehat{g}(\psi)|^2\Big)^{1/2}
\end{equation}
for every $\phi \in \widehat{A}$.  Note that
\begin{equation}
\label{sum_{psi in widehat{A}} |widehat{f}(phi overline{psi})|^2 = ...}
        \sum_{\psi \in \widehat{A}} |\widehat{f}(\phi \, \overline{\psi})|^2
                          = \sum_{\psi \in \widehat{A}} |\widehat{f}(\psi)|^2
\end{equation}
for every $\phi \in \widehat{A}$, using the change of variables $\psi
\mapsto \phi \, \overline{\psi}$.  Thus
\begin{equation}
\label{sum_{psi} |widehat{f}(phi overline{psi})| |widehat{g}(psi)| le ..., 2}
 \sum_{\psi \in \widehat{A}} |\widehat{f}(\phi \, \overline{\psi})| \,
                                                   |\widehat{g}(\psi)|
                                        \le \|f\|_{L^2(A)} \, \|g\|_{L^2(A)}
\end{equation}
for every $\phi \in \widehat{A}$, by (\ref{sum_{phi in widehat{A}}
  |widehat{f}(phi)|^2 le ||f||_2^2}).  Of course, $f \, g \in L^1(A)$
when $f, g \in L^2(A)$, so that
\begin{equation}
\label{widehat{(f g)}(phi) = int_A f(x) g(x) overline{phi(x)} dH(x)}
 \widehat{(f \, g)}(\phi) = \int_A f(x) \, g(x) \, \overline{\phi(x)} \, dH(x)
\end{equation}
is defined for every $\phi \in \widehat{A}$.  If $g \in
\mathcal{E}(A)$, then it is easy to see that (\ref{widehat{(f g)}(phi)
  = ...}) holds for every $\phi \in \widehat{A}$, where all but
finitely many terms in the sum on the right side of (\ref{widehat{(f
    g)}(phi) = ...})  are equal to $0$.  Using this and
(\ref{sum_{psi} |widehat{f}(phi overline{psi})| |widehat{g}(psi)| le
  ..., 2}), one can check that (\ref{sums for f(x), g(x)}) holds for
every $f, g \in L^2(A)$ and $\phi \in \widehat{A}$, by approximating
$g$ by elements of $\mathcal{E}(A)$ with respect to the $L^2$ norm.

        If $g \in L^1(A)$ and $\widehat{g} \in \ell^1(\widehat{A})$,
then
\begin{equation}
\label{sum_{phi in widehat{A}} widehat{g}(phi) phi(x)}
        \sum_{\phi \in \widehat{A}} \widehat{g}(\phi) \, \phi(x)
\end{equation}
can be defined for every $x \in A$, and is uniformly approximated by
finite subsums, as in Section \ref{compact groups}.  This implies that
(\ref{sum_{phi in widehat{A}} widehat{g}(phi) phi(x)}) defines a
continuous function on $A$, and we also have that $g(x)$ is equal to
(\ref{sum_{phi in widehat{A}} widehat{g}(phi) phi(x)}) for almost
every $x \in A$ with respect to $H$, by the remarks at the end of
Section \ref{compact groups}.  If $f \in L^1(A)$, then $\widehat{f}$
is bounded on $\widehat{A}$, as in (\ref{|widehat{f}(phi)| le int_A
  |f(x)| |phi(x)| dH(x) = int_A |f(x)| dH(x)}), and
\begin{eqnarray}
 \sum_{\psi \in \widehat{A}} |\widehat{f}(\phi \, \overline{\psi})| \,
                                                   |\widehat{g}(\psi)|
 & \le & \Big(\sup_{\psi \in \widehat{A}} |\widehat{f}(\psi)|\Big) \,
          \Big(\sum_{\psi \in \widehat{A}} |\widehat{g}(\psi)|\Big) \\
 & \le & \|f\|_{L^1(A)} \, \|\widehat{g}\|_{\ell^1(\widehat{A})} \nonumber
\end{eqnarray}
for every $\phi \in \widehat{A}$.  We also have that $f \, g \in
L^1(A)$ in this situation, so that (\ref{widehat{(f g)}(phi) = int_A
  f(x) g(x) overline{phi(x)} dH(x)}) is defined for every $\phi \in
\widehat{A}$.  As before, one can check that (\ref{widehat{(f g)}(phi)
  = ...}) holds for every $\phi \in \widehat{A}$ under these
conditions, by approximating $g$ by elements of $\mathcal{E}(A)$.

        Suppose that $f, g \in L^1(A)$ satisfy $\widehat{f}, \widehat{g}
\in \ell^1(\widehat{A})$, and let us check that $\widehat{(f \, g)}$ is
in $\ell^1(\widehat{A})$ too.  Of course,
\begin{equation}
\label{|widehat{(f g)}(phi)| le ...}
 |\widehat{(f \, g)}(\phi)| \le \sum_{\psi \in \widehat{A}}
             |\widehat{f}(\phi \, \overline{\psi})| \, |\widehat{g}(\psi)|
\end{equation}
for every $\phi \in \widehat{A}$, by (\ref{widehat{(f g)}(phi) =
  ...}), and hence
\begin{equation}
\label{sum_{phi in widehat{A}} |widehat{(f g)}(phi)| le ...}
        \sum_{\phi \in \widehat{A}} |\widehat{(f \, g)}(\phi)|
         \le \sum_{\phi \in \widehat{A}} \, \sum_{\psi \in \widehat{A}}
              |\widehat{f}(\phi \, \overline{\psi})| \, |\widehat{g}(\psi)|.
\end{equation}
Interchanging the order of summation, we get that
\begin{eqnarray}
\label{sum_{phi in widehat{A}} |widehat{(f g)}(phi)| le ..., 2}
 \sum_{\phi \in \widehat{A}} |\widehat{(f \, g)}(\phi)|
 & \le & \sum_{\psi \in \widehat{A}} \, \sum_{\phi \in \widehat{A}}
             |\widehat{f}(\phi \, \overline{\psi})| \, |\widehat{g}(\psi)| \\
 & = & \sum_{\psi \in \widehat{A}} \, \sum_{\phi \in \widehat{A}}
             |\widehat{f}(\phi)| \, |\widehat{g}(\psi)| \nonumber \\
 & = & \Big(\sum_{\phi \in \widehat{A}} |\widehat{f}(\phi)|\Big) \,
        \Big(\sum_{\psi \in \widehat{A}} |\widehat{g}(\psi)|\Big), \nonumber
\end{eqnarray}
using the change of variables $\phi \mapsto \phi \, \psi$ in the
second step.  One could also simply start with $f$ as in (\ref{f(x) =
  sum_{phi in widehat{A}} F(phi) phi(x)}), with coefficients $F(\phi)
\in \ell^1(\widehat{A})$, and similarly for $g$.  One could then
multiply the expansions for $f$ and $g$ directly, as in (\ref{f(x)
  g(x) = ...})  and (\ref{f(x) g(x) = ..., 2}), and estimate the
coefficients in the resulting expansion for $f \, g$ as before.

\section{Translation-invariant subspaces}
\label{translation-invariant subspaces}

        Let $A$ be a locally compact commutative topological group
with a Haar measure $H$, and let $1 \le p < \infty$ be given.
Also let $L^p(A)$ be the usual $L^p$ space associated to $H$,
with the $L^p$ norm $\|f\|_p$.  As before, put $f_a(x) = f(x - a)$
for each $f \in L^p(A)$ and $a \in A$, so that
\begin{equation}
\label{||f_a||_p = ||f||_p}
        \|f_a\|_p = \|f\|_p,
\end{equation}
because $H$ is invariant under translations.  If $f$ is a continuous
function on $A$ with compact support, then
\begin{equation}
\label{lim_{a to 0} ||f_a - f||_p = 0}
        \lim_{a \to 0} \|f_a - f\|_p = 0,
\end{equation}
where the limit is taken over $a \in A$ converging to $0$ with respect
to the given topology on $A$.  This is the same as (\ref{lim_{a to 0}
  int_A |f(x) - f(x - a)| dH(x) = 0}) when $p = 1$, and a similar
argument works for all $p$, using the fact that $f$ is uniformly
continuous on $A$, as in Section \ref{uniform continuity,
  equicontinuity}.  As in the $p = 1$ case, one can check that
(\ref{lim_{a to 0} ||f_a - f||_p = 0}) holds for every $f \in L^p(A)$
when $1 \le p < \infty$, by approximating $f$ by continuous functions
with compact support in $A$ with respect to the $L^p$ norm.  Note that
\begin{equation}
\label{||f_{a + b} - f_b||_p = ||f_a - f||_p}
        \|f_{a + b} - f_b\|_p = \|f_a - f\|_p
\end{equation}
for every $a, b \in A$, by translation-invariance of the $L^p$ norm,
as in (\ref{||f_a||_p = ||f||_p}).  It follows that the mapping from
$a \in A$ to $f_a \in L^p(A)$ is uniformly continuous, by (\ref{lim_{a
    to 0} ||f_a - f||_p = 0}) and (\ref{||f_{a + b} - f_b||_p = ||f_a
  - f||_p}).

        A closed linear subspace $\mathcal{L}$ of $L^p(A)$ is said to be
\emph{invariant under translations}\index{translation-invariant subspaces}
if for every $f \in \mathcal{L}$ and $a \in A$,
\begin{equation}
\label{f_a(x) = f(x - a) in mathcal{L}}
        f_a(x) = f(x - a) \in \mathcal{L}.
\end{equation}
If $\mathcal{L}$ has this property and $g \in L^1(A)$, then one can
check that
\begin{equation}
\label{f * g in mathcal{L}}
        f * g \in \mathcal{L}
\end{equation}
for every $f \in \mathcal{L}$.  Of course, $f * g$ is an integral of
translates of $f$, by definition, and so it suffices to show that this
integral can be approximated by linear combinations of translates of
$f$ with respect to the $L^p$ norm.  To do this, one can use the
uniform continuity of $a \mapsto f_a$ as a mapping from $A$ into
$L^p(A)$, as in the preceding paragraph.  The approximation is quite
straightforward when $g$ has compact support in $A$, and otherwise one
can approximate $g$ by integrable functions with compact support with
respect to the $L^1$ norm.

        Conversely, suppose that $\mathcal{L}$ is a closed linear subspace
of $L^p(A)$, $1 \le p < \infty$, such that (\ref{f * g in mathcal{L}})
holds for every $f \in \mathcal{L}$ and $g \in L^1(A)$.  Observe that
\begin{equation}
\label{f_a * g = f * g_a}
        f_a * g = f * g_a,
\end{equation}
for every $f \in L^p(A)$, $g \in L^1(A)$, and $a \in A$, by the
definition of the convolution in Section \ref{convolution of
  functions}.  It follows that
\begin{equation}
\label{f_a * g in mathcal{L}}
        f_a * g \in \mathcal{L}
\end{equation}
for every $f \in \mathcal{L}$, $g \in L^1(A)$, and $a \in A$, by our
hypothesis on $\mathcal{L}$.  Let $f \in \mathcal{L}$ and $a \in A$ be
given, and let $g$ be a nonnegative real-valued integrable function on
$A$, with
\begin{equation}
\label{int_A g(x) dH(x) = 1, 2}
        \int_A g(x) \, dH(x) = 1.
\end{equation}
Also let $U$ be an open set in $A$ that contains $0$, and suppose that
$g = 0$ almost everywhere on $A \setminus U$ with respect to $H$.  As
in Section \ref{some more estimates}, $f_a * g$ approximates $f_a$
with respect to the $L^p$ norm when $U$ is sufficiently small.  This
implies that $f_a \in \mathcal{L}$, because of (\ref{f_a * g in
  mathcal{L}}), and since $\mathcal{L}$ is supposed to be a closed set
in $L^p(A)$.  Thus $\mathcal{L}$ is invariant under translations
under these conditions.

        Let $\mathcal{L}$ be a closed linear subspace of $L^p(A)$,
$1 \le p < \infty$, which is invariant under translations again.
If $\nu$ is a regular complex Borel measure on $A$, then
\begin{equation}
\label{f * nu in mathcal{L}}
        f * \nu \in \mathcal{L}
\end{equation}
for every $f \in \mathcal{L}$, for essentially the same reasons as in
(\ref{f * g in mathcal{L}}).  More precisely, if $\nu$ has compact
support, then it is easy to approximate $f * \nu$ by linear
combinations of translates of $f$ with respect to the $L^p$ norm,
using the uniform continuity of $a \mapsto f_a$ as a mapping from $A$
into $L^p(A)$.  Otherwise, one can approximate $\nu$ by measures with
compact support, as before.

        Now let $\mathcal{L}$ be a linear subspace of the space $C_b(A)$
of bounded continuous functions on $A$, and suppose that $\mathcal{L}$
is a closed set with respect to the supremum norm on $C_b(A)$.  As
before, $\mathcal{L}$ is said to be invariant under
translations\index{translation-invariant subspaces} if (\ref{f_a(x) =
  f(x - a) in mathcal{L}}) for every $f \in \mathcal{L}$ and $a \in A$.
If the elements of $\mathcal{L}$ are uniformly continuous on $A$,
then it is easy to see that the previous arguments carry over to
this situation.  In particular, this applies to closed linear subspaces
of $C_0(A)$, and to closed linear subspaces of $C(A)$ when $A$ is compact.

        Note that $L^1(A)$ is a commutative Banach algebra, using
convolution of functions as multiplication.  A linear subspace
$\mathcal{L}$ of $L^1(A)$ is an \emph{ideal}\index{ideals in
  L^1(A)@ideals in $L^1(A)$} in $L^1(A)$ as a commutative Banach
algebra if (\ref{f * g in mathcal{L}}) holds for every $f \in
\mathcal{L}$ and $g \in L^1(A)$.  Thus the closed ideals in $L^1(A)$
are the same as the closed linear subspaces of $L^1(A)$ that are
invariant under translations, by the earlier discussion.

        If $E$ is any subset of $\widehat{A}$, then it is easy to see that
\begin{equation}
\label{mathcal{L}_E = {f in L^1(A) : widehat{f}(phi) = 0 for every phi in E}}
        \mathcal{L}_E = \{f \in L^1(A) : \widehat{f}(\phi) = 0
                                       \hbox{ for every } \phi \in E\}
\end{equation}
is a closed ideal in $L^1(A)$.  Remember that $\widehat{f}$ is
continuous with respect to the usual topology on $\widehat{A}$ for
each $f \in L^1(A)$, as in Section \ref{fourier transforms, 2}.
This implies that
\begin{equation}
\label{mathcal{L}_{overline{E}} = mathcal{L}_E}
        \mathcal{L}_{\overline{E}} = \mathcal{L}_E
\end{equation}
for every $E \subseteq \widehat{A}$, where $\overline{E}$ is the
closure of $E$ in $\widehat{A}$.  Thus we may as well restrict our
attention to closed sets $E \subseteq \widehat{A}$ here.  If
$\mathcal{L}$ is any subset of $L^1(A)$, then
\begin{equation}
\label{E(mathcal{L}) = ...}
        E(\mathcal{L}) = \{\phi \in \widehat{A} : \widehat{f}(\phi) = 0
                                         \hbox{ for every } f \in \mathcal{L}\}
\end{equation}
is a closed subset of $\widehat{A}$, because $\widehat{f}$ is
continuous on $\widehat{A}$ for each $f \in L^1(A)$, and
\begin{equation}
\label{mathcal{L} subseteq mathcal{L}_{E(mathcal{L})}}
        \mathcal{L} \subseteq \mathcal{L}_{E(\mathcal{L})},
\end{equation}
by construction.

        Let us restrict our attention to the case where $A$ is compact
for the rest of the section, and take Haar measure $H$ on $A$ to
be normalized so that $H(A) = 1$.  Let $E$ be any subset of
$\widehat{A}$, which is automatically a closed set, because the usual
topology on $\widehat{A}$ is discrete when $A$ is compact.  If $\phi,
\psi \in \widehat{A}$ and $\phi \ne \psi$, then $\phi$ is orthogonal
to $\psi$ in $L^2(A)$, as in Section \ref{compact, discrete groups},
which is equivalent to saying that
\begin{equation}
        \widehat{\phi}(\psi) = 0.
\end{equation}
If $\phi \in \widehat{A} \setminus E$, then it follows that $\phi \in
\mathcal{L}_E$, so that $\widehat{A} \setminus E \subseteq
\mathcal{L}_E$.  This implies that the closure of the linear span of
$\widehat{A} \setminus E$ in $L^1(A)$ is also contained in
$\mathcal{L}(E)$.  If $f$ is any integrable function on $A$, then $f$
can be approximated by linear combinations of $\phi \in \widehat{A}$
such that $\widehat{f}(\phi) \ne 0$, as in Section \ref{compact
  groups, continued}.  If $f \in \mathcal{L}_E$, then we get that
$f$ is in the closure of the linear span of $\widehat{A} \setminus E$,
by the definition of $\mathcal{L}_E$.  This shows that $\mathcal{L}_E$
is equal to the closure of the linear span of $\widehat{A} \setminus E$
in $L^1(A)$ when $A$ is compact.

        Remember that
\begin{equation}
\label{f * phi = widehat{f}(phi) phi}
        f * \phi = \widehat{f}(\phi) \, \phi
\end{equation}
for every $f \in L^1(A)$ and $\phi \in \widehat{A}$, as in (\ref{(phi
  * nu)(x) = widehat{nu}(phi) phi(x)}) in Section \ref{functions,
  measures}.  Let $\mathcal{L}$ be a closed ideal in $L^1(A)$, and let
$E(\mathcal{L})$ be as in (\ref{mathcal{L} subseteq
  mathcal{L}_{E(mathcal{L})}}).  Thus (\ref{f * phi = widehat{f}(phi)
  phi}) is an element of $\mathcal{L}$ for every $f \in \mathcal{L}$
and $\phi \in \widehat{A}$.  If $\phi \in \widehat{A} \setminus
E(\mathcal{L})$, then $\widehat{f}(\phi) \ne 0$ for some $f \in
\mathcal{L}$, and it follows that $\phi \in \mathcal{L}$.  This
implies that the closure of the linear span of $\widehat{A} \setminus
E(\mathcal{L})$ in $L^1(A)$ is contained in $\mathcal{L}$, because
$\mathcal{L}$ is a closed linear subspace of $L^1(A)$.  We also know
that the closure of the linear span of $\widehat{A} \setminus
E(\mathcal{L})$ in $L^1(A)$ is equal to
$\mathcal{L}_{E(\mathcal{L})}$, as in the previous paragraph.
Combining this with (\ref{mathcal{L} subseteq mathcal{L}_{E(mathcal{L})}}),
we get that
\begin{equation}
\label{mathcal{L} = mathcal{L}_{E(mathcal{L})}}
        \mathcal{L} = \mathcal{L}_{E(\mathcal{L})}
\end{equation}
in this case.

        If $E$ is any subset of $\widehat{A}$ and $1 < p < \infty$, then
it is easy to see that
\begin{equation}
\label{mathcal{L}_E cap L^p(A)}
        \mathcal{L}_E \cap L^p(A)
\end{equation}
is a closed linear subspace of $L^p(A)$ which is invariant under
translations.  One can also check that (\ref{mathcal{L}_E cap L^p(A)})
is the same as the closure of the linear span of $\widehat{A}
\setminus E$ in $L^p(A)$ under these conditions, for essentially the
same reasons as when $p = 1$.  If $\mathcal{L}$ is any closed linear
subspace of $L^p(A)$ that is invariant under translations, where $1 <
p < \infty$, and if $E(\mathcal{L})$ is as in (\ref{E(mathcal{L}) =
  ...}), then we have that
\begin{equation}
\label{mathcal{L} = mathcal{L}_{E(mathcal{L})} cap L^p(A)}
        \mathcal{L} = \mathcal{L}_{E(\mathcal{L})} \cap L^p(A),
\end{equation}
as in the preceding paragraph.  Similarly,
\begin{equation}
\label{mathcal{L}_E cap C(A)}
        \mathcal{L}_E \cap C(A)
\end{equation}
is a closed linear subspace of $C(A)$ which is invariant under
translations for each $E \subseteq \widehat{A}$, and
(\ref{mathcal{L}_E cap C(A)}) is the same as the closure of the linear
span of $\widehat{A} \setminus E$ in $C(A)$ with respect to the
supremum norm.  If $\mathcal{L}$ is any closed linear subspace of
$C(A)$ with respect to the supremum norm that is invariant under
translations, then
\begin{equation}
\label{mathcal{L} = mathcal{L}_{E(mathcal{L})} cap C(A)}
        \mathcal{L} = \mathcal{L}_{E(\mathcal{L})} \cap C(A),
\end{equation}
where $E(\mathcal{L})$ is as in (\ref{E(mathcal{L}) = ...}).

\section{Translation-invariant subalgebras}
\label{translation-invariant subalgebras}

        Let $A$ be a compact commutative topological group, with
Haar measure $H$ on $A$ normalized so that $H(A) = 1$.  As usual,
$\mathcal{E}(A)$ denotes the linear span of $\widehat{A}$ in $C(A)$,
which is a subalgebra of $C(A)$ that contains the constant functions
and is invariant under complex conjugation.  Similarly, let $B$ be a
subset of $\widehat{A}$, and let $\mathcal{E}_B(A)$ be the linear span
of $B$ in $C(A)$, which is interpreted as being $\{0\}$ when $B =
\emptyset$.  Also let $C_B(A)$ be the closure of $\mathcal{E}_B(A)$ in
$C(A)$, with respect to the supremum norm.

        If $f \in C_B(A)$, then it is easy to see that
\begin{equation}
\label{widehat{f}(phi) = 0}
        \widehat{f}(\phi) = 0
\end{equation}
for every $\phi \in \widehat{A} \setminus B$, by approximating $f$ by
elements of $\mathcal{E}_B(A)$.  Conversely, if $f \in C(A)$ satisfies
(\ref{widehat{f}(phi) = 0}) for every $\phi \in \widehat{A} \setminus
B$, then $f \in C_B(A)$, by the discussion in Section \ref{compact
  groups, continued}.  Thus
\begin{equation}
\label{C_B(A) = ...}
        C_B(A) = \{f \in C(A) : \widehat{f}(\phi) = 0 \hbox{ for every }
                                           \phi \in \widehat{A} \setminus B\},
\end{equation}
which is the same as (\ref{mathcal{L}_E cap C(A)}), with $E =
\widehat{A} \setminus B$.  This is a closed linear subspace of $C(A)$
that is invariant under translations, and every translation-invariant
closed linear subspace of $C(A)$ is of this form, as in the previous
section.  Note that
\begin{equation}
\label{B = widehat{A} cap C_B(A)}
        B = \widehat{A} \cap C_B(A).
\end{equation}

        Suppose that $B$ is a sub-semigroup of $\widehat{A}$, which
is to say that the product of any two elements of $B$ is also
contained in $B$.  This implies that $\mathcal{E}_B(A)$ is a
subalgebra of $\mathcal{E}(A)$, and hence that $C_B(A)$ is a
sub-algebra of $C(A)$.  One can also look at this in terms of
(\ref{widehat{(f g)}(phi) = ...})  and (\ref{C_B(A) = ...}).
Conversely, if $C_B(A)$ is a subalgebra of $C(A)$, then $B$ is a
sub-semigroup of $\widehat{A}$, by (\ref{B = widehat{A} cap C_B(A)}).

        Remember that the identity element of $\widehat{A}$ is the
constant function equal to $1$ at every point in $A$.  Using
(\ref{B = widehat{A} cap C_B(A)}) again, we get that $B$ contains
the identity element of $\widehat{A}$ if and only if $C_B(A)$
contains the constant functions on $A$.

        If $\phi \in \widehat{A}$, then $1/\phi$ is equal to the complex
conjugate $\overline{\phi}$ of $\phi$, because $\phi$ takes values in
${\bf T}$.  Put
\begin{equation}
\label{B^{-1} = {1/phi : phi in B} = {overline{phi} : phi in B}}
        B^{-1} = \{1/\phi : \phi \in B\} = \{\overline{\phi} : \phi \in B\},
\end{equation}
so that
\begin{equation}
\label{mathcal{E}_{B^{-1}}(A) = {overline{f} : f in mathcal{E}_B(A)}}
        \mathcal{E}_{B^{-1}}(A) = \{\overline{f} : f \in \mathcal{E}_B(A)\},
\end{equation}
and hence
\begin{equation}
\label{C_{B^{-1}}(A) = {overline{f} : f in C_B(A)}}
        C_{B^{-1}}(A) = \{\overline{f} : f \in C_B(A)\},
\end{equation}
by the definitions of $\mathcal{E}_{B^{-1}}(A)$ and $C_{B^{-1}}(A)$.
Similarly, one can use (\ref{C_B(A) = ...}) to get that
\begin{equation}
\label{C_B(A) cap C_{B^{-1}}(A) = C_{B cap B^{-1}}(A)}
        C_B(A) \cap C_{B^{-1}}(A) = C_{B \cap B^{-1}}(A).
\end{equation}
In particular, if $B^{-1} = B$, then it follows that $C_B(A)$ is
invariant under complex conjugation.  Conversely, if $C_B(A)$ is
invariant under complex conjugation, then $B = B^{-1}$, by (\ref{B =
  widehat{A} cap C_B(A)}).

        If $B$ is a subgroup of $\widehat{A}$, then it follows that
$C_B(A)$ is a subalgebra of $C(A)$ that contains the constant functions
and is invariant under complex conjugation.  Conversely, if $C_B(A)$
is a subalgebra of $C(A)$ that is invariant under complex conjugation,
then $B$ is a sub-semigroup of $\widehat{A}$ that satisfies $B^{-1} =
B$.  This implies that $B$ is a subgroup of $\widehat{A}$ when $B \ne
\emptyset$, so that $C_B(A) \ne \{0\}$.

        Let us now consider a basic example where $A = {\bf T}$.
Remember that
\begin{equation}
\label{z mapsto z^j, 2}
        z \mapsto z^j
\end{equation}
is a continuous homomorphism from ${\bf T}$ into itself for every $j
\in {\bf Z}$, and that every continuous homomorphism from ${\bf T}$
into itself is of this form.  Let $B$ be the collection of these
homomorphisms (\ref{z mapsto z^j, 2}) with $j \ge 0$, so that $B$ is a
sub-semigroup of the dual of ${\bf T}$ that contains the identity
element.  In this case, $C_B({\bf T})$ consists of the $f \in C({\bf
  T})$ that can be extended to a continuous function on the closed
unit disk $\overline{D}$ which is holomorphic on the open unit disk
$D$, as in Section \ref{holomorphic functions}.

        The case of subalgebras of $C(A)$ associated to subgroups of
$\widehat{A}$ is discussed further in the next section, and an extension
of the example in the preceding paragraph will be considered in
Section \ref{another class of subalgebras}.

\section{Subgroups and subalgebras}
\label{subgroups, subalgebras}

        Let $A$ be a compact commutative topological group, and let
$A_1$ be a closed subgroup of $A$.  Thus the quotient group $A / A_1$
is also a compact commutative topological group with respect to the
quotient topology, as in Section \ref{quotient groups}.  Every
continuous function on $A / A_1$ leads to a continuous function on $A$
that is constant on the translates of $A_1$ in $A$, by composition
with the quotient mapping $q_1$ from $A$ onto $A / A_1$.  Conversely,
if $f$ is any function on $A$ that is constant on the translates of
$A_1$ in $A$, then $f$ can be expressed as the composition of $q_1$
with a function $f_1$ on $A / A_1$.  If $f$ is also continuous on $A$,
then it is easy to see that $f_1$ is continuous on $A / A_1$, by the
definition of the quotient topology.

        Note that
\begin{equation}
\label{f_1 mapsto f_1 circ q_1}
        f_1 \mapsto f_1 \circ q_1
\end{equation}
defines an algebra homomorphism from $C(A / A_1)$ into $C(A)$, which
is also an isometric embedding with respect to the supremum norms on
$C(A / A_1)$ and $C(A)$.  As in the previous paragraph, the image of
$C(A / A_1)$ in $C(A)$ under (\ref{f_1 mapsto f_1 circ q_1}) is equal to
\begin{equation}
\label{{f in C(A) : f(a) = f(b) for every a, b in A such that a - b in A_1}}
        \quad \{f \in C(A) : f(a) = f(b) \hbox{ for every } a, b \in A
                                    \hbox{ such that } a - b \in A_1\}.
\end{equation}
It is easy to see that (\ref{{f in C(A) : f(a) = f(b) for every a, b
    in A such that a - b in A_1}}) is a closed subalgebra of $C(A)$
that contains the constant functions, is invariant under translations
on $A$, and is invariant under complex conjugation as well.

        Let $B_1$ be the set of $\phi \in \widehat{A}$ such that $A_1$
is contained in the kernel of $\phi$, which is the same as saying that
$\phi$ is constant on the translates of $A_1$ in $A$, because $\phi$
is a homomorphism.  If $\phi_1$ is a continuous homomorphism from $A / A_1$
into ${\bf T}$, then
\begin{equation}
\label{phi_1 circ q_1 in B_1}
        \phi_1 \circ q_1 \in B_1,
\end{equation}
and every element of $B_1$ is of this form.  This is analogous to the
discussion in the preceding paragraphs, and was also mentioned in
Section \ref{quotient groups}.  This defines a natural group
isomorphism from the dual $\widehat{(A / A_1)}$ of $A / A_1$ onto
$B_1$.

        Under these conditions, $C_{B_1}(A)$ is the same as
(\ref{{f in C(A) : f(a) = f(b) for every a, b in A such that a - b in A_1}}).
More precisely, $B_1$ is contained in (\ref{{f in C(A) : f(a) = f(b) for every 
a, b in A such that a - b in A_1}}) by construction, which implies that the
linear span $\mathcal{E}_{B_1}(A)$ of $B_1$ in $C(A)$ is contained in
(\ref{{f in C(A) : f(a) = f(b) for every a, b in A such that a - b in A_1}}),
and hence that the closure $C_{B_1}(A)$ of $\mathcal{E}_{B_1}(A)$ in $C(A)$
is contained in (\ref{{f in C(A) : f(a) = f(b) for every a, b in A such that 
a - b in A_1}}).  In the other direction, if $f_1 \in C(A / A_1)$, then
$f_1$ can be approximated by linear combinations of elements of
$\widehat{(A / A_1)}$ uniformly on $A / A_1$, because $A / A_1$ is
compact.  This implies that (\ref{{f in C(A) : f(a) = f(b) for every a, b in A 
such that a - b in A_1}}) is contained in $C_{B_1}(A)$, since every element of
(\ref{{f in C(A) : f(a) = f(b) for every a, b in A such that a - b in A_1}})
is of the form $f_1 \circ q_1$ for some $f_1 \in C(A / A_1)$.

        Now let $B$ be any subgroup of $\widehat{A}$, and consider
\begin{equation}
\label{A_1 = {a in A : phi(a) = 1 for every phi in B}}
        A_1 = \{a \in A : \phi(a) = 1 \hbox{ for every } \phi \in B\}.
\end{equation}
This is a closed subgroup of $A$, so that the quotient group $A / A_1$
is a compact commutative topological group with respect to the
quotient topology, as before.  If $q_1$ is the usual quotient mapping
from $A$ onto $A / A_1$, then every $\phi \in B$ can be expressed as
$\phi_1 \circ q_1$ for some $\phi_1 \in \widehat{(A / A_1)}$.  Put
\begin{equation}
\label{widetilde{B} = {phi_1 in widehat{(A / A_1)} : phi_1 circ q_1 in B}}
        \widetilde{B} = \{\phi_1 \in \widehat{(A / A_1)} :
                                           \phi_1 \circ q_1 \in B\},
\end{equation}
which is a subgroup of $\widehat{(A / A_1)}$ isomorphic to $B$.  Of
course, if $a \in A \setminus A_1$, then there is a $\phi \in B$ such
that $\phi(a) \ne 1$, by definition of $A_1$.  This implies that
$\widetilde{B}$ separates points in $A / A_1$ in this situation,
because of the way that $A_1$ is defined.  It follows that
\begin{equation}
\label{widetilde{B} = widehat{(A / A_1)}}
        \widetilde{B} = \widehat{(A / A_1)},
\end{equation}
because $A / A_1$ is compact, as in Section \ref{compact, discrete
  groups}.

        Let $B_1$ be the subgroup of $\widehat{A}$ associated to $A_1$
as before, consisting of the $\phi \in \widehat{A}$ such that the
kernel of $\phi$ contains $A_1$.  Thus $B \subseteq B_1$ in this case,
by the definition (\ref{A_1 = {a in A : phi(a) = 1 for every phi in
    B}}) of $A_1$.  If $\phi \in B_1$, then we have seen that $\phi =
\phi_1 \circ q_1$ for some $\phi_1 \in \widehat{(A / A_1)}$, and now
we have that $\phi_1 \in \widetilde{B}$, by (\ref{widetilde{B} =
  widehat{(A / A_1)}}).  This implies that $\phi \in B$, so that
\begin{equation}
\label{B = B_1}
        B = B_1.
\end{equation}
It follows that
\begin{equation}
\label{C_B(A) = C_{B_1}(A)}
        C_B(A) = C_{B_1}(A),
\end{equation}
and hence that $C_B(A)$ is equal to (\ref{{f in C(A) : f(a) = f(b) for
    every a, b in A such that a - b in A_1}}).

        Let $B_0$ be a nonempty subset of $\widehat{A}$, and let $B$
be the subgroup of $\widehat{A}$ generated by $B_0$.  Note that $B$
separates points in $A$ if and only if $B_0$ separates points in $A$.
In this case, $B = \widehat{A}$, because $A$ is compact, as in Section
\ref{compact, discrete groups}.  Otherwise, one can apply the previous
discussion to $B$, to get that $C_B(A)$ is equal to (\ref{{f in C(A) :
    f(a) = f(b) for every a, b in A such that a - b in A_1}}), with
$A_1$ as in (\ref{A_1 = {a in A : phi(a) = 1 for every phi in B}}).

\chapter{Additional topics}
\label{additional topics}

\section{Almost periodic functions}
\label{almost periodic functions}

        Let $A$ be a commutative topological group, and let $C_b(A)$
be the space of complex-valued functions $f$ on $A$ that are bounded
and continuous, equipped with the supremum norm $\|f\|_{sup}$.  Let $f
\in C_b(A)$ be given, and consider the set
\begin{equation}
\label{{f_a : a in A}}
        \{f_a : a \in A\}
\end{equation}
of all translates of $f$, where $f_a(x) = f(x - a)$, as in
(\ref{f_a(x) = f(x - a), 3}).  If (\ref{{f_a : a in A}}) is totally
bounded as a subset of $C_b(A)$, with respect to the supremum metric,
then $f$ is said to be \emph{almost periodic}\index{almost periodic
  functions} on $A$.  Equivalently, this means that the closure of
(\ref{{f_a : a in A}}) in $C_b(A)$ is compact, since $C_b(A)$ is
complete.  Let $\mathcal{AP}(A)$ denote the collection of almost
periodic functions on $A$.

        Remember that a subset $E$ of a metric space is totally bounded
if for each $r > 0$, $E$ can be covered by finitely many balls of
radius $r > 0$.  One can get an equivalent condition by requiring that
these balls be centered at elements of $E$.  More precisely, if $E$ is
convered by finitely many balls of radius $r/2$ that are not necessarily
centered at elements of $E$, then one can first throw away the balls
that do not intersect $E$.  The remaining balls are contained in balls
of radius $r$ centered at elements of $E$, by the triangle inequality.

        Constant functions on $A$ are obviously almost periodic, and it
is easy to see that continuous homomorphisms from $A$ into ${\bf T}$
are almost periodic as well.  If $f, g \in \mathcal{AP}(A)$, then one
can check their sum $f + g$ is almost periodic, because the sum of two
totally bounded subsets of $C_b(A)$ is totally bounded.  Similarly,
the product $f \, g$ is almost periodic on $A$, so that
$\mathcal{AP}(A)$ is a subalgebra of $C_b(A)$.  The complex-conjugate
of an almost periodic function on $A$ is almost periodic on $A$ too,
and in fact $f \in C_b(A)$ is almost periodic on $A$ if and only if
the real and imaginary parts of $f$ are almost periodic on $A$.
Note that $\mathcal{AP}(A)$ is a closed subset of $C_b(A)$ with
respect to the supremum metric.

        If $f \in \mathcal{AP}(A)$ and $\phi$ is a continuous
complex-valued function on ${\bf C}$, then one can check that their
composition $\phi \circ f$ is almost periodic on $A$ too.  This uses
the fact that $\phi$ is uniformly continuous on compact subsets of
${\bf C}$.  In particular, if $f \in \mathcal{AP}(A)$ and $|f(x)|$ has
a positive lower bound on $A$, then one can choose $\phi$ so that
\begin{equation}
\label{phi(f(x)) = 1/f(x)}
        \phi(f(x)) = 1/f(x)
\end{equation}
for every $x \in A$.  Thus $1/f \in \mathcal{AP}(A)$ when $f \in
\mathcal{AP}(A)$ and $|f(x)|$ has a positive lower bound on $A$.
Similarly, $|f| \in \mathcal{AP}({\bf R})$ for every $f \in
\mathcal{AP}({\bf R})$.

        If $f$ is any bounded continuous function on $A$ and $b \in A$,
then the translate $f_b(x) = f(x - b)$ of $f$ by $b$ is also bounded and
continuous on $A$.  Of course, $f_b$ generates the same set of translates
in $C_b(A)$ as $f$ does.  Thus $f_b \in \mathcal{AP}(A)$ for every
$b \in A$ when $f \in \mathcal{AP}(A)$.

        Uniform continuity of complex-valued functions on $A$
was defined in Section \ref{uniform continuity, equicontinuity},
and uniform continuity of functions on $A$ with values in a
metric space or a commutative topological group can be defined
similarly.  If $f \in C_b(A)$ is uniformly continuous on $A$, then
\begin{equation}
\label{a mapsto f_a}
        a \mapsto f_a
\end{equation}
is uniformly continuous as a mapping from $A$ into $C_b(A)$, with
respect to the supremum metric on $C_b(A)$.  In this case, if $A$ is
also totally bounded as a commutative topological group, then it
follows that (\ref{{f_a : a in A}}) is totally bounded in $C_b(A)$, so
that $f$ is almost periodic on $A$.  If $A$ is compact, then we have
seen that every continuous complex-valued function on $A$ is uniformly
continuous, and $A$ is automatically totally bounded as well.  It
follows that every continuous complex-valued function on $A$ is almost
periodic on $A$ when $A$ is compact.

        If $\mathcal{E} \subseteq C_b(A)$ is totally bounded with
respect to the supremum metric, then it is easy to see that
$\mathcal{E}$ is equicontinuous at every point in $A$.  If $f \in
\mathcal{AP}(A)$, then it follows that (\ref{{f_a : a in A}}) is
equicontinuous at every point in $A$.  The equicontinuity of
(\ref{{f_a : a in A}}) at any point in $A$ is equivalent to the
uniform continuity of $f$ on $A$.  In particular, every $f \in
\mathcal{AP}(A)$ is uniformly continuous on $A$.

        Let $B$ be another commutative topological group, and suppose
that $h$ is a continuous homomorphism from $A$ into $B$.  If $g \in
\mathcal{AP}(B)$, then one can check that $g \circ h \in \mathcal{AP}(A)$.
If $B$ is compact, then $g \circ h \in \mathcal{AP}(A)$ for every $g
\in C(B)$, because every continuous function on $B$ is almost periodic. 

        If $B$ is a compact commutative topological group, then it is
well known that the continuous homomorphisms from $B$ into ${\bf T}$
separate points on $B$.  If $g \in C(B)$, then it follows that $g$ can
be approximated uniformly on $B$ by finite linear combinations of
continuous homomorphisms from $B$ into ${\bf T}$, because of the
Stone--Weierstrass theorem, as mentioned in Section \ref{compact,
  discrete groups}.  If $h$ is a continuous homomorphism from $A$ into
$B$ and $f = g \circ h$, then it is easy to see that $f$ can be
approximated uniformly on $A$ by finite linear combinations of
continuous homomorphisms from $A$ into ${\bf T}$.  It will be shown in
the next section that every $f \in \mathcal{AP}(A)$ can be represented
in this way, which implies that every $f \in \mathcal{AP}(A)$ can be
approximated uniformly on $A$ by finite linear combinations of
continuous homomorphisms from $A$ into ${\bf T}$.

\section{Groups of isometries}
\label{groups of isometries}

        Let $(M, d(x, y))$ be a (nonempty) metric space, and let
$C(M, M)$ be the space of continuous mappings from $M$ into itself.
This is a semigroup with respect to composition of mappings, and with
the identity mapping on $M$ as the identity element.  Similarly, the
collection of homeomorphisms from $M$ onto itself is a group with
respect to composition.  Let $\mathcal{I}(M)$ be the collection of
mappings $\phi$ from $M$ onto itself that are
isometries,\index{isometric mappings} in the sense that
\begin{equation}
\label{d(phi(x), phi(y)) = d(x, y)}
        d(\phi(x), \phi(y)) = d(x, y)
\end{equation}
for every $x, y \in M$.  This is a subgroup of the group of
homeomorphisms from $M$ onto itself.

        Suppose now that $M$ is bounded, and put
\begin{equation}
\label{rho(phi, psi) = sup_{x in M} d(phi(x), psi(x))}
        \rho(\phi, \psi) = \sup_{x \in M} d(\phi(x), \psi(x))
\end{equation}
for every $\phi, \psi \in C(M, M)$, which is the supremum metric on
$C(M, M)$.  Note that (\ref{rho(phi, psi) = sup_{x in M} d(phi(x),
  psi(x))}) is invariant under composition of $\phi$ and $\psi$ on the
right by continuous mappings from $M$ onto itself, and under
composition of $\phi$ and $\psi$ on the left by isometries from $M$
into itself.  In particular, the restriction of the supremum metric to
$\mathcal{I}(M)$ is invariant under left and right translations on
$\mathcal{I}(M)$ as a group with respect to composition of mappings.
Using this, one can check that composition of mappings defines a
continuous mapping from $\mathcal{I}(M) \times \mathcal{I}(M)$ into
$\mathcal{I}(M)$, with respect to the topology on $\mathcal{I}(M)$
determined by the restriction of the supremum metric to
$\mathcal{I}(M)$.

        If $\phi$ is any homeomorphism from $M$ onto itself, then the
distance from $\phi^{-1}$ to the identity mapping on $M$ with respect to
the supremum metric is equal to the distance from $\phi$ to the identity
mapping.  If $\phi, \psi \in \mathcal{I}(M)$, then the distance from
$\phi^{-1}$ to $\psi^{-1}$ with respect to the supremum metric is equal
to the distance from $\phi$ to $\psi$.  This implies that $\phi \mapsto
\phi^{-1}$ is a continuous mapping from $\mathcal{I}(M)$ into itself
with respect to the restriction of the supremum metric to $\mathcal{I}(M)$.
It follows that $\mathcal{I}(M)$ is a topological group with respect
to composition of mappings and the topology determined by the restriction
of the supremum metric to $\mathcal{I}(M)$.

        It is easy to see that the collection of $f \in C(M, M)$
such that $f(M)$ is dense in $M$ is a closed subset of $C(M, M)$
with respect to the supremum metric.  If $M$ is compact, then $f(M)$
is a compact subset of $M$ for every $f \in C(M, M)$, which implies
that $f(M)$ is a closed set in $M$.  It follows that the collection
of $f \in C(M, M)$ such that $f(M) = M$ is a closed set in $C(M, M)$
with respect to the supremum metric when $M$ is compact.  Using
standard Arzela--Ascoli arguments, one can check that $\mathcal{I}(M)$
is compact with respect to the supremum metric when $M$ is compact.

        Now let $A$ be a commutative topological group, and let
$f \in \mathcal{AP}(A)$ be given.  Let $M$ be the closure of
(\ref{{f_a : a in A}}) in $C_b(A)$, which is a compact metric space
with respect to the restriction of the supremum metric on $C_b(A)$
to $M$.  As before, the group $\mathcal{I}(M)$ of isometries of $M$ onto
itself is a compact topological group with respect to the restriction of
the supremum metric on $C(M, M)$ to $\mathcal{I}(M)$.

        Put
\begin{equation}
\label{T_y(g)(x) = g(x + y)}
        T_y(g)(x) = g(x + y)
\end{equation}
for every $g \in C_b(A)$ and $x, y \in A$, so that $T_y$ defines a
linear mapping from $C_b(A)$ onto itself that preserves the supremum
norm.  Of course,
\begin{equation}
\label{T_y(f_a) = f_{a - y}}
        T_y(f_a) = f_{a - y}
\end{equation}
for every $a, y \in A$, which implies that $T_y$ maps (\ref{{f_a : a
    in A}}) onto itself for every $y \in A$, and hence that $T_y$ maps
$M$ onto itself.  We also have that
\begin{equation}
\label{T_y circ T_z = T_{y + z}}
        T_y \circ T_z = T_{y + z}
\end{equation}
for every $y, z \in A$, so that $y \mapsto T_y$ defines a homomorphism
from $A$ into the group of isometric linear mappings from $C_b(A)$
onto itself.  This leads to a homomorphism from $A$ into
$\mathcal{I}(M)$, which sends $y \in A$ to the restriction of $T_y$ to
$M$.

        Remember that $f$ is uniformly continuous on $A$, which implies
that (\ref{{f_a : a in A}}) is uniformly equicontinuous on $A$, and
hence that the elements of $M$ are uniformly equicontinuous on $A$.
Using this, one can check that the mapping from $y \in A$ to the
restriction of $T_y$ to $M$ is continuous with respect to the supremum
metric on $C(M, M)$.  Let $B_0$ be the subgroup of $\mathcal{I}(M)$
consisting of the restriction of $T_y$ to $M$ for each $y \in A$, and
let $B_1$ be the closure of $B_0$ in $\mathcal{I}(M)$ with respect to
the supremum metric on $C(M, M)$.  Thus $B_0$ is a commutative
subgroup of $\mathcal{I}(M)$, because $A$ is commutative, which
implies that $B_1$ is commutative as well.  Note that $B_1$ is
compact, because $\mathcal{I}(M)$ is compact.

        Put
\begin{equation}
\label{lambda_0(g) = g(0)}
        \lambda_0(g) = g(0)
\end{equation}
for each $g \in C_b(A)$, which defines a bounded linear functional on
$C_b(A)$ with respect to the supremum norm.  Also put
\begin{equation}
\label{F(T) = lambda_0(T(f)) = (T(f))(0)}
        F(T) = \lambda_0(T(f)) = (T(f))(0)
\end{equation}
for each $T \in C(M, M)$, where $T(f) \in M$ because $f \in M$, by
construction.  It is easy to see that (\ref{F(T) = lambda_0(T(f)) =
  (T(f))(0)}) is continuous as a function of $T$ with respect to the
supremum metric on $C(M, M)$.  In particular, the restriction of $F$
to $B_1$ defines a continuous function on $B_1$ such that
\begin{equation}
\label{F(T_y) = (T_y(f))(0) = f(y)}
        F(T_y) = (T_y(f))(0) = f(y)
\end{equation}
for every $y \in A$.  This shows that $f$ can be expressed as the
composition of a continuous homomorphism from $A$ into $B_1$ with a
continuous function on $B_1$, as mentioned in the previous section.

\section{Compactifications}
\label{compactifications}

        Let $A$, $B$ be commutative topological groups, and let
$\widehat{A} = \hom(A, {\bf T})$, $\widehat{B} = \hom(B, {\bf T})$
be the corresponding dual groups of continuous homomorphisms into
${\bf T}$.  If $h$ is a continuous homomorphism from $A$ into $B$,
then
\begin{equation}
\label{widehat{h}(phi) = phi circ h}
        \widehat{h}(\phi) = \phi \circ h
\end{equation}
defines a homomorphism from $\widehat{B}$ into $\widehat{A}$, as in
Section \ref{quotient groups}.  This is known as the dual homomorphism
associated to $h$, and it is also continuous with respect to the usual
topologies on $\widehat{A}$, $\widehat{B}$.  In this section, we shall
be especially interested in the case where $B$ is compact, so that the
usual topology on $\widehat{B}$ is the same as the discrete topology,
and the continuity of $\widehat{h}$ is trivial.

        We shall also be especially interested in the case where
$h(A)$ is dense in $B$.  Of course, one can always arrange for this
to be the case, by replacing $B$ with the closure of $h(A)$ in $B$.
If $h(A)$ is dense in $B$, then it is easy to see that $\widehat{h}$
is injective, as in Section \ref{quotient groups}.  Otherwise, the
kernel of $\widehat{h}$ is isomorphic as a group to the dual of the
quotient $B / \overline{h(A)}$, where $\overline{h(A)}$ is the closure
of $h(A)$ in $B$, and the quotient group $B / \overline{h(A)}$ is
equipped with the corresponding quotient topology.  If $B$ is locally
compact, then we have seen that the quotient of $B$ by any closed
subgroup is locally compact as well.  It is well known that the dual
of a locally compact commutative topological group separates points on
that group, and in particular that the dual of a nontrivial locally
compact commutative topological group is nontrivial.  Thus the dual of
$B / \overline{h(A)}$ is nontrivial when $B$ is locally compact and
$\overline{h(A)} \ne B$.  If $\widehat{h}$ is injective, then the
dual of $B / \overline{h(A)}$ is trivial, which implies that
$\overline{h(A)} = B$ when $B$ is locally compact.

        Let $\widehat{\widehat{A}}$, $\widehat{\widehat{B}}$ be the
duals of $\widehat{A}$, $\widehat{B}$, respectively, which are the
second duals of $A$, $B$.  Remember that there are natural
homomorphisms from $A$, $B$ into their second duals, as in Section
\ref{second dual}.  We have seen that these mappings are continuous on
compact subsets of their domains, and hence that they are continuous
under some additional conditions.  Let $\widehat{\widehat{h}}$ be the
dual homomorphism associated to $\widehat{h}$, which is a continuous
homomorphism from $\widehat{\widehat{A}}$ into $\widehat{\widehat{B}}$.
It is easy to see that the composition of the natural mapping from
$A$ into $\widehat{\widehat{A}}$ with $\widehat{\widehat{h}}$ is equal to
the composition of $h$ with the natural mapping from $B$ into
$\widehat{\widehat{B}}$.

        If $B$ is locally compact, then it is well known that the
natural mapping from $B$ into $\widehat{\widehat{B}}$ is an
isomorphism from $B$ onto $\widehat{\widehat{B}}$ as topological
groups.  Let $k$ be a continuous homomorphism from $\widehat{B}$ into
$\widehat{A}$, so that the dual $\widehat{k}$ of $k$ is a continuous
homomorphism from $\widehat{\widehat{A}}$ into
$\widehat{\widehat{B}}$.  If $B$ is locally compact, then
$\widehat{\widehat{B}}$ can be identified with $B$, and hence
$\widehat{k}$ can be identified with a continuous homomorphism from
$\widehat{\widehat{A}}$ into $B$.  In this case, the composition of
the natural mapping from $A$ into $\widehat{\widehat{A}}$ with
$\widehat{k}$ can be identified with a homomorphism from $A$ into $B$,
which is continuous when the natural mapping from $A$ into
$\widehat{\widehat{A}}$ is continuous.  Under these conditions, the
dual of this homomorphism from $A$ into $B$ is equal to $k$.

        Suppose now that $B$ is compact, so that the usual topology
on $\widehat{B}$ is the same as the discrete topology.  Let $k$ be a
homomorphism from $\widehat{B}$ into $\widehat{A}$, which is
automatically continuous, and let $g$ be the homomorphism from $A$
into $\widehat{\widehat{B}}$ which is the composition of the natural
mapping from $A$ into $\widehat{\widehat{A}}$ with the dual
$\widehat{k}$ of $k$.  Thus for each $\phi \in \widehat{B}$ and $a \in
A$, $k(\phi) \in \widehat{A}$, $g(a) \in \widehat{\widehat{B}}$, and
\begin{equation}
\label{g(a)(phi) = k(phi)(a)}
        g(a)(\phi) = k(\phi)(a),
\end{equation}
which is continuous as a function of $a \in A$ for each $\phi \in
\widehat{B}$.  In this situation, one can check that $g$ is continuous
as a mapping from $A$ into $\widehat{\widehat{B}}$, without additional
conditions on $A$.  This uses the fact that the compact subsets of
$\widehat{B}$ are finite, because $\widehat{B}$ is equipped with the
discrete topology.

        As before, the natural mapping from $B$ into $\widehat{\widehat{B}}$
is an isomorphism from $B$ onto $\widehat{\widehat{B}}$ as topological
groups when $B$ is compact, so that $g$ can be identified with a 
continuous homomorphism from $A$ into $B$, whose dual is equal to $k$.
If $k$ is injective, then $g(A)$ is dense in $B$, by the earlier arguments.

        Alternatively, let $C$ be a commutative group equipped with the
discrete topology, and let $k$ be a homomorphism from $C$ into $\widehat{A}$,
which is automatically continuous.  Thus the dual of $k$ is a continuous
homomorphism from $\widehat{\widehat{A}}$ into $\widehat{C}$, and $\widehat{C}$
is compact.  If $g$ is the composition of the natural mapping from $A$
into $\widehat{\widehat{A}}$ with $\widehat{k}$, then
\begin{equation}
\label{g(a)(c) = k(c)(a)}
        g(a)(c) = k(c)(a)
\end{equation}
for every $a \in A$ and $c \in C$.  One can check that $g$ is
continuous as a mapping from $A$ into $\widehat{C}$ under these
conditions, because the compact subsets of $C$ are finite, and
(\ref{g(a)(c) = k(c)(a)}) is continuous as a function of $a \in A$ for
each $c \in C$.  The dual of $g$ is a homomorphism from
$\widehat{\widehat{C}}$ into $\widehat{A}$, which is automatically
continuous, because the usual topology on $\widehat{\widehat{C}}$ is
the discrete topology, since $\widehat{C}$ is compact.  The
composition of the natural mapping from $C$ into
$\widehat{\widehat{C}}$ and the dual of $g$ is equal to $k$, and we
have seen that the natural mapping from $C$ into
$\widehat{\widehat{C}}$ is an isomorphism when $C$ is discrete.  Of
course, this means that $C$ is isomorphic to the dual of a compact
commutative topological group.

         In particular, one can apply this with $C$ equal to $\widehat{A}$
as a commutative group, equipped with the discrete topology.  The
identity mapping on $\widehat{A}$ can then be interpreted as a
homomorphism $k$ from $C$ into $\widehat{A}$.  This leads to a
continuous homomorphism $g$ from $A$ into $\widehat{C}$ as in the
previous paragraph, which is injective exactly when $\widehat{A}$
separates points in $A$.  Note that $g(A)$ is dense in $\widehat{C}$,
because $k$ is injective by construction, and $\widehat{C}$ is
compact.  This construction is known as the \emph{Bohr
  compactification}\index{Bohr compactification} of $A$.

\section{Averages on ${\bf R}$}
\label{averages on R}

        Let $f$ be a bounded continuous complex-valued function
on the real line that is almost periodic, and let $\epsilon > 0$
be given.  Thus there are finitely many elements $a_1, \ldots, a_n$
of ${\bf R}$ such that for each $a \in {\bf R}$,
\begin{equation}
\label{||f_a - f_{a_j}||_{sup} < epsilon}
        \|f_a - f_{a_j}\|_{sup} < \epsilon
\end{equation}
for some $j \in \{1, \ldots, n\}$, where $\|\cdot\|_{sup}$ is the
supremum norm on $C_b({\bf R})$.

        Let $I$ be a closed interval in ${\bf R}$, with length
$|I| > 0$.  Because $f$ is bounded on ${\bf R}$, there is an
$L(\epsilon) > 0$ such that
\begin{equation}
\label{difference of averages, 1}
 \biggl|\frac{1}{|I|} \, \int_I f(x - a_j) \, dx
                - \frac{1}{|I|} \, \int_I f(x) \, dx\biggr| < \epsilon
\end{equation}
for every $j = 1, \ldots, n$ when $|I| \ge L(\epsilon)$.  This implies
that
\begin{equation}
\label{difference of averages, 2}
 \biggl|\frac{1}{|I|} \, \int_I f(x - a) \, dx 
                - \frac{1}{|I|} \, \int_I f(x) \, dx\biggr| < 2 \, \epsilon
\end{equation}
for every $a \in {\bf R}$ when $|I| \ge L(\epsilon)$, by (\ref{||f_a -
  f_{a_j}||_{sup} < epsilon}).  Equivalently,
\begin{equation}
\label{difference of averages, 3}
 \biggl|\frac{1}{|I'|} \, \int_{I'} f(x) \, dx
                - \frac{1}{|I|} \, \int_I f(x) \, dx\biggr| < 2 \, \epsilon
\end{equation}
for every closed interval $I'$ in ${\bf R}$ with $|I'| = |I| \ge
L(\epsilon)$.

        Suppose now that $I'$ is a closed interval in ${\bf R}$ whose
length is equal to a positive integer $l$ times $|I|$.  This implies that
$I'$ can be expressed as
\begin{equation}
\label{I' = bigcup_{k = 1}^l I'_k}
        I' = \bigcup_{k = 1}^l I'_k,
\end{equation}
where $I'_1, \ldots, I'_l$ are closed intervals in ${\bf R}$ with
disjoint interiors, and $|I'_k| = |I|$ for $k = 1, \ldots, l$.  Thus
\begin{equation}
\label{frac{1}{|I'|} int_{I'} f(x) dx}
        \frac{1}{|I'|} \, \int_{I'} f(x) \, dx
 = \frac{1}{l} \, \sum_{k = 1}^l \frac{1}{|I'_k|} \, \int_{I'_k} f(x) \, dx,
\end{equation}
and
\begin{equation}
\label{difference of averages, 4}
 \biggl|\frac{1}{|I'_k|} \, \int_{I'_k} f(x) \, dx
                - \frac{1}{|I|} \, \int_I f(x) \, dx\biggr| < 2 \, \epsilon
\end{equation}
for each $k = 1, \ldots, l$ when $|I| \ge L(\epsilon)$, as before.  It
follows that (\ref{difference of averages, 3}) still holds in this
case when $|I| \ge L$.

        Suppose instead that $|I'|$ is a positive rational number times $|I|$.
This implies that there are positive integer $l$, $l'$ such that
\begin{equation}
\label{l' |I'| = l |I|}
        l' \, |I'| = l \, |I|,
\end{equation}
and we can take $I''$ to be a closed interval in ${\bf R}$ whose
length is equal to (\ref{l' |I'| = l |I|}).  If $|I|, |I'| \ge
L(\epsilon)$, then we get that
\begin{equation}
\label{difference of averages, 5}
 \biggl|\frac{1}{|I''|} \, \int_{I''} f(x) \, dx
                - \frac{1}{|I|} \, \int_I f(x) \, dx\biggr| < 2 \, \epsilon
\end{equation}
and
\begin{equation}
\label{difference of averages, 6}
 \biggl|\frac{1}{|I''|} \, \int_{I''} f(x) \, dx
                - \frac{1}{|I'|} \, \int_{I'} f(x) \, dx\biggr| < 2 \, \epsilon
\end{equation}
as in the preceding paragraph.  Hence
\begin{equation}
\label{difference of averages, 7}
 \biggl|\frac{1}{|I'|} \, \int_{I'} f(x) \, dx
                - \frac{1}{|I|} \, \int_I f(x) \, dx\biggr| < 4 \, \epsilon
\end{equation}
under these conditions.

        In particular, (\ref{difference of averages, 7}) holds when $I$, $I'$
are closed intervals in ${\bf R}$ whose lengths are positive rational
numbers greater than or equal to $L(\epsilon)$.  It follows that
\begin{equation}
\label{difference of averages, 8}
 \biggl|\frac{1}{|I'|} \, \int_{I'} f(x) \, dx
                - \frac{1}{|I|} \, \int_I f(x) \, dx\biggr| \le 4 \, \epsilon
\end{equation}
when $I$, $I'$ are any two closed intervals in ${\bf R}$ with $|I|,
|I'| \ge L(\epsilon)$, since we can approximate $I$ and $I'$ by
intervals with rational lengths.

        Let $I_1, I_2, I_3, \ldots$ be a sequence of closed intervals
in ${\bf R}$, where $|I_j| > 0$ for each $j$, and $|I_j| \to \infty$
as $j \to \infty$.  Using (\ref{difference of averages, 8}), it is
easy to see that the corresponding sequence of averages
\begin{equation}
\label{frac{1}{|I_j|} int_{I_j} f(x) dx}
        \frac{1}{|I_j|} \, \int_{I_j} f(x) \, dx
\end{equation}
is a Cauchy sequence in ${\bf C}$, which therefore converges in ${\bf
  C}$.  Put
\begin{equation}
\label{Lambda(f) = lim_{j to infty} frac{1}{|I_j|} int_{I_j} f(x) dx}
        \Lambda(f) = \lim_{j \to \infty} \frac{1}{|I_j|} \, \int_{I_j} f(x) \, dx,
\end{equation}
and observe that this does not depend on the particular choice of
sequence of $I_j$'s, because of (\ref{difference of averages, 8})
again.  If $I$ is any closed interval in ${\bf R}$ with $|I| \ge
L(\epsilon)$, then it is easy to see that
\begin{equation}
\label{|Lambda(f) - frac{1}{|I|} int_I f(x) dx| le 4 epsilon}
        \biggl|\Lambda(f) - \frac{1}{|I|} \, \int_I f(x) \, dx\biggr|
                                                \le 4 \, \epsilon,
\end{equation}
by (\ref{difference of averages, 8}).

        Of course, $\Lambda(f)$ is a linear functional on
$\mathcal{AP}({\bf R})$, and
\begin{equation}
\label{|Lambda(f)| le ||f||_{sup}}
        |\Lambda(f)| \le \|f\|_{sup}
\end{equation}
for every $f \in \mathcal{AP}({\bf R})$, because of the analogous
properties of the averages (\ref{frac{1}{|I_j|} int_{I_j} f(x) dx}).
Similarly, $\Lambda(f) \in {\bf R}$ when $f \in \mathcal{AP}({\bf R})$
is real-valued on ${\bf R}$, and
\begin{equation}
\label{Lambda(overline{f}) = overline{Lambda(f)}}
        \Lambda(\overline{f}) = \overline{\Lambda(f)}
\end{equation}
for every $f \in \mathcal{AP}({\bf R})$.  If $f \in \mathcal{AP}({\bf
  R})$ is real-valued and nonnegative on ${\bf R}$, then $\Lambda(f)
\ge 0$.  In this case, we shall show in a moment that $\Lambda(f) > 0$
when $f(x) > 0$ for some $x \in {\bf R}$.  If $f$ is any element of
$\mathcal{AP}({\bf R})$, then we have seen that $|f| \in
\mathcal{AP}({\bf R})$, and it is easy to check that
\begin{equation}
\label{|Lambda(f)| le Lambda(|f|)}
        |\Lambda(f)| \le \Lambda(|f|).
\end{equation}
We have also seen that $f_b(x) = f(x - b) \in \mathcal{AP}({\bf R})$
for every $b \in {\bf R}$ when $f \in \mathcal{AP}({\bf R})$, and one
can verify that
\begin{equation}
\label{Lambda(f_b) = Lambda(f)}
        \Lambda(f_b) = \Lambda(f),
\end{equation}
by construction.  Note that
\begin{equation}
\label{Lambda({bf 1}_{bf R}) = 1}
        \Lambda({\bf 1}_{\bf R}) = 1,
\end{equation}
where ${\bf 1}_{\bf R}(x)$ is the constant function on ${\bf R}$
equal to $1$ for every $x \in {\bf R}$.

        Suppose that $f \in \mathcal{AP}({\bf R})$ is real-valued and
nonnegative on ${\bf R}$, and that $f(x_0) > 0$ for some $x_0 \in {\bf R}$.
Put $\epsilon = f(x_0) / 3$, and let $\delta$ be a positive real
number such that
\begin{equation}
\label{f(x) ge 2 epsilon}
        f(x) \ge 2 \, \epsilon
\end{equation}
for every $x \in {\bf R}$ with $|x - x_0| \le \delta$, which exists because
$f$ is continuous at $x_0$.  Let $a_1, \ldots, a_n \in {\bf R}$ be as at
the beginning of the section, and observe that
\begin{equation}
\label{f_{a_j}(x) = f(x - a_j) ge 2 epsilon}
        f_{a_j}(x) = f(x - a_j) \ge 2 \, \epsilon
\end{equation}
when $|(x - a_j) - x_0| \le \delta$.  If $a \in {\bf R}$ and $1 \le j \le n$
satisfy (\ref{||f_a - f_{a_j}||_{sup} < epsilon}), then it follows that
\begin{equation}
\label{f_a(x) = f(x - a) ge epsilon}
        f_a(x) = f(x - a) \ge \epsilon
\end{equation}
when $|x - a_j - x_0| \le \delta$.

        Let $I_0$ be a closed interval in ${\bf R}$ such that
$a_j + x_0 \in I_0$ for each $j = 1, \ldots, n$ and $|I_0| \ge \delta$.
Put
\begin{equation}
\label{I_{0, j} = I_0 cap [a_j + x_0 - delta, a_j + x_0 + delta]}
        I_{0, j} = I_0 \cap [a_j + x_0 - \delta, a_j + x_0 + \delta]
\end{equation}
for each $j = 1, \ldots, n$, and observe that $|I_{0, j}| \ge \delta$
for every $j$.  If $a \in {\bf R}$ and $1 \le j \le n$ satisfy
(\ref{||f_a - f_{a_j}||_{sup} < epsilon}), then (\ref{f_a(x) = f(x -
  a) ge epsilon}) holds for each $x \in I_{0, j}$, and hence
\begin{equation}
\label{int_{I_0} f_a(x) dx ge int_{I_{0, j}} f_a(x) dx ge epsilon delta}
 \int_{I_0} f_a(x) \, dx \ge \int_{I_{0, j}} f_a(x) \, dx \ge \epsilon \, \delta.
\end{equation}
This implies that
\begin{equation}
\label{frac{1}{|I_0|} int_{I_0} f_a(x) dx ge epsilon delta |I_0|^{-1}}
 \frac{1}{|I_0|} \, \int_{I_0} f_a(x) \, dx \ge \epsilon \, \delta \, |I_0|^{-1}
\end{equation}
for every $a \in {\bf R}$, since for each $a \in {\bf R}$ there is a
$j$ such that (\ref{||f_a - f_{a_j}||_{sup} < epsilon}) holds.  This
is the same as saying that
\begin{equation}
\label{frac{1}{|I|} int_I f(x) dx ge epsilon delta |I_0|^{-1}}
        \frac{1}{|I|} \, \int_I f(x) \, dx \ge \epsilon \, \delta \, |I_0|^{-1}
\end{equation}
for every closed interval $I \subseteq {\bf R}$ with $|I| = |I_0|$.

        If $I$ is a closed interval in ${\bf R}$ with $|I| = l \, |I_0|$
for some positive integer $l$, then the average of $f$ over $I$ can be 
expressed as an average of the averages of $f$ over $l$ closed
intervals of length $|I_0|$, as before.  This implies that 
(\ref{frac{1}{|I|} int_I f(x) dx ge epsilon delta |I_0|^{-1}})
also holds in this situation.  Applying this to a sequence of closed
intervals with length $l \, |I_0|$ for each $l \in {\bf Z}_+$,
we get that
\begin{equation}
\label{Lambda(f) ge epsilon delta |I_0|^{-1}}
        \Lambda(f) \ge \epsilon \, \delta \, |I_0|^{-1}
\end{equation}
under these conditions.  This shows that $\Lambda(f) > 0$ when $f \in
\mathcal{AP}({\bf R})$ is real-valued and nonnegative on ${\bf R}$,
and $f(x_0) > 0$ for some $x_0 \in {\bf R}$.

\section{Some other perspectives}
\label{some other perspectives}

        Put
\begin{equation}
\label{e_t(x) = exp(i x t)}
        e_t(x) = \exp(i \, x \, t)
\end{equation}
for each $x, t \in {\bf R}$, so that $e_t$ is a continuous
homomorphism from ${\bf R}$ into ${\bf T}$ for each $t \in {\bf R}$,
and every continuous homomorphism from ${\bf R}$ into ${\bf T}$ is of
this form.  If $t \ne 0$, then
\begin{equation}
\label{int_a^b e_t(x) dx = - i t^{-1} (e_t(b) - e_t(a))}
        \int_a^b e_t(x) \, dx = - i \, t^{-1} \, (e_t(b) - e_t(a))
\end{equation}
for every $a, b \in {\bf R}$ with $a \le b$, which implies that
\begin{equation}
\label{|int_a^b e_t(x) dx| = |t|^{-1} |e_t(b) - e_t(a)| le 2 |t|^{-1}}
 \qquad  \frac{1}{b - a} \, \biggl|\int_a^b e_t(x) \, dx\biggr|
                         = |t|^{-1} \, |e_t(b) - e_t(a)| \, (b - a)^{-1}
                                       \le 2 \, |t|^{-1} \, (b - a)^{-1}.
\end{equation}
Note that $e_t \in \mathcal{AP}({\bf R})$ for every $t \in {\bf R}$,
so that $\Lambda(e_t)$ is defined as in the previous section for every
$t \in {\bf R}$.  It is easy to see that
\begin{equation}
\label{Lambda(e_t) = 0}
        \Lambda(e_t) = 0
\end{equation}
when $t \ne 0$, because (\ref{|int_a^b e_t(x) dx| = |t|^{-1} |e_t(b) -
  e_t(a)| le 2 |t|^{-1}}) tends to $0$ as $(b - a) \to \infty$.  Of
course, $e_0(x) = 1$ for every $x \in {\bf R}$, so that $\Lambda(e_0)
= 1$, as in (\ref{Lambda({bf 1}_{bf R}) = 1}).

        Let $I_1, I_2, I_3, \ldots$ be a sequence of closed intervals
in ${\bf R}$ such that $|I_j| > 0$ for each $j \in {\bf Z}_+$ and
$|I_j| \to \infty$ as $j \to \infty$.  Consider the set $E(\{I_j\}_{j
  = 1}^\infty)$ of $f \in C_b({\bf R})$ such that the corresponding
sequence of averages
\begin{equation}
\label{frac{1}{|I_j|} int_{I_j} f(x) dx, 2}
        \frac{1}{|I_j|} \, \int_{I_j} f(x) \, dx
\end{equation}
converges in ${\bf C}$.  It is easy to see that $E(\{I_j\}_{j =
  1}^\infty)$ is a linear subspace of $C_b({\bf R})$ that is invariant
under complex conjugation and contains the constant functions.  We
also have that $e_t \in E(\{I_j\}_{j = 1}^\infty)$ for every $t \in
{\bf R}$ with $t \ne 0$, since the corresponding sequence of averages
of $e_t$ converges to $0$ when $t \ne 0$, by (\ref{|int_a^b e_t(x) dx|
  = |t|^{-1} |e_t(b) - e_t(a)| le 2 |t|^{-1}}).

        Equivalently, $E(\{I_j\}_{j = 1}^\infty)$ consists of the
$f \in C_b({\bf R})$ such that the sequence of averages 
(\ref{frac{1}{|I_j|} int_{I_j} f(x) dx, 2}) is a Cauchy sequence
in ${\bf C}$.  Using this, one can check that $E(\{I_j\}_{j = 1}^\infty)$
is a closed set in $C_b({\bf R})$ with respect to the supremum metric.
This also uses the fact that
\begin{equation}
\label{|frac{1}{|I_j|} int_{I_j} f(x) dx| le ... le ||f||_{sup}}
        \biggl|\frac{1}{|I_j|} \, \int_{I_j} f(x) \, dx\biggr|
            \le \frac{1}{|I_j|} \, \int_{I_j} |f(x)| \, dx \le \|f\|_{sup}
\end{equation}
for every $f \in C_b({\bf R})$ and $j \ge 1$, where $\|f\|_{sup}$
denotes the supremum norm of $f$.  It follows that
\begin{equation}
\label{mathcal{AP}({bf R}) subseteq E({I_j}_{j = 1}^infty)}
        \mathcal{AP}({\bf R}) \subseteq E(\{I_j\}_{j = 1}^\infty),
\end{equation}
because $e_t \in E(\{I_j\}_{j = 1}^\infty)$ for every $t \in {\bf R}$,
as in the previous paragraph, and the linear span of the set of $e_t$
with $t \in {\bf R}$ is dense in $\mathcal{AP}({\bf R})$ with respect
to the supremum norm.

        Similarly, let $I_1', I_2', I_3', \ldots$ be another sequence
of closed intervals in ${\bf R}$ such that $|I_j'| > 0$ for each $j
\in {\bf Z}_+$ and $|I_j'| \to \infty$ as $j \to \infty$.  Consider
the set of $f \in C_b({\bf R})$ such that
\begin{equation}
\label{frac{1}{|I_j'|} int_{I_j'} f(x) dx - frac{1}{|I_j|} int_{I_j} f(x) dx}
        \frac{1}{|I_j'|} \, \int_{I_j'} f(x) \, dx
           - \frac{1}{|I_j|} \, \int_{I_j} f(x) \, dx
\end{equation}
converges to $0$ as $j \to \infty$.  It is easy to see that this is a
linear subspace of $C_b({\bf R})$ that is invariant under complex
conjugation and contains $e_t$ for each $t \in {\bf R}$.  One can also
check that this is a closed set in $C_b({\bf R})$ with respect to the
supremum norm, using (\ref{|frac{1}{|I_j|} int_{I_j} f(x) dx| le
  ... le ||f||_{sup}}) and its analogue for $I_j'$.  It follows that
$\mathcal{AP}({\bf R})$ is contained in this subspace of $C_b({\bf
  R})$, for the same reasons as before.

        This gives another way to show that the limit in
(\ref{Lambda(f) = lim_{j to infty} frac{1}{|I_j|} int_{I_j} f(x) dx})
exists for every $f \in \mathcal{AP}({\bf R})$, and that the
limit does not depend on the choice of sequence of $I_j$'s
when $f \in \mathcal{AP}({\bf R})$.  Once one has this, one can
derive properties of $\Lambda(f)$ for $f \in \mathcal{AP}({\bf R})$
like those in the previous section.

        Let $B$ be a compact commutative topological group, and suppose
that $h$ is a continuous homomorphism from ${\bf R}$ into $B$.
Suppose also that $h({\bf R})$ is dense in $B$, which implies that
\begin{equation}
\label{sup_{x in {bf R}} |f(h(x))| = sup_{y in B} |f(y)|}
        \sup_{x \in {\bf R}} |f(h(x))| = \sup_{y \in B} |f(y)|
\end{equation}
for every $f \in C(B)$.  If $f \in C(B)$, then $f \circ h \in
\mathcal{AP}({\bf R})$, as in Section \ref{almost periodic functions}.
Put
\begin{equation}
\label{Lambda_h(f) = Lambda(f circ h)}
        \Lambda_h(f) = \Lambda(f \circ h)
\end{equation}
for each $f \in C(B)$, which defines a nonnegative linear functional
on $C(B)$, because $\Lambda$ is a nonnegative linear functional on
$\mathcal{AP}({\bf R})$.  Similarly,
\begin{equation}
\label{|Lambda_h(f)| = ... le sup_{x in {bf R}} |f(h(x))| = sup_{y in B} |f(y)|}
        |\Lambda_h(f)| = |\Lambda(f \circ h)| \le \sup_{x \in {\bf R}} |f(h(x))|
                                               = \sup_{y \in B} |f(y)|
\end{equation}
for every $f \in C(B)$, by (\ref{|Lambda(f)| le ||f||_{sup}}).  If
${\bf 1}_B$ is the constant function on $B$ equal to $1$ at every
point, then
\begin{equation}
\label{Lambda_h({bf 1}_B) = Lambda({bf 1}_{bf R}) = 1}
        \Lambda_h({\bf 1}_B) = \Lambda({\bf 1}_{\bf R}) = 1,
\end{equation}
by (\ref{Lambda({bf 1}_{bf R}) = 1}).  If $f$ is real-valued and
nonnegative on $B$, and if $f > 0$ somewhere on $B$, then
$f \circ h > 0$ somewhere on ${\bf R}$, because $h({\bf R})$
is dense in $B$.  Under these conditions,
\begin{equation}
\label{Lambda_h(f) = Lambda(f circ h) > 0}
        \Lambda_h(f) = \Lambda(f \circ h) > 0,
\end{equation}
because of the analogous property of $\Lambda$, discussed in the
previous section.

        It is easy to see that (\ref{Lambda_h(f) = Lambda(f circ h)})
is invariant under translations of $f$ by elements of $h({\bf R})$,
because $\Lambda$ is invariant under translations on
$\mathcal{AP}({\bf R})$, as in (\ref{Lambda(f_b) = Lambda(f)}).  Using
this, one can check that (\ref{Lambda_h(f) = Lambda(f circ h)}) is
invariant under translations of $f$ by arbitrary elements of $B$.
More precisely, this also uses the hypothesis that $h({\bf R})$ be
dense in $B$, and the fact that $f$ is uniformly continuous on $B$,
because $B$ is compact.  This permits translations of $f$ by elements
of $B$ to be approximated by translations of $f$ by elements of
$h({\bf R})$, uniformly on $B$.  The remaining point is that
$\Lambda_h$ is bounded with respect to the supremum norm on $C(B)$, as
in (\ref{|Lambda_h(f)| = ... le sup_{x in {bf R}} |f(h(x))| = sup_{y
    in B} |f(y)|}), which implies that $\Lambda_h$ is continuous with
respect to the supremum metric on $C(B)$.

        Thus $\Lambda_h$ satisfies the requirements of a Haar integral
on $C(B)$, normalized as in (\ref{Lambda_h({bf 1}_B) = Lambda({bf
    1}_{bf R}) = 1}).  This implies that
\begin{equation}
\label{Lambda_h(f) = int_B f(y) dH_B(y)}
        \Lambda_h(f) = \int_B f(y) \, dH_B(y)
\end{equation}
for every $f \in C(B)$, where $H_B$ is Haar measure measure on $B$,
normalized so that $H_B(B) = 1$.  Alternatively, if $\phi$ is a
continuous homomorphism from $B$ into ${\bf T}$, then $\phi \circ h$
is a continuous homomorphism from ${\bf R}$ into ${\bf T}$, and hence
\begin{equation}
        \phi \circ h = e_t
\end{equation}
for some $t \in {\bf R}$.  If $\phi = {\bf 1}_B$, then $t = 0$, and
conversely $\phi = {\bf 1}_B$ when $t = 0$, because $h({\bf R})$ is
supposed to be dense in $B$.  If $\phi \ne {\bf 1}_B$, then
\begin{equation}
\label{int_B phi(y) dH_B(y) = 0}
        \int_B \phi(y) \, dH_B(y) = 0,
\end{equation}
as in (\ref{int_A phi(x) dH(x) = 0}) in Section \ref{compact, discrete
  groups}, and
\begin{equation}
\label{Lambda_h(phi) = Lambda(e_t) = 0}
        \Lambda_h(\phi) = \Lambda(e_t) = 0,
\end{equation}
as in (\ref{Lambda(e_t) = 0}), because $t \ne 0$.  Of course,
(\ref{Lambda_h(f) = int_B f(y) dH_B(y)}) follows from the
normalization $H_B(B) = 1$ when $f$ is a constant function.  This
gives another proof of (\ref{Lambda_h(f) = int_B f(y) dH_B(y)}) when
$f$ is a continuous homomorphism from $B$ into ${\bf T}$.  It follows
that (\ref{Lambda_h(f) = int_B f(y) dH_B(y)}) also holds when $f$ is
in the linear span $\mathcal{E}(B)$ of $\widehat{B}$ in $C(B)$, since
both sides of (\ref{Lambda_h(f) = int_B f(y) dH_B(y)}) are linear in
$f$.  More precisely, both sides of (\ref{Lambda_h(f) = int_B f(y)
  dH_B(y)}) determine bounded linear functionals on $C(B)$ with
respect to the supremum norm, so that (\ref{Lambda_h(f) = int_B f(y)
  dH_B(y)}) holds when $f$ is in the closure of $\mathcal{E}(B)$ in
$C(B)$ too.  This shows that (\ref{Lambda_h(f) = int_B f(y) dH_B(y)})
holds for every $f \in C(B)$, because $\mathcal{E}(B)$ is dense in
$C(B)$ with respect to the supremum norm.

\section{The associated inner product}
\label{the associated inner product}

        If $f, g \in \mathcal{AP}({\bf R})$, then we have seen that
$f \, \overline{g} \in \mathcal{AP}({\bf R})$ too, and so we can put
\begin{equation}
\label{langle f, g rangle = ... = Lambda(f overline{g})}
 \langle f, g \rangle = \langle f, g \rangle_{\mathcal{AP}({\bf R})}
                      = \Lambda(f \, \overline{g}),
\end{equation}
where $\Lambda$ is defined on $\mathcal{AP}({\bf R})$ as in Section
\ref{averages on R}.  By construction, (\ref{langle f, g rangle =
  ... = Lambda(f overline{g})}) is linear in $f$, conjugate-linear in
$g$, and satisfies
\begin{equation}
\label{langle g, f rangle = overline{langle f, g rangle}}
        \langle g, f \rangle = \overline{\langle f, g \rangle}
\end{equation}
for every $f, g \in \mathcal{AP}({\bf R})$.  We also have that
\begin{equation}
\label{langle f, f rangle = Lambda(|f|^2) ge 0}
        \langle f, f \rangle = \Lambda(|f|^2) \ge 0
\end{equation}
for every $f \in \mathcal{AP}({\bf R})$, and that this is equal to $0$
if and only if $f = 0$ on ${\bf R}$, as in Section \ref{averages on
  R}.  This shows that (\ref{langle f, g rangle = ... = Lambda(f
  overline{g})}) defines an inner product on $\mathcal{AP}({\bf R})$.
If $f, g \in \mathcal{AP}({\bf R})$ and $a \in {\bf R}$, then their
translates $f_a$, $g_a$ by $a$ are in $\mathcal{AP}({\bf R})$ too, and
\begin{equation}
\label{langle f_a, g_a rangle = langle f, g rangle}
        \langle f_a, g_a \rangle = \langle f, g \rangle,
\end{equation}
by the analogous property (\ref{Lambda(f_b) = Lambda(f)}) of
$\Lambda$.

        Let $e_t(x)$ be as in (\ref{e_t(x) = exp(i x t)}), and
observe that
\begin{equation}
\label{e_r(x) overline{e_t(x)} = e_r(x) e_{-t}(x) = e_{r - t}(x)}
        e_r(x) \, \overline{e_t(x)} = e_r(x) \, e_{-t}(x) = e_{r - t}(x)
\end{equation}
for every $x, r, t \in {\bf R}$.  It follows that
\begin{equation}
\label{langle e_r, e_t rangle = Lambda(e_{r - t}) = 0}
        \langle e_r, e_t \rangle = \Lambda(e_{r - t}) = 0
\end{equation}
when $r \ne t$, by (\ref{Lambda(e_t) = 0}), and that
\begin{equation}
\label{langle e_t, e_t rangle = Lambda(e_0) = 1}
        \langle e_t, e_t \rangle = \Lambda(e_0) = 1
\end{equation}
for every $t \in {\bf R}$, by (\ref{Lambda({bf 1}_{bf R}) = 1}).  If
$f \in \mathcal{AP}({\bf R})$ and $E$ is a finite subset of ${\bf R}$,
then
\begin{equation}
\label{sum_{t in E} |langle f, e_t rangle|^2 le langle f, f rangle}
        \sum_{t \in E} |\langle f, e_t \rangle|^2 \le \langle f, f \rangle,
\end{equation}
by standard results about inner product spaces.  Taking the supremum
over all finite subsets $E$ of ${\bf R}$, we get that
\begin{equation}
\label{sum_{t in {bf R}} |langle f, e_t rangle|^2 le langle f, f rangle}
        \sum_{t \in {\bf R}} |\langle f, e_t \rangle|^2 \le \langle f, f \rangle,
\end{equation}
and in particular that the sum on the left is finite.  Of course,
(\ref{sum_{t in {bf R}} |langle f, e_t rangle|^2 le langle f, f
  rangle}) implies that $\langle f, e_t \rangle \ne 0$ for only
finitely or countably many $t \in {\bf R}$.

        Note that
\begin{equation}
\label{langle f, f rangle = Lambda(|f|^2) le ||f||_{sup}^2}
        \langle f, f \rangle = \Lambda(|f|^2) \le \|f\|_{sup}^2
\end{equation}
for every $f \in \mathcal{AP}({\bf R})$, where $\|f\|_{sup}$ is the
supremum norm of $f$ on ${\bf R}$, by (\ref{|Lambda(f)| le
  ||f||_{sup}}).  As before, the linear span of the set of $e_t$ with
$t \in {\bf R}$ is dense in $\mathcal{AP}({\bf R})$ with respect to
the supremum norm.  This implies that the linear span of the set of
$e_t$ with $t \in {\bf R}$ is dense in $\mathcal{AP}({\bf R})$ with
respect to the norm associated to $\langle \cdot, \cdot \rangle$, by
(\ref{langle f, f rangle = Lambda(|f|^2) le ||f||_{sup}^2}).  It
follows that
\begin{equation}
\label{sum_{t in {bf R}} |langle f, e_t rangle|^2 = langle f, f rangle}
        \sum_{t \in {\bf R}} |\langle f, e_t \rangle|^2 = \langle f, f \rangle
\end{equation}
for every $f \in \mathcal{AP}({\bf R})$.

        If $f, g \in \mathcal{AP}({\bf R})$, then
\begin{eqnarray}
\label{sum_r |langle f, e_{t - r} rangle| |langle g, e_r rangle| le ...}
 \sum_{r \in {\bf R}} |\langle f, e_{t - r} \rangle| \, |\langle g, e_r \rangle|
 & \le & \Big(\sum_{r \in {\bf R}} |\langle f, e_{t - r} \rangle|^2\Big)^{1/2}
       \, \Big(\sum_{r \in {\bf R}} |\langle g, e_r \rangle|^2\Big)^{1/2} \\
  & = & \Big(\sum_{r \in {\bf R}} |\langle f, e_r \rangle|^2\Big)^{1/2}
      \, \Big(\sum_{r \in {\bf R}} |\langle g, e_r \rangle|^2\Big)^{1/2} \nonumber
\end{eqnarray}
for every $t \in {\bf R}$, using the Cauchy--Schwarz inequality in the
first step, and a change of variables in the second step.  Thus
\begin{equation}
\label{sum_r |langle f, e_{t - r} rangle| |langle g, e_r rangle| le ..., 2}
 \sum_{r \in {\bf R}} |\langle f, e_{t - r} \rangle| \, |\langle g, e_r \rangle|
        \le \langle f, f \rangle^{1/2} \, \langle g, g \rangle^{1/2}
         \le \|f\|_{sup} \, \|g\|_{sup}
\end{equation}
for every $t \in {\bf R}$, by (\ref{sum_{t in {bf R}} |langle f, e_t
  rangle|^2 le langle f, f rangle}) and (\ref{langle f, f rangle =
  Lambda(|f|^2) le ||f||_{sup}^2}).  Remember that $f \, g \in
\mathcal{AP}({\bf R})$, so that $\langle f \, g, e_t\rangle$ is
defined for every $t \in {\bf R}$.  Let us check that
\begin{equation}
\label{langle f g, e_t rangle = ...}
        \langle f \, g, e_t \rangle =
 \sum_{r \in {\bf R}} \langle f, e_{t - r} \rangle \, \langle g, e_r\rangle
\end{equation}
for every $t \in {\bf R}$, where the summability of the sum on the
right follows from (\ref{sum_r |langle f, e_{t - r} rangle| |langle g,
  e_r rangle| le ..., 2}).  If $g = e_w$ for some $w \in {\bf R}$,
then it is easy to see that (\ref{langle f g, e_t rangle = ...})
holds, using the definition of the inner product and the
orthonormality properties mentioned earlier.  This implies that
(\ref{langle f g, e_t rangle = ...}) also holds when $g$ is in the
linear span of the set of $e_w$ with $w \in {\bf R}$, by linearity.
In order to get the same conclusion for every $g \in \mathcal{AP}({\bf
  R})$, one can approximate $g$ by elements of the linear span of the
set of $e_w$ with $w \in {\bf R}$, as before.

\section{Some additional properties}
\label{some additional properties}

        Let $B$ be a compact commutative topological group, and suppose
that $h$ is a continuous homomorphism from ${\bf R}$ into $B$.  As
usual, we may as well ask that $h({\bf R})$ be dense in $B$, since
otherwise we can replace $B$ with the closure of $h({\bf R})$ in $B$.
Also let $H_B$ be Haar measure on $B$, normalized so that $H_B(B) =
1$.  If $f, g \in C(B)$, then we have seen that $f \circ h, g \circ h
\in \mathcal{A}({\bf R})$, so that the inner product of $f \circ h$
and $g \circ h$ can be defined as in the previous section.  Under these
conditions,
\begin{equation}
\label{langle f circ h, g circ h rangle_{mathcal{AP}({bf R})} = ...}
        \langle f \circ h, g \circ h \rangle_{\mathcal{AP}({\bf R})}
          = \int_B f(y) \, \overline{g(y)} \, dH_B(y),
\end{equation}
by (\ref{Lambda_h(f) = int_B f(y) dH_B(y)}) and (\ref{langle f, g
  rangle = ... = Lambda(f overline{g})}).  Note that
\begin{equation}
\label{{f circ h : f in C(B)}}
        \{f \circ h : f \in C(B)\}
\end{equation}
is a closed linear subspace of $\mathcal{AP}({\bf R})$ with respect to
the supremum norm.  This uses the fact that
\begin{equation}
\label{f mapsto f circ h}
        f \mapsto f \circ h
\end{equation}
is an isometric linear mapping from $C(B)$ into $\mathcal{AP}({\bf
  R})$ with respect to the corresponding supremum norms, as in
(\ref{sup_{x in {bf R}} |f(h(x))| = sup_{y in B} |f(y)|}), and the
completeness of $C(B)$.

        Let $\widehat{B}$ be the dual group of continuous homomorphisms
from $B$ into ${\bf T}$, which is equipped with the discrete topology,
because $B$ is compact.  If $\phi \in \widehat{B}$, then $\phi \circ
h$ is a continuous homomorphism from ${\bf R}$ into ${\bf T}$, and
hence there is a unique real number $t_h(\phi)$ such that
\begin{equation}
\label{phi(h(x)) = e_{t_h(phi)}(x) = exp(i x t_h(phi))}
        \phi(h(x)) = e_{t_h(\phi)}(x) = \exp(i \, x \, t_h(\phi))
\end{equation}
for every $x \in {\bf R}$.  It is easy to see that the mapping from
$\phi \in \widehat{B}$ to $t_h(\phi) \in {\bf R}$ is a homomorphism
with respect to the usual group structure on $\widehat{B}$, defined by
pointwise multiplication of functions from $B$ into ${\bf T}$, and
addition on ${\bf R}$.  This homomorphism corresponds exactly to the
dual homomorphism $\widehat{h}(\phi) = \phi \circ h$ from
$\widehat{B}$ into the dual of ${\bf R}$, using the identification
between ${\bf R}$ and its dual given by the mapping from $t \in {\bf
  R}$ to $e_t$.  Note that $\widehat{h}$ is injective, because $h({\bf
  R})$ is dense in $B$, so that $t_h$ is injective as a mapping from
$\widehat{B}$ into ${\bf R}$.

        If $f \in C(B)$ and $\phi \in \widehat{B}$, then
\begin{equation}
\label{widehat{f}(phi) = ... = langle f circ h, phi circ h rangle}
        \widehat{f}(\phi) = \int_B f(y) \, \overline{\phi(y)} \, dH_B(y)
           = \langle f \circ h, \phi \circ h \rangle_{\mathcal{AP}({\bf R})},
\end{equation}
by the definition of the Fourier transform $\widehat{f}(\phi)$
and (\ref{langle f circ h, g circ h rangle_{mathcal{AP}({bf R})} = ...}).
Equivalently,
\begin{equation}
\label{widehat{f}(phi) = langle f circ h, e_{t_h(phi)} rangle}
        \widehat{f}(\phi)
             = \langle f \circ h, e_{t_h(\phi)} \rangle_{\mathcal{AP}({\bf R})}
\end{equation}
for every $\phi \in \widehat{B}$, by the definition (\ref{phi(h(x)) =
  e_{t_h(phi)}(x) = exp(i x t_h(phi))}) of $t_h(\phi)$.  Remember that
$f$ can be approximated by finite linear combinations of elements of
$\widehat{B}$, uniformly on $B$.  This implies that $f \circ h$ can be
approximated by finite linear combinations of functions of the form
$e_{t_h(\phi)}$ with $\phi \in \widehat{B}$, uniformly on ${\bf R}$.
Using this and (\ref{langle e_r, e_t rangle = Lambda(e_{r - t}) = 0}),
it follows that
\begin{equation}
\label{langle f circ h, e_t rangle_{mathcal{AP}({bf R})} = 0}
        \langle f \circ h, e_t \rangle_{\mathcal{AP}({\bf R})} = 0
\end{equation}
for every $t \in {\bf R}$ which is not of the form $t_h(\phi)$ for
some $\phi \in \widehat{B}$.  

        More precisely, $f$ can be approximated by finite linear
combinations of $\phi \in \widehat{B}$ with $\widehat{f}(\phi) \ne 0$,
uniformly on $B$, as in Section \ref{compact groups, continued}.  This
implies that $f \circ h$ can be approximated by finite linear
combinations of functions of the form $e_{t_h(\phi)}$, where $\phi \in
\widehat{B}$ and (\ref{widehat{f}(phi) = langle f circ h, e_{t_h(phi)}
  rangle}) is not zero, uniformly on ${\bf R}$.  Because of
(\ref{langle f circ h, e_t rangle_{mathcal{AP}({bf R})} = 0}), this is
the same as saying that $f \circ h$ can be approximated by linear
combinations of functions of the form $e_t$, where $t \in {\bf R}$ and
\begin{equation}
\label{langle f circ h, e_t rangle_{mathcal{AP}({bf R})} ne 0}
        \langle f \circ h, e_t \rangle_{\mathcal{AP}({\bf R})} \ne 0,
\end{equation}
uniformly on ${\bf R}$.  If $F$ is any element of $\mathcal{AP}({\bf
  R})$, then $F$ can be expressed as $f \circ h$ for some $B$ and $h$
as before.  It follows that every $F \in \mathcal{AP}({\bf R})$ can be
approximated by linear combinations of functions of the form $e_t$,
where $t \in {\bf R}$ and
\begin{equation}
\label{langle F, e_t rangle_{mathcal{AP}({bf R})} ne 0}
        \langle F, e_t \rangle_{\mathcal{AP}({\bf R})} \ne 0,
\end{equation}
uniformly on ${\bf R}$.

\section{Bounded holomorphic functions}
\label{bounded holomorphic functions}

        Let $f(z)$ be a continuous complex-valued function on the
closed upper half-plane $\overline{U}$ that is holomorphic on the
open upper half-plane $U$.  Thus
\begin{equation}
\label{f_t(z) = f(z) exp(- i z t)}
        f_t(z) = f(z) \, \exp(- i \, z \, t)
\end{equation}
is a continuous function of $z$ on $\overline{U}$ that is holomorphic
on $U$ for every $t \in {\bf R}$.  Also let $I = [a, b]$ be a closed
interval in ${\bf R}$, and let $R > 0$ be given.  Using Cauchy's
theorem, we get that
\begin{eqnarray}
\label{int_a^b f_t(x) dx + ... = 0}
        \int_a^b f_t(x) \ dx \ + \, \int_0^R f_t(b + i \, y) \, i \, dy
                                & - & \int_a^b f_t(x + i \, R) \, dx \\
                        & - & \int_0^R f_t(a + i \, y) \, i \, dy = 0 \nonumber
\end{eqnarray}
for every $t \in {\bf R}$.  It follows that
\begin{eqnarray}
\label{|int_a^b f_t(x) dx| le ...}
 \biggl|\int_a^b f_t(x) \, dx\biggr| \le \int_0^R |f_t(b + i \, y)| \, dy
                             & + & \int_a^b |f_t(x + i \, R)| \, dx \\
                             & + & \int_0^R |f_t(a + i \, y)| \, dy \nonumber
\end{eqnarray}
for every $t \in {\bf R}$.

        Suppose now that $f(z)$ is also bounded on $\overline{U}$,
and let $\|f\|_{sup}$ be the supremum norm of $f$ on $\overline{U}$,
which is the same as the supremum norm of $f$ on ${\bf R}$, as in
Section \ref{bounded harmonic functions}.  This implies that
\begin{eqnarray}
\label{int_a^b |f_t(x + i R)| dx = ... le ||f||_{sup} exp(t R) (b - a)}
        \int_a^b |f_t(x + i \, R)| \, dx 
           & = & \int_a^b |f(x + i \, R)| \, \exp (t \, R) \, dx \\
          & \le & \|f\|_{sup} \, \exp(t \, R) \, (b - a) \nonumber
\end{eqnarray}
for every $t \in {\bf R}$.  If $t < 0$, then $\exp(t \, R) \to 0$ as
$R \to \infty$, and we get that
\begin{equation}
\label{|int_a^b f_t(x) dx| le ..., 2}
 \biggl|\int_a^b f_t(x) \, dx\biggr| \le \int_0^\infty |f_t(b + i \, y)| \, dy
                                        + \int_0^\infty |f_t(a + i \, y)| \, dy,
\end{equation}
by (\ref{|int_a^b f_t(x) dx| le ...}).  Similarly,
\begin{eqnarray}
\label{int_0^infty |f_t(x + i y)| dy le ... = |t|^{-1} ||f||_{sup}}
        \int_0^\infty |f_t(x + i \, y)| \, dy
           & = & \int_0^\infty |f(x + i \, y)| \, \exp(t \, y) \, dy \\
          & \le & \|f\|_{sup} \, \int_0^\infty \exp(t \, y) \, dy
             = |t|^{-1} \, \|f\|_{sup} \nonumber
\end{eqnarray}
for every $x \in {\bf R}$ when $t < 0$.  Combining this with
(\ref{|int_a^b f_t(x) dx| le ..., 2}), we get that
\begin{equation}
\label{|int_a^b f_t(x) dx| le 2 |t|^{-1} ||f||_{sup}}
        \biggl|\int_a^b f_t(x) \, dx\biggr| \le 2 \, |t|^{-1} \, \|f\|_{sup}
\end{equation}
when $t < 0$.

        Suppose that the restriction of $f$ to ${\bf R}$ is almost
periodic, in addition to the other conditions on $f$ already mentioned.
This implies that
\begin{equation}
\label{f_t(x) = f(x) exp(- i x t) = f(x) e_{-t}(x)}
        f_t(x) = f(x) \, \exp(- i \, x \, t) = f(x) \, e_{-t}(x)
\end{equation}
is almost periodic on ${\bf R}$ for every $t \in {\bf R}$ as well,
because $e_{-t}(x)$ is almost periodic on ${\bf R}$.  Thus
$\Lambda(f_t)$ is defined for every $t \in {\bf R}$, and it is easy to
see that
\begin{equation}
\label{Lambda(f_t) = 0}
        \Lambda(f_t) = 0
\end{equation}
when $t < 0$, using (\ref{|int_a^b f_t(x) dx| le 2 |t|^{-1} ||f||_{sup}}).
Equivalently, this means that
\begin{equation}
\label{langle f, e_t rangle_{mathcal{AP}({bf R})} = 0}
        \langle f, e_t \rangle_{\mathcal{AP}({\bf R})} = 0
\end{equation}
when $t < 0$, because $\overline{e_t(x)} = e_{-t}(x)$ for every $x, t
\in {\bf R}$.  Conversely, if $f \in \mathcal{AP}({\bf R})$ satisfies
(\ref{langle f, e_t rangle_{mathcal{AP}({bf R})} = 0}) for every $t
< 0$, then we would like to show that $f$ can be extended to a bounded
continuous function on $\overline{U}$ that is holomorphic on $U$.

        If $f \in \mathcal{AP}({\bf R})$ satisfies
(\ref{langle f, e_t rangle_{mathcal{AP}({bf R})} = 0})
for every $t < 0$, then $f(x)$ can be approximated by finite linear
combinations of functions of the form $e_t(x)$ with $t \ge 0$,
uniformly on ${\bf R}$, as in the previous section.  Of course,
$e_t(x) = \exp(i \, x \, t)$ has an obvious holomorphic extension to
$z \in {\bf C}$ for each $t \in {\bf R}$, which is $\exp(i \, z \,
t)$.  If $t \ge 0$, then $\exp(i \, z \, t)$ is also bounded on the
closed upper half-plane $\overline{U}$.  Thus our hypothesis on $f$
implies that there is a sequence $\{g_j(z)\}_{j = 1}^\infty$ of
bounded continuous functions on $\overline{U}$ that are holomorphic on
$U$ and whose restrictions to ${\bf R}$ converge uniformly to $f$.  As
in Section \ref{bounded harmonic functions},
\begin{equation}
\label{sup_{z in overline{U}} |g_j(z) - g_k(z)| = ...}
        \sup_{z \in \overline{U}} |g_j(z) - g_k(z)|
           = \sup_{x \in {\bf R}} |g_j(x) - g_k(x)|
\end{equation}
for every $j, k \ge 1$.  Because $\{g_j(x)\}_{j = 1}^\infty$ converges
to $f(x)$ uniformly on ${\bf R}$, $\{g_j(x)\}_{j = 1}^\infty$ is a
Cauchy sequence with respect to the supremum norm on $C_b({\bf R})$.
This implies that $\{g_j(z)\}_{j = 1}^\infty$ is a Cauchy sequence
with respect to the supremum norm on $C_b(\overline{U})$, because of
(\ref{sup_{z in overline{U}} |g_j(z) - g_k(z)| = ...}).  It is well
known that $\{g_j(z)\}_{j = 1}^\infty$ converges uniformly to a bounded
continuous function on $\overline{U}$ under these conditions, and that
the limit is holomorphic on $U$.  By construction, the limit is equal
to $f$ on ${\bf R}$, which shows that $f$ can be extended to a bounded
continuous function on $\overline{U}$ that is holomorphic on $U$.

\section{Some other averages}
\label{some other averages}

          Let $f$ be a bounded measurable complex-valued function
on ${\bf R}$, and put
\begin{equation}
\label{A_y(f)(x) = frac{1}{pi} int_{bf R} f(x - w) frac{y}{w^2 + y^2} dw}
 A_y(f)(x) = \frac{1}{\pi} \, \int_{\bf R} f(x - w) \, \frac{y}{w^2 + y^2} \, dw
\end{equation}
for every $x, y \in {\bf R}$ with $y > 0$.  This is equivalent to
(\ref{A_r(f)(w) = frac{1}{pi} int_{bf R} f(t) frac{r}{r^2 + (t - w)^2}
  dt}) in Section \ref{multiplication formula}, and to (\ref{h(z) =
  A_y(f)(x) = ...})  in Section \ref{harmonic extensions}.  Remember
that
\begin{equation}
\label{A_y({bf 1}_{bf R})(x) = 1}
        A_y({\bf 1}_{\bf R})(x) = 1
\end{equation}
for every $x \in {\bf R}$ and $y > 0$, by (\ref{frac{1}{pi} int_{bf R}
  frac{r}{r^2 + w^2} dw = 1}) in Section \ref{multiplication formula},
where ${\bf 1}_{\bf R}$ is the constant function equal to $1$ at every
point in ${\bf R}$.  Suppose that $f \in \mathcal{AP}({\bf R})$, and
let $\Lambda(f)$ be as in Section \ref{averages on R}.  Under these
conditions,
\begin{equation}
\label{A_y(f)(x) to Lambda(f) as y to infty}
        A_y(f)(x) \to \Lambda(f) \quad\hbox{as } y \to \infty,
\end{equation}
where the convergence is uniform over $x \in {\bf R}$.

        One way to see this is to express $A_y(f)(x)$ as an average
of averages of $f$ over intervals of the form $[x - r, x + r]$, where
$r > 0$.  This uses the fact that $y / (w^2 + y^2)$ is monotonically
decreasing as a function of $|w|$, and can be made precise by
integration by parts.  As $y$ increases, these averages of averages of
$f$ become concentrated on averages of $f$ over intervals of the form
$[x - r, x+ r]$ with $r$ large.  As in Section \ref{averages on R},
the average of $f$ over $[x - r, x + r]$ tends to $\Lambda(f)$ as
$r \to \infty$, uniformly over $x \in {\bf R}$.  This permits one to
show that (\ref{A_y(f)(x) to Lambda(f) as y to infty}) holds uniformly
over $x \in {\bf R}$, as desired.

        Alternatively, let $t \in {\bf R}$ be given, and consider
$A_y(e_t)$, where $e_t$ is as in (\ref{e_t(x) = exp(i x t)}) in
Section \ref{some other perspectives}.  Observe that
\begin{equation}
\label{A_y(e_t)(x) = e_t(x) A_y(e_t)(0)}
        A_y(e_t)(x) = e_t(x) \, A_y(e_t)(0)
\end{equation}
for every $x \in {\bf R}$ and $y > 0$, which is basically an instance
of (\ref{(phi * nu)(x) = widehat{nu}(phi) phi(x)}) in Section
\ref{functions, measures}, because $A_y(e_t)$ is the convolution of
$e_t$ with an integrable function on ${\bf R}$.  If $t = 0$, then
(\ref{A_y(e_t)(x) = e_t(x) A_y(e_t)(0)}) is equal to $1$ for every $x
\in {\bf R}$ and $y > 0$, as in (\ref{A_y({bf 1}_{bf R})(x) = 1}).  If
$t \ne 0$, then
\begin{equation}
\label{lim_{y to infty} A_y(e_t)(0) = 0}
        \lim_{y \to \infty} A_y(e_t)(0) = 0,
\end{equation}
and we shall say more about this in a moment.  This is consistent with
(\ref{A_y(f)(x) to Lambda(f) as y to infty}), because of
(\ref{Lambda(e_t) = 0}) in Section \ref{some other perspectives}.

         It is not too difficult to check (\ref{lim_{y to infty} A_y(e_t)(0) 
= 0}) directly, by integration by parts, for instance.  One can also
look at $A_y(e_t)(0)$ in terms of the Fourier transform of an explicit
function on ${\bf R}$.  This function was obtained as the Fourier
transform of another explicit function in Section \ref{some examples},
which can be used to determine $A_y(e_t)(0)$, and to verify
(\ref{lim_{y to infty} A_y(e_t)(0) = 0}).  As another approach,
$A_y(e_t)(x)$ corresponds to the unique extension of $e_t$ to a
continuous function on the closed upper half-plane $\overline{U}$ that
is harmonic on the open upper half-plane $U$, as in Sections
\ref{harmonic extensions} and \ref{bounded harmonic functions}.  In
this case, this extension of $e_t$ can be given explicitly in terms of
the complex exponential function, and satisfies (\ref{lim_{y to infty}
  A_y(e_t)(0) = 0}).  More precisely, $e_t(x)$ extends to a bounded 
continuous function on $\overline{U}$ that is holomorphic on $U$
when $t > 0$, and which is conjugate-holomorphic on $U$ when $t < 0$.
This permits $A_y(e_t)(x)$ to be evaluated using Cauchy's integral
formula as well.

        At any rate, once we have (\ref{lim_{y to infty} A_y(e_t)(0) = 0})
when $t \ne 0$, we get that (\ref{A_y(f)(x) to Lambda(f) as y to
  infty}) holds when $f$ is in the linear span of the set of functions
of the form $e_t$ with $t \in {\bf R}$.  This implies that
(\ref{A_y(f)(x) to Lambda(f) as y to infty}) holds for every $f \in
\mathcal{AP}({\bf R})$, because the linear span of the set of
functions of the form $e_t$ with $t \in {\bf R}$ is dense in
$\mathcal{AP}({\bf R})$.  This also uses (\ref{||A_r(f)||_p le
  ||f||_p}) in Section \ref{multiplication formula} with $p = \infty$,
to handle the relevant approximations.

        Let $g$ be a continuous complex-valued function on the unit
circle, let $t > 0$ be given, and put
\begin{equation}
        f(x) = g(e_t(x)) = g(\exp (i \, x \, t))
\end{equation}
for each $x \in {\bf R}$.  Thus $f$ is a continuous function on ${\bf
  R}$, that is periodic with period $2 \, \pi / t$, so that $f$ is an
element of $\mathcal{AP}({\bf R})$ in particular.  The restriction to
$t > 0$ is not too serious here, because one can always reduce to that
case using the fact that $\overline{e_t} = e_{-t}$.  An advantage of
takng $t > 0$ is that $e_t(x)$ extends to $\exp (i \, z \, t)$, as a
continuous function from the closed upper half-plane $\overline{U}$
onto the closed unit disk $\overline{D}$ which is holomorphic on $U$.
If $g(w)$ is the unique extension of $g$ to a continuous function on
$\overline{D}$ that is harmonic on $D$, then
\begin{equation}
\label{f(z) = g(exp (i z t))}
        f(z) = g(\exp (i \, z \, t))
\end{equation}
extends $f$ to a bounded continuous function on $\overline{U}$ that is
harmonic on $U$.  In this case,
\begin{equation}
\label{f(x + i y) to g(0) as y to infty}
        f(x + i \, y) \to g(0) \quad\hbox{as } y \to \infty
\end{equation}
uniformly over $x \in {\bf R}$, because $|\exp (i \, (x + i \, y) \,
t)| = \exp (-y \, t) \to 0$ as $y \to \infty$ and $g(w)$ is continuous
at $0$.  We also have that
\begin{equation}
        A_y(f)(x) = f(x + i \, y)
\end{equation}
for every $x \in {\bf R}$ and $y > 0$ in this situation, because of
the uniqueness of bounded harmonic extensions, as in Section
\ref{bounded harmonic functions}.  Of course, $g(0)$ is equal to the
average of $g$ on ${\bf T}$, which is the same as the average of $f$
on any interval of length $2 \, \pi / t$.  This is also the same as
the average of $f$ on any interval whose length is a positive integer
multiple of $2 \, \pi / t$, which implies that $\Lambda(f) = g(0)$.

        If $f_1, f_2 \in \mathcal{AP}({\bf R})$, then $f_1 \, f_2 \in
\mathcal{AP}({\bf R})$ too, so that $\Lambda(f_1 \, f_2)$ is defined.
If $f_1$, $f_2$ also have bounded continuous extensions to
$\overline{U}$ that are holomorphic on $U$, then $f_1 \, f_2$ has the
same property, and
\begin{equation}
\label{Lambda(f_1 f_2) = Lambda(f_1) Lambda(f_2)}
        \Lambda(f_1 \, f_2) = \Lambda(f_1) \, \Lambda(f_2).
\end{equation}
This can be derived from (\ref{langle f g, e_t rangle = ...}) in
Section \ref{the associated inner product}, with $t = 0$, using the
fact that
\begin{equation}
\label{langle f_1, e_r rangle = langle f_2, e_r rangle = 0}
        \langle f_1, e_r \rangle_{\mathcal{AP}({\bf R})}
         = \langle f_2, e_r \rangle_{\mathcal{AP}({\bf R})} = 0
\end{equation}
for every $r < 0$, as in (\ref{langle f, e_t rangle_{mathcal{AP}({bf
      R})} = 0}).  More precisely, the bounded continuous extension of
$f_1 \, f_2$ to $\overline{U}$ that is holomorphic on $U$ is the
product of the corresponding extensions of $f_1$ and $f_2$ to
$\overline{U}$.  Of course, holomorphic functions are harmonic, and
bounded continuous functions on $\overline{U}$ that are harmonic on
$U$ are uniquely determined by their boundary values on ${\bf R}$, as
in Section \ref{bounded harmonic functions}.  If $f$ is any bounded
continuous function on ${\bf R}$, then $A_y(f)(x)$ defines a harmonic
function on $U$ that determines a bounded continuous function on
$\overline{U}$ equal to $f$ on ${\bf R}$, as in Section \ref{harmonic
extensions}.  It follows that
\begin{equation}
\label{A_y(f_1 f_2)(x) = A_y(f_1) A_y(f_2)(x)}
        A_y(f_1 \, f_2)(x) = A_y(f_1) \, A_y(f_2)(x)
\end{equation}
for every $x \in {\bf R}$ and $y > 0$ under these conditions, so that
(\ref{Lambda(f_1 f_2) = Lambda(f_1) Lambda(f_2)}) can be derived from
(\ref{A_y(f)(x) to Lambda(f) as y to infty}) as well.

\section{Product spaces}
\label{product spaces}

        Let $A$, $B$ be topological spaces, and consider their Cartesian
product $A \times B$, equipped with the product topology.  Also let
$F(a, b)$ be a continuous complex-valued function on $A \times B$.
Thus for each $a \in A$, $b \in B$, and $\epsilon > 0$ there is an
open set $W \subseteq A \times B$ such that $(a, b) \in W$ and
\begin{equation}
\label{|F(a, b) - F(u, v)| < epsilon/2}
        |F(a, b) - F(u, v)| < \epsilon/2
\end{equation}
for every $(u, v) \in W$.  Because of the way that the product
topology is defined, we may as well take $W$ to be of the form $U
\times V$, where $U \subseteq A$ and $V \subseteq B$ are open sets, $a
\in U$, and $b \in V$.  If $u \in U$ and $v \in V$, then $(u, b), (u,
v) \in U \times V$, and hence
\begin{eqnarray}
\label{|F(u, b) - F(u, v)| le ... < epsilon/2 + epsilon/2 = epsilon}
 |F(u, b) - F(u, v)| & \le & |F(u, b) - F(a, b)| + |F(a, b) - F(u, v)| \\
                      & < & \epsilon/2 + \epsilon/2 = \epsilon, \nonumber
\end{eqnarray}
using (\ref{|F(a, b) - F(u, v)| < epsilon/2}) twice in the second step.

        Let $K$ be a nonempty compact subset of $A$, let $b$ be an element
of $B$, and let $\epsilon > 0$ be given again.  If $a \in K$, then
there is an open set $W(a)$ in $A \times B$ such that $(a, b) \in
W(a)$ and (\ref{|F(a, b) - F(u, v)| < epsilon/2}) holds for every $(u,
v) \in W(a)$.  As before, we can take $W(a)$ to be of the form $U(a)
\times V(a)$, where $U(a) \subseteq A$ and $V(a) \subseteq B$ are open
sets, $a \in U(a)$, and $b \in V(a)$.  It follows that (\ref{|F(u, b)
  - F(u, v)| le ... < epsilon/2 + epsilon/2 = epsilon}) holds for
every $u \in U(a)$ and $v \in V(a)$.  In particular, the open sets
$U(a)$ in $A$ with $a \in K$ form an open covering of $K$, and so
there are finitely many elements $a_1, \ldots, a_n$ of $K$ such that
\begin{equation}
\label{K subseteq bigcup_{j = 1}^n U(a_j)}
        K \subseteq \bigcup_{j = 1}^n U(a_j),
\end{equation}
because $K$ is compact.  If we put
\begin{equation}
\label{V = bigcap_{j = 1}^n V(a_j)}
        V = \bigcap_{j = 1}^n V(a_j),
\end{equation}
then $V$ is an open set in $B$ that contains $b$, and
\begin{equation}
\label{|F(u, b) - F(u, v)| < epsilon}
        |F(u, b) - F(u, v)| < \epsilon
\end{equation}
for every $u \in K$ and $v \in V$.  More precisely, if $u \in K$, then
$u \in U(a_j)$ for some $j = 1, \ldots, n$, $v \in V \subseteq V(a_j)$
for the same $j$, and (\ref{|F(u, b) - F(u, v)| < epsilon}) follows
from (\ref{|F(u, b) - F(u, v)| le ... < epsilon/2 + epsilon/2 =
  epsilon}) in the case where $U = U(a_j)$ and $V = V(a_j)$.

        Put
\begin{equation}
\label{F_b(a) = F(a, b)}
        F_b(a) = F(a, b)
\end{equation}
for every $a \in A$ and $b \in B$, which is a continuous function of
$a \in A$ for each $b \in B$.  Thus
\begin{equation}
\label{b mapsto F_b}
        b \mapsto F_b
\end{equation}
defines a mapping from $B$ into the space $C(A)$ of continuous
complex-valued functions on $A$.  The preceding argument shows that
this mapping is continuous with respect to the topology on $C(A)$
corresponding to uniform convergence on compact subsets of $A$, as in
Section \ref{spaces of continuous functions}.  If $A$ is compact, then
we can apply this argument to $K = A$, to get that (\ref{b mapsto
  F_b}) is continuous with respect to the topology on $C(A)$ defined
by the supremum norm.

        Conversely, if (\ref{b mapsto F_b}) is any mapping from $B$
into $C(A)$, then (\ref{F_b(a) = F(a, b)}) is a complex-valued
function on $A \times B$ that is continuous as a function of $a \in A$
for each $b \in B$.  If (\ref{b mapsto F_b}) is a continuous mapping
from $B$ into the space $C_b(A)$ of bounded continuous complex-valued 
functions on $A$, with respect to the topology on $C_b(A)$ determined
by the supremum norm, then one can check that (\ref{F_b(a) = F(a, b)})
is a continuous function on $A \times B$.  Suppose instead that
(\ref{b mapsto F_b}) is continuous as a mapping from $B$ into $C(A)$,
where $C(A)$ is equipped with the topology defined by the supremum
seminorms associated to nonempty compact subsets of $A$, as before.
If $A$ is locally compact, then (\ref{F_b(a) = F(a, b)}) is continuous
on $A \times B$, as in the previous case.

        Let us now take $A$ to be the real line, with the standard
topology, and let $F(x, b)$ be a bounded continuous complex-valued
function on ${\bf R} \times B$.  Let $U$ be the open upper half-plane
in ${\bf C}$, and consider the function $G(z, b)$ defined on $U \times B$ by
\begin{eqnarray}
\label{G(x + i y, b) = frac{1}{pi} int_{bf R} F(x - w, b) frac{y}{w^2 + y^2} dw}
        G(x + i \, y, b) & = & \frac{1}{\pi} \, \int_{\bf R} F(x - w, b) \,
                                             \frac{y}{w^2 + y^2} \, dw \\
 & = & \frac{1}{\pi} \, \int_{\bf R} \, F(w, b) \, \frac{y}{(x - w)^2 + y^2}
                                                    \, dw \nonumber
\end{eqnarray}
for every $x, y \in {\bf R}$ with $y > 0$ and $b \in B$.  This
corresponds to (\ref{A_r(f)(w) = frac{1}{pi} int_{bf R} f(t)
  frac{r}{r^2 + (t - w)^2} dt}) in Section \ref{multiplication
  formula}, (\ref{h(z) = A_y(f)(x) = ...}) in Section \ref{harmonic
  extensions}, and (\ref{A_y(f)(x) = frac{1}{pi} int_{bf R} f(x - w)
  frac{y}{w^2 + y^2} dw}) in Section \ref{some other averages},
applied to $F(x, b)$ as a function of $x \in {\bf R}$ for each $b \in
B$.  Note that $G(z, b)$ is a harmonic function of $z \in U$ for each
$b \in B$, as in Section \ref{harmonic extensions}.

        Put
\begin{equation}
\label{G(x, b) = F(x, b)}
        G(x, b) = F(x, b)
\end{equation}
for every $x \in {\bf R}$ and $b \in B$, so that $G(z, b)$ is defined
on $\overline{U} \times B$.  Thus $G(z, b)$ is a bounded continuous
function of $z \in \overline{U}$ for each $b \in U$, as in Section
\ref{harmonic extensions} again.  More precisely,
\begin{equation}
\label{|G(z, b)| le sup_{w in {bf R}} |F(w, b)|}
        |G(z, b)| \le \sup_{w \in {\bf R}} |F(w, b)|
\end{equation}
for every $z \in \overline{U}$ and $b \in B$, as in (\ref{|h(z)| le
  ||f||_infty}).  This implies that $G(z, b)$ is bounded on
$\overline{U} \times B$, because $F(x, b)$ is supposed to be bounded
on ${\bf R} \times B$.  Similarly,
\begin{equation}
\label{|G(z, b) - G(z, b')| le sup_{w in {bf R}} |F(w, b) - F(w, b')|}
        |G(z, b) - G(z, b')| \le \sup_{w \in {\bf R}} |F(w, b) - F(w, b')|
\end{equation}
for every $z \in \overline{U}$ and $b, b' \in B$, for essentially
the same reasons as in (\ref{|G(z, b)| le sup_{w in {bf R}} |F(w, b)|}).

        Suppose that (\ref{b mapsto F_b}) defines a continuous mapping
from $B$ into $C_b({\bf R})$, with respect to the supremum norm on
$C_b({\bf R})$.  Using (\ref{|G(z, b) - G(z, b')| le sup_{w in {bf R}}
  |F(w, b) - F(w, b')|}), it follows easily that $G(z, b)$ corresponds
in the same way to a continuous mapping from $B$ into
$C_b(\overline{U})$, with respect to the supremum norm on
$C_b(\overline{U})$.  In particular, this implies that $G(z, b)$ is a
continuous function on $\overline{U} \times B$ with respect to the
product topology, since we already know that $G(z, b)$ is continuous
as a function of $z \in \overline{U}$ for each $b \in B$, as mentioned
earlier.  One can also check that $G(z, b)$ is continuous on
$\overline{U} \times B$ when $F(x, b)$ is bounded and continuous on
${\bf R} \times B$.

        Suppose now that
\begin{equation}
\label{F_b(x) = F(x, b) in mathcal{AP}({bf R})}
        F_b(x) = F(x, b) \in \mathcal{AP}({\bf R})
\end{equation}
as a function of $x \in {\bf R}$ for every $b \in B$, so that
\begin{equation}
\label{lambda(b) = Lambda(F_b)}
        \lambda(b) = \Lambda(F_b)
\end{equation}
is defined for every $b \in B$, as in Section \ref{averages on R}.
Observe that
\begin{equation}
\label{lim_{y to infty} G(x + i y, b) = lambda(b)}
        \lim_{y \to \infty} G(x + i \, y, b) = \lambda(b)
\end{equation}
for every $x \in {\bf R}$ and $b \in B$, by (\ref{A_y(f)(x) to
  Lambda(f) as y to infty}) in Section \ref{some other averages}, and
where the convergence is uniform over $x \in {\bf R}$ for each $b \in
B$.  If (\ref{b mapsto F_b}) is a continuous mapping from $B$ into
$\mathcal{AP}({\bf R})$ with respect to the supremum norm on
$\mathcal{AP}({\bf R})$, then $\lambda(b)$ is a continuous function on
$B$, because $\Lambda$ is a bounded linear functional on
$\mathcal{AP}({\bf R})$ with respect to the supremum norm.  In this
case, $G(z, b)$ corresponds in the same way to a continuous mapping
from $B$ into $C_b(\overline{U})$, with respect to the supremum norm
on $C_b(\overline{U})$, as in the preceding paragraph.  Using this,
one can check that the convergence in (\ref{lim_{y to infty} G(x + i
  y, b) = lambda(b)}) is uniform over $x \in {\bf R}$ and $b$ in any
compact subset of $B$ under these conditions.

\section{Upper half-spaces}
\label{upper half-spaces}

        Let $B$ be a compact commutative topological group, and let $h$
be a continuous homomorphism from ${\bf R}$ into $B$, where ${\bf R}$
is considered as a topological group with respect to addition and the
standard topology.  Also let $f$ be a continuous complex-valued function
on $B$, and put
\begin{equation}
\label{F(x, b) = f(h(x) + b)}
        F(x, b) = f(h(x) + b)
\end{equation}
for every $x \in {\bf R}$ and $b \in B$.  Note that $f$ is bounded and
uniformly continuous on $B$, because $B$ is compact.  This implies
that $F(x, b)$ is bounded and uniformly continuous on ${\bf R} \times
B$, where ${\bf R} \times B$ is considered as a topological group with
respect to coordinatewise addition and the product topology.  Of course,
\begin{equation}
\label{F(x + a, b) = F(x, h(a) + b)}
        F(x + a, b) = F(x, h(a) + b)
\end{equation}
for every $a, x \in {\bf R}$ and $b \in B$, because $h : {\bf R} \to B$
is a homomorphism.

        Let $G(z, b)$ be defined on $\overline{U} \times B$ as in the
previous section, so that $G(z, b)$ is given by (\ref{G(x + i y, b) =
  frac{1}{pi} int_{bf R} F(x - w, b) frac{y}{w^2 + y^2} dw}) when $z
\in U$, and (\ref{G(x, b) = F(x, b)}) holds for every $x \in {\bf R}$
and $b \in B$.  It is easy to see that (\ref{b mapsto F_b}) defines a
continuous mapping from $B$ into $C_b({\bf R})$, with respect to the
supremum norm on $C_b({\bf R})$, because $F(x, b)$ is uniformly continuous
on ${\bf R} \times B$.  This implies that $G(z, b)$ determines a continuous
mapping from $B$ into $C_b(\overline{U})$, with respect to the supremum
norm on $C_b(\overline{U})$, as in the preceding section.  In particular,
it follows that $G(z, b)$ is a continuous function on $\overline{U} \times B$,
as before.  We also have that $G(z, b)$ is bounded on $\overline{U} \times B$,
because $F(x, b)$ is bounded on ${\bf R} \times B$, and that
\begin{equation}
\label{G(z + a, b) = G(z, h(a) + b)}
        G(z + a, b) = G(z, h(a) + b)
\end{equation}
for every $a \in {\bf R}$, $z \in \overline{U}$, and $b \in B$, by
(\ref{F(x + a, b) = F(x, h(a) + b)}) and the definition of $G(z, b)$.

        In this situation, (\ref{F_b(x) = F(x, b) in mathcal{AP}({bf R})})
holds for every $b \in B$, as in Section \ref{almost periodic
  functions}.  Thus $\lambda(b)$ can be defined as in (\ref{lambda(b)
  = Lambda(F_b)}) for each $b \in B$, and is a continuous function of
$b$, because (\ref{b mapsto F_b}) is a continuous mapping from $B$
into $\mathcal{AP}({\bf R})$ with respect to the supremum norm.  Note
that the convergence in (\ref{lim_{y to infty} G(x + i y, b) =
  lambda(b)}) is uniform over $x \in {\bf R}$ and $b \in B$ under
these conditions, because $B$ is compact.  It is easy to see that
\begin{equation}
\label{lambda(h(a) + b) = lambda(b)}
        \lambda(h(a) + b) = \lambda(b)
\end{equation}
for every $a \in {\bf R}$ and $b \in B$, because $\Lambda$ is
invariant under translations, as in (\ref{Lambda(f_b) = Lambda(f)}) in
Section \ref{averages on R}, and using (\ref{F(x + a, b) = F(x, h(a) +
  b)}).  This implies that $\lambda$ is constant on the translates of
the closure of $h({\bf R})$ in $B$, since $\lambda$ is continuous on
$B$.  In particular, $\lambda(b)$ is constant on $B$ when $h({\bf R})$
is dense in $B$.  Let $H_B$ be Haar measure on $B$, normalized so that
$H_B(B) = 1$.  If $h({\bf R})$ is dense in $B$, then
\begin{equation}
\label{lambda(b) = Lambda(F_b) = int_B f(u + b) dH_B(u) = int_B f(u) dH_B(u)}
        \lambda(b) = \Lambda(F_b) = \int_B f(u + b) \, dH_B(u)
                                  = \int_B f(u) \, dH_B(u)
\end{equation}
for every $b \in B$, using (\ref{Lambda_h(f) = int_B f(y) dH_B(y)}) in
Section \ref{some other perspectives} in the second step.

        Put
\begin{equation}
\label{g(b, y) = G(i y, b)}
        g(b, y) = G(i \, y, b)
\end{equation}
for every $b \in B$ and $y \in {\bf R}$ with $y \ge 0$, which defines
a bounded continuous function on $B \times [0, \infty)$.  Thus
\begin{equation}
\label{g(b, 0) = G(0, b) = F(0, b) = f(h(0) + b) = f(b)}
        g(b, 0) = G(0, b) = F(0, b) = f(h(0) + b) = f(b)
\end{equation}
for every $b \in B$, and
\begin{equation}
\label{G(x + i y, b) = g(h(x) + b, i y)}
        G(x + i \, y, b) = g(h(x) + b, i \, y)
\end{equation}
for every $x \in {\bf R}$, $y \ge 0$, and $b \in B$, by (\ref{G(z + a,
  b) = G(z, h(a) + b)}).  We also have that
\begin{equation}
\label{lim_{y to infty} g(b, y) = lambda(b)}
        \lim_{y \to \infty} g(b, y) = \lambda(b)
\end{equation}
for every $b \in B$, as in (\ref{lim_{y to infty} G(x + i y, b) =
  lambda(b)}).  The convergence in (\ref{lim_{y to infty} g(b, y) =
  lambda(b)}) is uniform over $b \in B$, because $B$ is compact,
as in the previous paragraph.

\section{Pushing measures forward}
\label{pushing measures forward}

        Let $(X, \,\mathcal{A})$ and $(Y, \mathcal{B})$ be measurable
spaces, so that $X$, $Y$ are sets, and $\mathcal{A}$, $\mathcal{B}$
are $\sigma$-algebras of subsets of $X$, $Y$, respectively.  We may
also refer to elements of $\mathcal{A}$, $\mathcal{B}$ as measurable
subsets of $X$, $Y$, respectively.  Suppose that $h : X \to Y$ is
measurable, in the sense that $h^{-1}(E) \in \mathcal{A}$ for every $E
\in \mathcal{B}$.  If $\mu$ is a nonnegative countably-additive
measure on $(X, \mathcal{A})$, then
\begin{equation}
\label{nu(E) = mu(h^{-1}(E))}
        \nu(E) = \mu(h^{-1}(E))
\end{equation}
defines a nonnegative countably-additive measure on $(Y,
\mathcal{B})$.  If $f$ is a measurable function on $Y$, then $f \circ
h$ is a nonnegative measurable function on $X$, because $h$ is
measurable.  If $f$ is nonnegative, then it is easy to see that
\begin{equation}
\label{int_X f(h(x)) d mu(x) = int_Y f(y) d nu(y)}
        \int_X f(h(x)) \, d\mu(x) = \int_Y f(y) \, d\nu(y),
\end{equation}
by approximating $f$ by nonnegative measurable simple functions on
$Y$.  If $f$ is a measurable real or complex-valued function on $Y$,
then it follows that $f$ is integrable on $Y$ with respect to $\nu$ if
and only if $f \circ h$ is integrable on $X$ with respect to $\mu$, in
which case (\ref{int_X f(h(x)) d mu(x) = int_Y f(y) d nu(y)}) still
holds.

        Similarly, if $\mu$ is a complex-valued countably-additive
measure on $(X, \mathcal{A})$, then (\ref{nu(E) = mu(h^{-1}(E))})
defines a complex-valued countably-additive measure on $(Y,
\mathcal{B})$.  Let $|\mu|$, $|\nu|$ be the total variation measures
on $X$ and $Y$ associated to $\mu$, $\nu$, respectively, which are
finite nonnegative measures on $X$ and $Y$.  Thus
\begin{equation}
\label{|nu(E)| = |mu(h^{-1}(E))| le |mu|(h^{-1}(E))}
        |\nu(E)| = |\mu(h^{-1}(E))| \le |\mu|(h^{-1}(E))
\end{equation}
for every $E \in \mathcal{B}$, which implies that
\begin{equation}
\label{|nu|(E) le |mu|(h^{-1}(E))}
        |\nu|(E) \le |\mu|(h^{-1}(E))
\end{equation}
for every $e \in \mathcal{B}$.  If $f$ is a bounded complex-valued
measurable function on $Y$, then $f \circ h$ is a bounded measurable
function on $X$, and the integrals on both sides of (\ref{int_X
  f(h(x)) d mu(x) = int_Y f(y) d nu(y)}) are defined.  It is easy to
see that (\ref{int_X f(h(x)) d mu(x) = int_Y f(y) d nu(y)}) holds
under these conditions as well.

        Suppose now that $X$, $Y$ are topological spaces, and that
$h : X \to Y$ is continuous.  Thus we can take $\mathcal{A}$, $\mathcal{B}$
to be the corresponding collections of Borel sets, for which $h$ is
measurable.  Let us also ask that $X$ and $Y$ be Hausdorff, so that
compact subsets of $X$ and $Y$ are closed sets, and hence Borel sets.
Let $\mu$ be a nonnegative Borel measure on $X$, and let $\nu$ be
defined on $Y$ as in (\ref{nu(E) = mu(h^{-1}(E))}).  If $K \subseteq
X$ is compact, then $h(K)$ is a compact subset of $Y$, $K \subseteq
h^{-1}(h(K))$, and so
\begin{equation}
\label{mu(K) le mu(h^{-1}(h(K))) = nu(h(K))}
        \mu(K) \le \mu(h^{-1}(h(K))) = \nu(h(K)).
\end{equation}
Using this, one can check that $\nu$ is inner regular on $Y$ when
$\mu$ is inner regular on $X$.  This implies that $\nu$ is outer
regular on $Y$ when $\mu(X) = \nu(Y) < \infty$.  The regularity
properties of a complex Borel measure $\mu$ on $X$ are defined in
terms of the corresponding regularity properties of $|\mu|$.  If
$|\mu|$ is inner regular on $X$, then
\begin{equation}
\label{|mu|(h^{-1}(E))}
        |\mu|(h^{-1}(E))
\end{equation}
defines an inner regular Borel measure on $Y$, as before.  If $\nu$ is
defined on $Y$ as in (\ref{nu(E) = mu(h^{-1}(E))}), then it follows
that $|\nu|$ is inner regular on $Y$ as well, by (\ref{|nu|(E) le
  |mu|(h^{-1}(E))}).  This implies that $|\nu|$ is outer regular on
$Y$ too, since $|\nu|(Y) < \infty$ holds automatically for a complex
measure $\nu$ on $Y$.

        Let $A$ and $B$ be commutative topological groups, and let $h$
be a continuous homomorphism from $A$ into $B$.  Also let $\mu$ be a 
complex Borel on $A$, and let $\nu$ be the complex Borel measure defined
on $B$ as in (\ref{nu(E) = mu(h^{-1}(E))}).  If $\phi \in \widehat{B}$,
then
\begin{equation}
\label{widehat{nu}(phi) = ... = widehat{mu}(phi circ h)}
        \widehat{\nu}(\phi) = \int_B \overline{\phi(y)} \, d\nu(y)
 = \int_A \overline{\phi(h(x))} \, d\mu(x) = \widehat{\mu}(\phi \circ h),
\end{equation}
as in (\ref{int_X f(h(x)) d mu(x) = int_Y f(y) d nu(y)}).

\section{Compact subgroups}
\label{compact subgroups}

        Let $B$ be a commutative topological group, and let $A$
be a compact subgroup of $B$.  Also let $f$ be a continuous complex-valued
function on $B$, and put
\begin{equation}
\label{F(a, b) = f(a + b)}
        F(a, b) = f(a + b)
\end{equation}
for each $a \in A$ and $b \in B$, which is a continuous function on $A
\times B$.  As in Section \ref{product spaces}, $F_b(a) = F(a, b)$
determines a continuous mapping from $B$ into $C(A)$, with respect to
the supremum norm on $C(A)$.  This could also be derived from the fact
that $f$ is uniformly continuous along the translates of $A$ in $B$,
because $A$ and hence its translates are compact subsets of $B$, as in
Section \ref{uniform continuity, equicontinuity}.

        Let $H_A$ be Haar measure on $A$, normalized so that $H_A(A) = 1$,
and put
\begin{equation}
\label{f_A(b) = int_A f(a + b) dH_A(a) = int_A F_b(a) dH_A(a)}
        f_A(b) = \int_A f(a + b) \, dH_A(a) = \int_A F_b(a) \, dH_A(a)
\end{equation}
for every $b \in B$.  This defines a continuous function on $B$,
because $b \mapsto F_b$ is a continuous mappping from $B$ into $C(A)$,
as in the preceding paragraph.  Of course,
\begin{equation}
\label{f_A(a' + b) = f_A(b)}
        f_A(a' + b) = f_A(b)
\end{equation}
for every $a' \in A$, because $H_A$ is invariant under translations
on $A$.

        Let $B / A$ be the quotient of $B$ by $A$, as a commutative
topological group with respect to the quotient topology, and let $q$
be the usual quotient mapping from $B$ onto $B / A$.  Because $f_A$ is
constant on the translates of $A$ in $B$, by (\ref{f_A(a' + b) =
  f_A(b)}), there is a function $\widetilde{f}_A$ on $B / A$ such that
\begin{equation}
\label{widetilde{f}_A(q(b)) = f_A(b)}
        \widetilde{f}_A(q(b)) = f_A(b)
\end{equation}
for every $b \in B$.  More precisely, $\widetilde{f}_A$ is a
continuous function on $B / A$, because $f_A$ is continuous on $B$,
and by the way that the quotient topology on $B / A$ is defined.

        Suppose now that $B$ is locally compact, which implies that
$B / A$ is locally compact.  Let $H_{B / A}$ be a choice of Haar measure
on $B / A$.  If $f$ is a continuous function on $B$ with support
contained in a compact set $K$, then the support of $\widetilde{f}_A$
is contained in $q(K)$, which is a compact subset of $B / A$.  This
implies that
\begin{equation}
\label{int_{B / A} widetilde{f}_A dH_{B / A}}
        \int_{B / A} \widetilde{f}_A \, dH_{B / A}
\end{equation}
is defined, which may be considered as a nonnegative linear functional
on the space $C_{com}(B)$ of continuous complex-valued functions on
$A$.  It is easy to see that this linear functional is invariant under
translations on $B$, because $H_{B / A}$ is invariant under
translations on $B / A$.  If $f$ is a nonnegative real-valued
continuous function on $B$ with compact support such that $f(b_0) > 0$
for some $b_0 \in B$, then $f_A(b_0) > 0$, by the positivity property
of Haar measure on $A$.  This implies that $\widetilde{f}_A(q(b_0)) >
0$, and hence that (\ref{int_{B / A} widetilde{f}_A dH_{B / A}}) is
strictly positive as well, by the positivity property of Haar measure
on $B / A$.

        Thus (\ref{int_{B / A} widetilde{f}_A dH_{B / A}}) satisfies
the requirements of a Haar integral on $B$.  It follows that there
is a choice of Haar measure $H_B$ on $B$ such that
\begin{equation}
\label{int_B f dH_B = int_{B / A} widetilde{f}_A dH_{B / A}}
        \int_B f \, dH_B = \int_{B / A} \widetilde{f}_A \, dH_{B / A}
\end{equation}
for every $f \in C_{com}(B)$.  If $B$ is compact, then $B / A$ is
compact, and one can choose $H_{B / A}$ so that $H_{B / A}(B / A) =
1$, which implies that $H_B(B) = 1$.

        Let $U$ be an open set in $B$ such that $0 \in U$ and the
closure $\overline{U}$ of $U$ in $B$ is compact.  This implies that
\begin{equation}
\label{A + overline{U}}
        A + \overline{U}
\end{equation}
is a compact subset of $B$, because $A$ is compact by hypothesis, and
using the continuity of addition on $B$.  If $K$ is a compact subset
of $B / A$, then $K$ is contained in the union of finitely many
translates of $q(U)$, because $q(U)$ is an open set in $B / A$.  Thus
$q^{-1}(K)$ is contained in the union of finitely many translates of
(\ref{A + overline{U}}), and $q^{-1}(K)$ is also a closed set in $B$,
because $K$ is a closed set in $B / A$, since it is compact.  It
follows that $q^{-1}(K)$ is a compact subset of $B$ under these
conditions, because the union of finitely many translates of (\ref{A +
  overline{U}}) is compact.

        If $g$ is a continuous function on $B / A$ with compact
support, then $g \circ q$ is a continuous function on $B$ with
compact support, by the remarks in the preceding paragraph.  Let
$H_B$ be a choice of Haar measure on $B$, so that
\begin{equation}
\label{int_B g circ q dH_B}
        \int_B g \circ q \, dH_B
\end{equation}
is defined, and determines a nonnegative linear functional on
$C_{com}(B / A)$.  It is easy to see that this linear functional is
invariant under translations and strictly positive when $g \in
C_{com}(B / A)$ is real-valued, nonnegative, and strictly positive at
some point in $C_{com}(B / A)$, because of the corresponding
properties of $H_B$.  This means that (\ref{int_B g circ q dH_B})
satisfies the requirements of a Haar integral on $C_{com}(B / A)$, so
that there is a choice of Haar measure $H_{B / A}$ on $B / A$ such
that
\begin{equation}
\label{int_{B / A} g dH_{B / A} = int_B g circ q dH_B}
        \int_{B / A} g \, dH_{B / A} = \int_B g \circ q \, dH_B
\end{equation}
for every $g \in C_{com}(B / A)$.  Of course, this is basically the
same as (\ref{int_B f dH_B = int_{B / A} widetilde{f}_A dH_{B / A}})
when $f$ is constant on the translates of $A$ in $B$.

        Let $H_B$ be a choice of Haar measure on $B$ again,
and let $f$ be a continuous function on $B$ with support contained
in a compact set $K \subseteq B$.  It is easy to see that the
support of $f_A$ is contained in $A + K$, which is also a compact
subset of $B$, because $A$ is compact.  Under these conditions,
\begin{eqnarray}
\label{int_B f_A(b) dH_B = ... = int_B f(b) dH_B(b)}
 \int_B f_A(b) \, dH_B & = & \int_B \int_A f(a + b) \, dH_A(a) \, dH_B(b) \\
 & = & \int_A \int_B f(a + b) \, dH_B(b) \, dH_A(a) \nonumber \\
 & = & \int_A \int_B f(b) \, dH_B(b) \, dH_A(a) = \int_B f(b) \, dH_B(b),
                                                    \nonumber
\end{eqnarray}
using Fubini's theorem in the second step, invariance of $H_B$ under
translations in the third step, and the normalization $H_A(A) = 1$ in
the last step.

\section{Compact subgroups, continued}
\label{compact subgroups, continued}

        Let $A$ be a compact commutative topological group, and let $H_A$
be Haar measure on $A$, normalized so that $H_A(A) = 1$.  Also let 
${\bf 1}_A$ be the constant function on $A$ equal to $1$ at every point,
which is the identity element of $\widehat{A}$.  If $\phi \in \widehat{A}$
and $\phi \ne {\bf 1}_A$, then
\begin{equation}
\label{int_A overline{phi(a)} dH_A(a) = 0}
        \int_A \overline{\phi(a)} \, dH_A(a) = 0,
\end{equation}
as in (\ref{int_A phi(x) dH(x) = 0}) in Section \ref{compact, discrete
  groups}.  Let $g$ be a complex-valued function on $A$ which is
integrable with respect to $H_A$, and let $\widehat{g}(\phi)$ denotes
the Fourier transform of $g$ evaluated at $\phi \in \widehat{A}$, as
usual.  If $g$ is a constant function on $A$, then
\begin{equation}
\label{widehat{g}(phi) = 0}
        \widehat{g}(\phi) = 0
\end{equation}
for every $\phi \in \widehat{A}$ such that $\phi \ne {\bf 1}_A$, by
(\ref{int_A overline{phi(a)} dH_A(a) = 0}).

          If $g$ is any integrable function on $A$, then put
\begin{equation}
\label{g_1(x) = g(x) - int_A g(a) dH_A(a)}
        g_1(x) = g(x) - \int_A g(a) \, dH_A(a)
\end{equation}
for each $x \in A$, so that $g_1$ is also an integrable function on
$A$, and
\begin{equation}
\label{int_A g_1(a) dH_A(a) = 0}
        \int_A g_1(a) \, dH_A(a) = 0.
\end{equation}
Thus $\widehat{g_1}({\bf 1}_A) = 0$, and
\begin{equation}
\label{widehat{g_1}(phi) = widehat{g}(phi)}
        \widehat{g_1}(\phi) = \widehat{g}(\phi)
\end{equation}
for every $\phi \in \widehat{A}$ such that $\phi \ne {\bf 1}_A$, by
(\ref{int_A overline{phi(a)} dH_A(a) = 0}).  If (\ref{widehat{g}(phi)
  = 0}) holds for every $\phi \in \widehat{A}$ such that $\phi \ne
{\bf 1}_A$, then it follows that $\widehat{g_1}(\phi) = 0$ for every
$\phi \in \widehat{A}$.  This implies that $g_1 = 0$ almost everywhere
on $A$ with respect to $H_A$, as in Section \ref{compact groups}.  Of
course, this means that $g$ is equal to a constant almost everywhere
on $A$ with respect to $H_A$, where the constant is the integral of
$g$ with respect to $H_A$.

        Let $B$ be a commutative topological group, and suppose that
$A$ is a compact subgroup of $B$.  Also let $f$ be a continuous
complex-valued function on $B$, and put
\begin{equation}
\label{F_b(a) = F(a, b) = f(a + b)}
        F_b(a) = F(a, b) = f(a + b)
\end{equation}
for every $a \in A$ and $b \in B$, as in the previous section.
Suppose that $f$ is constant each translate of $A$ in $B$, so that
$F_b(a)$ is a constant function on $A$ for each $b \in B$.  Let $\phi
\in \widehat{B}$ be given, so that the restriction of $\phi$ to $A$
defines an element of $\widehat{A}$.  If $\phi(a) \ne 1$ for some $a
\in A$, then $\phi$ satisfies (\ref{int_A overline{phi(a)} dH_A(a) =
  0}), and hence
\begin{equation}
\label{int_A f(a + b) overline{phi(a)} dH_A(a) = ... = 0}
        \int_A f(a + b) \, \overline{\phi(a)} \, dH_A(a)
          = \int_A F_b(a) \, \overline{\phi(a)} \, dH_A(a) = 0
\end{equation}
for each $b \in B$, because $F_b(a)$ is constant on $A$, by
hypothesis.  Thus
\begin{equation}
\label{int_A f(a + b) overline{phi(a + b)} dH_A(a) = ... = 0}
        \int_A f(a + b) \, \overline{\phi(a + b)} \, dH_A(a)
 = \overline{\phi(b)} \, \int_A f(a + b) \, \overline{\phi(a)} \, dH_A(a) = 0
\end{equation}
for every $b \in B$.  If $B$ is locally compact, $H_B$ is a choice of
Haar measure on $B$, and $f$ has compact support in $B$, then we can
integrate (\ref{int_A f(a + b) overline{phi(a + b)} dH_A(a) = ... =
  0}) over $b \in B$ with respect to $H_B$ as in (\ref{int_B f_A(b)
  dH_B = ... = int_B f(b) dH_B(b)}), to get that
\begin{equation}
\label{widehat{f}(phi) = int_B f(b) overline{phi(b)} dH_B(b) = 0}
        \widehat{f}(\phi) = \int_B f(b) \, \overline{\phi(b)} \, dH_B(b) = 0.
\end{equation}
If $B$ is compact, then this follows from the discussion in Sections
\ref{translation-invariant subalgebras} and \ref{subgroups,
  subalgebras} as well.

        Conversely, suppose that $f$ is a continuous complex-valued
function on $B$ that satisfies (\ref{int_A f(a + b) overline{phi(a)}
  dH_A(a) = ... = 0}) for every $b \in B$ and $\phi \in
\widehat{B}$ such that $\phi(a) \ne 1$ for some $a \in A$.  Let $R$ be
the natural restriction mapping from $\widehat{B}$ into $\widehat{A}$,
which sends each element of $\widehat{B}$ to its restriction to $A$.
Thus our hypothesis on $f$ can be reformulated as saying that
\begin{equation}
\label{widehat{F_b}(psi) = 0}
        \widehat{F_b}(\psi) = 0
\end{equation}
for every $b \in B$ and $\psi \in R(\widehat{B})$ such that $\psi \ne
{\bf 1}_A$, where $\widehat{F_b}$ is the Fourier transform of $F_b$ as
a function on $A$.  If $R(\widehat{B}) = \widehat{A}$, then it follows
that $F_b$ is a constant function on $A$ for every $b \in B$, so that
$f$ is constant on the translates of $A$ in $B$.  Note that
$R(\widehat{B}) = \widehat{A}$ when the elements of $\widehat{B}$
separate points in $A$, as in Section \ref{compact, discrete groups}.
In particular, it is well known that this holds when $B$ is compact,
as mentioned previously.  In this case, the fact that $f$ is constant
on the translates of $A$ in $B$ could also be derived from
(\ref{widehat{f}(phi) = int_B f(b) overline{phi(b)} dH_B(b) = 0}) and
the discussion in Section \ref{translation-invariant subalgebras}.

        Let $B$ be a compact commutative topological group, and let
$A$ be a closed subgroup of $B$, so that $A$ is compact as well.
Also let $E$ be a nonempty subset of $\widehat{B}$, and let
$\mathcal{E}_E(B)$ be the linear span of $E$ in $C(B)$, as in Section
\ref{translation-invariant subalgebras}.  If $f \in \mathcal{E}_E(B)$
and $F_b(a)$ is as in (\ref{F_b(a) = F(a, b) = f(a + b)}), then
it is easy to see that
\begin{equation}
\label{F_b in mathcal{E}_{R(E)}(A)}
        F_b \in \mathcal{E}_{R(E)}(A)
\end{equation}
for every $b \in B$.  Here $R : \widehat{B} \to \widehat{A}$ is the
restriction mapping mentioned in the preceding paragraph, so that
$R(E)$ is the image of $E$ under $R$ in $\widehat{A}$, and
$\mathcal{E}_{R(E)}(A)$ is the linear span of $R(E)$ in $C(A)$.  Let
$C_E(B)$ be the closure of $\mathcal{E}_E(B)$ in $C(B)$ with respect
to the supremum norm, as in Section \ref{translation-invariant
  subalgebras}.  Similarly, let $C_{R(E)}(A)$ be the closure of
$\mathcal{E}_{R(E)}(A)$ in $C(A)$ with respect to the supremum norm.
If $f \in C_E(B)$, then it follows that
\begin{equation}
\label{F_b in C_{R(E)}(A)}
        F_b \in C_{R(E)}(A)
\end{equation}
for every $b \in B$, by approximating $f$ by elements of
$\mathcal{E}_E(B)$, and using (\ref{F_b in mathcal{E}_{R(E)}(A)}) for
the approximations.

        Conversely, suppose that $f \in C(B)$ satisfies
(\ref{F_b in C_{R(E)}(A)}) for every $b \in B$.  This implies that
\begin{equation}
\label{widehat{F_b}(psi) = 0, 2}
        \widehat{F_b}(\psi) = 0
\end{equation}
for every $\psi \in \widehat{A} \setminus R(E)$, as in
(\ref{widehat{f}(phi) = 0}) in Section \ref{translation-invariant
  subspaces}, where $\widehat{F_b}$ is the Fourier transform of $F_b$
as a function on $A$.  In particular, if $\phi \in \widehat{B}$ and
$R(\phi) \not\in R(E)$, then we can take $\psi = R(\phi)$ in
(\ref{widehat{F_b}(psi) = 0, 2}), to get that
\begin{equation}
\label{int_A f(a + b) overline{phi(a)} dH_A(a) = ... = 0, 2}
        \int_A f(a + b) \, \overline{\phi(a)} \, dH_A(a)
          = \int_A F_b(a) \, \overline{\phi(a)} \, dH_A(a) = 0
\end{equation}
for every $b \in B$.  As before, it follows that
\begin{equation}
\label{int_A f(a + b) overline{phi(a + b)} dH_A(a) = ... = 0, 2}
        \int_A f(a + b) \, \overline{\phi(a + b)} \, dH_A(a)
 = \overline{\phi(b)} \, \int_A f(a + b) \, \overline{\phi(a)} \, dH_A(a) = 0
\end{equation}
for every $b \in B$.  Integrating this over $b \in B$ with respect to
Haar measure $H_B$ on $B$ as in (\ref{int_B f_A(b) dH_B = ... = int_B
  f(b) dH_B(b)}), we get that
\begin{equation}
\label{widehat{f}(phi) = int_B f(b) overline{phi(b)} dH_B(b) = 0, 2}
        \widehat{f}(\phi) = \int_B f(b) \, \overline{\phi(b)} \, dH_B(b) = 0,
\end{equation}
where $\widehat{f}$ is the Fourier transform of $f$ as a function on $B$.

        Put
\begin{equation}
\label{E_1 = R^{-1}(R(E))}
        E_1 = R^{-1}(R(E)),
\end{equation}
so that $E \subseteq E_1 \subseteq \widehat{B}$ and hence $R(E)
\subseteq R(E_1) \subseteq R(E)$, which is to say that $R(E) =
R(E_1)$.  If $f \in C(B)$ satisfies (\ref{F_b in C_{R(E)}(A)}) for
each $b \in B$, then (\ref{widehat{f}(phi) = int_B f(b)
  overline{phi(b)} dH_B(b) = 0, 2}) holds for every $\phi \in
\widehat{B} \setminus E_1$, as in the previous paragraph.  This
implies that
\begin{equation}
\label{f in C_{E_1}(B)}
        f \in C_{E_1}(B),
\end{equation}
as in (\ref{C_B(A) = ...}) in Section \ref{translation-invariant
  subalgebras}.  Of course, (\ref{F_b in C_{R(E)}(A)}) is the same as
\begin{equation}
\label{F_b in C_{R(E_1)}(A)}
        F_b \in C_{R(E_1)}(A),
\end{equation}
because $R(E) = R(E_1)$, and every $f \in C_{E_1}(B)$ satisfies
(\ref{F_b in C_{R(E_1)}(A)}) for each $b \in B$, by the same argument
as before.  Thus we get that $f \in C(B)$ is an element of
$C_{E_1}(B)$ if and only if (\ref{F_b in C_{R(E_1)}(A)}) holds for
every $b \in B$.

\section{Another class of subalgebras}
\label{another class of subalgebras}

        Let $B$ be a compact commutative topological group, and let
$h$ be a continuous homomorphism from ${\bf R}$ into $B$, where ${\bf R}$
is considered as a topological group with respect to addition and the
standard topology.  If $\phi \in \widehat{B}$, then
\begin{equation}
\label{widehat{h}(phi) = phi circ h, 2}
        \widehat{h}(\phi) = \phi \circ h
\end{equation}
is a continuous homomorphism from ${\bf R}$ into ${\bf T}$, and
$\widehat{h}$ defines a homomorphism from $\widehat{B}$ into the dual
of ${\bf R}$.  As before, this corresponds to a homomorphism $t_h$
from $\widehat{B}$ into ${\bf R}$, where
\begin{equation}
\label{phi(h(x)) = exp(i x t_h(phi))}
        \phi(h(x)) = \exp(i \, x \, t_h(\phi))
\end{equation}
for every $\phi \in \widehat{B}$ and $x \in {\bf R}$.  Under these
conditions,
\begin{equation}
\label{E(h) = {phi in widehat{B} : t_h(phi) ge 0}}
        E(h) = \{\phi \in \widehat{B} : t_h(\phi) \ge 0\}
\end{equation}
is a sub-semigroup of $\widehat{B}$ that contains the identity element
in $\widehat{B}$.  This leads to a closed subalgebra $C_{E(h)}(B)$ of
$C(B)$ that contains the constant functions and is invariant under
translations, as in Section \ref{translation-invariant subalgebras}.

        Let $f \in C(B)$ be given, and put
\begin{equation}
\label{F(x, b) = f(h(x) + b), 2}
        F(x, b) = f(h(x) + b)
\end{equation}
for every $x \in {\bf R}$ and $b \in B$, as in Section \ref{upper
  half-spaces}.  Thus $F(x, b)$ is bounded uniformly continuous on
${\bf R} \times B$, and an element of $\mathcal{AP}({\bf R})$ as a
function of $x \in {\bf R}$ for every $b \in B$, as before.  If $f \in
C_{E(h)}(B)$, then $f$ can be approximated by linear combinations of
elements of $E(h)$ uniformly on $B$.  This implies that $F(x, b)$ can
be approximated by linear combinations of functions of the form
\begin{equation}
\label{phi(h(x) + b) = phi(h(x)) phi(b)  = exp(i x t_h(x)) phi(b)}
 \phi(h(x) + b) = \phi(h(x)) \, \phi(b)  = \exp(i \, x \, t_h(x)) \, \phi(b)
\end{equation}
with $\phi \in E(h)$, uniformly on ${\bf R} \times B$.

        In particular, for each $b \in B$, $F_b(x) = F(x, b)$ can be
approximated by linear combinations of functions of the form
\begin{equation}
\label{e_t(x) = exp(i x t), 2}
        e_t(x) = \exp(i \, x \, t)
\end{equation}
with $t \ge 0$, uniformly on ${\bf R}$.  This implies that for each $b
\in B$, $F_b(x)$ can be extended to a bounded continuous function on
the closed upper half-plane $\overline{U}$ that is holomorphic on the
open upper half-plane $U$, by the argument given at the end of Section
\ref{bounded holomorphic functions}.  Of course, holomorphic functions
on $U$ are harmonic, so that this extension also corresponds to the
function $G(z, b)$ on $\overline{U} \times B$ discussed in Section
\ref{product spaces}, as in Section \ref{bounded harmonic functions}.

        Conversely, suppose that $f \in C(B)$, and that $F_b(x) = F(x, b)$
has a bounded continuous extension to $\overline{U}$ that is
holomorphic on $U$ for each $b \in B$.  Remember that $F_b \in
\mathcal{A}({\bf R})$ for every $b \in B$, and that $\langle \cdot,
\cdot \rangle_{\mathcal{AP}({\bf R})}$ is the inner product on
$\mathcal{AP}({\bf R})$ defined in Section \ref{the associated inner
product}.  Under these conditions, we have that
\begin{equation}
\label{langle F_b, e_t rangle_{mathcal{AP}({bf R})} = 0}
        \langle F_b, e_t \rangle_{\mathcal{AP}({\bf R})} = 0
\end{equation}
for every $b \in B$ and $t < 0$, as in (\ref{langle f, e_t
  rangle_{mathcal{AP}({bf R})} = 0}) in Section \ref{bounded
  holomorphic functions}.  We would like to use this to show that
\begin{equation}
\label{widehat{f}(phi) = int_B f(b) overline{phi(b)} dH_B(b) = 0, 3}
        \widehat{f}(\phi) = \int_B f(b) \, \overline{\phi(b)} \, dH_B(b) = 0
\end{equation}
when $\phi \in \widehat{B}$ satisfies $t_h(\phi) < 0$, which is to say
that $\phi \in \widehat{B} \setminus E(h)$.  This would imply that $f$
is an element of $C_{E(h)}(B)$, by (\ref{C_B(A) = ...}).

        If $h({\bf R})$ is dense in $B$, then
\begin{equation}
\label{widehat{f}(phi) = langle f circ h, e_{t_h(phi)} rangle, 2}
        \widehat{f}(\phi)
             = \langle f \circ h, e_{t_h(\phi)} \rangle_{\mathcal{AP}({\bf R})}
\end{equation}
for every $\phi \in \widehat{B}$, by (\ref{widehat{f}(phi) = langle f
  circ h, e_{t_h(phi)} rangle}) in Section \ref{some additional
  properties}.  By construction, $f \circ h = F_0$, so that
(\ref{widehat{f}(phi) = int_B f(b) overline{phi(b)} dH_B(b) = 0, 3})
follows from (\ref{langle F_b, e_t rangle_{mathcal{AP}({bf R})} = 0})
and (\ref{widehat{f}(phi) = langle f circ h, e_{t_h(phi)} rangle, 2})
when $t_h(\phi) < 0$, as desired.

        Otherwise, let $A$ be the closure of $h({\bf R})$ in $B$,
which is a compact subgroup of $B$, because $V$ is compact.  Put
\begin{equation}
\label{F'_b(a) = F'(a, b) = f(a + b)}
        F'_b(a) = F'(a, b) = f(a + b)
\end{equation}
for every $a \in A$ and $b \in B$, and let $H_A$ be Haar measure on
$A$, normalized so that $H_A(A) = 1$.  If $\phi \in \widehat{B}$, then
\begin{equation}
\label{int_A F'_b(a) overline{phi(a)} dH_A(a) = ...}
 \int_A F'_b(a) \, \overline{\phi(a)} \, dH_A(a)
            = \langle F'_b \circ h, \phi \circ h \rangle_{\mathcal{AP}({\bf R})}
            = \langle F_b, e_{t_h(\phi)} \rangle_{\mathcal{AP}({\bf R})}
\end{equation}
for every $b \in B$.  This is analogous to (\ref{widehat{f}(phi) =
  langle f circ h, e_{t_h(phi)} rangle, 2}), and can be derived from
(\ref{langle f circ h, g circ h rangle_{mathcal{AP}({bf R})} = ...})
as in (\ref{widehat{f}(phi) = ... = langle f circ h, phi circ h
  rangle}) and (\ref{widehat{f}(phi) = langle f circ h, e_{t_h(phi)}
  rangle}) in Section \ref{some additional properties}.  Combining
this with (\ref{langle F_b, e_t rangle_{mathcal{AP}({bf R})} = 0}),
we get that
\begin{equation}
\label{int_A f(a + b) overline{phi(a)} dH_A(a) = ... = 0, 3}
        \int_A f(a + b) \, \overline{\phi(a)} \, dH_A(a)
          = \int_A F'_b(a) \, \overline{\phi(a)} \, dH_A(a) = 0
\end{equation}
for every $b \in B$ and $\phi \in \widehat{B}$ with $t_h(\phi) < 0$.
It follows that
\begin{equation}
\label{int_A f(a + b) overline{phi(a + b)} dH_A(a) = ... = 0, 3}
        \int_A f(a + b) \, \overline{\phi(a + b)} \, dH_A(a)
 = \overline{\phi(b)} \, \int_A f(a + b) \, \overline{\phi(a)} \, dH_a(a) = 0
\end{equation}
for every $b \in B$ and $\phi \in \widehat{B}$ with $t_h(\phi) < 0$.
This implies that (\ref{widehat{f}(phi) = int_B f(b) overline{phi(b)}
  dH_B(b) = 0, 3}) holds for every $\phi \in \widehat{B}$ such that
$t_h(\phi) < 0$, by integrating (\ref{int_A f(a + b) overline{phi(a +
    b)} dH_A(a) = ... = 0, 3}) over $b \in B$ with respect to $H_B$,
and using (\ref{int_B f_A(b) dH_B = ... = int_B f(b) dH_B(b)}).

        Note that
\begin{equation}
\label{E(h)^{-1} = {phi in widehat{B} : t_h(phi) le 0}}
        E(h)^{-1} = \{\phi \in \widehat{B} : t_h(\phi) \le 0\},
\end{equation}
where $E(h)^{-1}$ is as in (\ref{B^{-1} = {1/phi : phi in B} =
  {overline{phi} : phi in B}}) in Section \ref{translation-invariant
  subalgebras}, so that
\begin{equation}
\label{E(h) cap E(h)^{-1} = {phi in widehat{B} : t_h(phi) = 0}}
        E(h) \cap E(h)^{-1} = \{\phi \in \widehat{B} : t_h(\phi) = 0\}.
\end{equation}
This is the same as the kernel of $\widehat{h}$, which is trivial
when $h({\bf R})$ is dense in $B$.

        Suppose that $C$ is a commutative group equipped with the
discrete topology, and that $k$ is a homomorphism from $C$ into ${\bf
  R}$, as a commutative group with respect to addition.  In
particular, one can simply take $C$ to be a subgroup of ${\bf R}$,
and $k$ to be the natural inclusion mapping.  Remember that the
dual $\widehat{C}$ of $C$ is a compact commutative topological group,
with respect to the usual topology.  The dual homomorphism $\widehat{k}$
associated to $k$ is a continuous homomorphism from the dual of ${\bf R}$
into $\widehat{C}$, and the dual of ${\bf R}$ is isomorphic to ${\bf R}$
as a commutative topological group.  This takes us back to the
same type of situation as at the beginning of the section, with
\begin{equation}
\label{B = widehat{C}}
        B = \widehat{C},
\end{equation}
and where $h$ corresponds to $\widehat{k}$, and $\widehat{h}$ can be
identified with $k$.

\backmatter

\newpage

\addcontentsline{toc}{chapter}{Index}

\printindex

\end{document}